Н. М. Кіяновська, Н. В. Рашевська,
С. О. Семеріков

# Теоретико-методичні засади використання інформаційно-комунікаційних технологій у навчанні вищої математики студентів інженерних спеціальностей у Сполучених Штатах Америки

Спецвипуск містить монографію Н. М. Кіяновської, Н. В. Рашевської, С. О. Семерікова, у якій виокремлено етапи розвитку теорії та методики використання ІКТ у навчанні вищої математики студентів інженерних спеціальностей у США, теоретично обґрунтовано та розроблено дидактичні моделі використання ІКТ навчання вищої математики у технічних ВНЗ на виокремлених етапах, визначено основні підходи для застосування досвіду США щодо використання ІКТ у навчанні вищої математики студентів інженерних спеціальностей в Україні. Практична частина дослідження містить рекомендації для викладачів з використання ІКТ навчання вищої математики студентів інженерних спеціальностей.

Для науковців, викладачів та студентів вищих навчальних закладів, аспірантів та всіх тих, кого цікавлять сучасні теорія та методика використання ІКТ в освіті.

## ЗМІСТ

# ПЕРЕЛІК УМОВНИХ ПОЗНАЧЕНЬ

| | |
|---|---|
| ABET | Accreditation Board for Engineering and Technology (Рада з акредитації інженерних розробок і технологій) |
| LMS | Learning Management System (система управління навчанням) |
| MOOC | Massive open online course (масовий відкритий дистанційний курс – МВДК) |
| ВНЗ | вищий навчальний заклад |
| ДВНЗ | Державний вищий навчальний заклад |
| ДН | дистанційне навчання |
| ЕОР | електронні освітні ресурси |
| ІКТ | інформаційно-комунікаційні технології |
| ІТ | інформаційні технології |
| КНУ | Криворізький національний університет |
| КОНМК | комп'ютерно орієнтований навчально-методичний комплекс |
| ММС | мобільне математичне середовище |
| МТІ | Массачусетський технологічний інститут (Massachusetts Institute of Technology – MIT) |
| ООН | Організація Об'єднаних Націй |
| ППЗ | програмно-педагогічний засіб |
| СКМ | систем комп'ютерної математики |
| СПН | системи підтримки навчання |
| США | Сполучені Штати Америки |
| ЮНЕСКО | United Nations Educational, Scientific and Cultural Organization, UNESCO (Організація Об'єднаних Націй з питань освіти, науки і культури) |

# ВСТУП

На сучасному етапі розвитку інформаційного суспільства використання засобів інформаційно-комунікаційних технологій (ІКТ) сприяє глобалізації освіти, розвитку міжнародного ринку праці, зростанню різних видів мобільності особистості. Важливим наслідком глобалізації є підвищення мобільності студентів, абітурієнтів та випускників університетів: особа, що має високий рівень мобільності, можете вчитися, працювати, співпрацювати та бути конкурентоздатною в будь-якій країні. Зростання академічної мобільності, уведення міжнародних норм і стандартів, за допомогою яких академічні кваліфікації з різних країн можуть бути порівняні та визнані, призводить до збільшення конкуренції між ВНЗ та сприяє підвищенню якості вищої освіти.

Необхідною умовою суспільного й економічного розвитку будь-якої країни є інвестиції в освіту населення. У цьому контексті глобалізація освіти сприяє особистісному та професійному розвитку фахівців, які займаються розробкою та впровадженням нових технологій – інженерів.

Вищі технічні навчальні заклади США мають значні педагогічні досягнення і розвинену систему підготовки фахівців інженерних напрямів на основі системного використання засобів ІКТ. У глобалізованому просторі вищої освіти проблему підвищення якості підготовки фахівців у вітчизняних ВНЗ доцільно розв’язувати через інтеграцію з кращими здобутками світової педагогічної думки і творче використання досвіду передових ВНЗ інженерного профілю. Проте різниця у технологічних укладах, що домінують у США (п’ятий та шостий) та в Україні (четвертий та п’ятий), вимагає звернення не тільки до сучасного досвіду використання засобів ІКТ у інженерній освіті США, а й до аналізу історичного досвіду.

Для України інформатизація інженерної освіти є надзвичайно актуальним у контексті її економічного, соціального та культурного розвитку: як зазначено у Законі України «Про основні засади розвитку інформаційного суспільства в Україні на 2007-2015 роки», основним напрямом використання ІКТ є створення системи освіти, орієнтованої на використання новітніх ІКТ у формуванні всебічно розвиненої особистості, що надає можливість кожній людині самостійно здобувати знання, уміння та навички під час навчання, виховання та професійної підготовки [263].

Метою Національної стратегії розвитку освіти в Україні на 2012-2021 роки є: оновлення змісту, форм, методів і засобів навчання шляхом широкого впровадження у навчально-виховний процес сучасних ІКТ та електронного контенту. А пріоритетом розвитку освіти є впровадження сучасних ІКТ, що забезпечують удосконалення навчально-виховного процесу, доступність та ефективність освіти, підготовку молодого покоління до життєдіяльності в інформаційному суспільстві [241].

Однією із складових системи професійної підготовки сучасного інженера є фундаментальна підготовка, основним завданням якої є удосконалення професійної підготовки і всебічного розвитку студента як особистості та яке включає в себе: оволодіння науково-дослідницькими методами розв'язання виробничих задач; розробку раціоналізаторських пропозицій і участь у винахідницькій роботі; врахування технічного прогресу і еволюцію потреб, щоб керуватись не лише усталеною практикою, а схилятись до новаторської позиції в інженерній діяльності; знання технології і техніки із сфери своєї спеціалізації (спеціальності) та оволодіння різноманітними формами самоосвіти, що неможливе без ґрунтовних знань з вищої математики та умінь застосовувати набуті знання на практиці та професійній діяльності.

Використання ІКТ у процесі навчання вищої математики студентів інженерних спеціальностей створює умови для самореалізації студента, що сприяє підвищенню його пізнавальної активності, розвитку критичного мислення, формуванню у студентів навичок організації самостійної роботи, розвитку творчих здібностей та лідерських якостей, підвищенню відповідальності за результати своєї праці, а також вдосконаленню процесу навчання та підвищенню його якості.

Тому виникає необхідність дослідження історії та сучасного стану розвитку засобів ІКТ навчання вищої математики студентів інженерних спеціальностей у технічних ВНЗ США, що займають найвищі позиції у рейтингу найкращих ВНЗ світу [140], з метою модернізації системи вищої інженерної освіти України та її спрямування на підготовку фахівців, здатних до швидкого просування науково-технічного прогресу.

Аналіз наукової літератури показав, що проблеми розвитку теорії й методики використання ІКТ в освіті традиційно перебувають у полі зору вітчизняних науковців. Теоретичні основи застосування ІКТ у процесі навчання досліджувались у роботах М. І. Жалдака [172], Ю. І. Машбиця [232], Н. В. Морзе [234], Ю. С. Рамського [267], С. О. Семерікова [292], Ю. В. Триуса [316] та інших дослідників. У роботах В. Ю. Бикова [169], Ю. Г. Запорожченко [180], М. П. Лещенко [227], О. М. Спіріна [306], О. В. Овчарук [243], Н. В. Сороко [302], Б. І. Шуневича [328] та інших здійснені порівняльно-педагогічні дослідження щодо зарубіжного досвіду застосування ІКТ в освіті. Теорія та методики використання ІКТ у навчанні вищої математики розроблялись у роботах К. В. Власенко [158], Ю. В. Горошка [162], В. І. Клочка [214], С. А. Ракова [266],

О. В. Співаковського [303] та інших.

Серед робіт зарубіжних науковців важливими є дослідження:

– К. Блертона (C. Blurton) [14], Л. Ларсона (L. Larson) [71] з історії та сучасного стану використання ІКТ у процесі навчання;

– Дж. Вавріка (J. Wavrik) [136], Дж. Енгельбрехта (J. Engelbrecht) [39], Дж. Панкіна (J. Pankin) [96], Р. Пеа (R. Pea) [98], Дж. Харві (J. Harvey) [55] з використання ІКТ у процесі навчання математичних дисциплін студентів інженерних спеціальностей;

– Д. Меріно (D. Merino) [78], Б. Хана (B. Khan) [68], Дж. Ітмазі (J. Itmazi) [64], Дж. Гамільтона (J. Hamilton) [52], Дж. Прадоса (J. Prados) [105] з використання ІКТ у підготовці студентів інженерних спеціальностей у США.

Існують певні суперечності, зокрема, між сучасними вимогами до фахівця інженера та реальним рівнем їх підготовки у ВНЗ, прагненням підвищувати кваліфікацію викладачів математичних дисциплін та рівнем їх обізнаності у засобах ІКТ навчання математичних дисциплін. Залишаються недослідженими загальні тенденції розвитку засобів ІКТ навчання вищої математики студентів інженерних спеціальностей у США у контексті їх еволюції та конвергенції, що визначило необхідність проведення даного дослідження, хронологічні межі якого охоплюють період з 1965 р. по сьогодення.

Монографія складається з трьох розділів.

У першому розділі охарактеризовано понятійно-термінологічний дискурс проблеми дослідження та здійснено історико-педагогічний аналіз розвитку ІКТ навчання вищої математики студентів інженерних спеціальностей США.

У другому розділі проаналізовано провідні ІКТ навчання вищої математики студентів інженерних спеціальностей США на сучасному етапі, Національний план використання освітніх ІКТ (National Education Technology Plan, Transforming American Education: Learning Powered by Technology), розглянуто тенденції та моделі використання ІКТ у навчанні вищої математики студентів інженерних спеціальностей США.

У третьому розділі зроблено аналіз вітчизняних ІКТ навчання вищої математики, обґрунтовано можливості застосування досвіду США для удосконалення процесу навчання вищої математики у технічних ВНЗ України, розроблено структурно-функціональну схему та рекомендації щодо використання ІКТ у навчанні вищої математики студентів інженерних спеціальностей України, базуючись на позитивному досвіді США, визначено ІКТ-компетентності викладача вищої математики, описано розроблений спецкурс, спрямований на розвиток ІКТ-компетентності викладача вищої математики.

Автори щиро дякують ініціаторам створення монографії – завідувачу кафедри комп'ютерних технологій Черкаського державного технологічного університету д. пед. н., проф. Ю. В. Триусу, завідувачу відділу інформатизації навчально-виховних закладів Інституту інформаційних технологій і засобів навчання НАПН України к. філос. н., с. н. с. М. П. Шишкіній, та всім співробітникам спільної науково-дослідної лабораторії з питань використання хмарних технологій в освіті ДВНЗ «Криворізький національний університет» та Інституту інформаційних технологій і засобів навчання НАПН України, які вносили пропозиції щодо структури та змісту цієї роботи.

# РОЗДІЛ 1
# ТЕОРЕТИЧНІ ОСНОВИ ВИКОРИСТАННЯ ІНФОРМАЦІЙНО-КОМУНІКАЦІЙНИХ ТЕХНОЛОГІЙ У НАВЧАННІ ВИЩОЇ МАТЕМАТИКИ СТУДЕНТІВ ІНЖЕНЕРНИХ СПЕЦІАЛЬНОСТЕЙ У СПОЛУЧЕНИХ ШТАТАХ АМЕРИКИ

## 1.1 Характеристика інженерної освіти в Сполучених Штатах Америки

*1.1.1 Зміст математичної складової інженерної освіти в Сполучених Штатах Америки.* Стрімкий розвиток технологій, спричинений науково-технічним прогресом, вимагає підготовки спеціалістів високого рівня, наділеного аналітичними здібностями, здатного шукати і знаходити необхідні дані, точно формулювати проблеми, для певної сукупності даних виводити закономірності, вміти розв'язувати складні міждисциплінарні задачі. Технічна освіта спрямована на розв'язання саме таких задач, готуючи спеціалістів – інженерів.

За С. І. Ожеговим [244] *інженер* – це спеціаліст з вищою технічною освітою.

*Інженер* – особа, що професійно займається інженерією, тобто на основі поєднання прикладних наукових знань, математики та винахідництва знаходить нові розв'язання технічних проблем. Зміст творчості інженера дає вагомі підстави визнавати інженерів одними з основних творців ноосфери в частині матеріальної культури та прикладної науки, «відповідальних» за науково-технічний прогрес (загально) людської цивілізації та, відповідно, «технологічний добробут» людства [186].

Під *інженерною освітою* будемо розуміти систему підготовки фахівців для роботи в галузях технічних наук: металургії, гірництві, ливарній справі, машинознавстві, електротехніці, теплотехніці, гідротехніці, радіотехніці, будівництві й інших.

Інженерна освіта має достатньо давню історію (принаймні з XVI ст.), проте як усталений термін (engineering education) у англомовних джерелах з'явилась лише у 1830-х рр. (рис. 1.1).

Так, у перший рік видання журналу «The Civil Engineer and Architect's Journal» (1837-1838 рр.) опубліковано редакційну статтю про інженерну освіту цивільних інженерів [93], у якій аналізується сучасний на той момент стан підготовки інженерів у Європі. Зокрема, автор указує, що французька «політехнічна школа – це не просто цивільна інженерна школа, вона готує інженерів артилерії, морської піхоти, армії і флоту для військово-інженерної діяльності, виробництво пороху і селітри, і для всіх суспільних ситуацій, що вимагають великих знань у галузі фізико-математичних наук. ... Основними навчальними предметами є аналіз, механіка, геодезія і топографія, статистика, фізика, хімія, архітектура, французька, німецька та англійська мови. ... Повна освіта триває два роки». У статті наведено приклад Женевської та Стокгольмської технологічних шкіл з трирічним терміном навчання. У статті повідомляється, що англійська інженерна освіта є більш досконалою, тому що спонукає до самостійної роботи студентів, виховує відповідальність та ініціативність: «важливе джерело нашої переваги полягає в готовності прийняти все, що є новим і корисним, відкидаючи все, що застаріло або абсурдне. Це саме можна сказати і про інженерну освіту в Сполучених Штатах», однією з головних проблем якої є необхідність підвищення престижності професії інженера.

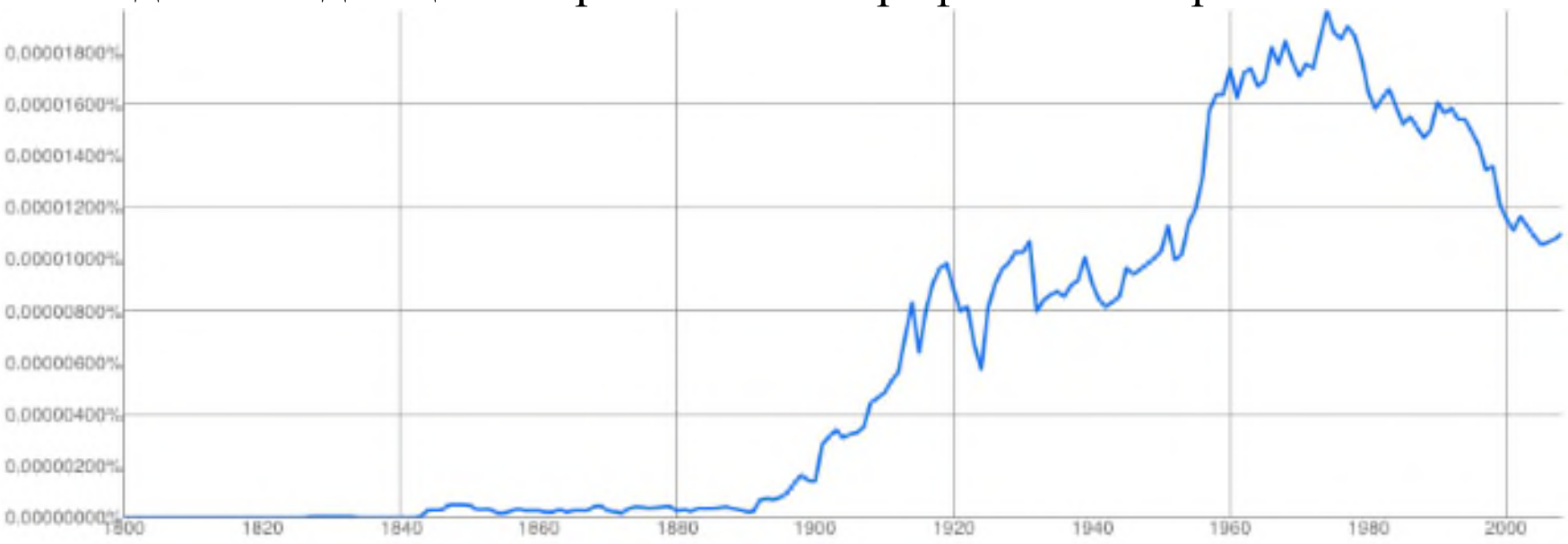

Рис. 1.1. Використання терміну «інженерна освіта» (engineering education) у англомовних джерелах

У середині 1990-х рр. у США було випущено 60000 бакалаврів інженерії – 5% від усього випуску бакалаврів. Для порівняння, 90000 бакалаврів інженерії в Японії склали 18% від усього випуску бакалаврів, а 120000 у Китаї – 46% [52, 22]. Це відповідає загальним тенденціям кінця XX ст. (рис. 1.2) і «не є проблемою лише Великобританії, Німеччини та США: в більшості західних країн повідомляють про подібні труднощі у спрямуванні молоді на науково-технічні спеціальності» [52, 24], проте і через 160 років «низький соціальний статус інженера у Великобританії залишається проблемою» [52, 32].

Базовий рівень інженерної освіти у США – бакалавр (4 роки). На відміну від Великобританії та України, у США відсутні державні галузеві стандарти вищої освіти: натомість існує потужна система акредитації на чолі з ABET (Accreditation Board for Engineering and Technology) – неурядовою організацією, що, зокрема, оцінює якість підготовки на інженерних спеціальностях на основі «Інженерних критеріїв 2000 року» (Engineering Criteria 2000 – EC2000) [72].

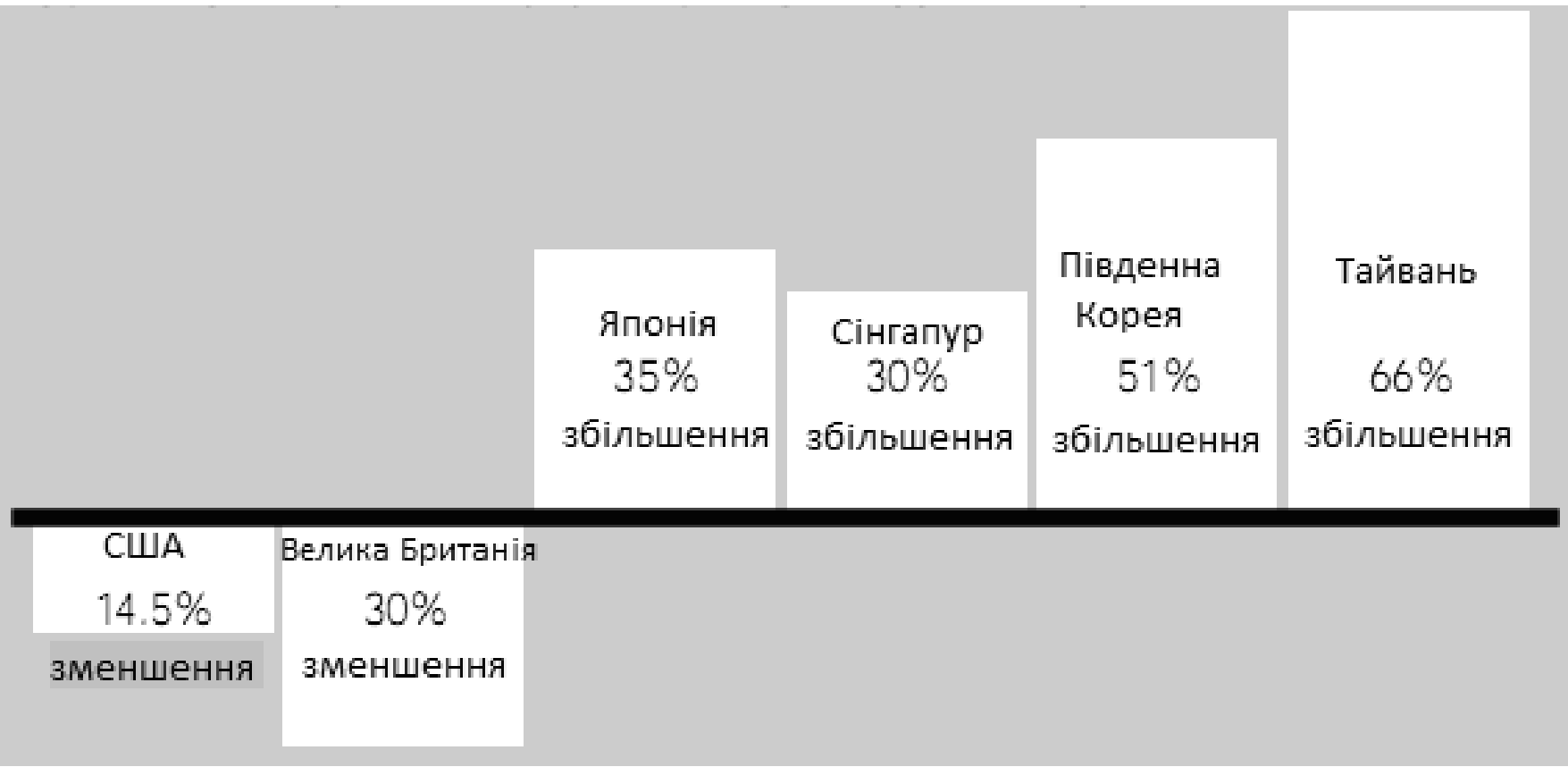

Рис. 1.2. Зміна кількості випускників інженерних спеціальностей у 1985-1998 рр. (за Дж. Гамільтоном (Sir James Hamilton) [52, 24])

Серед 8 критеріїв у EC2000 (1. Студенти. 2. Освітні цілі програми. 3. Основні задачі програми і оцінка ступеня їх виконання. 4. Професійна підготовка. 5. Викладацький склад. 6. Матеріальна база. 7. Підтримка з боку ВНЗ та фінансові ресурси. 8. Спеціальні вимоги до програма) найбільшу увагу привертає третій. Цей критерій визначає вимоги до знань та вмінь випускників інженерних ВНЗ, серед яких найвищу вагу має «здатність застосовувати прикладні знання з вищої математики, науки та інженерії у професійній діяльності» [72, 18].

Як зазначає Дж. Прадос (John W. Prados), у інженерній освіті США за останні 60 років відбулись суттєві зміни, пов'язані з післявоєнним розвитком технологій. «У 1947 році інженерія була високоприкладною галуззю із незначним застосуванням математики на рівні елементарних обчислень. ... У період з 1950 по 1960 рр. відбулася зміна парадигми інженерної освіти з прикладної на математичну, академічну, інженерно-наукову. Курси будови машин ... та подібні до них зникли та були замінені диференціальними рівняннями, теорією управління...» [105, 1].

У 1985 р. у звіті з інженерної освіти США зазначалось, що «погана підготовка з математики ... перешкоджає інженерній підготовці» [40, 65]. Для виправлення цього пропонувалось увести в курс середньої школи (К-12) інтегровану природничо-наукову програму, що отримала назву STEM – наука (science), технологія (technology), інженерія (engineering) та математика (mathematics) [11]. На подолання кризи інженерної освіти спрямована також програма Національної академії інженерії США (National Academy of Engineering) «Інженер 2020 року» (Engineer of 2020), метою якої є підготовка технологічно досконалих та інноваційних інженерів, здатних працювати у швидкозмінному середовищі [34].

Наприкінці XX століття посилився вплив інформаційних технологій на інженерну освіту і практику у зв'язку із тенденціями глобалізації виробництва та переходом до суспільства сталого розвитку, тому новими компетенціями інженера стали навички комунікації, спільної роботи, навчання протягом всього життя [105, 2]. Для їх реалізації була необхідна нова освітня парадигма активного проектно-орієнтованого навчання: посилення прикладної спрямованості навчання вищої математики, тісний зв'язок із виробництвом, широке використання ІКТ.

На початку XXI століття у США студенти всіх інженерних спеціальностей на першому курсі навчання вивчають вищу математику, загальну хімію, англійську мову, загальну та сучасну фізику, комп'ютерні науки (насамперед програмування), вступ до інженерії. Загальноінженерні курси також включають інженерну графіку, інженерію матеріалів, інженерну механіку, опір матеріалів, електроінженерію, термодинаміку, механіку рідин. Із другого курсу починають вивчати дисципліни спеціалізації.

Аналіз джерел з проблеми дослідження надав можливість зробити висновок, що сучасна інженерна освіта США має такі основні особливості: відсутність державних галузевих стандартів; недержавна система акредитації; математизація та комп'ютеризація загальноінженерних та спеціальних дисциплін; прикладна спрямованість навчання вищої математики; широке використання засобів ІКТ у навчанні вищої математики.

Вказані особливості зумовлюють необхідність узагальнення та систематизації різних акредитованих програм підготовки з вищої математики з метою виділення провідних напрямів навчання вищої математики майбутніх інженерів у США.

Незважаючи на недержавну форму акредитації та традиційне різноманіття пропонованих математичних курсів (як обов'язкових, так і факультативних), навчання вищої математики майбутніх інженерів у США здійснюється за схожими навчальними програмами. Розглянемо, наприклад, зміст навчання вищої математики студентів інженерних спеціальностей одного із провідних ВНЗ США – Массачусетського технологічного інституту (МТІ, Massachusetts Institute of Technology – MIT).

У 2012 році MIT приймав на заняття в осінньому семестрі (Fall 2012) за 44 напрямами підготовки [83]. Усталені напрями підготовки є пронумерованими від Курсу 1 (Course 1) до Курсу 24 (Course 24): саме цей номер і використовується для позначення того, студенти якого напряму мають відвідувати той чи інший навчальний курс. Наприклад, 18.03 Differential Equations – це курс диференціальних рівнянь для напряму підготовки Course 18 – Mathematics.

Споріднені напрями підготовки об’єднані у школи, найбільш популярною з яких у 2011 р. була Школа інженерії МТІ (MIT School of Engineering), у якій навчалось 62,7% студентів МТІ за 19 напрямами інженерної підготовки, провідними з яких є такі:

Курс 1 – Цивільна інженерія та охорона навколишнього середовища (Civil and Environmental Engineering);

Курс 2 – Машинобудування (Mechanical Engineering);

Курс 3 – Матеріалознавство (Materials Science and Engineering);

Курс 6 – Електротехніка та комп’ютерні науки (Electrical Engineering and Computer Science);

Курс 10 – Хімічна інженерія (Chemical Engineering);

Курс 16 – Аеронавтика та астронавтика (Aeronautics and Astronautics);

Курс 20 – Біотехнології (Biological Engineering);

Курс 22 – Ядерна інженерія (Nuclear Science and Engineering).

Крім указаних напрямів, починаючи з 2012 р., Школа інженерії МТІ [80] пропонує гнучку програму підготовки інженерів (The MIT Flexible Engineering Degree Program), що відображає найновіші досягнення у галузі інженерії та відповідає концепціям, викладеним Національною академією інженерії США у програмі адаптації інженерної освіти до нового століття [34].

Всі першокурсники (freshman) мають опанувати набір обов’язкових навчальних дисциплін (core curriculum) – загальноінститутських вимог (General Institute Requirements – GIRs) [79]. Нормативна підготовка з вищої математики включає в себе елементи математичного аналізу функції однієї змінної та багатьох змінних (Calculus I та Calculus II відповідно): границі, функції, похідна, інтеграли, ряди.

Числення однієї змінної (Calculus I, шифр курсу 18.01) та числення багатьох змінних (Calculus II, шифр курсу 18.02) пропонуються у декількох версіях: основний курс (шифри 18.01, 18.02), додаткові розділи (шифри 18.014, 18.022), факультативний курс (шифри 18.01A, 18.02A).

Студенту необхідно опанувати два блоки: (18.01 або 18.01A, або 18.014) та (18.02 або 18.02A, або 18.022, або 18.023, або 18.024). Із навчальних дисциплін, перерахованих у дужках, необхідно обрати лише одну – якщо вибрати дві та більше, кредити із «зайвих» дисциплін не враховуватимуться. Так, у 18.014 – числення з теорією (Calculus with Theory) пропонується той самий матеріал, що й у 18.01, проте на більш глибокому рівні, а 18.01A є оглядовим шеститижневим курсом числення функції однієї змінної, призначеним для тих, хто вже опанував однорічний курс математичного аналізу на достатньому та високому рівнях (наприклад, у коледжі), з метою швидкого перезарахування кредитів.

До складу курсу «Числення однієї змінної» включено такі основні теми: диференціальне та інтегральне числення функцій однієї змінної, їх застосування; неформальний вступ до границь та нескінченності; диференціювання: означення, основні правила, застосування до побудови графіків функцій, наближення, екстремуми; невизначений інтеграл; диференціальні рівняння першого порядку з відокремлюваними змінними; визначений інтеграл, основна теорема аналізу; застосування інтегралів у геометрії; елементарні функції; методи інтегрування; полярні координати; правило Лопіталя; невласні інтеграли; ряди: геометричні, гармонійні, ознаки порівняння рядів, степеневі ряди для деяких елементарних функцій.

До складу курсу «Числення багатьох змінних» включено такі основні теми: числення декількох змінних; векторна алгебра в тривимірному просторі, визначники, матриці; вектор-функції однієї змінної, рух у просторі; скалярні функції декількох змінних: частинні похідні, градієнт, методи оптимізації; подвійні інтеграли та криволінійні інтеграли на площині; точні диференціали та потенційне векторне поле; теорема Гріна та її застосування, потрійні інтеграли, лінійні та поверхневі інтеграли у просторі, теорема дивергенції, теорема Стокса; застосування числення декількох змінних.

Крім обов’язкових навчальних дисциплін (core curriculum) на різних напрямах підготовки пропонується додатково певний список дисциплін для кожного напряму підготовки.

Так, наприклад, на кафедрі Електротехніка та комп’ютерні науки (Electrical Engineering & Computer Science – EECS), що є найбільшою в MTI, студентам пропонується чотири освітні бакалаврські програми (6-1, 6-2, 6-3, 6-7) і дві програми магістра інженерії (6 – A або 6 – P) [130]. Опишемо ці курси.

**Курс 6-1**: Electrical Science and Engineering. На курсі 6-1 студенти опановують схеми і пристрої, матеріали і нанотехнології, комунікації, управління та обробку сигналів, а також прикладну фізику.

Студенти курсу 6-1 Electrical Science and Engineering повинні пройти курси із загальноінститутських вимог (General Institute Requirements – GIRs), а також кафедральні курси. Порядок появи дисциплін в курсі 6-1 наведено на рис. 1.3.

Як видно з рис. 1.3, крім курсу «Числення багатьох змінних» студенту необхідно додатково опанувати курс «Диференціальні рівняння» (Differential Equations).

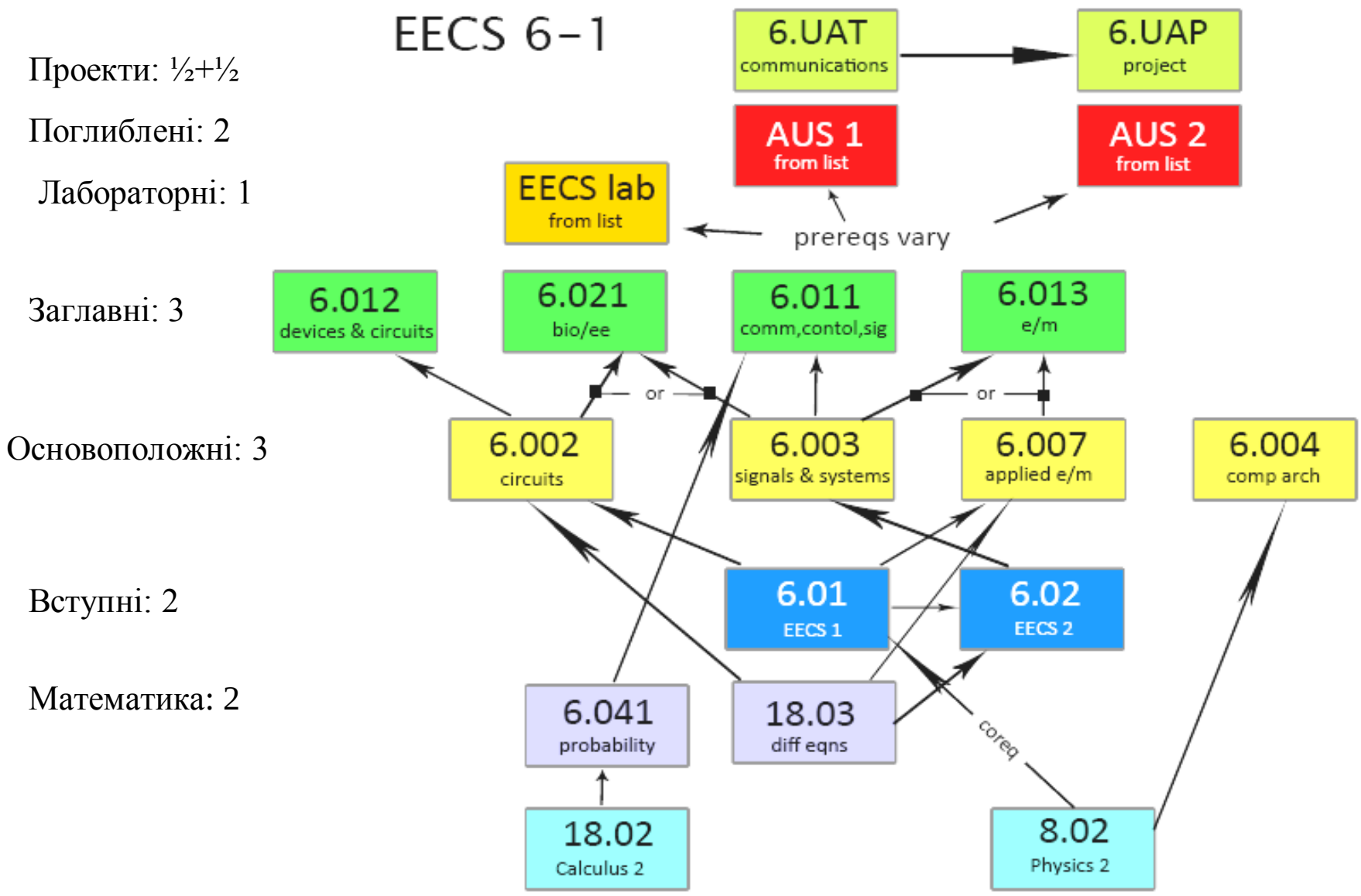

Рис. 1.3. Послідовність вивчення курсів на Курсі 6-1 [26]

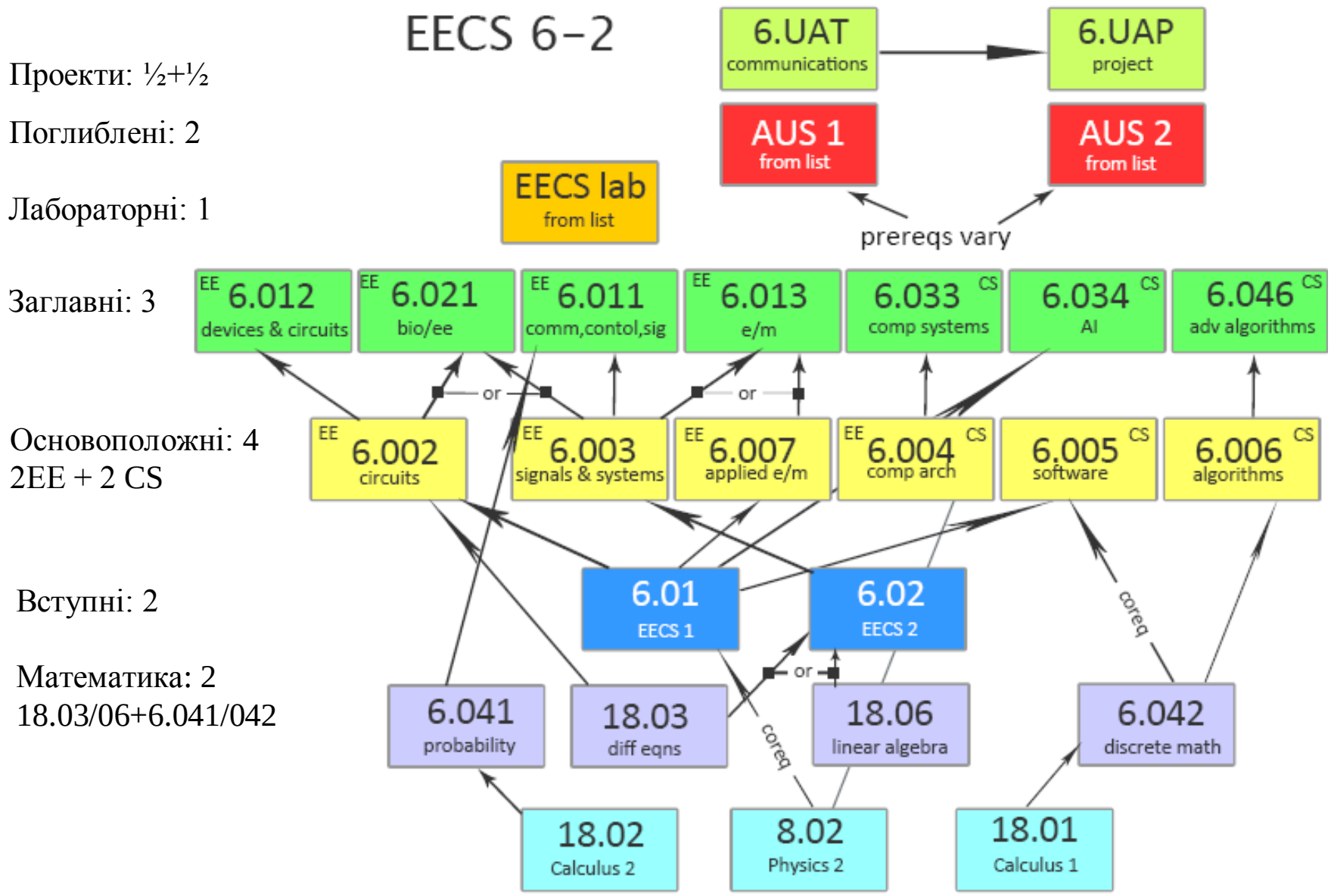

Рис. 1.4. Послідовність вивчення курсів на Курсі 6-2 [27]

**Курс 6-2**: Electrical Eng. & Computer Science. Програма 6-2 готує студентів для інженерної кар'єри та науково-дослідної роботи, де розуміння апаратного та програмного забезпечення має важливе значення. На курсі 6-2 студенти повинні пройти курси із загальноінститутських вимог, а також кафедральні курси, що відображені на рис. 1.4.

З рис. 1.4 видно, що студенту, крім курсів «Числення однієї змінної» та «Числення багатьох змінних», необхідно опанувати курси «Диференціальні рівняння» та «Лінійна алгебра» (Linear Algebra).

**Курс 6-3**: Computer Science and Engineering. Курс 6-3 зосереджується на обчислювальних структурах, штучному інтелекті, розробці програмного забезпечення, комп'ютерних алгоритмах і комп'ютерних системах. На курсі 6-3 студенти повинні пройти курси із загальноінститутських вимог, а також кафедральні курси, що відображені на рис. 1.5.

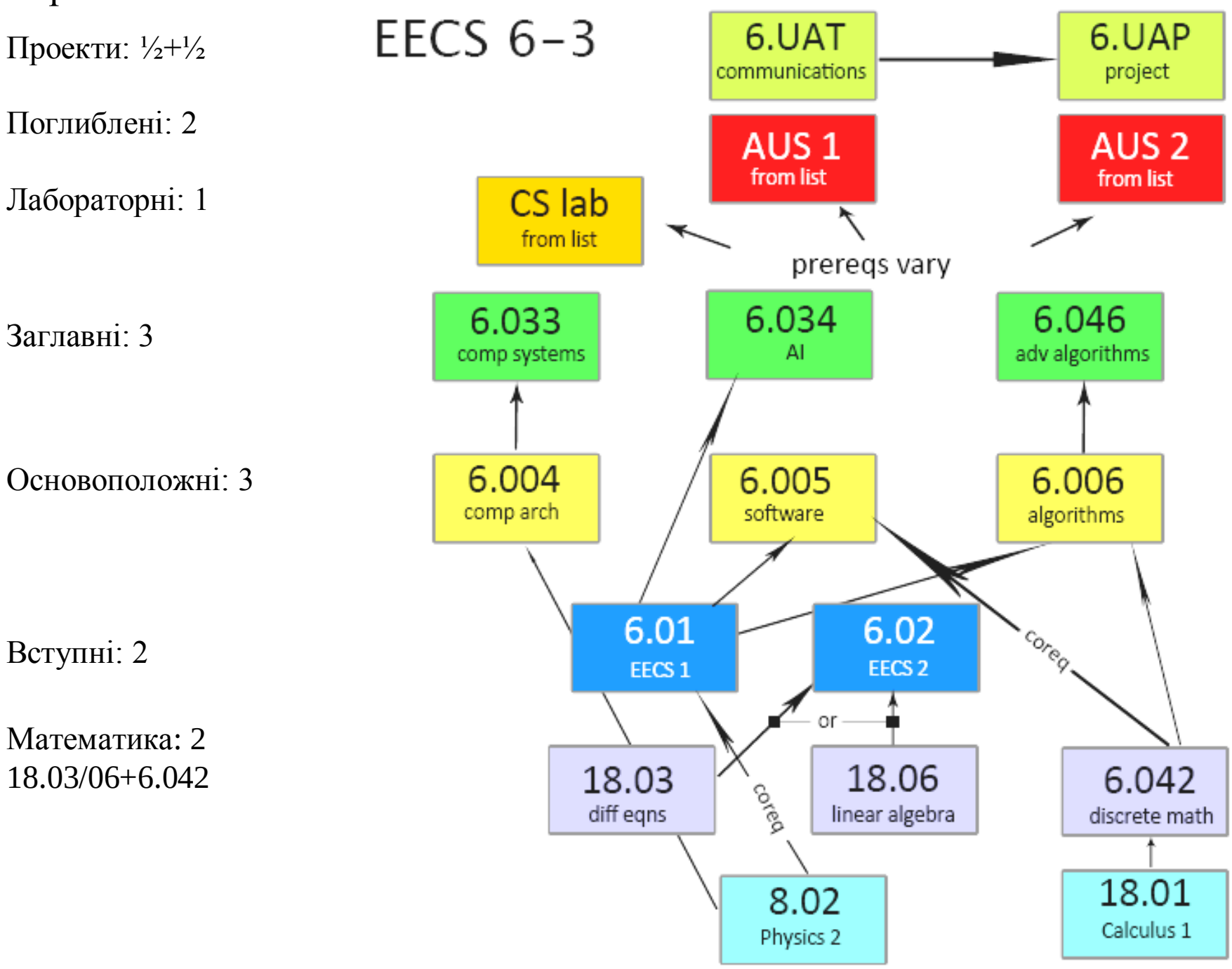

Рис. 1.5. Послідовність вивчення курсів на Курсі 6-3 [28]

Згідно з даними рис. 1.5, студенту необхідно опанувати курс «Числення однієї змінної», а також курси «Диференціальні рівняння» та «Лінійна алгебра».

**Курс 6-7**: Computer Science and Molecular Biology. Курс 6-7 забезпечує навчання як з молекулярної біології, так і з комп'ютерних наук. Навчальна програма пропонується спільно EECS та біологічним факультетом MIT, що надає можливість готувати студентів до роботи в області біології та техніки – в тому числі лікарських препаратів, фармації та обчислювальної молекулярної біології. На курсі 6-7 студенти повинні пройти курси з GIRs, а також кафедральні курси. Порядок появи дисциплін наведено на рис. 1.6.

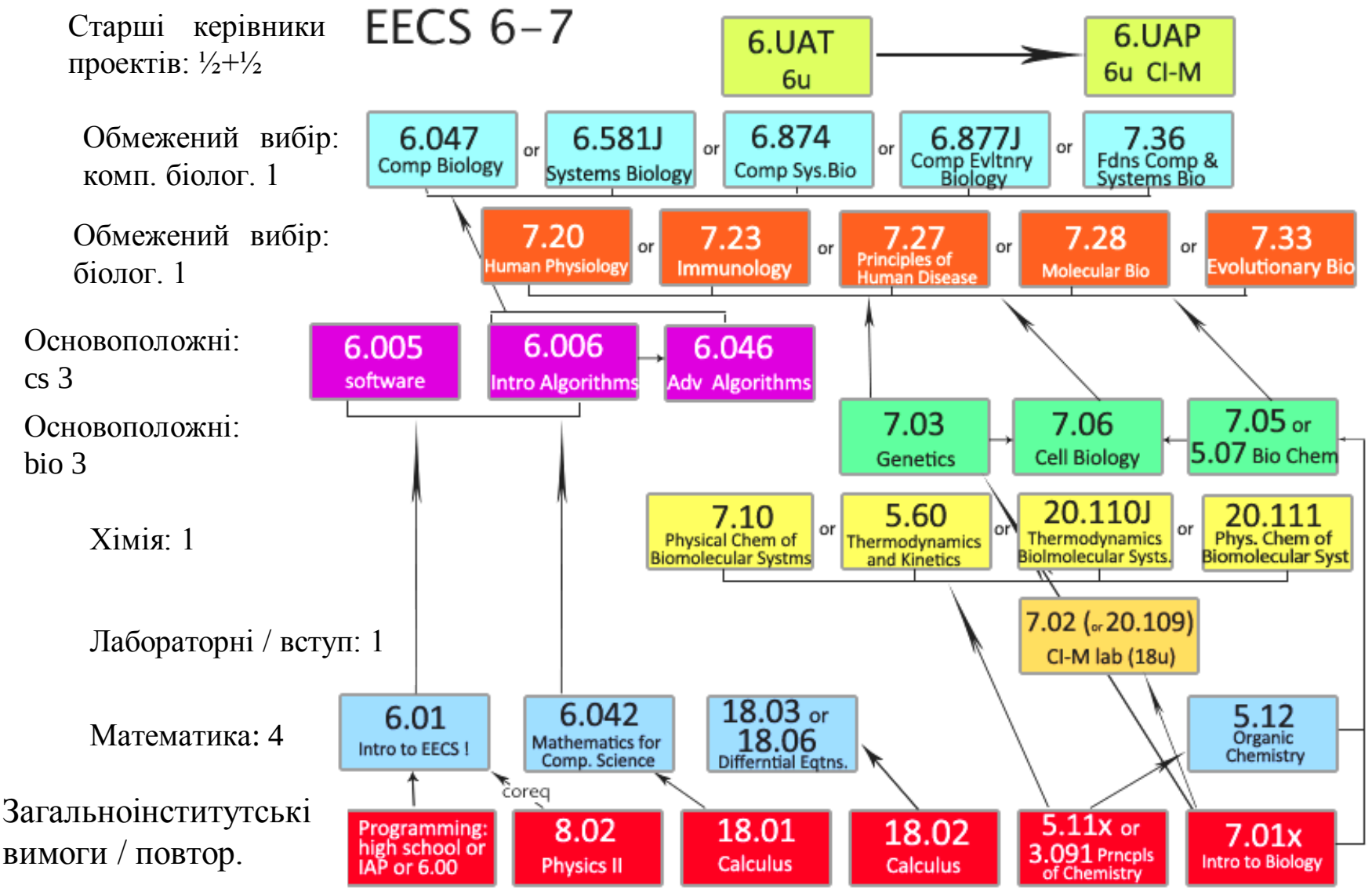

Рис. 1.6. Послідовність вивчення курсів на Курсі 6-7 [29]

Як видно з рис. 1.6, разом з «Числення однієї змінної» та «Числення багатьох змінних» студенту необхідно також опанувати курс «Диференціальні рівняння».

Порівняння розподілу кредитів на вивчення вищої математики на інженерних спеціальностях Массачусетського технологічного інституту (таблиця А.1 додатку А), а також на інженерних спеціальностях інших ВНЗ США (рис. 1.7) та проведений аналіз навчальних курсів інженерних спеціальностей у ВНЗ США надав можливість зробити такі висновки:

1. Незважаючи на недержавну форму акредитації та традиційне розмаїття пропонованих математичних курсів (як обов'язкових, так і факультативних), навчання вищої математики майбутніх інженерів у США здійснюється за схожими навчальними програмами.

2. Системи підготовки інженерів у ВНЗ США та у ВНЗ України мають такі спільні риси: а) високий рівень математизації та

комп'ютеризації загальноінженерних та спеціальних дисциплін; б) навчання фундаментальних дисциплін, зокрема вищої математики (лінійної алгебри, математичного аналізу та диференціальних рівнянь), відбувається переважно на молодших курсах, загальнопрофесійних дисциплін – на середніх курсах та спеціальних професійних – на старших; в) зміст навчання вищої математики є професійно орієнтованим та диференційованим за рівнями початкової підготовки студентів; г) у навчанні вищої математики широко використовуються засоби ІКТ.

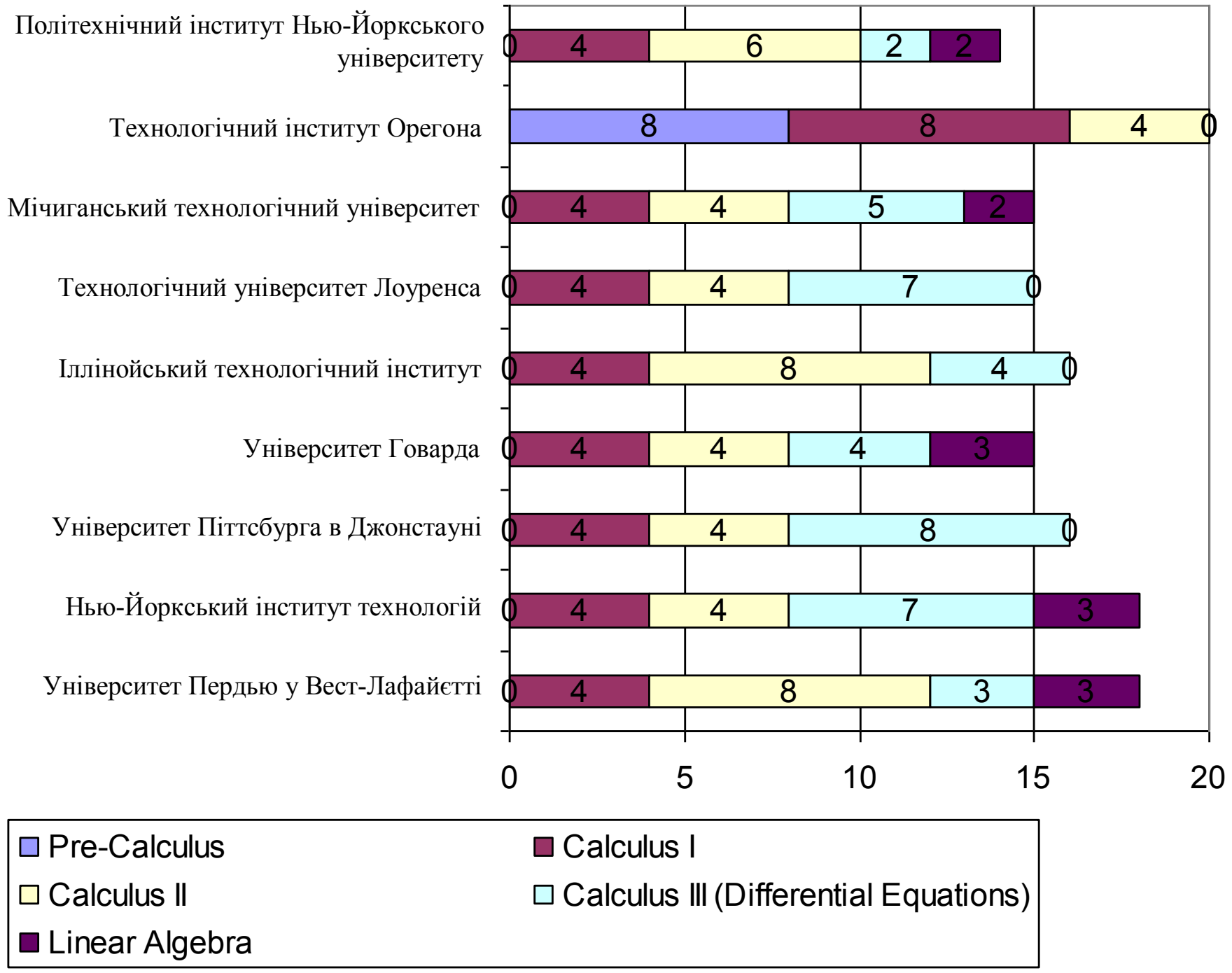

Рис. 1.7. Порівняльний аналіз кредитів на вивчення вищої математики за курсом Computer Engineering (комп'ютерна інженерія)

Аналіз програм з підготовки студентів інженерних спеціальностей у Сполучених Штатах Америки може бути корисним для визначення напрямів модернізації вітчизняної системи інженерної освіти.

*1.1.2 Національна стратегія Сполучених Штатів Америки з використання інформаційно-комунікаційних технологій у інженерній освіті*. Стрімкий розвиток і розповсюдження ІКТ набуває сьогодні глобального характеру: широко використовуються в політиці, економіці, управлінні, фінансах, науці, культурі та інших сферах людської життєдіяльності у рамках національних кордонів і світі у

цілому. ІКТ сьогодні є одним з найважливіших факторів, що впливають на формування суспільства XXI століття, так як надавши унікальні можливості у сфері пересування капіталу, товарів і послуг, інформаційні і комунікаційні технології стали основою формування нового типу економіки – «економіки знань», «інформаційної економіки», «кіберекономіки». Існуючі економічні системи, покликані пристосовуватися до інформаційної і комп'ютерної реальності.

Саме тому сьогодні в суспільстві виникають такі поняття, як «електронний уряд», «електронне громадянство», «кіберполітика», «кібердемократія», «комп'ютерно-опосередкована політична комунікація». У процесі формування суспільства поступово стираються кордони між країнами і людьми, радикально змінюється структура світової економіки, значно більш динамічним і конкурентоспроможним стає ринок. Дані і знання стають одним із стратегічних ресурсів держави, масштаби використання якого стали одним із факторів соціально-економічного розвитку. У зв'язку з цим до числа найважливіших задач кожної держави відносять формування розвитку інформаційної інфраструктури та інтеграції у глобальне інформаційне суспільство. Вирішення цих задач стає сьогодні необхідною умовою стійкого розвитку держави та її повноцінного входження у світову економіку [250].

*Інформаційно-комунікаційні технології* утворюють інформаційні технології, що базуються на використанні комп'ютерів, комп'ютерних мереж і засобів зв'язку. Інформаційно-комунікаційні технології (ІКТ, від англ. Information and communications technology – ICT) – іноді вживають як синонім до інформаційних технологій (IT), хоча ІКТ є більш загальним терміном, що підкреслює роль уніфікованих технологій та інтеграцію телекомунікацій (телефонних ліній, бездротових з'єднань), комп'ютерів, програмного забезпечення, накопичувальних та аудіовізуальних систем, що надають можливість користувачам створювати та зберігати данні, змінювати їх, передавати ці данні іншим користувачам.

Терміни «інформаційні технології» та «інформаційні та комунікаційні технології» з'являються в російськомовних джерелах з 80-х рр. XX століття (рис. 1.8).

Одним із перших видань, в якому зустрічається поняття «нові інформаційні технології», є журнал «Научно-техническая информация: Информационные процессы и системы» за 1983 рік [240].

У свою чергу, поняття «інформаційні та комунікаційні технології» вперше зустрічається в журналі «Проблемы теории и практики управления» в №1-6 за 1995 рік [264].

Значно раніше ці терміни з'являються в англомовних джерелах (рис. 1.9).

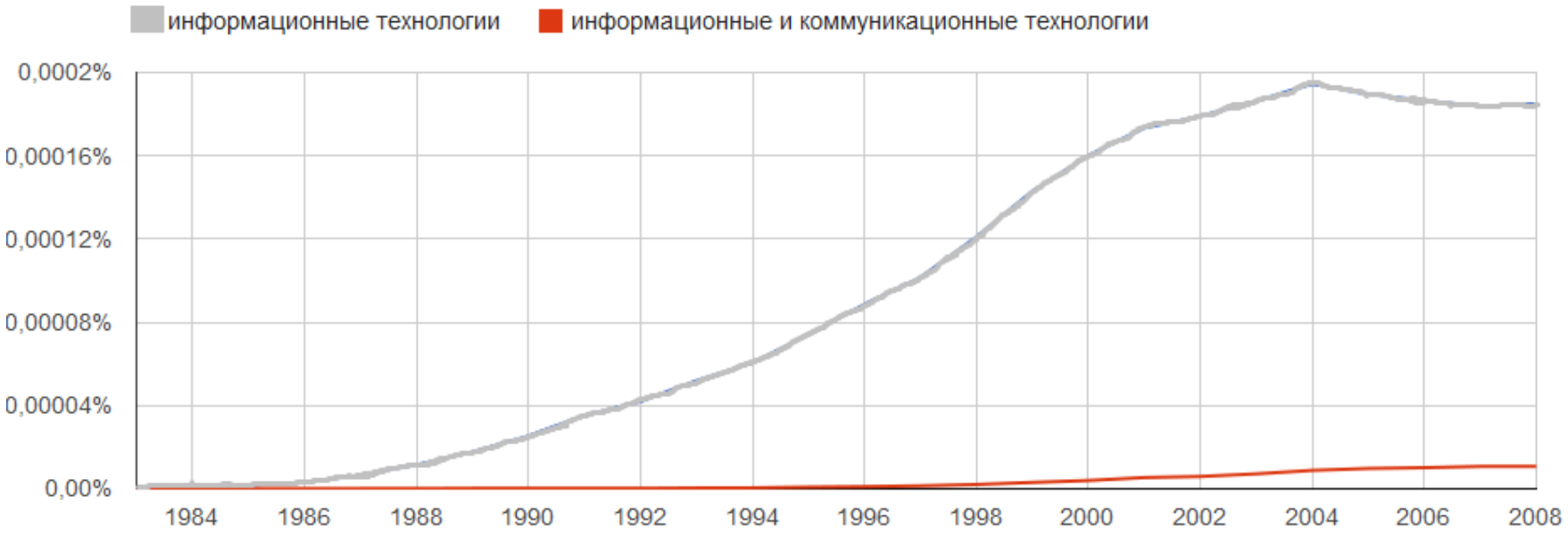

Рис. 1.8. Використання термінів «*інформаційні технології*» та «*інформаційні та комунікаційні технології*» у російськомовних джерелах

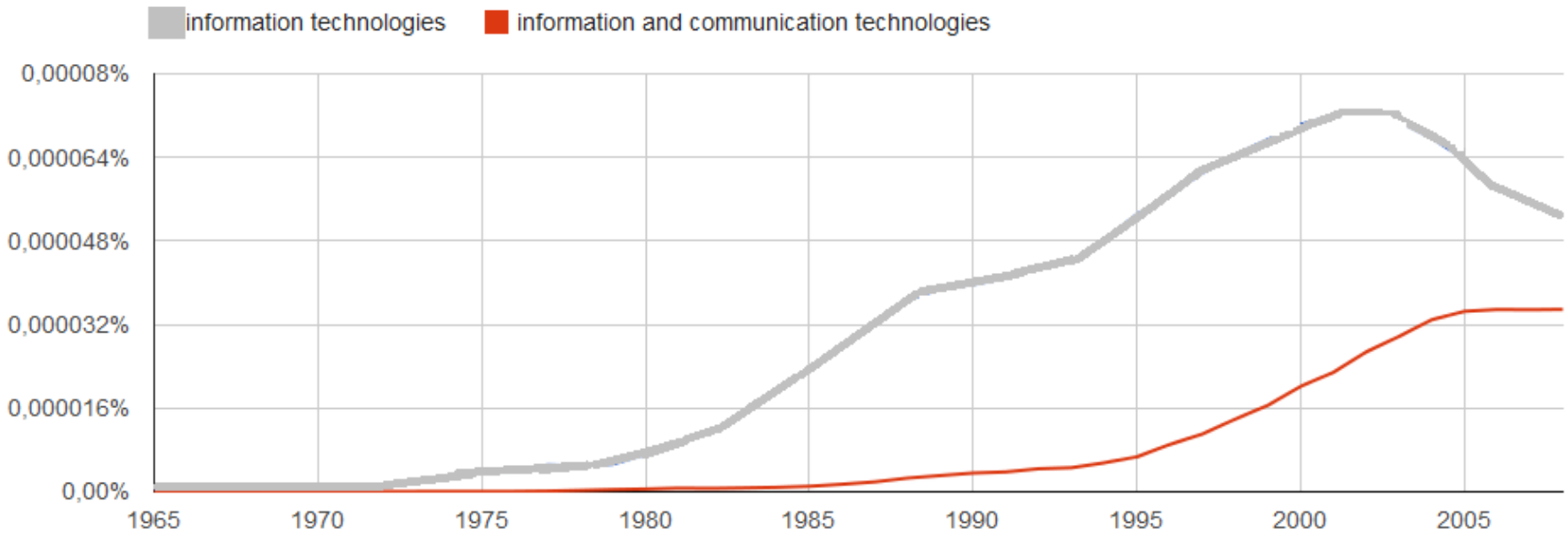

Рис. 1.9. Використання терміну «*інформаційні технології*» (information technologies) та «*інформаційні та комунікаційні технології*» ( information and communication technologies) у англомовних джерелах

Вираз «information technologies» з'являється у зарубіжній літературі в 1964 р. у статті Дж. Дайболда (John Diebold) [121].

Одним із перших закордонних видань, в якому з'являється вираз «information and communication technologies», є газета «Computerworld» за 19 вересня 1977 року [59].

У науково-педагогічній літературі зустрічається декілька тлумачень поняття «інформаційно-комунікаційні технології».

За визначенням К. Блертона (Craig Blurton) [14] (ЮНЕСКО) *інформаційно-комунікаційні технології* – це різноманітні технологічні інструменти та ресурси, що використовується для спілкування і для створення, поширення, зберігання, та управління даними.

*Інформаційно-комунікаційні технології*, з точки зору системи освіти США, включають в себе використання комп'ютера та різноманітних комунікацій, що сприяють перетворюванню навчальних відомостей у знання [62]. Використання ІКТ навчання підтримує традиційну освіту і допомагає людям впоратися з навчанням протягом всього життя, оскільки робить процес навчання легким, вмотивованим, індивідуальним та гнучким [35].

Повний аналіз категоріального поля дослідження розглянуто у додатку Б.

Стратегія США з інформатизації економіки полягає в освоєнні і використанні ІКТ для забезпечення сталого розвитку суспільства.

У промові 2009 року міністр освіти США А. Дункан (Arne Duncan) закликав лідерів освіти зосередитися на чотирьох напрямах реформи освіти [63]:

– підвищення рівня середньої освіти;

– прийняття жорстких стандартів, що готують студентів до успіху в навчанні у ВНЗ та на робочому місці;

– створення інформаційних систем для відстеження навчальних досягнень студентів та роботи викладачів;

– залучення до організації процесу навчання висококваліфікованих викладачів, особливо на заняттях, де використовуються сучасні ІКТ.

Для вирішення цих завдань Департаментом освіти США був виданий Національний план технічної освіти 2010 року. Національний план використання освітніх технологій «Трансформація американської освіти: навчання за допомогою технологій» (National Education Technology Plan, Transforming American Education: Learning Powered by Technology) вимагає застосування передових технологій, що використовуються в повсякденній роботі всією системою освіти для покращення навчання студентів, прискорення прийняття ефективних практик і використання даних та відомостей для безперервного вдосконалення [127]. У ньому представлені п’ять основних рекомендацій для держави, районів, федеральних урядів та інших зацікавлених сторін. Ці рекомендації розглядають одну з п’яти основних компонентів навчання на базі технологій: навчання, оцінювання, викладацька діяльність, засоби і продуктивність.

ІКТ, з точки зору системи освіти США, включають в себе використання комп’ютера та комунікацій, надаючи можливість перетворювати дані на відомості та знання [62].

У сучасній високо конкурентній і глобальній економіці, освіта не обмежується класами та аудиторіями. Використання ІКТ навчання сприяє традиційній освіті і допомагає людям у навчанні протягом всього життя. ІКТ-орієнтоване навчання забезпечує мотивації навчання всередині й за межами аудиторії. Ці нові системи навчання забезпечують індивідуальний підхід до потреб студента з точки зору змісту (те, чому навчається) і методу (як і коли навчається) [35].

*1.1.3 Ресурси освітніх мереж.* Ресурси Інтернет освітніх мереж та взагалі ІКТ є тими інструментами, що не тільки допомагають студентам отримати знання з різних предметів, але й сформувати в них творчий підхід до процесу навчання, використовуючи при цьому інноваційні технології та методи. У цьому їм допомагають існуючи освітні портали і мережі: BECTA, Globalschoolnet, Educared, Happychild, Teachers Network та інші. Створені на основі Інтернет, вони мають свої особисті напрями розвитку та розвивальні методики навчання з різних предметів [230; 231].

Із самого початку створення мережі Інтернет її ресурси були використані у процесі навчання саме предметів гуманітарного циклу, до яких відносять такі дисципліни як антропологія, філософія, історія, культурологія, філологія, педагогіка, мистецтвознавство, літературознавство, психологія, журналістика, етика, естетика, археологія, графологія, лінгвістика, соціологія, політологія, юриспруденція, економіка, право, етнографія, наукознавство, криміналістика [231, 7].

Одними з перших були розроблені ресурси з вивчення іноземних мов (переважно англійської), історії та права.

Такий стан є закономірним результатом історії створення мережі Інтернет, «батьками» якої є американські вчені – Б. М. Лейнер (Barry M. Leiner), В. Дж. Серф (Vinton G. Cerf), Д. Д. Кларк (David D. Clark), Р. Е. Кан (Robert E. Kahn), Л. Клейнрок (Leonard Kleinrock), Д. С. Лінч (Daniel C. Lynch), Дж. Постел (Jon Postel), Л. Дж. Робертс (Larry G.Roberts), С. Вулф (Stephen Wolff) [17], які свої розробки і програми створювали і створюють англійською мовою. Тому, всім бажаючим працювати або спілкуватися через мережу було необхідно володіти англійською мовою. Починаючи з найпростіших програм та сервісів того часу, що були задіяні (електронна пошта, BBS, новини, чати тощо) [231].

Створення сервісів Веб 2.0 пов’язують з іменем С. Пейперта, який ще у 1995 році сформулював у своїх роботах новий підхід щодо нових напрямів педагогічних інновацій, стверджуючи, що створення студентами своїх особистих інтелектуальних структур, що допомагають їм не тільки шукати та отримувати навчальні відомості, але й аналізувати їх, креативно та творчо підходити до них з метою отримання знань, значно ефективніше ніж просто отримувати знання від викладача. Особливо це важливо, коли студент безпосередньо зацікавлений в отриманих знаннях та самостійно опрацьовує навчальні відомості, постійно обговорює та обмінюється думками з іншими [97]. Більш ефективним стає навчання, коли студенти займаються справою, що є особисто-мотивованою. У цьому випадку дуже важливими є загальні принципи відкритості, конструктивної діяльності і особистої відповідальності. Саме тому велика увага на той час приділяється розвитку сервісів Веб 2.0, так званих соціальних сервісів мережі Інтернет, що допомагають створити свій особистий навчальний простір, свою особисту систему навчання.

*Веб* 2.0 (англ. Web 2.0) – поняття, яким користуються для позначення ряду технологій та послуг всесвітньої павутини, відомої також як Веб (англ. WWW – World Wide Web). Окрім цього цим поняттям описують зміну сприйняття Інтернету користувачами [246].

До україномовних проектів Веб 2.0 відноситься ITEACH WIKI – майданчик для спілкування учасників програми Intel «Навчання для майбутнього» створеного для спільного планування, обміну досвідом, ведення мережних проектів учасниками програми [231].

Розробка різних програм та розвиток сервісів Веб 2.0 надали можливість, завдяки Інтернет, не тільки розміщувати різні відомості й спілкуватися, але й створювати та працювати у спільних проектах з різних предметів, розвивати та наповнювати сайти, портали, освітні мережі розробками навчальних матеріалів до занять та методичними рекомендаціями з викладання різних предметів.

Можна зазначити, що на цей час практично всі відомі освітні мережі використовують такі сервіси, що задіяні для вивчення та навчання різних предметів і створення віртуальних товариств.

Загальний план розвитку системи освіти Великобританії «Використання технологій: навчання наступного покоління, 2008-14 роки» [54], розміщений на сайті освітнього порталу *BECTA (British Educational Communications and Technology Agency)* визначає необхідність створення нових інформаційних технологій, які б відповідали основним цілям та потребам навчання й викладання в школах. Найбільш відомий освітній портал BECTA сприяє інтеграції ІКТ у процес навчання різних шкільних предметів, що сприяє покращенню освітніх стандартів та створенню баз даних з прикладами використання ІКТ. На сторінках BECTA надаються рекомендації, пропозиції, можливі плани занять з таких предметів як: мистецтво та дизайн, англійська мова, географія, громадянська освіта, сучасні мови, математика, фізична освіта, музика, історія, ІКТ, бізнес-освіта, релігійна освіта, що уможливлюється завдяки використанню сервісів Веб 2.0.

Одна з найбільш поширених глобальних освітніх мереж *Globalschoolnet* (globalschoolnet.org – США, Сан-Дієго) створює в мережі сервіси, що спрямовані на розвиток комунікативних навичок та креативного мислення як учнів, так і вчителів, використовуючи при цьому новітнє програмне забезпечення, інноваційні комп'ютерні програми та підходи. Створення віртуальних середовищ для навчання завдяки проведенню віртуальних проектів стає повсякденною роботою мережі. Технології Веб 2.0 надають можливість створювати віртуальні навчальні середовища, що допомагають більш тісному спілкуванню та обміну думками між учнями, студентами, викладачами, школами, ВНЗ. Розвитку мережі, що сприяє покращенню організації процесу навчання, надається підтримка урядовими органами систем освіти різних країн [49].

Іспаномовна глобальна освітня мережа *Educared* (educared.org) містить навчальні ресурси для дітей, починаючи з дошкільного віку, включаючи учнів, педагогів, батьків та науковців, охоплюючи школярів різних вікових категорій загальної середньої школи. Ця освітня мережа надає можливість вивчати іспанську мову, географію, історію, літературу, природничі науки, освоювати комп'ютерні програми. За допомогою інформаційних технологій створюються різноманітні розвивальні та навчальні програми-ігри.

Основними напрями роботи мережі є [46]:

– сприяння інтеграції ІКТ у навчання та викладання різних предметів (переважно гуманітарних);

– розміщення методичних рекомендацій;

– консультації для викладачів, студентів, учнів та їх батьків;

– сприяння інноваційним процесам та розповсюдженню педагогічних інновацій у навчальному процесі;

– онлайн навчання з опанування ІКТ.

У мережі розміщені матеріали з таких предметів як: іспанська мова та діалекти, іноземні мови, географія, музика, біологія, мистецтво, математика, інформаційні технології, історія, історія мистецтва, природничі науки тощо.

Глобальна освітня мережа *Happychild* (happychild.org – Летерхед, Великобританія), що розміщує на своїх сторінках навчальні онлайн ресурси для дітей початкової та середньої школи, особливу увагу приділяє навчанню письма та читання дітей дошкільного віку початкової школи, вивченню іноземних мов (особливо англійської мови), освоєнню матеріалу з деяких курсів історії та географії з використанням сервісів Веб 2.0 [53].

Створення віртуальних спільнот поступово формує єдиний освітній простір у глобальному вимірі. Такими прикладами можуть слугувати такі освітні мережі:

1) *мережа вчителів Teachers network* (teachers.net – *США, Каліфорнія)*, що охоплює всі рівні освіти починаючи з початкової (середню та вищу) та допомагає в організації процесу навчання вчителям та викладачам. На сайтах розміщені конспекти занять, рекомендації щодо проведення занять і покращення процесу навчання за допомогою онлайн ресурсів [122];

2) *мережа творчих вчителів (Росія)* створена для педагогів, які зацікавлені у можливостях підвищення якості навчання за допомогою використання ІКТ. На порталі мережі (www.it-n.ru) розміщені різноманітні матеріали та ресурси, що стосуються використання ІКТ у процесі навчання, створюються можливості для спілкування колег. До таких ресурсів слід віднести:

– бібліотеку готових навчальних проектів з використанням ІКТ, а також різні проектні ідеї, за допомогою яких можна створити свій власний проект;

– бібліотеку методичних рекомендацій щодо проведення занять із використанням різних електронних ресурсів;

– посібники та корисні поради з використання програмного забезпечення в процесі навчання;

– посилання на аналітичні та тематичні статті педагогів;

– можливість взяти участь у роботі викладацьких спільнот та творчих груп, обговорити важливі питання на форумі.

Орієнтація України на європейські та світові освітні стандарти, що вимагають підготовки фахівців високого рівня, здатних у конкурентних умовах до прийняття самостійних нестандартних рішень, також спонукає вітчизняних педагогів до науково-методичного та технологічного пошуку на шляху вдосконалення, насамперед, самостійної роботи в навчанні.

Відповідність технологій самостійного навчання студентів таким критеріям, як [228]: відтворюваності, надійності, оптимальності, диференційованості на основі індивідуально-психологічних властивостей студентів сприятиме організації продуктивного самонавчання для студентів з різним рівнем пізнавальної активності. Такі варіативні технології навчання важливо доповнювати психолого-педагогічними технологіями, окремими методами й прийомами вдосконалення різних видів пам'яті, мислення, уваги, розвитку різнопланових здібностей, загальноосвітніх навичок. Це створить можливості для подолання таких суттєвих недоліків технологічного підходу, як ігнорування особистісних якостей студентів, орієнтація на навчання репродуктивного типу, несформованість мотивації до навчання тощо.

Уведення технологій самонавчання студентів до організації процесу навчання буде сприяти:

– поглибленню та закріпленню знань, умінь і навичок;

– формуванню навичок роботи з технологічним інструментарієм;

– розвитку технологічного мислення, умінь самостійно планувати,

алгоритмізувати, стандартизувати своє учіння;

– формуванню спрямованості та навичок раціонально організовувати самонавчання, чітко дотримуватися вимог технологічної дисципліни в різних видах діяльності [229, 168].

Таким чином, зростання ролі ІКТ у багатьох видах людської діяльності цілком природно спричинює зміни в системі освіти, спрямовані на переорієнтацію навчально-виховного процесу. Використання ІКТ у процесі навчання може забезпечити передачу знань і доступ до різноманітних навчальних відомостей нарівні, а іноді й інтенсивніше та ефективніше, ніж за традиційного навчання.

*1.1.4 Інформаційно-комунікаційні технології навчання.* Інформаційні технології зазнали значної еволюції протягом останніх двох десятиліть. Це, в свою чергу, привело до змін у системі вищої освіти. Зокрема, з метою вивчення шляхів використання технологічних інновацій та поліпшення освітніх програм все більше університетів пропонує Інтернет-курси і програми, гібридні або змішані (комбіновані) курси.

Нововведення, такі як віртуальні навчальні середовища (virtual learning environments, VLE), забезпечують розвиток навчання за допомогою Інтернет і мультимедіа, дистанційного та мобільного навчання, надаючи студентам можливості доступу до курсів та освітніх програм у зручний для них час та у зручному місці, створюючи умови для індивідуалізації та персоналізації процесу навчання [75].

У всій системі вищої освіти, навчання математичних дисциплін також слідує чітким тенденція до збільшення використання таких Інтернет-ресурсів, як:

– віртуальні навчальні середовища VLE – Moodle, Sakai або Blackboard;

– математичне програмне забезпечення спеціального призначення – наприклад, Mathematica, Maple, Statistica, SPSS або R.

Використання зазначених технологій сприяє і надає можливість досліджувати і розробляти нові підходи до організації процесу навчання, що, в свою чергу, приводить до розвитку нової стратегії та методології навчання вищої математики у вищій технічній школі [75].

*1.1.4.1 Дистанційне навчання* (ДН) В. Ю. Биков [169, 191] визначає як форму організації і реалізації навчально-виховного процесу, за якою його учасники здійснюють навчальну взаємодію принципово й переважно екстериторіально (на відстані, що не надає можливості і не передбачає безпосередню навчальну взаємодію учасників віч-на-віч, коли учасники територіально перебувають за межами можливої безпосередньої навчальної взаємодії і коли у процесі навчання їх особиста присутність у певних навчальних приміщеннях навчального закладу не є обов’язковою).

На сучасному етапі розвитку освітянського процесу стає все більш актуальним проектування, розвиток та впровадження новітніх інформаційних технологій у дистанційному навчанні. За останнє десятиріччя дистанційне навчання закріпило свої позиції у корпоративному, післядипломному навчанні, і поступово поширюється у вищій освіті, доповнюючи традиційні форми організації навчання і закладаючи фундамент для розвитку комбінованого (змішаного) навчання.

Головною метою дистанційного навчання є надання однакових освітніх можливостей населенню за допомогою інформаційних і телекомунікаційних засобів, а також підвищення якісного рівня освіти за рахунок більш активного використання наукового та освітнього потенціалу провідних університетів, академій, інститутів, наукових центрів та інших освітніх установ [308].

Протягом останнього десятиріччя дистанційне навчання стало одним з найважливіших елементів системи вищої освіти США. З одного боку, це зумовлено бурхливим розвитком інформаційних технологій, з іншого – політикою уряду США та інших розвинених країн у галузі освіти, їхнім прагненням зробити навчання та освіту будь-якого рівня максимально доступними для всіх верств населення [32].

За В. Ю. Биковим [146], *дистанційна освіта* – це «різновид освітньої системи, в якій використовуються переважно дистанційні технології навчання та організації освітнього процесу, або одна з форм отримання освіти, за якою опанування тим чи іншим її рівнем за тою чи іншою спеціальністю (напрямом підготовки, перепідготовки або підвищення кваліфікації) здійснюється в процесі дистанційного навчання».

Серед факторів, що спричинили появу наприкінці XX ст. сучасних форм дистанційної освіти, можна вказати такі: глобалізація; підвищення динаміки соціально-економічного розвитку суспільства; поява нових потреб у тих, хто навчається; розвиток ІКТ та їх всебічне впровадження практично у всі сфери життєдіяльності людини; необхідність широкого застосовування в освітній практиці як засобу навчання і предмета вивчення.

У наукових працях зустрічається велика кількість тлумачень поняття «дистанційне навчання», що характеризує різні підходи до його розуміння.

П. В. Стефаненко з дидактичної точки зору вирізняє два підходи до цих тлумачень [308]:

– у *першому* під дистанційним навчанням розуміють обмін

навчальними відомостями між педагогами і тими, хто навчається за допомогою електронних мереж чи інших засобів телекомунікацій. Студент при цьому одержує навчальні матеріали і завдання щодо їх засвоєння, а потім результати своєї самостійної роботи надсилає педагогу, який оцінює якість та рівень засвоєння матеріалу. При цьому особиста діяльність студента із здобуття знань майже не організовується;

– домінантою дистанційного навчання у *другому* підході виступає особистісна продуктивна діяльність того, хто навчається, яка вибудовується за допомогою сучасних засобів телекомунікацій. Цей підхід передбачає інтеграцію інформаційних і педагогічних технологій, що забезпечують інтерактивність взаємодії педагога і студента, а також продуктивність навчального процесу. Навчання в даному разі здійснюється в реальному часі (чат, відеозв'язок тощо), а також асинхронно (за допомогою електронної пошти). Особистісний і телекомунікаційний характер навчання – основні ознаки дистанційного навчання цього типу.

Аналіз наукових джерел [161; 169; 314; 323] надав можливість зробити висновок, що *дистанційне навчання* – це форма навчання за якої взаємодія викладача та студента відбувається на відстані, але при цьому відображає всі властиві процесу навчання компоненти (мету, зміст, методи, організаційні форми, засоби навчання) і реалізується за допомогою технологій, що забезпечюють інтерактивність процесу навчання.

Залежно від характеру організації навчальної комунікацій між учасниками навчально-виховного процесу та організаторами освіти і способу побудови комунікаційного каналу навчального середовища В. Ю. Биков розрізняє традиційне дистанційне навчання та е-ДН [169, 191].

*Традиційне дистанційне навчання* – різновид дистанційного навчання, за яким учасники і організатори навчального процесу здійснюють взаємодію переважно асинхронно у часі, значною мірою використовуючи як транспортну систему доставки навчальних матеріалів пошту, телефонний або телеграфний зв'язок [169].

*Електронне дистанційне навчання* надає можливість підтримувати індивідуалізовану взаємодію між учасниками процесу навчання як асинхронно так і синхронно, використовуючи електронні мережні засоби комунікації [146].

В. М. Кухаренко [308] вказує, що в межах дистанційного навчання реалізуються всі існуючі традиційні *дидактичні принципи* навчання та з'являються нові: принцип педагогічної доцільності використання нових інформаційних технологій; принцип забезпечення безпеки даних, що циркулює при дистанційному навчанні; принцип відповідності технологій навчання; принцип мобільності навчання.

Дистанційне навчання має низку переваг. До *основних переваг* дистанційного навчання Б. А. Демида [166] відносить:

– *свобода і гнучкість* – можливість навчатися одночасно в різних місцях, на різних курсах не тільки в одному, а й у декількох університетах чи навіть країнах;

– *індивідуальність* – самостійний вибір студентами власного темпу та графіка навчання у звичній для них обстановці і в зручний час;

– отримання *освіти інвалідами та людьми з різними фізичними вадами;*

– набуття студентами таких якостей, як *самостійність, мобільність і відповідальність*;

– *навчання більшої кількості людей* різних вікових груп порівняно з іншими формами навчання.

*1.1.4.2 Навчання за допомогою Інтернет і мультимедіа.* Аналіз тлумачень «навчання за допомогою Інтернет і мультимедіа» у системі вищої освіти США надав можливість стверджувати, що це поняття не відокремлюється від поняття «дистанційне навчання», а є деяким етапом на шляху розвитку всієї дистанційної освіти США.

Наведемо декілька тлумачень поняття «навчання за допомогою Інтернет і мультимедіа».

М. Розенберг (Marc Rosenberg) дав таке тлумачення терміну e-Learning: *e-Learning* – використання Інтернет-технологій для надання широкого спектра для вирішення питань, що забезпечують підвищення знань та продуктивності праці; e-Learning базується на трьох основних принципах: робота здійснюється в мережі; доставка навчальних матеріалів студенту здійснюється за допомогою комп’ютера з використанням стандартних Інтернет-технологій [107, 10].

Е. Роззетт (Allison Rossett) визначає e-Learning так: *Web-навчання (WBT) або навчання за допомогою Інтернет і мультимедіа, або онлайн навчання* – це є підготовка кадрів, що знаходиться на сервері або на комп’ютері, який підключений до мережі Інтернет (World Wide Web) [32].

Фахівці ЮНЕСКО вважають, що *e-Learning* – це навчання за допомогою Інтернет і мультимедіа [10, 7].

Існує велика кількість тлумачень, що роблять акцент на інших аспектах e-Learning. Наведемо декілька із них:

*e-Learning* – навчання, що ґрунтується на використанні інформаційних і телекомунікаційних технологій, які надають можливість підтримувати весь процес навчання: від доставки навчальних матеріалів студентам до їх контролю рівня їх засвоєння [37].

*Електронне навчання (навчання за допомогою Інтернет і мультимедіа)* – інноваційна технологія, спрямована на професіоналізацію та підвищення мобільності тих, хто навчається, і на сучасному етапі розвитку засобів ІКТ може розглядатися як технологічна основа фундаменталізації вищої освіти [292, 109–110].

Виступаючи в якості повної заміни або як доповнення до традиційного навчання, асинхронне навчання за допомогою Інтернет і мультимедіа є, скоріше, найбільш швидко зростаючим сегментом у сфері вищої освіти США. Останні дослідження в США показують, що навчання за допомогою Інтернет і мультимедіа, виступаючи в якості повної заміни традиційного навчання, має в середньому щорічне збільшення чисельності студентів і охоплює трохи менше 20% всіх студентів у період між 2002 і 2008 роками, приблизно 300000 викладачів займаються навчанням за допомогою Інтернет і мультимедіа (у тому числі у США в 2008 році від 20 до 25 % студентів реєструвалися хоча б в одному онлайн-класі) [75, 240].

Кількість студентів, що пройшли принаймні один онлайн-курс з 1 602 970 у 2002 р. і 1 971 397 у 2003 році збільшилася до 6 700 000 у 2012 році [6].

Розвиток навчання за допомогою Інтернет і мультимедіа відбувався на трьох етапах [292, 103–105]. *Перший етап* (20–50-ті роки XX століття) охоплює період з моменту появи електромеханічних комп’ютерів до широкого впровадження електронних комп’ютерів. Цей етап характеризується застосуванням різних механічних, електромеханічних та електронних індивідуалізованих пристроїв, за допомогою яких подавався навчальний матеріал та виконувався контроль і самоконтроль знань.

*Другий етап* охоплює період 50–80-х років минулого століття та пов’язаний з широким впровадженням електронних обчислювальних машин у практику. Ключовими термінами цього періоду стали інтелектуальні навчаючі системи, комп’ютерно-орієнтовані системи навчання, комп’ютерна підтримка навчального процесу, комп’ютерні системи контролю знань. У цей період була створена велика кількість спеціалізованого програмного забезпечення – автоматизованих навчальних систем PLATO, Coursewriter, Tutor та інші. Цьому сприяли очевидні переваги електронних комп’ютерів над електромеханічними – наявність пам’яті для зберігання навчальних матеріалів, висока швидкість опрацювання та розрахунків, більш широкі засоби для перегляду навчальних матеріалів та багато інших. Головним недоліком розробок цього періоду була їх стаціонарність та автономність, пов’язана з використанням «великих» обчислювальних машин або, в кращому випадку, зв’язаних з ними терміналів. Також було важко реалізувати обмін освітніми ресурсами та послугами між великою кількістю користувачів.

*Третій етап* (з 80-х років минулого століття) розпочався з появою комп’ютерних мереж та персональних комп’ютерів. Виключно потужний імпульс у розвитку освітніх технологій пов’язаний з використанням глобальної мережі Інтернет. Використання спільних та розподілених ресурсів, Web-технологій, віддалений доступ до навчальних матеріалів забезпечив суттєве підвищення ефективності професійної підготовки, її доступності та масовості. Мережні технології, висока якість та підвищення апаратного забезпечення уможливили створення професійних середовищ та систем для надання освітніх послуг і реалізації різних видів формальної (організованої) та неформальної (спеціально не організованої) освіти. Ключовими термінами цього періоду є Інтернет, Web-курси, гіпертекст, віртуальне навчання, віртуальний університет, неперервна освіта, навчання протягом усього життя, дистанційне навчання, навчання за допомогою Інтернет і мультимедіа та мобільне навчання.

Останнім часом все більшого поширення набуває термін e-Learning 2.0. Цей термін відображає тенденції в сфері організації навчання за допомогою Інтернет і мультимедіа, пов’язані з використанням технологій Веб 2.0. На відміну від навчання за допомогою Інтернет і мультимедіа, що передбачає використання дистанційних курсів, які пропонуються студентам з метою організації процесу навчання, e-Learning 2.0 передбачає використання засобів Веб 2.0: блогів, вікі, підкастів, соціальних мереж тощо [37].

Спираючись на зазначені характерні риси і принципи побудови навчання за допомогою Інтернет і мультимедіа, В. М. Кухаренко вказує такі його специфічні *якісні властивості* [224]:

1) *гнучкість і адаптивність* навчального процесу до потреб і можливостей студентів, які, в основному не відвідують регулярних занять, а працюють у зручний (як для викладача, так і для студента) для такої роботи час у зручному місці й зручному темпі;

2) *модульність побудови* навчальних програм;

3) *нова роль викладача*: викладач координує навчально-пізнавальний процес, коригує курс, який викладає, керує навчальними проектами, перевіряє поточні завдання, консультує при складанні індивідуального навчального плану, управляє навчальними групами взаємопідтримки;

4) *спеціалізовані форми контролю* якості навчальних досягнень: традиційні формами контролю якості освіти та дистанційні (співбесіди, практичні, курсові та проектні роботи, екстернат, робота в середовищі комп'ютерних інтелектуальних тестових систем тощо);

5) *використання спеціалізованих засобів навчання*.

У зв'язку з тим, що навчання за допомогою Інтернет і мультимедіа в останні роки набуває все більшої популярності, виникає необхідність в стандартизації підходів до створення курсів навчання за допомогою Інтернет і мультимедіа. Саме тому Міністерство оборони США та Департамент політики в галузі науки і технології Адміністрації Президента США в листопаді 1997 р. оголосили про створення ініціативи ADL (Advanced Distributed Learning – Розвиток розподіленого навчання). Метою створення даної ініціативи є розвиток стратегії, що проводиться Міністерством оборони і урядом в області модернізації навчання і тренінгу, а також для об'єднання вищих навчальних закладів та комерційних підприємств для створення стандартів у сфері навчання за допомогою Інтернет і мультимедіа [307].

Створення стандарту SCORM (Sharable Content Object Reference Model – «Зразкова модель об'єкта вмісту для спільного використання») є першим кроком на шляху до розвитку концепції ADL, оскільки даний стандарт визначає структуру навчальних матеріалів і інтерфейс середовища виконання. Завдяки цьому навчальні об'єкти можуть бути використані в різних системах електронної дистанційної освіти. SCORM описує структуру такої освіти за допомогою декількох основних принципів, специфікацій і стандартів, ґрунтуючись при цьому на інших, раніше створених, специфікаціях і стандартах електронної та дистанційної освіти.

У процесі роботи над SCORM були сформульовані кілька вимог до всіх систем, що розробляються відповідно до цього стандарту:

– *доступність*: здатність визначати місцезнаходження та отримати доступ до навчальних компонентів з точки віддаленого доступу і поставити їх багатьом іншим точкам віддаленого доступу;

– *адаптованість*: здатність адаптувати навчальну програму відповідно до індивідуальних потреб студента, або потребам організацій;

– *ефективність*: здатність збільшувати ефективність і продуктивність, скорочуючи час і витрати на доставку відомостей;

– *довговічність*: здатність відповідати новим технологіям без додаткового та дорогого доопрацювання;

– *інтероперабельність*: здатність використовувати навчальні матеріали незалежно від платформи, на якій вони створені;

– *можливість багаторазового використання*: здатність використовувати матеріали в різних додатках і контекстах [307].

Усі ці принципи успішно можуть бути дотримані в тому випадку, якщо спочатку орієнтуватися на використання освітнього контенту у

Web-середовищі [307].

Для телекомунікаційного середовища (зокрема, мережі Інтернет) характерна клієнт-серверна модель. Така модель використовується і в стандарті SCORM. Сервером у даному випадку є LMS (Learning Management System – система управління навчанням).

Термін «LMS», який використовується в SCORM, позначає набір функціональних можливостей, розроблених для розповсюдження, контролю та управління освітнім контентом і навчальним процесом. Цей термін відноситься як до простих систем управління, так і до складних організаційних систем [108].

У контексті SCORM широко використовуються LMS програми. SCORM зосереджується на інтерфейсі, що використовується в освітньому контенті та LMS, але не стосується особливостей внутрішньої реалізації LMS. У SCORM термін LMS означає середовище сервера. Іншими словами, згідно специфікації стандарту SCORM, LMS визначає які дані і куди треба надати, і відстежує роботу користувача з матеріалом.

Опишемо структуру стандарту [307]:

**SCORM Sequencing & Navigation** (SN) описує правила і методи для здійснення впорядкування навчального матеріалу (курсів, окремих уроків, завдань тощо). Під упорядкуванням у даному випадку розуміється розташування окремих частин курсу при проходженні курсу навчання. В SN реалізується процес переходу між різними частинами курсу – за якими правилами, якими діями може бути ініційований перехід, які форми переходу дозволені, а які ні.

**SCORM Content Aggregation Model** (CAM) описує компоненти, що використовуються в освітніх системах, відповідних стандартам SCORM, способи обміну цими компонентами та їх описи для пошуку і запуску, і правила впорядкування компонентів. CAM описує яким чином здійснювати зберігання змісту курсів, його маркування, обмін і розкриття змісту. SCORM CAM також визначає вимоги до створення змісту (наприклад, курсів, уроків, модулів тощо). CAM містить відомості про створення пакетів змісту, застосування методичних даних до компонентів у середині пакета і застосуванні правил впорядкування і навігації всередині певного пакета. CAM взаємопов’язана зі SCORM Run-Time Environment.

**SCORM Run-Time Environment** (RTE) описує вимоги до LMS у частині управлінням часу виконання (тобто, процесом запуску та обміном даних). У ній описується, які методи і властивості повинна підтримувати LMS і оточення виконання. Ці методи можуть використовуватися як для обміну даними між LMS та освітніми об’єктами, так і для управління процесом навчання.

При використанні моделі SCORM можливе створення електронних курсів, незалежних від самої системи, що легко переносяться і багаторазово використовуються в інших системах управління навчанням [108].

На теперішній час навчання за допомогою Інтернет і мультимедіа є невід'ємною частиною освітнього процесу. До його складу можна віднести електронні курси, електронні бібліотеки, нові програми та системи навчання.

*Електронний навчальний курс* (дистанційний курс) – комплекс навчально-методичних матеріалів та освітніх послуг, створених для організації індивідуального та групового навчання з використанням технологій дистанційного навчання під керівництвом викладача, який реалізується засобами Інтернет-технологій, відео конференцій, телебачення, інших засобів і вимагає активного спілкування викладачів зі студентами, студентів між собою, у якому навчальний матеріал подається у структурованому електронному вигляді та зберігається на спеціальному навчальному порталі [233].

Елементи навчання за допомогою Інтернет і мультимедіа, його переваги, недоліки та основні проблеми наведені у додатку В.

Сьогодні навчання за допомогою Інтернет і мультимедіа в Україні може повноцінно розвиватися при наявності нормативно-правової бази, навчальних закладів навчання за допомогою Інтернет і мультимедіа, контингенту студентів, кваліфікованих викладачів, навчальних програм і курсів, відповідної матеріально-технічної бази, фінансової підтримки тощо.

Дані про стан навчання за допомогою Інтернет і мультимедіа в нашій країні та в усьому світі свідчать про нагальну необхідність його стимулювання, щоб забезпечити динамічний і прогресивний розвиток та впровадження на всіх рівнях освіти, перш за все, – вищої [294, 84].

*1.1.4.3 Мобільне навчання.* Характерною рисою останнього десятиріччя стало активне використання засобів мобільного зв'язку та різноманітних апаратних мобільних пристроїв. Сучасний мобільний телефон має функціональність, що не поступається комп'ютерам початкового рівня, а в деяких випадках і середньої потужності. В першу чергу це стосується смартфонів та персональних комунікаторів (КПК із засобами зв'язку). Поширеність серед користувачів мобільного зв'язку смартфонів та персональних комунікаторів, на думку фахівців, складає близько 10% і має чітку тенденцію до зростання [187, 291; 114].

*Мобільне навчання* (mobile learning, M-Learning) – сучасний напрямок розвитку систем дистанційної освіти із застосуванням мобільних телефонів, смартфонів, КПК, електронних книжок [294].

*Мобільне навчання* – це технологія навчання, що базується на інтенсивному застосуванні сучасних мобільних засобів та технологій. Мобільне навчання тісно пов'язане з навчальною мобільністю в тому сенсі, що студенти повинні мати можливість брати участь в освітніх заходах без обмежень у часі та просторі. Використання мобільних технологій відкриває нові можливості для навчання, особливо для тих, хто живе ізольовано або у віддалених місцях чи стикається з труднощами в навчанні. Можливість навчання будь-де та будь-коли, що властиво мобільному навчанню, сьогодні є загальною тенденцією інтенсифікації життя в інформаційному суспільстві [294, 79].

В. О. Куклєв [221] розглядає *мобільне навчання* як навчання за допомогою Інтернет і мультимедіа за допомогою мобільних засобів, незалежно від часу та місця, з використанням спеціального програмного забезпечення на педагогічній основі міждисциплінарного та модульного підходів.

К. Л. Бугайчук виділяє такі апаратні пристрої для мобільного навчання [149]:

– *телефони*: звичайні мобільні телефони, смартфони, комунікатори;

– *портативні комп'ютери*: ноутбуки, нетбуки, Інтернет-планшети;

– *пристрої зберігання і відтворення даних*: електронні «рідери» (Pocket Book, Amazon Kindle), MP3/MP4 плеєри.

Існуючі мобільні пристрої суттєво відрізняються за своєю функціональністю, розмірами та ціною як між собою, так і в порівнянні зі звичайними технічними пристроями навчання за допомогою Інтернет і мультимедіа (стандартний ПК та периферійні пристрої). Основними рисами, що їх об'єднують, є мобільність та придатність для бездротового з'єднання [294].

До особливостей мобільного навчання М. Шарплес (M. Sharples) відносить [109; 110]: спільну онлайн роботу над проектом, моблоггінг (мобільний блоггінг), персоналізоване навчання, роботу у групах, онлайнові дослідження, рівний доступ до навчання.

На відміну від дистанційного навчання, мобільне навчання є більш доступним для більшості студентів, а мобільні ІКТ навчання мають достатній потенціал за гнучкістю навчання для використання та підтримки традиційного навчання [141].

Оскільки в основі мобільного навчання лежить використання мобільних пристроїв і бездротових комунікаційних технологій, розглянемо технологічну складову мобільного навчання (рис. 1.10).

Іноді окремо виділяють *віртуальне навчання,* під яким розуміють всі форми та підходи до навчання з використанням Інтернет, що надає можливість об'єднати мобільне та навчання за допомогою Інтернет і мультимедіа [294].

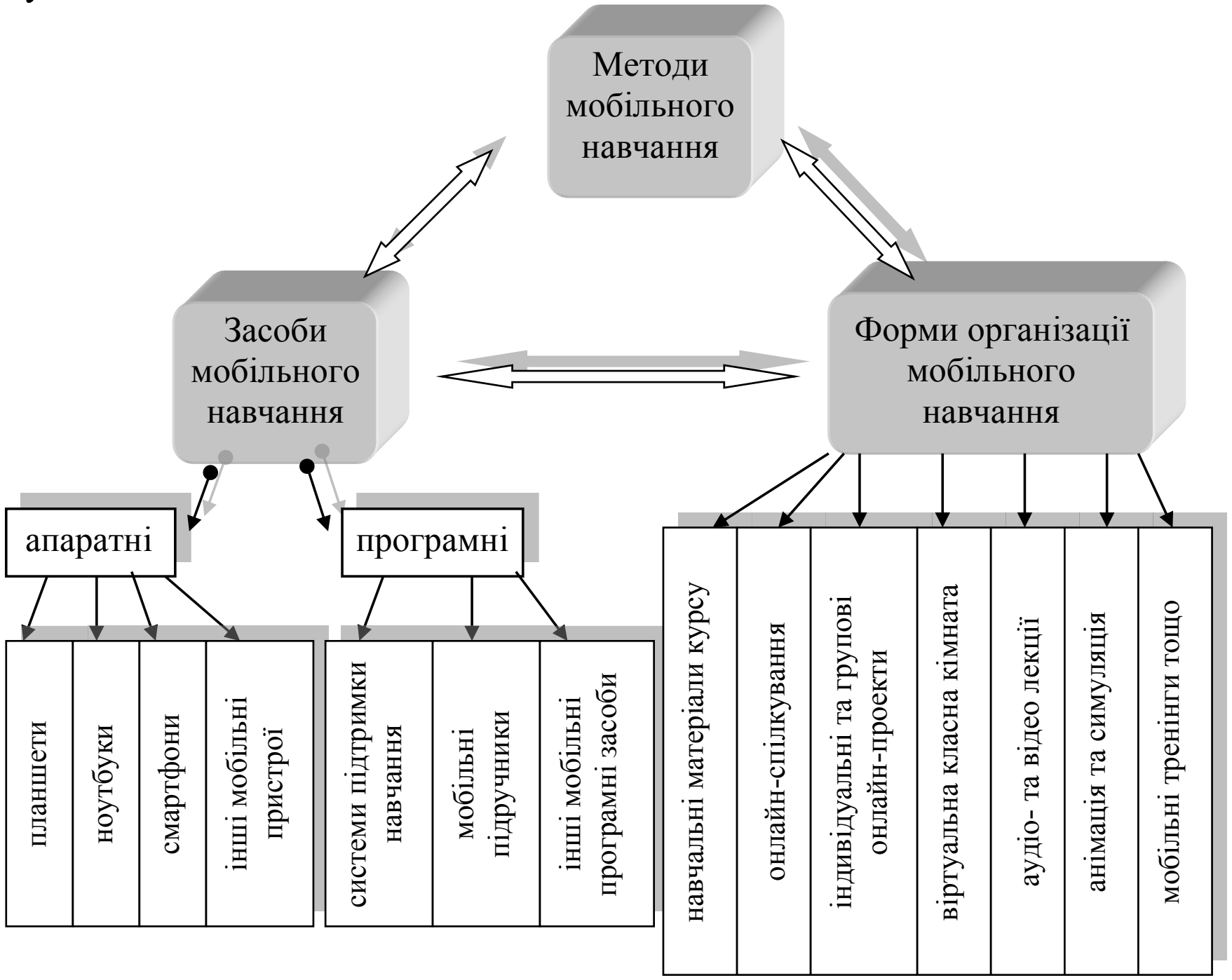

Рис. 1.10. Технологічна складова мобільного навчання

У роботах вітчизняних авторів поєднання мобільного, дистанційного, традиційного навчання та навчання за допомогою Інтернет і мультимедіа розглядається як комбіноване (змішане) навчання – процес навчання, за якого традиційні технології навчання поєднуються з інноваційними технологіями дистанційного, мобільного навчання та навчання за допомогою Інтернет і мультимедіа з метою створення гармонійного поєднання теоретичної та практичної складових процесу навчання.

А. М. Стрюк [311] *комбіноване навчання* тлумачить як цілеспрямований процес здобування знань, умінь та навичок в умовах інтеграції аудиторної та позааудиторної навчальної діяльності суб’єктів освітнього процесу на основі взаємного доповнення технологій традиційного, дистанційного, мобільного навчання та навчання за допомогою Інтернет і мультимедіа.

Ю. В. Триусом [319] сутність змішаного навчання з позицій вітчизняних авторів подано у схематичному вигляді (рис. 1.11).

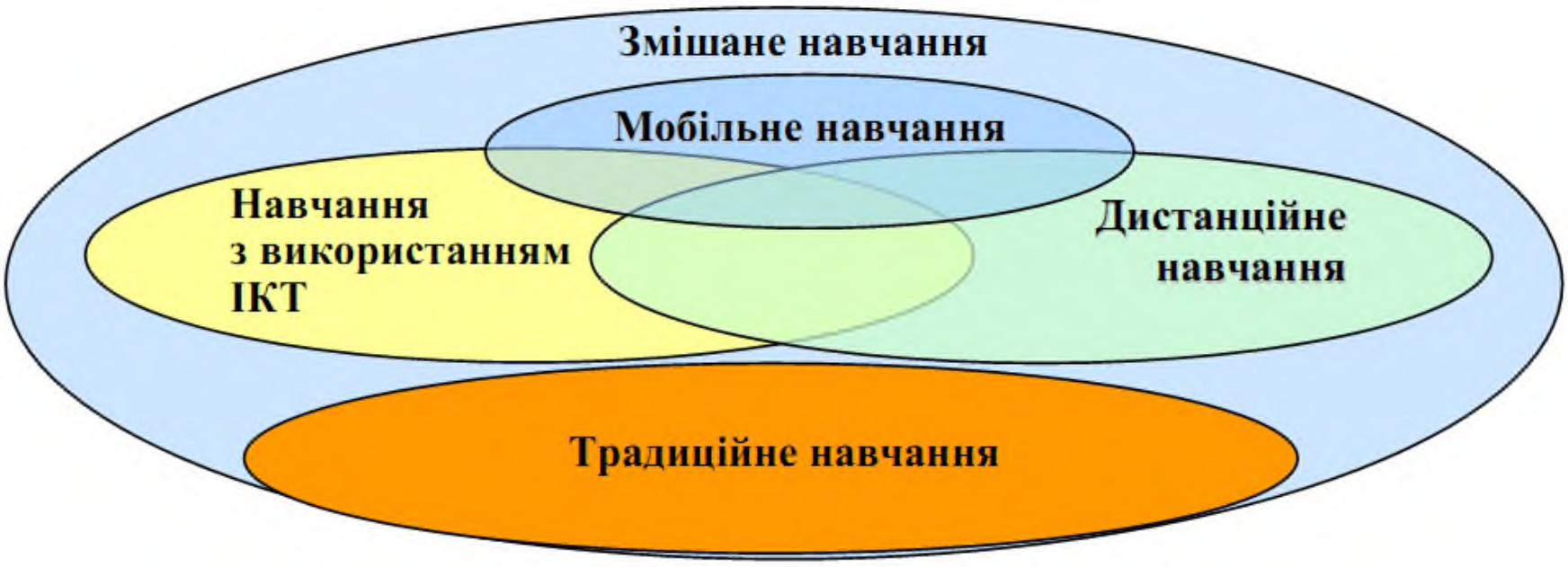

Рис. 1.11. Співвідношення понять Змішане навчання = Традиційне + З використанням ІКТ + Дистанційне + Мобільне навчання

Дж. Панкін (Jeff Pankin) [96] з США визначає *комбіноване навчання* як структуровані можливості вчитися, що використовують більше одного типу або методу навчання всередині або за межами аудиторії. Це тлумачення включає в себе різні форми та методи навчання (лекції, дискусії, практичні екскурсії, читання, ігри, тематичні дослідження, моделювання), різні способи доставки навчального матеріалу (в аудиторію або опосередковано за допомогою комп'ютера), взаємодія (синхронна або асинхронна) і різний рівень управління (особистий інструктор або експерт, або груповий / соціальний викладач). Використання комбінованого навчання відкриває можливості для створення ефективного навчання, заощаджуючи час і гроші установ, роблячи навчання більш цікавим і зручним для студентів, і пропонує викладачам можливості використовувати інновації.

На підставі зв'язку понять традиційного, дистанційного, мобільного навчання та навчання за допомогою Інтернет і мультимедіа, К. Л. Бугайчук виділяє *типові ознаки* мобільного навчання [149]:

– використання мобільних апаратних пристроїв (мобільні телефони, ноутбуки, ПК, планшети, електронні читачі «рідери» тощо);

– взаємодія учасників навчального процесу забезпечується і підтримується за допомогою бездротових мереж;

– навчання здійснюється незалежно від часу і місця знаходження його учасників;

– мобільне навчання надає можливість мобільне персональне навчальне середовище.

Опишемо моделі мобільного навчання, що виокремленні у системі освіти США [85]:

**Web-модель мобільного навчання**. У зв'язку з появою сучасних мобільних пристроїв і Інтернет-технологій, користувачі мають доступ до

Інтернет-ресурсів з мобільних пристроїв у будь-якому місці, поки у них є бездротове підключення до Інтернету. У результаті мобільний пристрій починає функціонувати як персональний комп'ютер. Якщо студенти використовують свої мобільні пристрої для доступу до навчальних відомостей у блозі або відвідують Web-сторінку, то таку діяльність можна назвати мобільною навчальною діяльністю.

**Прикладна модель мобільного навчання**. Завантажуючи додатки на мобільний телефон – прикладне програмне забезпечення, системи підтримки навчанням, системи комп'ютерної математики, довідкові додатки (словники, енциклопедії, перекладачі), що призначені для мобільних пристроїв, користувач отримує різні можливості, такі як доступ до навчальних матеріалів, до графіків проведення занять та консультацій, можливості виконувати мобільні онлайн вправи (Student Responce Systems (SRS), quick responce exercises, mobile quizzes, mobile polls), розв'язувати картки (flashcards) з математики для навчання та закріплення отриманих знань та багато іншого. Для доступу до додатків різних вправ використовують графічні QR-коди.

**Стільникова модель мобільного навчання**. Стільникова модель мобільного навчання використовує телефонні можливості мобільного пристрою. Кожен телефон має можливості передачі голосу і даних. При навчанні за допомогою текстових повідомлень передаються, наприклад, відповіді на тести. У той час як Web-модель та прикладна модель вимагають підключення до Інтернет, при стільниковій моделі мобільного навчання потрібно тільки підключення до мобільної мережі.

Дж. Тракслер (John Traxler) виділяє кілька напрямів реалізації мобільного навчання [84; 128]:

– технологічно орієнтоване мобільне навчання – окремі конкретні технологічні інновації, впроваджені у навчальний процес для демонстрації технічних переваг та педагогічних можливостей;

– мінінавчання за допомогою Інтернет і мультимедіа – мобільні, бездротові і портативні технології, які використовуються для повторного впровадження рішень і підходів, що вже використовуються у традиційних електронних засобах навчання, можливо, перенесення деяких технологій навчання за допомогою Інтернет і мультимедіа, таких, як віртуальні навчальні середовища (VLE), на мобільні платформи (MLE);

– комбіноване навчання – це навчання, що поєднує традиційне навчання з мобільним навчанням з метою створення гармонійного поєднання теоретичної та практичної складових процесу навчання;

– неформальне, персоналізоване, ситуативне мобільне навчання – мобільні технології з додатковою функціональністю, наприклад, залежні

від місця розташування;

– мобільні тренінги – технології, що використовуються для підвищення продуктивності та ефективності мобільних працівників шляхом надання матеріалів для підтримки «точно у термін» і в контексті їхніх першочергових пріоритетів;

– віддалене (сільське) розвивальне мобільне навчання – мобільні технології використовуються для вирішення інфраструктурних та екологічних проблем та підтримки освіти там, де традиційні технології навчання малоефективні.

Основне призначення мобільного навчання полягає в тому, щоб покращити знання людини в тій галузі, в якій вона бажає, і в той момент, коли їй це потрібно. Завдяки сучасним технологіям мобільного зв’язку (взаємодія «студент-викладач» здійснюється в високошвидкісному середовищі обміну повідомленнями) через мобільне навчання забезпечується високий ступінь інтерактивності, що має вирішальне значення для навчання [294].

## 1.2 Етапи розвитку теорії та методики використання інформаційно-комунікаційних технологій у навчанні вищої математики студентів інженерних спеціальностей у Сполучених Штатах Америки

Існує багато різних технологій, доступних для навчання вищої математики, що базуються на інформаційних технологіях. Деякім вже кілька десятиліть, тоді як інші є більш пізніми за походженням, для деяких вже створено широку базу для використання, в той час як використання інших пропонує більший потенціал у майбутньому.

На рис. 1.12 показано динаміку впровадження ІКТ у систему освіти (е-готовність, е-інтенсивність, е-вплив) з подальшою зміною рівнів наповнення систем освіти ІКТ. При впровадженні в освіту засобів ІКТ важливо, щоб викладачі та студенти мали доступ до обладнання та програмного забезпечення, а також мали основні навички роботи з цими засобами [51].

Н. Сінклер (Nathalie Sinclair) [66, 235-253] пропонує класифікувати використання освітніх ІКТ не за математичним змістом, а за способом їх застосування. Автором виділено три основні «хвилі» ІКТ:

– *хвиля* 1: використання ІКТ у процесі вивчення вищої математики. Технології «першої хвилі» зосереджують увагу на прямому взаємозв’язку між окремими студентами і дисципліною «Вища математика». Використання цих технологій у процесі навчання вищої математики полегшують процес навчання вищої математики, але методична частина їх використання є нерозробленою. До першої хвилі ІКТ віднесено комп'ютерне навчання. Засоби ІКТ першої хвилі звертають увагу на пряму зустріч певного користувача з математичним змістом, який погодився на заглиблення в математичне навчання, в якому викладачі, групи, навчальні програми та інші студенти не мали б істотної ролі;

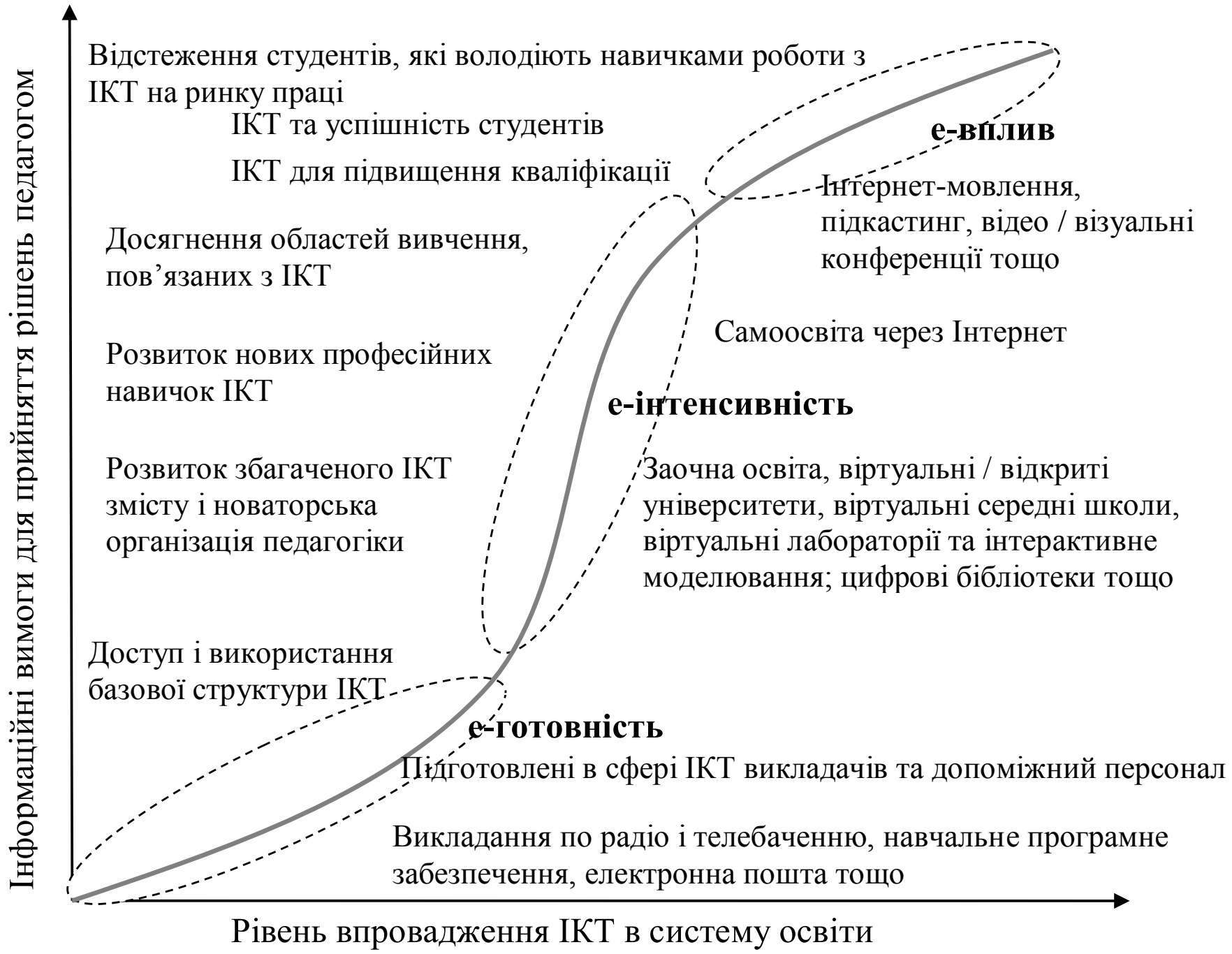

Рис. 1.12. Інформаційні потреби на різних рівнях впровадження ІКТ в систему освіти з часом [51]

– *хвиля* 2: розробка змісту навчання. Поява великої кількості новітніх ІКТ навчання (таких як графічні калькулятори, таблиці, засоби динамічної геометрії, статистичні програми, а також предметні мікросвіти) сприяла зміні методів та форм організації навчання вищої математики. У процесі навчання увага почала акцентуватися на використанні засобів ІКТ студентами для підтримку процесу навчання;

– *хвиля 3*: ІКТ майбутнього. Хвиля новітніх ІКТ характеризується розширенням можливостей викладачів та студентів. Третя хвиля використання ІКТ ґрунтується на визнанні того факту, що використання технологій має підтримувати можливість спілкування студентів із викладачем, один з одним, з усією групою і за межами аудиторії. До технологій «третьої хвилі» відносять: графічні калькулятори, мультимедійні дошки, засоби Інтернет (зокрема, відео-конференції).

Перевагою використання таких інструментів є те, що студенти отримають технічні навички в їх використанні, при цьому розширюються фізичні та інтелектуальні межі навчальної аудиторії.

У додатку Г розглянуто періодизацію розвитку та впровадження ІКТ у процес навчання, запропоновану С. А. Лисенком.

Розглядаючи впровадження ІКТ у навчання вищої математики та розвиток теорії та методики використання ІКТ у навчанні вищої математики студентів інженерних спеціальностей у США, можна виділити такі етапи:

Перший етап – 1965–1973 рр. – пов’язаний із появою достатньої кількості комп’ютерних засобів різного рівня, оснащених мовами високого рівня та специфікою апаратного забезпечення ІКТ (використання мейнфреймів з обмеженим мережним доступом).

До 1965 року використання ІКТ у навчанні вищої математики не мало системного характеру, незважаючи на вдалі, але поодинокі спроби. В Україні до таких робіт відноситься розробка спеціалізованого програмно-апаратного комплексу в 1960 році в Інституті кібернетики АН УРСР під керівництвом академіка В. М. Глушкова, що у Київському Вищому радіотехнічному училищі військ протиповітряної оборони (КВІРТУ, начальник училища генерал Г. Т. Ростунов) використовували для подання навчальних відомостей та контролю знань (тестування) курсантів з використанням ЕОМ «Промінь» [173].

Показовим прикладом спрямування використання ІКТ у навчанні вищої математики в США є оголошення в липневому номері журналу «New Scientist» за 1957 рік, у якому серед кваліфікаційних вимог до посади викладача (помічника лектора – асистента) були знання мов програмування (насамперед FORTRAN) та комп’ютерної техніки.

До речі, назва цієї мови походить від FORmula TRANslation («переклад формул» мовою, зрозумілою для комп’ютера) та відображає прикладний аспект інженерної математичної підготовки. Проте пряме перенесення програмування цією мовою у навчальний процес вищої школи зазнало утруднень, пов’язаних із:

1) станом розвитку засобів ІКТ: формування мейнфреймів – «великих» комп’ютерів високої вартості, – обслуговування яких було спрямоване на зменшення витрат від пристроїв, тому провідним режимом роботи таких комп’ютерів був пакетний (неінтерактивний режим виконання певної послідовності програм з відкладеним аналізом результатів їх роботи), у той час як для навчання провідним режимом роботи мав бути діалоговий;

2) обмеженими фінансовими можливостями закладів освіти, наслідком яких було використання насамперед застарілих

мікрокомп'ютерів, що часто не мали засобів розробки мовою FORTRAN;

3) нерозробленістю психолого-педагогічних засад використання засобів ІКТ у навчанні.

На розв'язання останньої із проблем й була спрямована робота Б. Скіннера і Н. А. Кроудера (Norman Allison Crowder) [112] з так званого «програмованого навчання». У сучасній зарубіжній психології запропонований ним підхід носить назву «інструкціоністського навчання» – гілки біхевіоризму.

Як зазначає М. І. Жалдак у [172], незважаючи на хибність цієї концепції, її розробка та застосування сприяли розробці перших «неінтуїтивно створених» навчальних програм, що надавали можливість обмеженої реалізації індивідуальної навчальної траєкторії.

У концепції програмованого навчання передбачалась така організація процесу засвоєння знань, умінь і навичок, що на кожному ступені навчального процесу чітко обумовлювались ті знання, уміння і навички, що мають бути засвоєні, і контролювався процес засвоєння [331].

Головна ідея цієї концепції – управління учінням, навчально-пізнавальними діями студентів за допомогою навчальної програми.

Основним поняттям даної концепції була навчальна програма, під якою розуміли алгоритм пізнавальних дій, що містить послідовні мікроетапи опанування одиницею знань або дій. Алгоритм складався з трьох частин:

– доза відомостей про предмет, що вивчається;

– завдання (операції) щодо роботи з даними та їх засвоєнням;

– контрольні завдання і вказівки про повторення вправи або перехід до наступного етапу.

Розв'язання другої проблеми було виконано у 1964 році Дж. Кемені (John Kemeny) та Т. Курцем (Thomas Kurtz) [86]. Створена ними мова BASIC була чи не найпершою спробою із реалізації інструкціоністського навчання на рівні мови програмування:

– по-перше, створена ними мова передбачала активне використання операторів уведення / виведення, тобто діалоговий режими роботи;

– по-друге, компактність ядра мови надавала можливість реалізації інтерпретаторів з неї на застарілому комп'ютерному обладнанні;

– по-третє, можливість роботи інтерпретатора BASIC на діалогових терміналах надало можливість застосування мейнфреймів із віддаленим доступом до них.

За таких умов університетські комп'ютерні комплекси могли ставати регіональними осередками із використання ІКТ у навчанні.

Не слід вважати, що лише комп'ютери у вузькому сенсі були об'єктом інформатизації в системі освіти США – так, у 1964 році відбулась публічна дискусія у коледжі Хоупа на тему надання високошвидкісних комп'ютерів та електронних калькуляторів для потреб навчання [118].

Електронні калькулятори (за сутністю спеціалізовані мікрокомп'ютери) у навчанні вищої математики у той час були провідними, проте далеко не єдиними засобами – так, у [113] наводиться перелік засобів ІКТ для навчання математичної логіки: комп'ютерний термінал з можливістю візуального та аудіального подання навчальних матеріалів, клавіатура для введення письмових відповідей, мікрофон для аудіовідповідей та світлове перо для вибору об'єктів. Автори [113] характеризують це не як набір засобів, а як навчальний комплекс, що складається з керуючого мінікомп'ютера та шести мережних «станцій», кожна з яких була розташована у кімнаті малого розміру (2х2,5 м).

Таким чином, на початок 1965 року у системі освіти США була достатня кількість комп'ютерних засобів різного рівня, оснащених мовами високого рівня, що надає можливість вважати цей рік умовною нижньою межею першого етапу розвитку теорії та методики використання ІКТ у навчанні вищої математики студентів інженерних спеціальностей у США.

У 1965 році фірмою DEC (Digital Equipment Corporation) був випущений перший комерційно успішний мінікомп'ютер – PDP-8.

Як згадується у [50; 56], комп'ютер PDP-8 використовували на заняттях з математики на математичних факультетах. У статті [56] показано як два вчителі математики запропонували створення окремого комп'ютерного відділу, частково, для навчання дітей, які не мають математичних схильностей, і частково тому, що вони не встигали надавати консультаційну допомогу учням за комп'ютером, якщо вони цього потребували, одночасно виконуючі іншу роботу.

Сімейство мінікомп'ютерів PDP часто зустрічається у розвідках з історії ІКТ саме через їх пошireність. Оптимальне на той час поєднання засобів ІКТ за поміркованих цін сприяло поширенню цього сімейства у наукових та освітніх установах. Одним із дослідницьких проектів із створення діалогової математичної системи стала Reduce (www.reduce-algebra.com) – система комп'ютерної алгебри загального призначення.

Розроблені згідно концепції Б. Ф. Скіннера засоби навчання мали спільний суттєвий недолік – їх використання не сприяло розвитку особистості у процесі навчання через неврахування розробок Л. С. Виготського із «зони найближчого розвитку» [159, 115].

Дослідження Ж. Піаже (Jean William Fritz Piaget) із психології раннього дитинства надали можливість його учню С. Пейперту (Seymour

Papert) у 1967 році у Массачусетському технологічному інституті, не відходячи повністю від інструкціонізму, запропонувати новий засіб навчання – мову LOGO [86], що базується на конструктивістському підході до навчальної діяльності. Сутність цього підходу найкраще можна визначити фразою «вчуся, навчаючи». Дійсно, у середовищі LOGO її користувач – програміст – виступав у ролі «вчителя» для головного об'єкта мікросвіту LOGO – черепахи, «навчаючи» її через програмування виконувати певні дії.

Як зазначено в [86], мова програмування LOGO призначена для заохочування до строгого мислення в математиці. С. Пейперт хотів, щоб для дітей було доступно і просто висловити процедуру розв'язання простих завдань. Він використовував LOGO для навчання математики. С. Пейперт наполягав, що не потрібно викладати математику, але треба вчити дітей бути математиками. LOGO незабаром стала мовою комп'ютерної грамотності у молодшій школі.

Дизайн середовища LOGO справив значний вплив на подальший розвиток засобів навчання та навчальні концепції. Так, один із співробітників С. Пейперта – А. Кей (Alan Curtis Kay), якого сьогодні вважають батьком об'єктно-орієнтованого програмування, запропонував у 1968 році Dynabook [31] – «персональний комп'ютер для дітей будь-якого віку», оснащений середовищем мови Smalltalk (www.smalltalk.org/main).

Це була перша мова високого рівня, що підтримувала експериментування із широким набором математичних об'єктів – від чисел з різних множин до геометричних об'єктів.

Сьогодні об'єктно-орієнтований підхід є основою побудови інтерфейсів користувача, а Dynabook знайшов своє втілення у програмних (середовище Squeak – www.squeak.org) та апаратних (ноутбуки та планшетні комп'ютери) засобах.

Таким чином, із виходом LOGO та Smalltak-72 програмоване навчання перестало бути домінуючою концепцією, що зумовило вибір верхньої межі першого етапу.

Знову звернувшись у якості прикладу до історії використання ІКТ у коледжі Хоупа, можна прослідкувати їх еволюцію протягом першого етапу:

1967 рік – створення обчислювальної лабораторії та обладнання її калькуляторами IME-86 [57];

1968 рік – застосування IBM 1130 для генерування псевдовипадкових чисел [5];

1969 рік – проект Національного наукового фонду США «Використання комп'ютерів у навчанні статистики» [116];

1970 рік – дворічний інтегрований курс «Прикладна статистика та програмування» мовою FORTFAN [118];

1972 рік – видане керівництво із виконання лабораторних робіт з теорії ймовірностей та математичної статистики [120].

Традиційно не залишилися осторонь й провідні ВНЗ: так, у MTI у 1968 році у рамках проекту MAC (Mathematics and Computation, пізніші тлумачення – Multiple Access Computer, Machine Aided Cognitions, або Man and Computer) було створено систему комп’ютерної математики Macsyma (maxima.sourceforge.net), а у 1971 році IBM під керівництвом Р. Дженкса (Richard Jenks) – систему Axiom (www.axiom-developer.org), що в той час мала назву Scratchpad.

Винайдення у 1971 році флоппі-диску сприяло персоналізації використання мейнфреймів та мінікомп’ютерів: 81,6 Кб даних, що уміщувались на 8-дюймовому диску, було цілком достатньо для зберігання текстових документів (статей, програм тощо), найбільш поширених у академічному середовищі [126].

Для підтримки навчання математики у цей період були розроблені ряд спеціалізованих пристроїв:

– KENBAK-1 (розробник – Дж. Бланкенбейнер (John Blankenbaker), 1972 рік) – перший комп’ютер для навчання, розрахований на непрофесійних користувачів. Як зазначається у [13], цей комп’ютер поєднував гнучкість та доступність;

– HP-35 (розробник – Hewlett Packard, 1972 рік) – не перший калькулятор HP, але один з перших, що містив мікропроцесор. Це сприяло його компактності та зручності використання у навчанні математики. Дж. К. Хорн (Joseph K. Horn), який використовував калькулятори HP як учитель математики [60], є автором багатьох статей з їх застосування. Як спеціалізований мікрокомп’ютер, HP-71 (подальший розвиток HP-35) надавав можливість програмування мовою BASIC, а його наступник – мовою СКМ Derive.

Як зазначає Дж. Г. Харві (John G. Harvey) [55], у навчанні вищої математики калькулятори доцільно застосовувати, зокрема, для тестування.

HP-35 був розроблений «для інженерів та студентів інженерних спеціальностей» [55, 140]. Наукові та графічні калькулятори HP (як програмовані, так й непрограмовані) підтримують й нетрадиційні для калькуляторів дії над матрицями та статистичні функції.

Ціна на чотирьох-функціональний калькулятор HP-35 в той час знаходилась в діапазоні від 4 до 7 доларів США, порівняно з ціною на простий (не програмований) калькулятор TI Data Math, що становила від 10 до 15 доларів США. В результаті калькулятори стають більш доступними для студентів навіть із сімей з невисоким матеріальним статком. Стає можливим придбання власних калькуляторів усіма студентами для використання їх на заняттях з математики [55, 140].

Як приклади застосування калькуляторів при навчанні вищої математики можна навести такі: калькулятори використовувались в Університеті штату Огайо при вивченні алгебри (надавали можливість виконувати операції з матрицями, проводити дослідження змінних величин, їх застосування при побудові таблиць) (1988 р.); крім того, в Університеті штату Огайо було розроблено тестові завдання з математики, в яких необхідно було використовувати ручні калькулятори (1984 р.); виконання завдань, в яких необхідно було застосовувати графічні інструменти, такі як графічні калькулятори, розроблені фірмами Casio та Sharp (1989 р.); два викладача розробили модулі з математики для закладів освіти в Техасі, що передбачають використання калькуляторів TI Math Explorer – дробового калькулятора для навчання операцій з дробами [55, 141].

Головною проблемою, на думку Дж. Г. Харві [55, 145], була проблема нерозробленості методики ефективного використання калькуляторів у навчанні математики, з одного боку, та нерозробленості методики навчання математики, орієнтованої на використання калькуляторів.

Протекціонізм комп'ютерних фірм до окремих ВНЗ привів до суттєвого зростання забезпеченості навчальних закладів засобами ІКТ: так, якщо у 1965 році менше 5 % всіх навчальних закладів були забезпечені комп'ютерами для навчальних потреб [124, 51], то у 1972 році передавання даних через мережу набуло актуальності за рахунок суттєвого зростання забезпеченості комп'ютерною технікою та комунікаційними засобами [135, 11].

Проведений аналіз надає можливість визначити основні особливості використання засобів ІКТ у навчанні вищої математики студентів інженерних спеціальностей на першому етапі їх розвитку:

1) діалоговий режим роботи з навчальними програмами;

2) поява перших систем підтримки математичної діяльності без програмування мовами загального призначення;

3) розмаїття апаратного та програмного забезпечення;

4) домінування біхевіоризму в обґрунтуванні використання ІКТ та розробці навчальних програм;

5) уведення програмування в курси вищої математики.

Указані особливості породили ряд протиріч:

1) між недостатньою адекватністю програмованого навчання для опису навчальної діяльності та особистим розвитком студента у процесі

навчання;

2) між необхідністю посилення прикладного аспекту навчання вищої математики майбутніх інженерів та недостатністю навчального часу на одночасне опанування математики та програмування;

3) між розмаїттям засобів ІКТ та необхідністю уніфікації засобів навчання.

Часткове розв'язання вказаних протиріч було досягнуто ще на першому етапі – з'явились нові, більш адекватні підходи до моделювання процесу навчання, перші системи комп'ютерної математики та у 1969 році – операційна система UNIX [67], спрямована на об'єднання різних засобів ІКТ у єдиному мережному середовищі.

Тому другий етап розвитку теорії та методики використання ІКТ у навчанні вищої математики студентів інженерних спеціальностей у США пов'язаний саме із цією операційною системою.

<u>Другий етап</u> – 1973–1981 рр. – пов'язаний із поширенням в університетах США мережної операційної системи UNIX, використанням міні- та мікрокомп'ютерних систем.

Із спогадів викладача вищої математики з Каліфорнійського університету в Берклі [56]: «було вирішено створити середовище в середній школі, що було б найбільш близьким до лабораторій MTI і Стенфордського університету, де навчався я. Це означало, що необхідна була потужна комп'ютерна система, з великою кількістю програмних засобів, готовність аудиторії навчатися не традиційно і організація процесу навчання, що відрізнялося від шкільного. Було встановлено на PDP-11/70 версію UNIX 7; ми були альфа-тестувальниками BSD 2.9, версії Berkeley Unix для PDP-11. Установка, тестування і налагодження цієї нової системи було проведено виключно за рахунок студентів».

Цей фрагмент містить згадку про PDP-11 – подальший розвиток застосовуваної на першому етапі PDP-8 та показує зацікавленість фірми-виробника (DEC) у наданні власних засобів ІКТ університетам із ОС UNIX. Версія UNIX – BSD – відображає внесок університету в її розробку (Berkeley Software Distribution) [132].

Головна особливість UNIX – мобільність (насамперед, програмна) – створила умови для поширення програмного забезпечення під управлінням цієї системи на різні засоби ІКТ – від мейнфреймів до міні- (а надалі і мікро-) комп'ютерів.

Багато відомих на теперішній час систем комп'ютерної математики були створені саме у цій ОС: MATLAB (наприкінці 1970-х), Maple (1979) та інші були розроблені як мови програмування для навчання студентів вищої математики, використовуючи різні математичні бібліотеки, не вивчаючи мову FORTRAN.

У зв'язку із тим, що ОС UNIX та її програмне забезпечення були мобільними, з'явилась можливість для об'єднання не лише обчислювальних ресурсів різних комп'ютерів, а й користувачів, їх програм та даних у мережному середовищі: «отримана система UNIX надавала користувачам можливості віддаленого доступу до терміналу і спільно використовувати файлові системи; вихідний код поставлявся з системою, користувачі могли обмінюватися даними та програмами безпосередньо і неофіційно; оскільки UNIX працювала на відносно недорогий міні-ЕОМ, малі групи дослідників могли вільно експериментувати з нею» [134].

Потенціал UNIX для навчального процесу був відзначений на IV Конференції з використання комп'ютерів у навчальному процесі, що відбулася в 1970 році в Університеті Північної Кароліни [118].

Початок другого етапу відзначався узагальненням досвіду використання засобів ІКТ навчання вищої математики; зокрема, на восьмій конференції з використання комп'ютерів у навчанні студентів [115] (1977), конференції НАТО з комп'ютерно-орієнтованого навчання [118] (1976) та інших.

У середині 1970-х рр. з'явився Журнал мічиганської асоціації користувачів комп'ютера у навчанні, у якому в 1978 році було опубліковано фундаментальну роботу по навчанню теорії ймовірностей та математичної статистиці із застосуванням комп'ютерів [117].

Розпочата винаходом дискет персоналізація засобів ІКТ була продовжена розробками 1977 року – комп'ютерами Apple II та TRS-80.

Застосування конструктивістського підходу до навчання вищої математики сприяли нові засоби Apple II для пересічних користувачів: графічний дисплей та маніпулятор «миша». Це надавало можливість створювати навчальні ігри, системи комп'ютерної графіки та динамічної геометрії. А убудована мова Apple II – BASIC – була доповнена графічними командами. Конструкція виявилась настільки вдалою, що стала родоначальником не лише комп'ютерів Apple Macintosh, а й вітчизняного комп'ютера для системи освіти «Агат» [15].

Як зазначено в [25], в рамках конструктивізму, методичною метою розробленого проекту навчання математичних дисциплін на комп'ютерах Apple II було створення навчального середовища, що сприяло побудові у студентів математичних понять через повторювані цикли розроблених завдань, вивчаючи завдання і розмірковуючи про способи їх розв'язання.

У [91] згадується про використання комп'ютерів TRS-80 як засобів навчання в аудиторії. Викладачі використовували програму на Level II BASIC для комп'ютера TRS-80, що імітувала машину Тюрінга і демонструвала природу пристрою. Програма запускалася в динамічному режимі і була призначена для використання в якості навчального посібника з інформатики або математики для студентів, які вивчали теорію числення або теорію автоматів.

Поширення персоналізованих комп'ютерів з графічним інтерфейсом, орієнтованих на ігрову діяльність (Atari 800 та інші) також сприяло розвитку ігрових засобів навчання математики [38].

Опубліковані у 1978 році «Рекомендації з виконання лабораторних робіт з теорії ймовірностей та математичної статистики» [119] зіграло значну роль у розвитку теорії та методики навчання математичних дисциплін із застосуванням комп'ютерів. Через 15 років після перших спроб уведення комп'ютерів у навчання вищої математики майбутніх інженерів відчувається суттєве зміщення: з навчання математики разом з програмуванням до обґрунтованого використання програмних засобів у процесі навчання.

Індустрія програмних засобів навчального призначення на другому етапі була орієнтована переважно на Apple II – й досі у Інтернет магазині Amazon можна придбати програмне забезпечення динамічної геометрії для історичної платформи Apple II [47].

У 1975 році інженерами Xerox PARC (тієї самої лабораторії, де були винайдені миша, графічний інтерфейс та Smalltalk) був отриманий патент США на технологію Ethernet [89], а вже у 1979 році з'явилась перша відкрита інформаційна онлайн-служба Compuserve. І хоча ядром її були UNIX-системи, об'єднувала вона й персональні комп'ютери.

Таким чином, можна виокремити такі характерні риси другого етапу розвитку теорії та методики використання ІКТ у навчанні вищої математики студентів інженерних спеціальностей у США:

1) перехід від використання мов програмування у навчанні до використання математичних бібліотек, систем комп'ютерної математики та мов високого рівня;

2) застосування комп'ютерної графіки у навчальних програмах;

3) поява нових класів навчальних програм – навчальних ігор, систем динамічної геометрії та електронних таблиць;

4) поява та поширення комп'ютерних мереж, що об'єднували викладачів та студентів;

5) розвиток засобів ІКТ навчання вищої математики – графічних та символьних калькуляторів.

У розвитку засобів ІКТ на цьому етапі можна виділити ряд протиріч:

1) між потенціалом використання мультимедійних засобів комп'ютерних систем та не розробленістю психолого-педагогічних основ їх використання;

2) між потребою студентів та викладачів у персональних засобах ІКТ

та недостатністю пропозицій виробників;

3) між потребою виробників персональних комп'ютерів у адаптованому для них варіанті UNIX та недостатнім апаратним забезпеченням її функціонування.

Вказані протиріччя визначили початок третього етапу, на якому вони були розв'язані.

<u>Третій етап</u> – 1981–1989 рр. – пов'язаний із поширенням персональних комп'ютерів.

Нижня межа етапу (1981 рік) відповідає появі ОС MS DOS та комп'ютера IBM PC. Не будучи ані бажаним клоном UNIX, ані першим персональним засобом ІКТ, вони відіграли надзвичайно значну роль, сформувавши сучасний ринок програмних засобів навчання вищої математики.

Персоналізація (користувача, комп'ютера, програм, даних) була провідною концепцією цього етапу. Навчання з використанням комп'ютера є одним з найбільш поширених підходів, призначених підвищити навчальні досягнення студентів вищої школи.

У ці роки значна кількість робіт була присвячена психолого-педагогічному обґрунтуванню навчання з використанням персональних комп'ютерів і розробці навчальних курсів з дисциплін. Аналізуючи публікації того часу, можна зазначити, що найбільш поширеними мовами програмування були BASIC, Pascal, FORTRAN, Algol. Так, у Массачусетському технологічному інституті на заняттях з вищої математики студенти займались написанням комп'ютерних програми для проведення досліджень з різних тем математики. Тенденції і досягнення в області численних методів, обчислення і новаторські дослідження в галузі прикладної математики сприяли впровадженню комп'ютерного викладання прикладної математики [133].

Дж. Енгельбрехт (Johann Engelbrecht) [39] виявляє ряд переваг застосування комп'ютерів у процесі викладанні вищої математики: відпрацювання практичних умінь та навичок; розмаїття подання навчального матеріалу; використання у процесі моделювання та програмування.

Викладач кафедри вищої математики університету Каліфорнії Дж. Дж. Ваврік (John J. Wavrik) [136] – автор курсу «Комп'ютерна алгебра» для студентів у символьних обчисленнях (середина 80-х років XX століття) один із перших стверджував, що використання комп'ютерів є необхідним у процесі навчання вищої математики, оскільки їх використання сприяє підвищенню рівня навчальних досягнень студентів [65].

Стаття Дж. Дж. Вавріка [136] «Комп'ютери і кратні корені полінома»

«…не описує різні повороти долі, що ведуть працівника в чистій галузі математики такій, як обчислювальна геометрія, пов'язаної з комп'ютерами. Цілком ймовірно, що персональні комп'ютери стають більш поширеними серед алгебраїстів, які роблять акцент на використання їх для науково-дослідної роботи… Природньо сподіватися, що комп'ютери можуть бути використані для полегшення складності обчислень. Як виявилося, цей процес не такий вже й і простий, як спочатку може здатися… Ми хотіли б мати машину, що може допомогти, наприклад, в обчисленнях інваріантів для специфічних станів об'єктів дослідження…» У статті автор, розглядаючи питання знаходження кратних коренів полінома, зауважує, що «проблеми ефективного обчислення найбільшого спільного дільника двох многочленів з цілими коефіцієнтами отримала велику увагу… Комп'ютерні системи, призначені для алгебраїчних обчислень, надали можливість «з нескінченною точністю» робити арифметичні обчислення та використовувати для алгоритмів систем цього типу. Машини, що виконують обчислення з заданою точністю часто допускають точність лише між 6 і 16 цифрами наближень... Ця стаття містить комп'ютерну програму, що служить прикладом реалізації на конкретній машині. Програма, написана на BASIC – найбільш поширеній мові для мікрокомп'ютерів». Автор вказує, що «метою даної статті є представлення читачеві деяких цікавих аспектів алгебраїчних обчислень. Програма, включена в статтю, не подається читачам як частина готового програмного забезпечення, а скоріше, щоб дати їм ідеї для розробки своїх власних програм».

У 1984 році Р. Д. Пеа (Roy D. Pea) [98] використовував програму AlgebraLand, застосування якої надало студентам можливість автоматизувати алгебричні обчислення і зосередити увагу на розв'язанні завдань більш високого рівня складності. AlgebraLand розроблена була для того, щоб з її допомогою студенти «вивчали питання приросту швидкості», закріплюючи навички при розв'язанні задач.

Якщо звернутися у якості прикладу до історії використання ІКТ у коледжі Хоупа, можна прослідкувати їх еволюцію протягом третього етапу [118]:

1982 рік – перша Міжнародна конференція з викладання статистики (ICOTS I) у «Використання мікрокомп'ютерів для розуміння понять теорії ймовірностей та математичної статистики»;

1983 рік – щорічні збори Математичної асоціації Америки в Денвері «Використання мікрокомп'ютерів надає можливість ілюструвати основні поняття теорії ймовірностей та математичної статистики»;

1984 рік – конференція з використання мікрокомп'ютерів у статистиці Американської асоціації статистиків Університету штату Делавер «Використання мікрокомп'ютерів у галузі викладання теорії ймовірностей та математичної статистики»;

1985-1996 роки – створення обчислювальної лабораторії «Комп'ютерна лабораторія вступу до статистики» та обладнання її ПК IBM з BASIC;

1987 рік – щорічні збори Американської статистичної асоціації у Сан-Франциско, Каліфорнія «Комп'ютерне моделювання для підтвердження теоретичних понять».

Значного розвитку на цьому етапі набувають системи комп'ютерної математики, а саме: Cayley (1982), MathCAD (1985), Fermat (1985), GAP (1986), Mathematica (1986), Derive (1988), MuPad (1989).

У звіті 1989 р. Національної академії наук США про майбутнє математичної освіти [41] (так, як воно бачилось наприкінці третього етапу) вказується, що символьні комп'ютерні системи вимагають фундаментального переосмислення методики навчання математики з використанням комп'ютера: якщо до цього використовувались обчислювальні та графічні програми, що мали похибки округлення та візуального подання, то символьні обчислення є точними та відповідають діям, що виконує людина. «Пріоритети математичної освіти повинні змінюватися з метою відображення шляхів використання комп'ютерів у математиці» [41, 78]: провідними засобами навчання вищої математики укладачі звіту вважали електронні таблиці, пакети числового аналізу, символьні комп'ютерні системи, засоби графічного подання відомостей, перспективними – електронні підручники, віддалені класні кімнати, інтегровані навчальні середовища.

Аналізуючи зростаючу роль комп'ютерів у математичній освіті, укладачі звіту прогнозували наступне [41, 62-63]:

1. Шкільна математика може стати більш схожою на математику, яку люди використовують у наукових дослідженнях та практичній діяльності. Використовуючи машини для прискорення розрахунків, студенти можуть отримати досвід реальної математичної діяльності як дослідницької.

2. Слабкі алгебраїчні навички більше не будуть перешкодою для розуміння студентами математичних понять: так само, як комп'ютерні системи перевірки правопису надають можливість зняти психологічний блок у письменників із жахливим правописом, так і нові калькулятори надають можливість мотивувати студентів, слабких в курсі алгебри або тригонометрії, опанувати математичний аналіз чи статистику. Використання калькуляторів в класі може допомогти зробити вищу математику більш доступною.

3. Навчання математики може стати більш активним і динамічним, а,

отже, і більш ефективними. Перекладаючи більшу частину обчислювального навантаження на домашні завдання, калькулятори і комп'ютери дозволять студентам дослідити більш широке коло прикладів, пересвідчитись у динамічній природі більшості процесів, побудувати математичні моделі з обчислювальним експериментом над реалістичними – а не спеціально дібраними і спрощеними – даними, а також зосередити увагу на важливих концепціях, а не рутинних розрахунках.

4. Студенти можуть вивчати математику самостійно, ставити і відповідати на незліченні «а що буде, якщо?..». Хоча калькулятори і комп'ютери не обов'язково спонукають студентів думати самостійно, вони можуть забезпечити середовище, в якому згенеровані студентом математичні ідеї можуть процвітати.

5. Час, інвестований у дослідження в процесі навчання математики, повертається довгостроковою інтуїцією та осяяннями, а не тільки методами обчислень, що швидко забуваються. Інноваційне навчання на основі нового симбіозу машинних обчислень і людського мислення може зсунути баланс навчання у бік розуміння, проникливості та математичної інтуїції.

«Підручники, програмне забезпечення, комп'ютерні мережі у найближчі роки скомбінуються у новий гібридний освітньо-інформаційний ресурс» [41, 82] – саме цей прогноз ознаменував кінець третього та початок четвертого етапу розвитку теорії та методики використання ІКТ у навчанні вищої математики студентів інженерних спеціальностей у США.

Таким чином, до характерних рис третього етапу розвитку теорії та методики використання ІКТ у навчанні вищої математики студентів інженерних спеціальностей у США відносяться:

1) широке використання математичних бібліотек, систем комп'ютерної математики та проблемно орієнтованих мов;

2) широке впровадження персональних та персоналізованих засобів ІКТ у навчання математичних дисциплін;

3) використання ІКТ загального призначення (текстові редактори, електронні таблиці, бази даних тощо) для підтримки навчання математичних дисциплін.

У розвитку засобів ІКТ на цьому етапі можна виділити ряд протиріч:

1) між потенціалом об'єднання використання глобальних комп'ютерних мереж та недостатньою розробленістю персональних засобів доступу до них;

2) між доцільністю перенесення гіпертекстових та гіпермедіальних систем навчального призначення у мережне середовище та недостатньою

розробленістю засобів доступу до них;

3) між потребою студентів та викладачів у перенесенні навчальних матеріалів у мережу та не розробленістю психолого-педагогічних основ її використання.

Вказані протиріччя визначили початок четвертого етапу, на якому вони були розв'язані.

Четвертий етап – 1989–1997 рр. – пов'язаний із створенням World Wide Web та використанням технологій Web 1.0.

У статті «Навчання бакалаврів математики засобами Інтернет» [39] Дж. Енгельбрехт (Johann Engelbrecht) і А. Хардінг (Ansie Harding) зазначають, що Інтернет-освіта з математики розвивається як новий метод навчання зі своїми особливостями і можливостями, відрізняючись від традиційних методів навчання. Інтернет викладання настільки відрізняється від будь-якої з категорій дистанційного навчання, що йому передували, що це, по суті практика без дослідження основ. Інтернет навчання має багато аспектів, що не притаманні традиційному навчанню.

Автори [39] вважають, що при розробці онлайн-курсів з математики, слід звернути увагу на ряд педагогічних проблем та попередній досвід, що передбачає тонкий і складний процес організації курсу, де викладачі і розробники курсу створюють та удосконалюють навчальне середовище та інфраструктуру, що сприяє навчальній діяльності студентів, заснованої на існуючих знаннях і заохочує та формує цілеспрямованість студентів.

Автори статті підкреслюють, що навчальне середовище повинно бути розроблене таким чином, щоб заохочувати студентів до навчання, підтримувати процес навчання. На початку створення Інтернет-освіти виникало багато проблем, оскільки методика навчання в мережі була нерозробленою – викладачі просто переносили свої традиційні курси в Інтернет. Перше покоління Web-курсів критикують за відсутність основ педагогіки та за використання сумнівних навчальних стратегій. Навіть хороша лекція в традиційній освіті стає неефективною при перетворенні її в електронну сторінку, оскільки в результаті такого перетворення втрачається вербальна складова процесу навчання та можливість моніторингу сприйняття навчального матеріалу студентами.

Аналіз літератури даного періоду надав можливість визначити основні причини неефективного використання технологій Веб 1.0 [39]:

– нездатність залучати студентів до створення ефективних навчальних матеріалів;

– заплутаність інтерактивної взаємодії;

– акцент викладачами був спрямований на зміст, а не на результат;

– неадаптоване застосування традиційних дидактичних підходів до нових технологій;

– невизнання соціального характеру навчання.

Навчальні матеріали з вищої математики для студентів технічних ВНЗ повинні бути доступними та включати в себе використання різноманітних засобів ІКТ, а також бути відкритим. Посилання на інші Інтернет-ресурси збагачують курс. Студенти повинні вільно орієнтуватися в інтерфейсі курсу, а навчальний матеріал повинен бути максимально наближений до «живого» викладання.

Автори [39] до переваг онлайн навчання того часу відносять широкий спектр доступних освітніх ресурсів. Студенти отримують доступ до глобальних ресурсів, відомості стають динамічними, а навчальні матеріали можна зберігати і переглядати в зручний для себе час. Асинхронна природа онлайн-курсів надає студентам можливість самостійно вирішувати коли і де виконувати завдання, забезпечуючи можливість студентам працювати у своєму власному темпі.

У книзі Б. Х. Хана (Badrul H. Khan) [68] 1997 року Web-орієнтоване навчання (Web-based instruction – WBI) розглядаються як інноваційний підхід для надання навчальних інструкцій для віддаленої аудиторії, використовуючи Інтернет в якості навчального середовища.

Серед компонентів WBI виділяють [68]:

1. *Розробка навчального змісту дисципліни*: навчальні теоретичні матеріали; педагогічне проектування; розробка навчальних програм.

2. *Мультимедійна складова*: текст і графіка; потокове аудіо; потокове відео; графічний інтерфейс користувача – використовуються іконки, графіка, вікна і вказівні пристрої замість текстового режиму інтерфейсу (Microsoft Windows і MacOS є прикладами графічних користувацьких інтерфейсів); технологія стиснення.

3. *Інтернет-інструменти*: засоби зв’язку (асинхронні: електронна пошта, списки розсилки, групи новин тощо; синхронні: інструменти текстових та аудіо-відео конференцій); засоби віддаленого доступу (Telnet, передача файлів); інструменти Інтернет-навігації (Gopher, Lynx); пошукові системи; лічильники відвідувань.

4. *Комп’ютери та пристрої зберігання навчальних даних*: комп’ютерні платформи Unix, DOS, Windows і операційна система Macintosh; сервери, жорсткі диски, компакт-диски.

5. *З’єднання і послуги*: модеми; послуги Інтернет-провайдерів.

6. *Авторські програми*: мови програмування (HTML – мова гіпертекстової розмітки, VRML – мова моделювання віртуальної реальності, Java, JavaScript); Authoring Tools (простіше у використанні, ніж мови програмування); HTML-конвертери і редактори.

7. *Сервера*: HTTP сервери (HTTPD), інтерфейс Common Gateway Interface (CGI) тощо.

8. *Браузери та інші програми*: текстові браузери, графічні браузери, VRML браузери; посилання (наприклад, гіпертекстові посилання, гіперпосилання, гіпермедіа посилання, 3D посилання, карти посилань); додатки, які можуть бути додані в Web-браузери, такі як плагіни.

Широкого поширення на заняттях з вищої математики в ці часи зазнають графічні калькулятори TI-83. Серед викладачів були побоювання з приводу того, що студенти втрачають основні навички проведення досліджень функцій, та проведені дослідження підтвердили думку про те, що при використанні студентами графічних калькуляторів, їх навчальні досягнення суттєво зростають. Студенти, які використовують графічні калькулятори при вивченні вищої математики мають більш глибоке і широке розуміння поняття функції, ніж при використанні олівця та паперу [18, 39].

У цей час розроблено такі системи комп'ютерної математики: MuPad (1989), Magma (1990), Python (1991), Maxima (1998).

З 1997 року починає широко застосовуватися система MAPLE V у навчанні математичних дисциплін. Так, у коледжі Хоупа в цей час проходять конференції та семінари, присвячені питанню використання MAPLE V при навчанні теорії чисел, теорії ймовірностей та математичної статистики [118].

Таким чином, можна виокремити такі характерні риси четвертого етапу розвитку теорії та методики використання ІКТ у навчанні вищої математики студентів інженерних спеціальностей у США:

1) широке впровадження комп'ютерних мереж, що об'єднували викладачів та студентів;

2) розвиток теорії та методики використання Web-технологій у навчанні математичних дисциплін;

3) розвиток систем управління навчанням.

У розвитку засобів ІКТ на цьому етапі можна виділити ряд протиріч:

1) між потребою студентів та викладачів у автоматизації адміністрування процесу навчання із залученням Інтернет-технологій та недостатньою пропозицією з боку виробників;

2) між потенціалом використання електронних засобів навчання та нерозробленістю стандартів у сфері навчання за допомогою Інтернет і мультимедіа.

Вказані протиріччя визначили початок п'ятого етапу, на якому вони були розв'язані.

<u>П'ятий етап</u> – 1997–2003 рр. – пов'язаний із появою та розробкою систем управління навчанням (Learning Management System – LMS). Ці системи виникають наприкінці 90-х років XX століття, забезпечуючи не тільки організацію і контроль використання комп'ютерних тренінгів, але й адміністрування процесу навчання в цілому, зокрема, його традиційних форм. Для того, щоб навчальні курси різних розробників, у тому числі і навчальні ресурси, що розроблені викладачем самостійно, були сумісні з різними LMS-платформами, необхідна розробка їх стандартизації.

Як зазначалося раніше, на початку п'ятого етапу в листопаді 1997 р. Міністерством оборони США та Департаментом політики в галузі науки і технології Адміністрації Президента США було оголошено про створення ініціативи ADL (Advanced Distributed Learning), метою якої було створення стандартів у сфері навчання за допомогою Інтернет і мультимедіа [307].

Дж. А. Ітмазі (Itmazi Jamil Ahmad) [64] визначає LMS як програму, що автоматизує управління навчальними ресурсами. LMS надає можливість контролювати вхід у систему зареєстрованих користувачів, керує каталогами, відслідковує діяльність студента і його результати роботи, а також надає звіти викладачу. Особливістю LMS є можливість розробки авторського змісту.

Дж. А. Ітмазі до основних функцій LMS відносить [64]:

– управління курсами і програмами;

– надання й адміністрування реєстрації;

– відстеження реєстрації студентів, доступу і роботи;

– управління навчальною звітністю;

– контроль за процесом навчання;

– забезпечення планування курсів і адміністрації.

Кількість LMS того періоду зростала дуже швидко: 27 LMS у 1997-1998, 59 LMS у 2000 році і 70 LMS у 2002 році [64]. Про виникнення LMS та їх вагу на ринку навчальних установах можна судити за рис. 1.13.

Найбільш поширеними системами для організації та підтримки навчання станом на 2012 рік є Moodle; Edmodo; BlackBoard; SumTotal System; Skillsoft; Cornerstone; Desire2Learn; Schoology; NetDimensions; Collaborize Classrom; Interactyx; Docebo; Instructure та інші (рис. 1.14) [125].

Деякі LMS, що є розповсюдженими в ВНЗ США, розглянуто в додатку Д.

Залежно від того, яка LMS використовується у процесі навчання, змінюються і засоби комунікації. LMS надають викладачам можливість спілкуватися зі своїми студентами в чотирьох конкретних напрямах:

1) оцінювання і виставляння балів у реальному часі з динамічними, поточними розрахунками;

2) отримувати безпосередні відповіді на питання або завдання за допомогою віртуального робочого дня;

3) проводити довгострокові діалогів (наприклад, протягом усього

семестру або навчального року) за допомогою електронної пошти;

4) організовувати віртуальні дискусії у групі за допомогою звичайних дискусійних форумів або інструментів онлайн роботи.

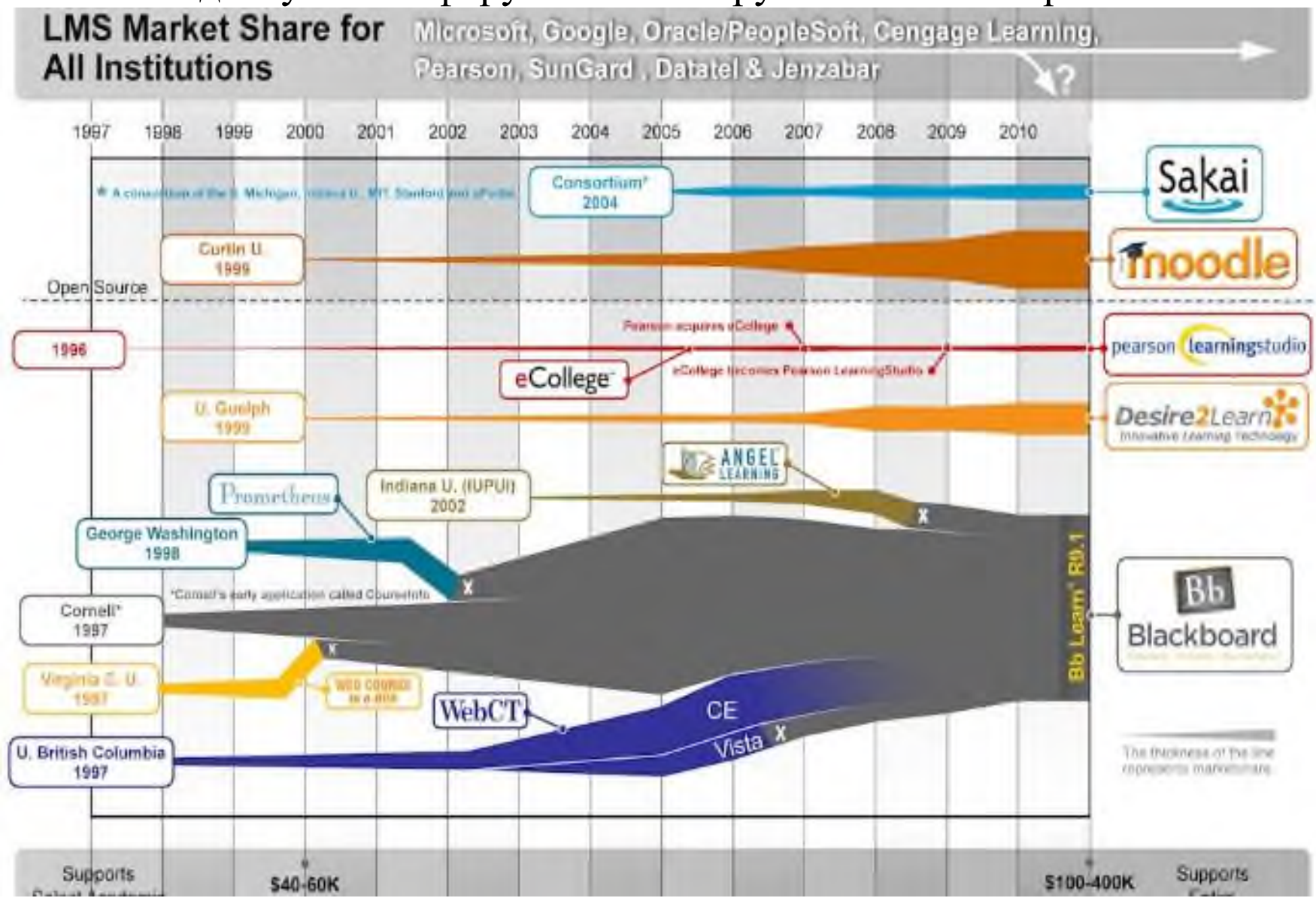

Рис. 1.13. Схема виникнення LMS та їх вага в навчальних установах [58]

Серед розмаїття сайтів, що використовують у системі вищої інженерної освіти США є сайти, наприклад, як RCampus, що пропонують безкоштовно всі інструменти для повного LMS, у тому числі динамічні звіти групи, електронну пошту та обмін повідомленнями. Gmail є одним з багатьох сайтів, де можна створити безкоштовно облікові записи електронної пошти. Віртуальні зустрічі можуть бути організовані за допомогою відкритого вихідного коду інструменту під назвою DimDim, у той час як форуми та віртуальні чати можна створити за допомогою інструменту Chatzy [100].

Таким чином, можна виокремити такі особливості п’ятого етапу розвитку теорії та методики використання ІКТ у навчанні вищої математики студентів інженерних спеціальностей у США:

1) виникнення систем управління навчанням та швидке впровадження їх у процес навчання, в зв’язку з чим дуже швидке зростання кількості навчальних курсів, що мали Web-підтримку;

2) Web-навчання вдало застосовується майже до кожного навчального курсу навчальної програми вищої математики;

3) під час проходження Web-курсів студенти можуть пройти такий

самий обсяг навчального матеріалу, як і при традиційних курсах; але при цьому, кількість студентів, які не закінчили Web-курс більше, ніж при традиційних курсах, як і при дистанційному навчанні;

4) велика кількість викладачів розробляє Web-додатки для своїх навчальних курсів, при цьому студенти, які використовують ці ресурси мають кращі результати в навчанні, ніж ті, які цього не роблять;

5) на розробку та проходження Web-курсів необхідно витратити значно більше часу, ніж на традиційні курси.

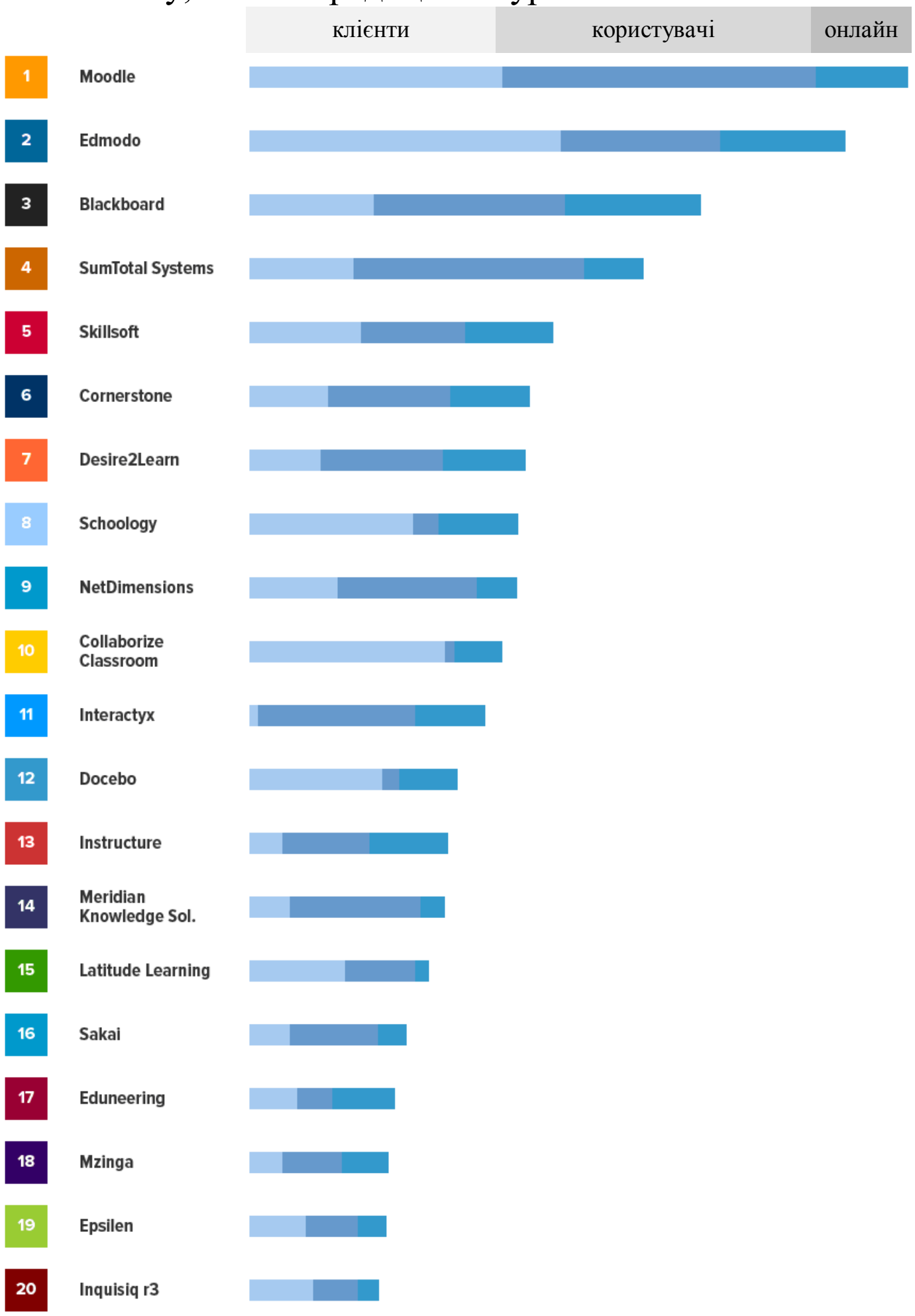

Рис. 1.14. Рейтинг найбільш популярних систем управління навчанням (2012 р.) [125]

У розвитку засобів ІКТ на цьому етапі можна виділити ряд протиріч:

1) між високою вартістю на програмне забезпечення та не достатнім фінансуванням навчальних закладів;

2) між великими потребами в ресурсах для збереження навчальних матеріалів та відсутністю можливостей у навчальних закладів встановлювати достатню кількості серверів;

3) між необхідністю постійного пересування та не можливістю забезпечення віддаленого доступу до даних, що зберігаються на персональному комп'ютері.

Вказані протиріччя визначили початок шостого етапу.

Шостий етап – з 2003 р. по теперішній час – пов'язаний із перенесенням у Web-середовище засобів підтримки математичної діяльності та становленням і розвитком хмарних технологій навчання.

За останнє десятиліття розроблено та впроваджено в процес навчання велика кількість різноманітних засобів підтримки математичної навчальної діяльності. Серед них можна виділити такі:

– засоби для передачі аудіо- та відеоданих;

– інструменти для спільної роботи над проектами;

– засоби підтримки предметно-орієнтованої практичної діяльності.

Серед *засобів для передачі аудіо- та відеоданих* найбільшою популярністю у викладачів користуються YouTube, Viddler, Voki, VoiceThread та інші [100].

У заснованій у 2005 р. Інтернет-службі *YouTube* (www.youtube.com) надаються послуги розміщення відеоматеріалів. Користувачі можуть додавати, продивлятися і коментувати ті чи інші відеозаписи. YouTube став одним із найпопулярніших місць для розміщення викладачами навчальних відеофайлів завдяки легкості завантаження та швидкому пошуку необхідних навчальних матеріалів.

У 2005 р. командою розробників таких компаній, як RackSpace, IBM і Macromedia було розроблено відео онлайн-платформу *Viddler* (http://www.viddler.com/), що є комерційним проектом та призначена для завантаження, перегляду та коментування відеоматеріалів. Готові відеоматеріали можна розміщувати на сайтах або в соціальних мережах.

*Voki* (http://www.voki.com) є вільно розповсюджуваним освітнім інструментом, що надає можливість користувачам створювати свої власні аватари, які вміють говорити голосом користувача. Символам Voki можна налаштовувати різну зовнішність, за допомогою мікрофона зробити аудіозапис, використовувати особистий номер телефону для підтримки зв'язку зі студентами, завантажувати аудіо файли. Символи Voki можна пересилати електронною поштою, ділилися в соціальних мережах і вбудовувати у Web-сайти.

Д. С. Пердью (Diana S. Perdue) [100] зазначає, що використовуючи

Voki як засіб підтримки процесу навчання, можна забезпечити різноманітне подання навчальних матеріалів. Студенти починають приділяти більше уваги навчальній дисципліні, коли чують голос викладача, який пояснює матеріал, а не коли просто читають його у друкованій формі.

Інтернет-ресурс *VoiceThread* (www.voicethread.com) призначений для проведення групових бесід та обміну думками з будь-якого міста. В VoiceThread можна проводити мультимедійні слайд-шоу, що містять зображення, документи і відео, і забезпечує можливість користувачам переміщуватися по сторінках і залишати коментарі одним із способів: голосові повідомлення (запис з мікрофону або телефону), текстові, аудіо-файли або відео (через Web-камеру) [100].

Використовуючи *інструменти для спільної роботи над проектами*, викладачі можуть підтримувати процес навчання на відстані, поза аудиторією.

Використовуючи документи Google *GoogleDocs* (http://www.google.com/ google-d-s/intl/ru/tour1.html), викладач може створювати різні документів, працювати з ними разом зі студентами в режимі реального часу і зберігати документи та інші файли в Інтернеті, надаючи доступ до цих документів іншим. Крім того надається можливість отримати доступ до своїх документів і файлів з будь-якого комп'ютера в будь-якому місці, якщо він підключений до Інтернету.

У статті Д. С. Пердью [100] зазначено, що GoogleDocs є кращим інструментом для роботи з невеликою групою студентів над спільним навчальним проектом. Використання GoogleDocs надає можливість забезпечити динамічну, постійно доступну платформу, де студенти і викладач можуть досягти освітніх цілей.

Подібні можливості надають сервіси Microsoft *Office Web Apps* (http://office.microsoft.com/): перегляд і спільне використання файлів не залежно від місця розташування користувачів; можливість одночасної роботи з іншими людьми над документами на різних платформах, пристроях і в різних версіях Office і навіть якщо у них не встановлений Office; збереження документів в Інтернеті для спільного використання.

У 2008 році розроблено хмарне сховище даних *DropBox* (*https://www.dropbox.com/*), використання якого надає можливість користувачам зберігати свої дані на серверах у хмарі і ділитися ними з іншими користувачами в Інтернеті. Робота побудована на синхронізації даних.

Комплексним рішенням для Інтернет-уроків є платформи для онлайн навчання та проведення Web-конференцій. *Elluminate* (*http://www.elluminate.com/*), *iLinc* (*http://www.ilinc.com/*), *WiZiQ*

(*http://www.wiziq.com/*) є прикладами ресурсів для відео-конференцій та створення онлайн-класів. За допомогою цих ресурсів можна організовувати проведення вебінарів (Інтернет-конференцій) в онлайн режимі або у записі. Опублікованими курсами в будь-якому форматі в бібліотеці можна поділитися у віртуальному класі або поширити їх безпосередньо серед студентів. Доступно для завантаження найбільш часто використовувані типи файлів: PDF, Word, PowerPoint, Excel документи, відео і аудіо файлів. Існує можливість користуватися дошкою. Можна завантажувати навчальні матеріали перед заняттям або під час заняття. Завантажені файли розміщуються в хмарі [92, 86].

Існує велика кількість Інтернет-ресурсів для підтримки інтерактивної позааудиторної взаємодії викладачів та студентів: платформа Piazza, Skype, Google+ та багато інших.

Використання платформи *Piazza* (*https://piazza.com/*) надає можливість підтримувати навчальну діяльність студентів за межами аудиторії. Так, наприклад, на платформі Piazza класи з математичної тематики мають в університеті Колумбії – «Числення однієї змінної», в університеті штату Іллінойс в Урбана-Шампейн – «Алгебра», в Корнеллі – «Числення багатьох змінних», в університеті Торонто – «Вища математика», в університеті Британської Колумбії – «Дійсні числа», в університеті штату Нью-Мексико – «Тригонометрія» та «Елементарне числення» та інші [103].

Із використанням *Skype* (*www.skype.com*) викладач може не лише підтримувати спілкування зі студентами, а й проводити конференцій та відео-заняття. За допомогою програми MP3 Skype Recorder можна зробити запис розмови в Skype в аудіоформаті та розмістити його як навчальний матеріал в Інтернеті.

Комп'ютерна підтримка *предметно-орієнтованої практичної діяльності* полягає у наданні користувачу (студенту або викладачу) набору засобів і інструментів, що автоматизують і надають можливість перевірити процес розв'язування практичної задачі. Така система повинна бути забезпеченою повним комплектом методичної підтримки. Для вищої математики – це підручник, задачник, довідник, робочий зошит студента, збірник контрольних робіт та тестів, методичні рекомендації викладача тощо.

Основною метою розробки та впровадження у навчальний процес педагогічних програмних середовищ підтримки практичної діяльності студентів є якісне підвищення ефективності проведення занять у лабораторно-практичній частині навчального курсу, для самостійної роботи дослідницького характеру [297; 298].

Проблема вибору систем комп'ютерної математики (СКМ) та підтримки великої інсталяційної бази розв'язана через застосування мережних технологій, коли користувач за допомогою спеціалізованого клієнтського програмного забезпечення звертається до серверної частини СКМ, де виконуються команди користувача та повертається результат до клієнтського програмного забезпечення. Такі послуги надаються, зокрема, *MATLAB Web Server, webMathematica* та *wxMaxima*. І хоча далеко не у всі СКМ включені вбудовані мережні засоби, для тих з них, в яких поряд з візуальним підтримується командний інтерфейс, можливе створення мережної надбудови [294].

Новим перспективним напрямом розвитку СКМ є мобільні математичні середовища; перші представники таких систем з'явились лише на початку XXI століття. Відбувається інтеграція СКМ у єдине мережне середовище, тобто перенесення прикладного ПЗ (навіть «робочих столів») у Web-середовища. Важливим є те, що застосування мобільних Web-середовищ у навчальному процесі надає можливість інтегрувати аудиторну й позааудиторну роботу у безперервний навчальний процес.

За допомогою Web-СКМ можна:

1) виконувати будь-які обчислення, як аналітичні (дії з алгебраїчними виразами, розв'язування рівнянь, диференціювання, інтегрування тощо), так і чисельні (точні – з будь-якою розрядністю, наближені – з будь-якою, наперед заданою точністю);

2) подавати результати обчислень у зручній для сприйняття формі, будувати дво- та тривимірні графіки кривих та поверхонь, гістограми та будь-які інші зображення (в тому числі анімаційні);

3) поєднувати обчислення, текст та графіку на робочих листах з можливістю їх друкування, оприлюднення в мережі та спільної роботи над ними;

4) створювати за допомогою вбудованих мов програмування моделей для виконання навчальних досліджень.

До основних напрямів застосування Web-СКМ у процесі навчання вищої математики відноситься [294]:

– графічна інтерпретація математичних моделей та теоретичних понять;

– автоматизація рутинних обчислень;

– підтримка самостійної роботи;

– організація математичних досліджень.

Яскравим представником Web-СКМ є *Sage*, розробником якої у 2005 році є математик Університету Вашингтона У. Стейн (William Stein). На основі Sage побудовано мобільне математичне середовище, що виступає в якості інтегратора різних математичних пакетів із забезпеченням спільного Web-інтерфейсу.

У 2007 році було запущено проект *Wolfram Demonstrations Project* (*http://demonstrations.wolfram.com/*), що містить чимало демонстрації з вищої математики та має відкритий код доступу. В демонстраціях можна надавати різні значення параметрам і проводити дослідження об'єктів [21, 297].

Проект *Wolfram|Alpha,* що почав розроблятися С. Вольфрамом (Stephen Wolfram) у 2009 році на мові Mathematica, являє собою базу даних та набір обчислювальних алгоритмів з різних дисциплін, зокрема з математики. Сайт є зручним у використанні. Зайшовши на головну сторінку за посиланням www.wolframalpha.com та ввівши питання, після натискання символу «=» з'явиться не список посилань, як у стандартних пошукових системах, а розгорнута відповідь на запитання, що генерується, спираючись на власну базу знань та алгоритмів [21, 286-287].

8 лютого 2012 відбувся випуск нової версії Wolfram|Alpha Pro, ключовою особливістю якої є можливість завантаження багатьох типів файлів і даних для автоматичного аналізу, враховуючи первинні табличні дані, зображення, аудіо, XML, а також десятки спеціалізованих наукових, медичних та математичних форматів. Серед інших функціональних можливостей – наявність розширеної клавіатури, динамічність з CDF (формат обчислюваних документів), завантаження даних і можливість індивідуального налаштування та збереження графічних і табличних результатів [90].

Освітній ресурс *GeoGebra* (2010) є простим у використанні аплетом, що може бути використаний безпосередньо у Web-браузері або завантажений на комп'ютер. Користувачам надаються порожні двовимірні площині, на яких вони можуть проводити побудови та виконувати маніпуляції з ними. У відкритому полі користувачі можуть задавати рівняння. Ці динамічні особливості GeoGebra роблять її універсальним і корисним інструментом для вивчення широкого спектру тем двовимірної плоскої геометрії [21, 296]. GeoGebra є ПЗ з вільним відкритим вихідним кодом. На заняттях з вищої математики засобами GeoGebra зручно ілюструвати основні математичні поняття і конструкції, проводити обчислення [4].

Розроблений у 2011 році КПК TI-Nspire нагадує графічний калькулятор, наприклад, TI-84, але має великі екран і клавіатуру. По суті, TI-Nspire є ручним комп'ютером з управляючим пакетом під назвою TI-Nspire, що може бути доступною у вигляді пакету програмного забезпечення для ОС Windows. TI-Nspire можна використовувати: як звичайний калькулятор (у тому числі з символьними алгебраїчними маніпуляціями); для геометричних побудов; для роботи із списками та таблицями; для роботи та обробки статистичних даних [92, 63].

Концепція навчання впродовж всього життя стає основою системи вищої освіти. У зв'язку зі стрімким розвитком технологій за час навчання у ВНЗ підготувати інженера, готового до роботи без отримання додаткових знань, неможливо. Для кожної людини безперервна освіта має бути основою для формування і задоволення її пізнавальних запитів, розвитку здібностей у різних навчальних закладах за допомогою різних форм навчання, а також шляхом самоосвіти і самовиховання.

Основою безперервної освіти стають масові відкриті дистанційні курси (МВДК: Massive open online course – MOOC), спрямовані на велику кількість учасників та відкритий доступ через мережу Інтернет.

Термін «масовий відкритий дистанційний курс» запропонували Б. Александр (Bryan Alexander) та Д. Кормьє (Dave Cormier) під час роботи над курсом «Connectivism and Connective Knowledge» («Коннективізм та сполучні знання»). Даний курс був оснований на курсі, розробниками якого у 2008 році були Дж. Сименс (George Siemens) і С. Даунс (Stephen Downes) та базувався на ідеї коннективізму [77].

В основу масових відкритих дистанційних курсів, що виникли у 2008 році, було покладено філософію коннективізму (Connectivist MOOC – cMOOC), що підтримувала ідею безперервної освіти на основі мережі Інтернет. Особливістю такого курсу є відкритість у спілкуванні, можливість проводити дискусії та діалоги. Такий курс – це соціальне середовище, де мета навчання визначається тільки самим студентом.

Групою дослідників з США у звіті [6] було зроблено аналіз поняття «масові відкриті дистанційні курси» (табл. 1.1).

До *традиційних* віднесено курси, що не використовують Інтернет-технологій – зміст курсу доступний у письмовій або усній формі.

До курсів з *Інтернет-підтримкою* було віднесено такі курси, в яких від 1 до 29 % змісту доступно онлайн. Це, по суті, традиційні курси з використання систем управління навчанням (LMS) або Web-сторінок з навчальними планами і завданнями.

До *дистанційних* було віднесено такі курси, в яких щонайменше 80 % змісту курсу розміщено в Інтернеті. Як правило, у таких курсах не буває аудиторних занять.

У *змішаних* (*комбінованих*) курсах поєднано дистанційне та традиційне навчання – від 30 до 80 % змісту курсу доступно в Інтернеті. Як правило, використовуються форуми; зменшена кількість аудиторних занять.

У звіті проаналізовано стан впровадження MOOC у ВНЗ США станом на 2012 рік [6]:

*Таблиця* 1.1

**Аналіз проведення дистанційних курсів у ВНЗ США**

| **Частка навчального змісту курсу, доступного через мережу Інтернет** | **Тип курсу** | **Опис курсу** |
|---|---|---|
| 0 % | Традиційні | Курс, де немає технологій, що використовуються онлайн – зміст доступний у письмовій або усній формі. |
| 1-29 % | З Інтернет підтримкою | У курсі використовуються Інтернет-технології для полегшення доступу до змісту курсу. Це, по суті, традиційний курс. Можливе використання систем управління навчанням (LMS) або Web-сторінок з навчальними планами і завданнями. |
| 30-79 % | Змішані (комбіновані) | У курсі поєднано дистанційне та традиційне навчання. Значна частка змісту курсу доступна через Інтернет. Як правило, використовуються форуми; знижена кількість аудиторних занять. |
| > 80 % | Дистанційні | Курси, в яких більшість або весь навчальний матеріал доступний через Інтернет. Як правило не буває аудиторних занять. |

1. Зазначається, що лише в незначній частині ВНЗ проводяться МООС, у більшості ВНЗ МООС у стадії планування:

– зараз проводять МООС у 2,6 % ВНЗ; у 9,4 % МООС знаходяться в стадії планування;

– у більшості ВНЗ (55,4 %) ще не визначилися з МООС, у той час як у менше ніж у третині ВНЗ (32,7 %) говорять, що у них немає планів стосовно проведення МООС; у провідних ВНЗ закріпилася думка, що МООС – це надійний спосіб для проведення онлайн-курсів і є важливим засобом для установ, щоб розвивати методику Інтернет-навчання.

2. У 2002 році у менше ніж у половині всіх ВНЗ повідомлялося, що Інтернет-освіта має вирішальне значення в їх довгострокових планах. У 2012 році це число зросло до 70 %: у провідних ВНЗ говорять, що МООС мають вирішальне значення в їх довгострокових планах, зараз таких ВНЗ 69,1 % – найвищий відсоток протягом десятирічного періоду 2002-2012 рр.; частка ВНЗ, в яких МООС не включено в довгострокові плани впала до рекордно низького рівня – 11,2 %.

3. У 2003 році у 57,2 % провідних ВНЗ оцінили результати навчання в МООС вище ніж при традиційному навчанні. Це число в 2012 р. становить 77,0 %.

4. Викладачі провідних ВНЗ висловлюють свою стурбованість з приводу існуючих перешкод для проведення МООС. Кількість таких ВНЗ збільшилася з 80 % у 2007 році до 88,8 % у 2012 році.

За останні роки проведено кілька десятків сМООС – МООС, побудованих на коннективістському підході до навчання. сМООС ґрунтуються на активній участі сотень і тисяч студентів, які самі організовують свою участь у курсі відповідно до особистих цілей навчання, попередніх знань і навичок, а також спільних інтересів [101]. Велика частина діяльності студентів у сМООС відбувається за межами LMS [77], в інших вузлах мережі, наприклад, в особистих блогах, особистих портфоліо, Web-сайтах, Twitter, YouTube, віртуальних світах та інше.

Починаючи з осені 2011 року провідні ВНЗ США починають розробляти та впроваджувати МООС. До програми залучаються університети Стенфорда, Пенсільванії, Мічигану та Прінстону [101].

На початку 2012 року з’являються нові курси – хМООС. Такі курси починає впроваджувати Массачусетський технологічний інститут на платформі MITx на основі розміщених у мережі відкритих освітніх ресурсів. хМООС базуються на когнітивно-біхевіоріському підході. Основою таких курсів є традиційні університетські навчальні курси, в яких мету навчання визначає викладач. Такі курси мають високу ступінь фінансування та автоматизовану систему контролю успішності студентів. хМООС Стенфордського університету та МТІ мають загальну фінансову підтримку понад 80 млн. доларів.

У 2012 – 2013 роках у ВНЗ США пропонувалися такі курси:

– на базі платформи MITx: аеродинаміка літальних апаратів; вступ до аеродинаміки;

– на базі платформи UC BerkeleyX: вступ до статистики; теорія ймовірностей; описова статистика; квантова механіка і квантові обчислення;

– на базі платформи HarvardX: здоров’я в цифрах – чисельні методи в клінічних та громадських досліджень у галузі охорони здоров’я та інші.

На рис. 1.15 показано схему виникнення та розподіл МООС у провідних ВНЗ США, розроблену Дж. Сименсом [111].

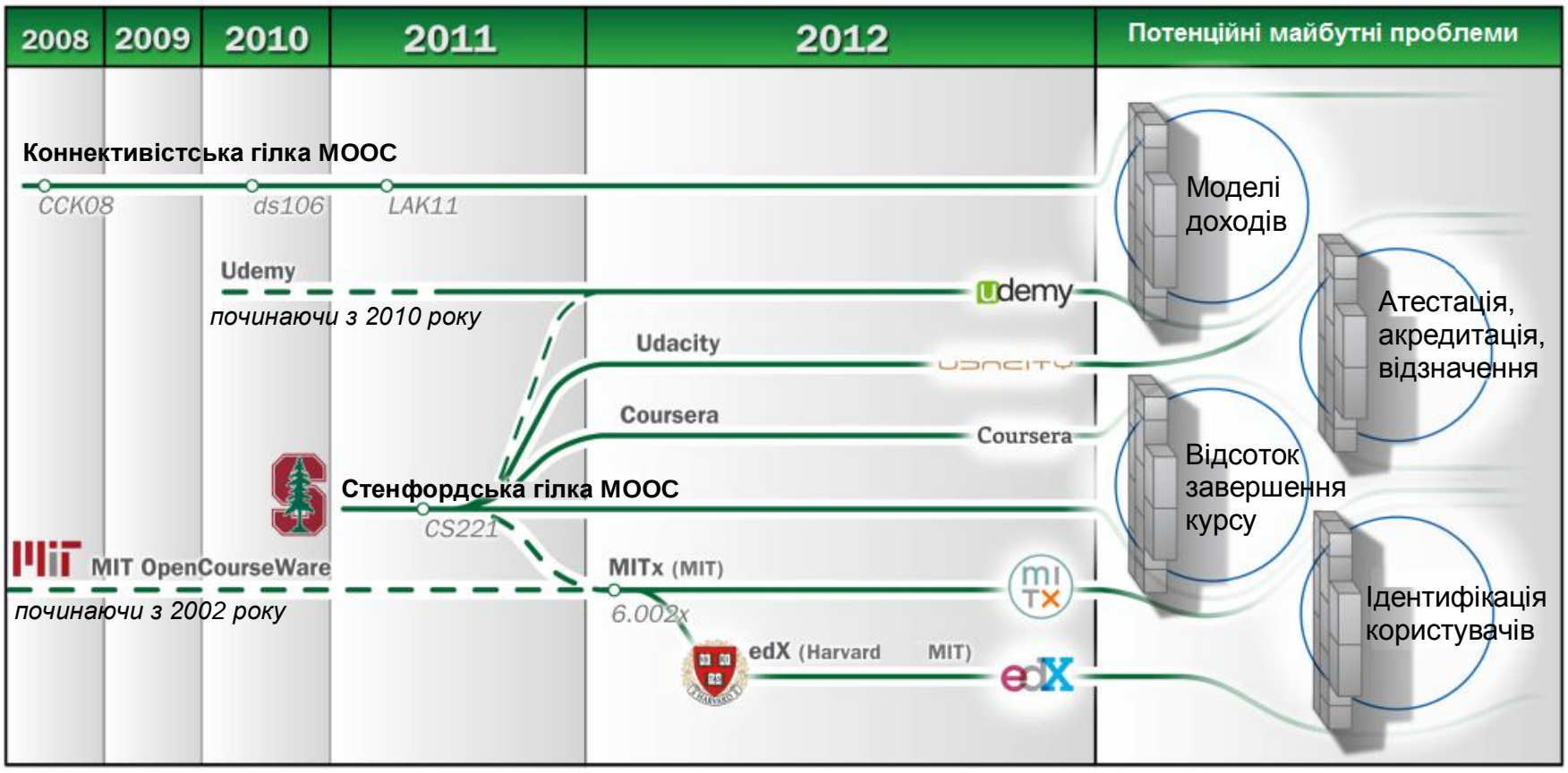

Рис. 1.15. Схема виникнення МООС за Дж. Сименсом

Таким чином, згідно з [275] МООС є:

– навчальний курс з дисципліни, що розміщений у вільно поширюваній системі підтримки навчання з використанням хмарних технологій;

– доступ до навчальних матеріалів є вільним і не залежить від місця розташування викладача і студента та кількості студентів, що здобувають знання;

– навчальні матеріали можуть бути опрацьовані на довільному апаратному, зокрема, мобільному, інформаційно-комунікаційному засобі;

– навчальні матеріали курсу добре структуровані, мають гіперпосилання, містять відео та аудіо додатки;

– система тестування та оцінювання знань є простою у використанні.

До особливостей МООС відносять [12]:

– для роботи з матеріалами курсу необхідно лише мати доступ до мережі Інтернет;

– при створенні курсу можна використовувати усі можливі технології, залежно від можливостей платформи дистанційного навчання чи цільової групи;

– навчальні матеріали, розміщені в курсі є загальними для всієї цільової групи і можуть бути швидко змінені, адаптовані та доповнені;

– процес навчання відбувається у зручній для студентів та викладачів час та не залежить від географічного розташування студента та викладача;

– процес навчання, організований за допомогою МООС відповідає сучасній парадигмі освіти «навчання протягом усього життя»;

– процес навчання здійснюється за допомогою неформальних знань, що виникають у курсі в наслідок обміну навчальними відомостями між його учасниками;

– після закінчення студентами навчального курсу усі відомості залишаються в мережі – тобто створені навчальні мережі є більш стійкими;

– можливість створення власного персонального навчального середовища і власної персональної навчальної мережі за допомогою інших учасників курсу.

На думку В. М. Кухаренка, відкритий дистанційний курс базується на чотирьох основних видах діяльності [223]:

– *співпраця*: у курсі даються посилання на різні навчальні матеріали, що необхідні для опрацювання та подальшого обговорення. Серед наданих навчальних матеріалів, студент вибирає тільки ті, що є необхідними йому в даний момент часу. Процес навчання у співпраці надає студенту можливість отримати узагальнене уявлення про відповідний розділ курсу;

– *ремікс*: навчальні матеріали курсу органічно поєднані між собою, тому студент, після опрацювання навчальних матеріалів, повинен обговорити отримані знання на вебінарах та на форумі, поділитися своїм контентом з іншими людьми;

– *перепрофілювання*: основним завданням курсу є допомога учаснику набути знання, необхідні для нього в подальшій діяльності. І це є найскладнішою частиною процесу навчання, оскільки студент починає навчання не з нуля (в курсі використовується термін «перепрофілювання» замість «створення»), а аналізує отримані знання та синтезує своє розуміння навчального матеріалу. Відкритий курс вчить, як читати, розуміти і працювати зі змістом інших людей і як створити своє власне нове розуміння. У курсі, як правило, надаються інструменти, що можна використовувати для створення власного контенту;

– *повідомлення*: завдання викладача полягає в організації спільної роботи з іншими студентами в даному процесі навчання. Але студент може працювати самостійно, не обговорюючи проблеми курсу в групі.

Організація процесу навчання вищої математики за допомогою масового відкритого дистанційного курсу створює умови для надання масового характеру індивідуалізованому навчанню.

Таким чином, можна виокремити такі особливості шостого етапу розвитку теорії та методики використання ІКТ у навчанні вищої математики студентів інженерних спеціальностей у США:

1) перенесення математичної діяльності викладачів та студентів у мережне середовище;

2) застосування засобів хмарних технологій для підтримки навчальної діяльності;

3) становлення Web-орієнтованих методичних систем навчання вищої математики;

4) розвиток масових відкритих дистанційних курсів з математичних дисциплін.

Можна виділити ряд проблем, що виникають при впровадженні засобів ІКТ шостого періоду:

1) наявність потрібного обладнання та програмного забезпечення у всіх учасників навчального процесу в аудиторії та за її межами;

2) необхідність постійного вивчення нових засобів та їх можливостей викладачами та студентами, що займає достатню кількість часу;

3) необхідність достатнього фінансування на закупівлю ліцензій на придбання комерційних засобів підтримки навчальної діяльності;

4) можливість самостійно розібратися студентам з постійно змінним інтерфейсом деяких Web-засобів.

**Висновки до розділу 1**

1. Сучасна інженерна освіта США має такі основні особливості: недержавна система акредитації; відсутність державних галузевих стандартів; математизація та комп'ютеризація загальноінженерних та спеціальних дисциплін; прикладна спрямованість навчання вищої математики; широке використання засобів ІКТ у навчанні вищої математики.

2. Аналіз фундаментальної підготовки студентів інженерних спеціальностей у ВНЗ США показав, що, незважаючи на недержавну систему акредитації, відсутність державних галузевих стандартів та традиційне різноманіття пропонованих математичних курсів (як обов'язкових, так і факультативних), навчання вищої математики майбутніх інженерів у США здійснюється за схожими навчальними програмами.

3. Системи підготовки інженерів у ВНЗ США та у ВНЗ України мають такі спільні риси: а) високий рівень математизації та комп'ютеризації загальноінженерних та спеціальних дисциплін; б) навчання фундаментальних дисциплін, зокрема вищої математики (лінійної алгебри, математичного аналізу та диференціальних рівнянь), відбувається переважно на молодших курсах, загальнопрофесійних дисциплін – на середніх курсах та спеціальних професійних – на старших; в) зміст навчання вищої математики є професійно орієнтованим та диференційованим за рівнями початкової підготовки студентів; г) у навчанні вищої математики широко використовуються засоби ІКТ. Таким чином, вивчення досвіду навчання вищої математики студентів інженерних спеціальностей у США може бути корисним для визначення напрямів модернізації вітчизняної системи інженерної освіти.

4. Проаналізувавши джерела з проблеми дослідження, можна стверджувати, що на кожному етапі розвитку теорії та методики використання ІКТ у навчанні вищої математики студентів інженерних спеціальностей у США адекватно розвивається і педагогічна наука, що враховує науково-технічні досягненням свого часу. Потреби педагогічної практики визначають еволюцію засобів навчання, їх розвиток спрямовується на задоволення цих потреб. Вплив наукових і технічних досягнень людства на зміст, структуру і організацію процесу навчання опосередковується і має матеріальний вираз в інформаційно-комунікаційних засобах навчання як знаряддях навчальної діяльності. Оскільки сучасний стан освіти в світі, й зокрема в Україні, не відповідає зростаючим потребам суспільства, вирішити цю проблему можна із залученням в освітній процес ІКТ, що ґрунтуються на сучасних досягненнях педагогічних технологій.

5. Аналіз розвитку засобів ІКТ надав можливість виокремити такі етапи розвитку теорії та методики використання ІКТ у навчанні вищої математики студентів інженерних спеціальностей у США:

– перший етап – 1965–1973 рр. – пов’язаний із появою достатньої кількості комп’ютерних засобів різного рівня, оснащених мовами високого рівня та специфікою апаратного забезпечення ІКТ (використання мейнфреймів з обмеженим мережним доступом). Характерні риси етапу: діалоговий режим роботи з навчальними програмами; поява перших систем підтримки математичної діяльності без програмування мовами загального призначення; розмаїття апаратного та програмного забезпечення; домінування біхевіоризму в обґрунтуванні використання ІКТ та розробці навчальних програм; уведення програмування в курси вищої математики;

– другий етап – 1973–1981 рр. – пов’язаний з поширенням в університетах США мережної операційної системи UNIX, використанням міні- та мікрокомп’ютерних систем. Характерні риси етапу: перехід від використання мов програмування у навчанні до використання математичних бібліотек, систем комп’ютерної математики та мов високого рівня; застосування комп’ютерної графіки у навчальних програмах; поява нових класів навчальних програм – навчальних ігор, систем динамічної геометрії та електронних таблиць; поява та поширення комп’ютерних мереж, що об’єднували викладачів та студентів; розвиток засобів ІКТ навчання вищої математики – графічних та символьних калькуляторів;

– третій етап – 1981–1989 рр. – пов'язаний із поширенням персональних комп'ютерів. Характерні риси етапу: широке використання математичних бібліотек, систем комп'ютерної математики та проблемно орієнтованих мов; широке впровадження персональних та персоналізованих засобів ІКТ у навчання математичних дисциплін; використання ІКТ загального призначення (текстові редактори, електронні таблиці, бази даних тощо) для підтримки навчання математичних дисциплін;

– четвертий етап – 1989–1997 рр. – пов'язаний із створенням World Wide Web та використанням технологій Web 1.0. Характерні риси етапу: широке впровадження комп'ютерних мереж, що об'єднували викладачів та студентів; розвиток теорії та методики використання Web-технологій у навчанні математичних дисциплін; розвиток систем управління навчанням;

– п'ятий етап – 1997–2003 рр. – пов'язаний із появою та розробкою систем управління навчанням. Характерні риси етапу: виникнення систем управління навчанням та швидке впровадження їх у процес навчання, в зв'язку з чим дуже швидке зростання кількості навчальних курсів, що мали Web-підтримку; Web-навчання вдало застосовується майже до кожного навчального курсу навчальної програми вищої математики; під час проходження Web-курсів студенти можуть пройти такий самий обсяг навчального матеріалу, як і при традиційних курсах, але при цьому, кількість студентів, які не закінчили Web-курс більше, ніж при традиційних курсах, як і при дистанційному навчанні; велика кількість викладачів розробляє Web-додатки для своїх навчальних курсів, при цьому студенти, які використовують ці ресурси мають кращі результати в навчанні, ніж ті, які цього не роблять; на розробку та проходження Web-курсів необхідно витратити значно більше часу, ніж на традиційні курси;

– шостий етап – з 2003 р. по теперішній час – пов'язаний із перенесенням у Web-середовище засобів підтримки математичної діяльності та становленням і розвитком хмарних технологій навчання. Характерні риси етапу: перенесення математичної діяльності викладачів та студентів у мережне середовище; застосування засобів хмарних технологій для підтримки навчальної діяльності; становлення Web-орієнтованих методичних систем навчання вищої математики; розвиток масових відкритих дистанційних курсів з математичних дисциплін.

6. Систематизовано напрями, функції та засоби використання ІКТ у навчанні вищої математики студентів інженерних спеціальностей США з позицій американської педагогічної думки другої половини ХХ століття – початку ХХІ століття. Проведений аналіз етапів розвитку теорії та методики використання ІКТ у навчанні вищої математики студентів інженерних спеціальностей у США надав можливість зробити такий висновок: поява нового типу апаратних чи програмних засобів ІКТ впливає на процес організації навчання вищої математики і на сучасному етапі створює умови для реалізації Web-орієнтованого навчання вищої математики.

# РОЗДІЛ 2
# МЕТОДИЧНІ ОСНОВИ ВИКОРИСТАННЯ ІНФОРМАЦІЙНО-КОМУНІКАЦІЙНИХ ТЕХНОЛОГІЙ У НАВЧАННІ ВИЩОЇ МАТЕМАТИКИ СТУДЕНТІВ ІНЖЕНЕРНИХ СПЕЦІАЛЬНОСТЕЙ У СПОЛУЧЕНИХ ШТАТАХ АМЕРИКИ

## 2.1 Використання інформаційно-комунікаційних технологій навчання вищої математики у інженерній освіті в Сполучених Штатах Америки

Один з найкращих інженерних ВНЗ США за даними британського видання Times Higher Education (THE) [140] – Массачусетський технологічний інститут – відрізняється не лише рівнем фундаментальної підготовки, а й підтримкою засобами ІКТ традиційного процесу навчання. Так, у 1966 р. лабораторія комп'ютерних наук та штучного інтелекту МТІ (MIT Computer Science and Artificial Intelligence Laboratory) налічувала близько 100 комп'ютерних терміналів, розташованих як у студмістечку, так і у приватних будинках, що одночасно могли бути використаними 30 студентами для розв'язання навчальних задач, моделювання процесів та явищ, а також спільної роботи у мережі. У 1968 р. у рамках проекту **MAC** (Mathematics and Computation, пізніше backronymed to Multiple Access Computer, Machine Aided Cognitions, або Man and Computer) лабораторії була створена перша у світі система комп'ютерної алгебри Maxima (MACSYMA – Project MAC's SYmbolic MAnipulator) – родоначальник усіх сучасних систем комп'ютерної математики [284, 4-5].

Визначним етапом розвитку засобів ІКТ навчання у МТІ став спільний із DEC (Digital Equipment Corporation) та IBM (International Business Machines Corporation) восьмирічний (1983-1991) проект **Athena**, спрямований на інтеграцію різних засобів ІКТ у Массачусетському технологічному інституті і за його межами з метою створення освітнього середовища [88].

**Athena** – академічне обчислювальне середовище МТІ, що обслуговує обчислювальні кластери (лабораторії), особисті робочі місця, сервери віддаленого доступу і особисті машини по всьому університетському містечку.

Основними цілями проекту були: 1) створення комп'ютерно-орієнтованих засобів навчання, придатних для роботи у різних навчальних середовищах; 2) розробка методик впровадження ІКТ у навчальний процес; 3) створення мобільного ІКТ середовища.

Проект стимулював створення таких додатків, як обмін миттєвими повідомленнями, активні каталоги і X Window System. Але реальна мета була зробити потужну програму для роботи студентів, викладачів і співробітників, що включала б майже всі потреби, від відправки електронної пошти та написання листа до аналізу даних і створення нових додатків. Отриманні можливості від Athena були оцінені всіма членами спільноти MIT. Хоча потужні бездротові мережі забезпечують майже універсальний Wi-Fi доступ по всій території студмістечка, студенти надають перевагу використанню кластерів Athena [82].

Створені у рамках проекту технології широко використовуються сьогодні не лише у навчанні: так, X Window System є основою переважної більшості графічних інтерфейсів сучасних операційних систем [252]. Athena є основою комп'ютерно-орієнтованих засобів навчання у MTI, надаючи її користувачам такі переваги: легкість адміністрування та використання, стійкість до збоїв, швидкий та повсюдний доступ до навчальних матеріалів і засобів навчання.

Основна концепція проекту Athena – інтеграція та конвергенція засобів ІКТ навчання з навчальними планами MTI традиційної освіти (з традиційними технологіями). У рамках проекту вперше було виконано інтеграцію стороннього програмного забезпечення (зокрема, Mathematica, MATLAB та Maple) в обчислювальне середовище, що могло бути використано у великій кількості навчальних курсів з метою надання студентам та викладачам вільного доступу до нього, забезпечуючи легкий доступ до файлової системи AFS для особистого і групового зберігання файлів. Поточною версією Athena є Debathena.

*2.1.1 Проект Массачусетського технологічного інституту OpenCourseWare*. Надання вільного доступу до навчальних матеріалів, створених провідними фахівцями MTI, є головною метою проекту MTI **OpenCourseWare** (OCW). Опубліковані на сайті проекту [43] матеріали включають плани курсів, конспекти лекцій, домашні завдання, екзаменаційні питання, відеозаписи лекцій тощо (рис. 2.1).

MTI OCW виник із ряду ініціатив, що проводилися радою інституту по освітнім технологіям (MIT Council on Educational Technology). Метою ради було підвищення якості освіти MTI через належне застосування технологій в студмістечках для підтримки комбінованого навчання. Функції ради включали в себе координацію розподілу централізованого фінансування для ініціативи технологій в освіту, моніторинг ефективності проведених програм та визначення пріоритетів для інвестицій в нові освітні технології. У 1999 році до складу ради входила група розробників стратегії під керівництвом проректора Р. Брауна (Robert A. Brown), яка розглядала можливості використання Інтернет в освітніх цілях. У цей час дистанційна освіта вважалася перспективним бізнесом і компанії розглядали можливості співпраці з університетами, в тому числі і з MTI, щоб зайняти місце на даному ринку. Для оцінки можливостей інституту в мінливій обстановці Інтернет-освіти, група з розробки стратегії вирішила звернутися до консалтингового агентства McKinsey & Company [2].

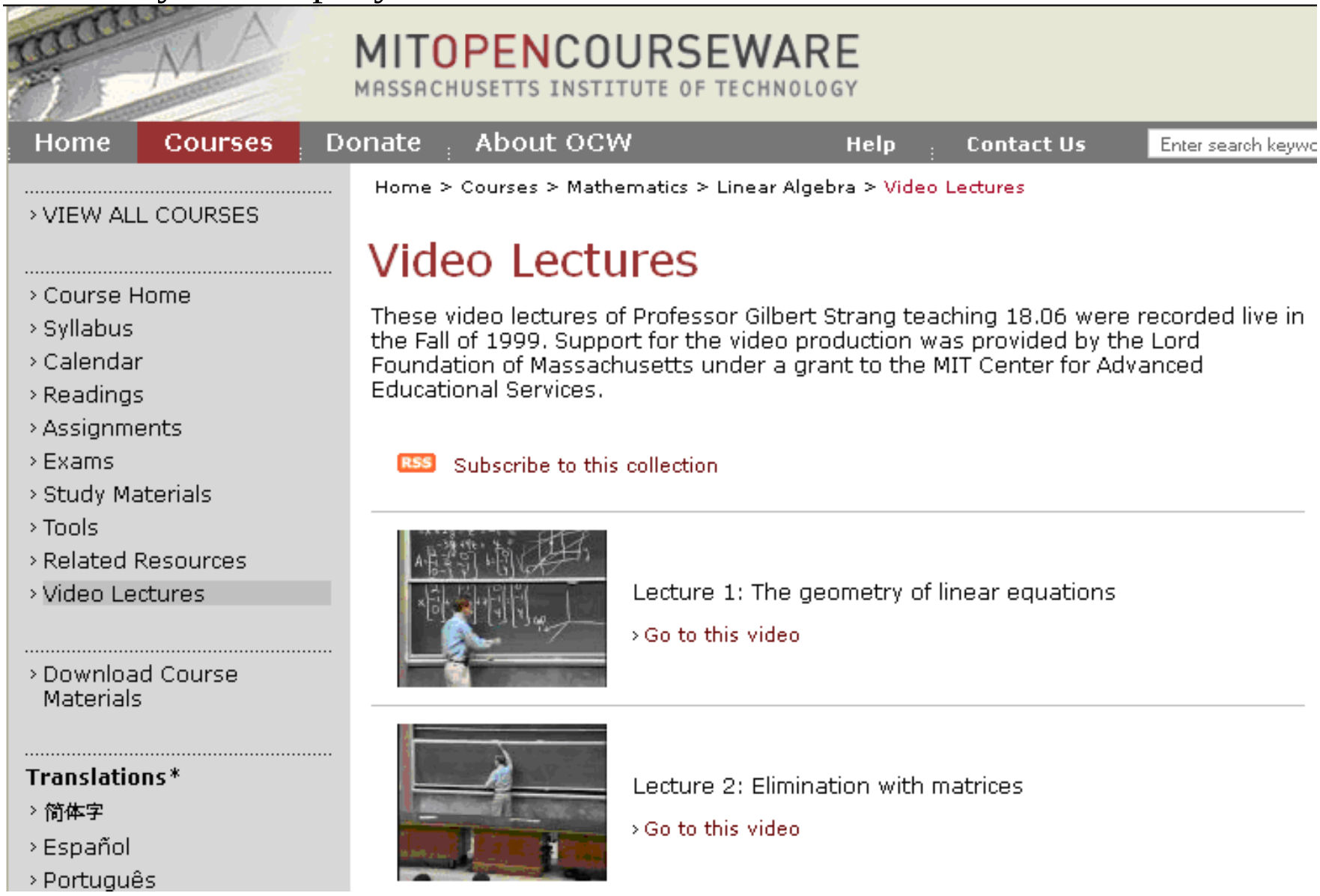

Рис. 2.1. Сайт відкритого курсу лінійної алгебри на сайті MIT OCW

Проведене MTI та McKinsey & Company опитування співробітників інституту показало, що вони не є прихильниками роботи на масовому ринку безперервної (дистанційної) освіти. Для даної аудиторії в 2000 році групою був розроблений план проекту Knowledge Updates MIT, що передбачав створення міні-курсів за новими напрямами в технічних і міждисциплінарних областях. Однак, дослідження ринку показали, що для того, щоб окупити проект, потрібен довгий термін і більша, в порівнянні із спочатку запланованою, аудиторія, адаптувати курси для якої більшість викладачів інституту не могло. Це поставило під сумнів комерційний потенціал проекту. Члени групи також побоювалися, що, навіть при успішній реалізації, проект не зможе істотно підвищити вплив інституту.

Опитування співробітників показало, що близько 20 % викладачів інституту вже створили сайти для своїх курсів за власною ініціативою, і це спонукало групу розробників проекту розглянути можливість безкоштовного розповсюдження освітніх матеріалів (OCW MIT). Різницю між OCW і тим, що вже було зроблено викладачами вбачали в

тому, що [74]:

– МТІ мав намір систематично створювати Web-сайти для всіх навчальних курсів інституту;

– план централізованої організаційної підтримки мав допомогти розробляти Web-сайти, не вимагаючи великих зусиль від викладачів (командою, що займалася наповненням сайту, передбачалося спланувати роботу таким чином, щоб викладачі витрачали мінімум часу на надання всіх необхідних навчальних матеріалів з курсу);

– створення єдиної пошукової системи, що надасть змогу охопити всі курси;

– забезпечення послідовного викладення матеріалів курсів;

– проект надасть можливість вільного і відкритого повторного використання OCW матеріалів у всіх некомерційних освітніх і дослідницьких цілях.

У складеному за результатами виконаної роботи звіті групи говорилося, що безкоштовні матеріали OCW будуть залучати клієнтів до платних Knowledge Updates, проте незабаром ці проекти були розділені (курси Knowledge Updates готувалися з 2003 по 2005 рік і не змогли залучити достатньої кількості клієнтів) [2].

Представлення президентом МТІ Ч. Вестом (Charles Vest) проекту OCW у 2001 році визвало велике здивування стосовно того, що МТІ віддає всі свої навчальні матеріали в Інтернет. Основна ідея проекту будувалася на тому, що коштами та засобами інституту будуть видаватися методичні рекомендації з організації традиційного навчання (аудиторного), і при цьому інститут буде залишатися лідером у наданні доступних освітніх ресурсів по всьому світу для некомерційного використання і модифікації на безкоштовній основі.

Проект спочатку розраховувався на 10 років. У 2002 році запущена пілотна версія сайту проекту, а в 2003 відбулося офіційне відкриття. З 2004 року матеріали почали поширюватися за ліцензією Creative Commons Attribution-NonCommercial-ShareAlike. Дану ліцензію запропонував використовувати Х. Абельсон (Hal Abelson) – один із засновників як МТІ OCW, так і ліцензії Creative Commons. Обидві ці ініціативи з’явилися приблизно в один і той же час, переслідували схожі цілі (надання доступу до публікацій, їх подальше поширення та переробка), і вихідна ліцензія OCW 2002 року була подібна до тієї, на яку проект перейшов у 2004. МТІ OCW став першим великим зібранням робіт, опублікованих під ліцензією Creative Commons [74].

У 2007 і 2011 роках на сайті проекту з’явилися два додаткових розділи: Highlights for High School з ресурсами для викладання природничих дисциплін у середній школі, і OCW Scholar з матеріалами для самоосвіти [81].

У процесі розвитку проекту з’ясувалося, що, крім підвищення авторитету інституту, він приносить і додаткову користь. Близько 35% студентів, що вступили в МТІ у 2007 році, погодились з тим, що на їх вибір ВНЗ вплинув OCW. Також вважається, що використання OCW у процесі навчання підвищує якість викладання, оскільки розміщені в освітньому середовищі матеріали, доступні багатьом користувачам, є якісними [8].

Сьогодні практично всі курси МТІ можна знайти на OCW, багато з них переведено на кілька мов. Понад 70 мільйонів людей відвідали Web-сайт. Найбільшою популярністю користуються курси відеолекцій лінійної алгебри Г. Стренга (Gilbert Strang) та відеолекції з демонстраціями У. Левіна (Walter Lewin), зокрема закони Ньютона, що мають найбільшу кількість переглядів на YouTube [94].

Пояснюючи причини створення проекту, його ініціатори наголошували, що OCW продемонструє навчальну програму інституту, вплине на викладання в інших навчальних закладах, і покаже, що МТІ ставить знання вище фінансової вигоди. Цю точку зору поділяла більшість викладачів, що в подальшому було підтверджено добровільною участю у проекті 78% професорсько-викладацького складу МТІ [1].

Опубліковані на серпень 2012 року 2100 курсів складають близько 85 % навчальної програми МТІ. Але, як зазначено в звіті компанією для ректора МТІ про OCW [81], в планах OCW продовжувати публікувати нові курси і оновлювати кожного року існуючі курси.

До загальних цілей OCW розробники відносять [81]:

– публікації матеріалів курсів МТІ високої якості;

– збільшення використання OCW для викладання та навчання;

– максимально використовувати переваги OCW для спільноти МТІ;

– підтримку по всьому світу відкритих освітніх ресурсів (ВОР) і руху OCW;

– підтримку програми МТІ OCW.

Зауважимо, що OCW не є програмою дистанційного навчання в МТІ, а є лише зрізом того, як певний предмет викладався у певний період. На сайті подано навчальні матеріали минулих років: відеолекції, мультимедійні програми, а також зразки екзаменаційних завдань [75].

У таблиці 2.1 вказано статистику та призначення використання користувачами OCW (до користувачів в основному включають самоосвітян, студентів і педагогів).

Оскільки на OCW містяться навчальні матеріали загального плану, що не пов’язані з навчальними планами інституту, тому кожна кафедра має свій сайт, де розміщує відомості, актуальні у відповідному семестрі.

*Таблиця* 2.1

**Статистика використання користувачами OpenCourseWare [45]**

| **Користувачі** | **Напрями використання** | **% використання** |
|---|---|---|
| Педагоги | Для поліпшення особистих знань | 31 % |
| | Вивчення нових методів навчання | 23 % |
| | Включення OCW матеріалів у курси | 20 % |
| | Як довідковий матеріал для своїх студентів | 15 % |
| | Розробка навчальних програм для свого відділу або школи | 8 % |
| Студенти | Для підвищення особистих знань | 46 % |
| | Доповнення поточного курсу | 34 % |
| | Перегляд плану навчального курсу | 16 % |
| Самоосвітяни | Дослідити простір за межами своєї професійної області | 40 % |
| | Огляд основних понять у своїй професійній сфері | 18 % |
| | Підготовка до майбутнього дослідження | 18 % |
| | Слідкування за розвитком подій в своїй області | 17 % |
| | Заповнення, пов’язані з роботою проектів або завдань | 4 % |

На кафедральних сайтах Шкіл МТІ представлено поточний стан курсів. Так, на сайті кафедри вищої математики МТІ (http://math.mit.edu) містяться такі відомості про засоби ІКТ навчання математичних дисциплін у осінньому семестрі 2012-2013 н. р., а саме:

– персональні дані лектора, адміністратора курсу та викладачів, що проводять практичні заняття (http://math.mit.edu/people/directory_faculty.php);

– платформа Piazza для інтерактивної навчальної позааудиторної взаємодії (https://piazza.com/mit/fall2012/1801/home): форум для спілкування учасників курсу, на якому будь-хто може розмістити питання чи коментарі з матеріалів курсу, домашніх завдань тощо;

– дистанційна та мобільна система управління навчанням Stellar, розроблена у MIT (https://stellar.mit.edu/courseguide/course/18/fa12/18.01/);

– допоміжні навчальні матеріали, приклади та розв’язання задач у форматі PDF (http://math.mit.edu/classes/18.01/1801_Supplementary%20Notes.html);

– студентський Центр навчання математики (Math Learning Center)

для надання консультативної підтримки з курсу (http://math.mit.edu/learningcenter);

– версія курсу у OCW Scholar, призначена для самостійного опрацювання (http://ocw.mit.edu/courses/mathematics/18-01sc-single-variable-calculus-fall-2010).

Офіційний сайт кафедри математики МТІ (рис. 2.2) зазначає, що кафедра вищої математики інституту є визнаною у світі у галузі теоретичної та прикладної математичної науки та освіти інженерів.

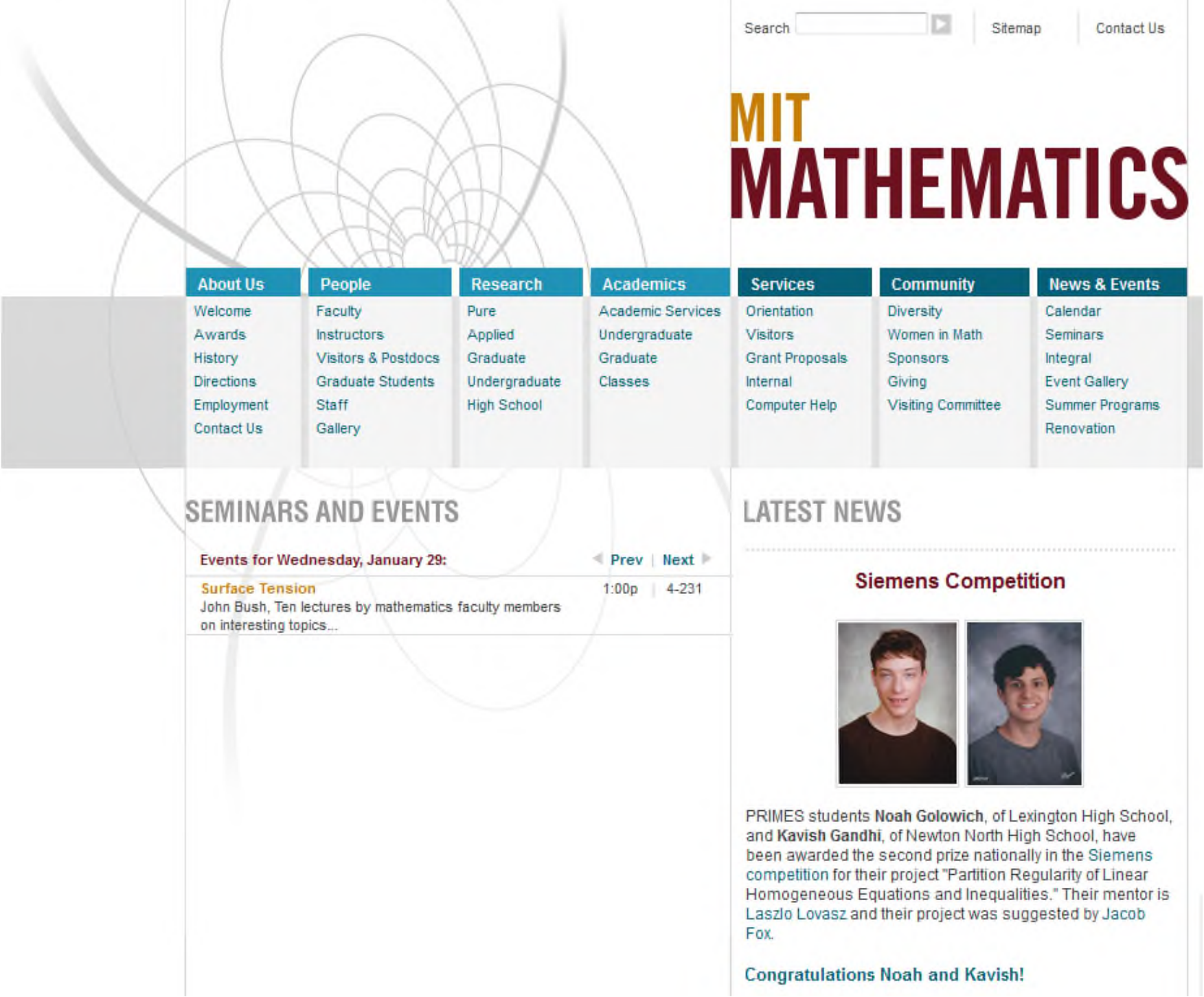

Рис. 2.2. Сайт кафедри вищої математики МТІ (http://math.mit.edu)

На кафедрі ведуться дослідження з провідних напрямів математики, особлива увага приділяється саме прикладній математиці, де інноваційні математичні міркування приводять до нового розуміння і доповнення. Група математиків у прикладній області фокусується на біології, комбінаториці, комп’ютерних науках, наукових обчисленнях, математичному аналізі та механіці.

На 2012-2013 навчальний рік на кафедрі працює 50 викладачів, серед яких є стипендіати: премії Абеля; Національної медалі науки; премій Мак-Артура (MacArthur Awards), Бохера (Bôcher), Коула (Cole), Веблена (Veblen), Фалкерсона (Fulkerson) та Стіла (Steele). Шістнадцять викладачів є членами Національної академії науки і техніки і двадцять три є членами Американської академії мистецтв і наук.

Інструктори з математики та прикладної математики, якими є аспіранти, мають можливість проживати в МТІ протягом двох або трьох років, і працювати на факультеті. Вони мають невеликі навчальні навантаження, що забезпечує їм сприятливі умови для дослідницької роботи. На кафедра навчається постійно близько 30 аспірантів.

Опишемо зазначені засоби ІКТ навчання вищої математики в МТІ.

*Персональні дані* лектора, адміністратора курсу та викладачів, що проводять практичні заняття (рис. 2.3).

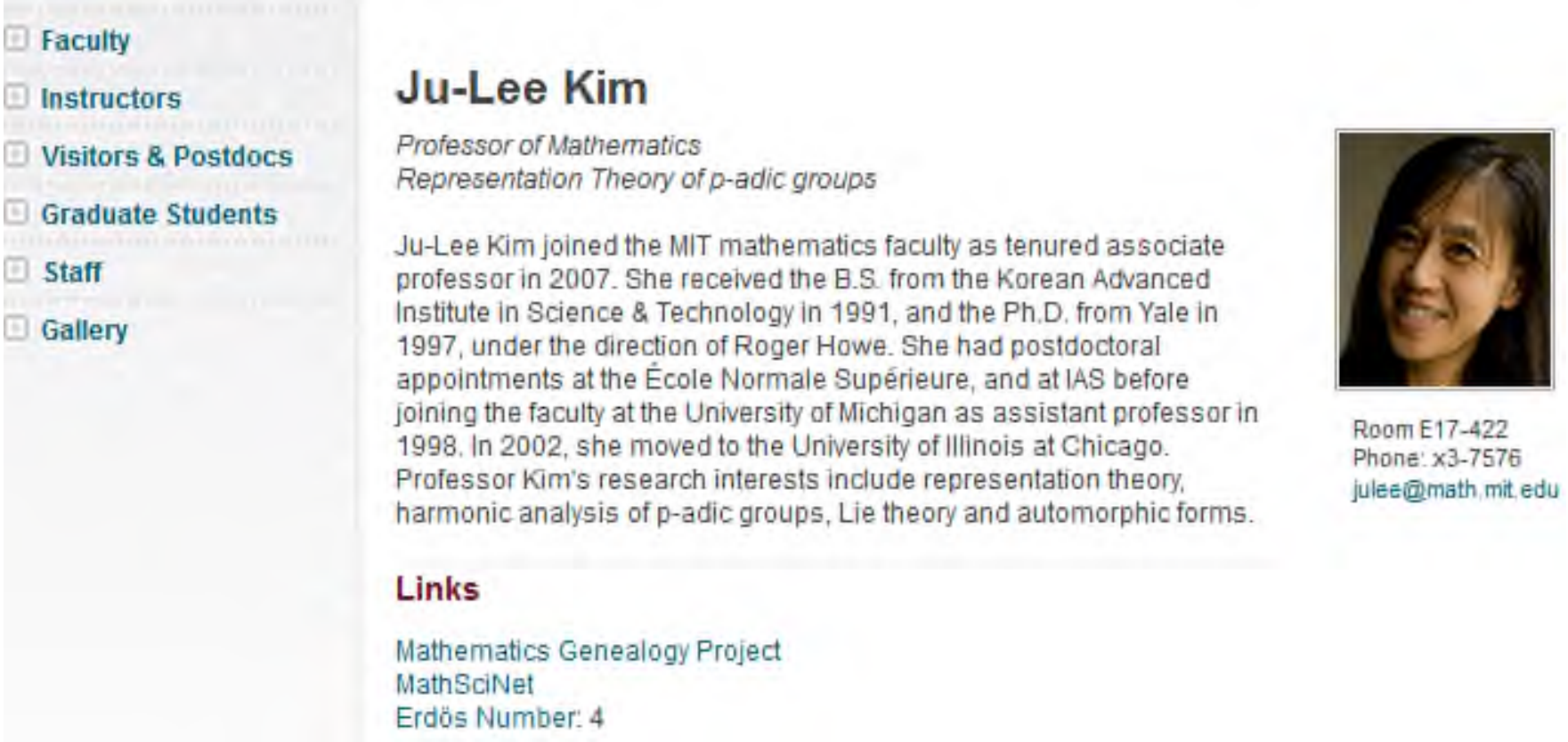

Рис. 2.3. Сайт кафедри вищої математики з відомостями про персональні дані лектора, адміністратора курсу та викладачів, що проводять практичні заняття

На сайті містяться відомості про наукове звання, сферу наукової діяльності, наукові досягнення; номер аудиторії, в якій можна знайти викладача та номер його телефону, електронну адресу. По кожному викладачу можна дізнатися про його науковий та викладацький шлях як у МТІ, так і на попередніх місцях роботи. Для деяких викладачів є посилання на їх персональні сайти, на електронну бібліотеку публікацій MathSciNet із зазначенням їх числа Ердьоша (Erdős number) – кількості статей у співавторстві, сторінку в Wikipedia.

Професорський склад кафедри містить посилання на сайт Mathematics Genealogy Project (www.genealogy.ams.org), де зібрано відомості про всіх осіб, які отримали докторський ступінь у галузі математики. Для кожного доктора відведено сторінку, на якій показано наступне: повне ім’я одержувача ступеня; назва університету, що надав ступінь; рік, в якому був удостоєний ступеня; повна назва дисертації; повне ім'я (імена) радника (радників).

На сайті кафедри існує підрозділ з назвою «Математичні академічні послуги» (Math Academic Services – MAS). До академічних послуг відносять надання підтримки студентам і викладачам щодо курсів, навчання, консультування, освітніх програм та студентських досліджень. В обов'язки підрозділу входить складання та контроль за навчальним графіком курсів (тобто розклад лекцій, екзаменів, консультацій), складання програм та вимог з курсів, адміністрування та розрахунок заробітної плати, адміністрування підготовки викладачів, адміністрування студентського Центру навчання математики.

Студентський *Центр навчання математики* (Math Learning Center) (http://math.mit.edu/learningcenter) (рис. 2.4) спрямований надавати консультативну підтримку з курсів. У Центр навчання математики кафедрою вищої математики запрошуються на роботу досвідчені студенти й аспіранти, де вони проводять консультації (не безкоштовні) за курсами Числення однієї змінної, Числення багатьох змінних, диференціальних рівнянь та лінійної алгебри. В обов'язки консультантів входить забезпечення сприятливого середовища для підтримки самостійної роботи студентів, підготовки до іспитів, при цьому в разі потреби надається допомога експерта. Серед вимог, що ставлять до консультантів Центру навчання математики, є приязність, терплячість і готовність працювати по дві години один раз на тиждень.

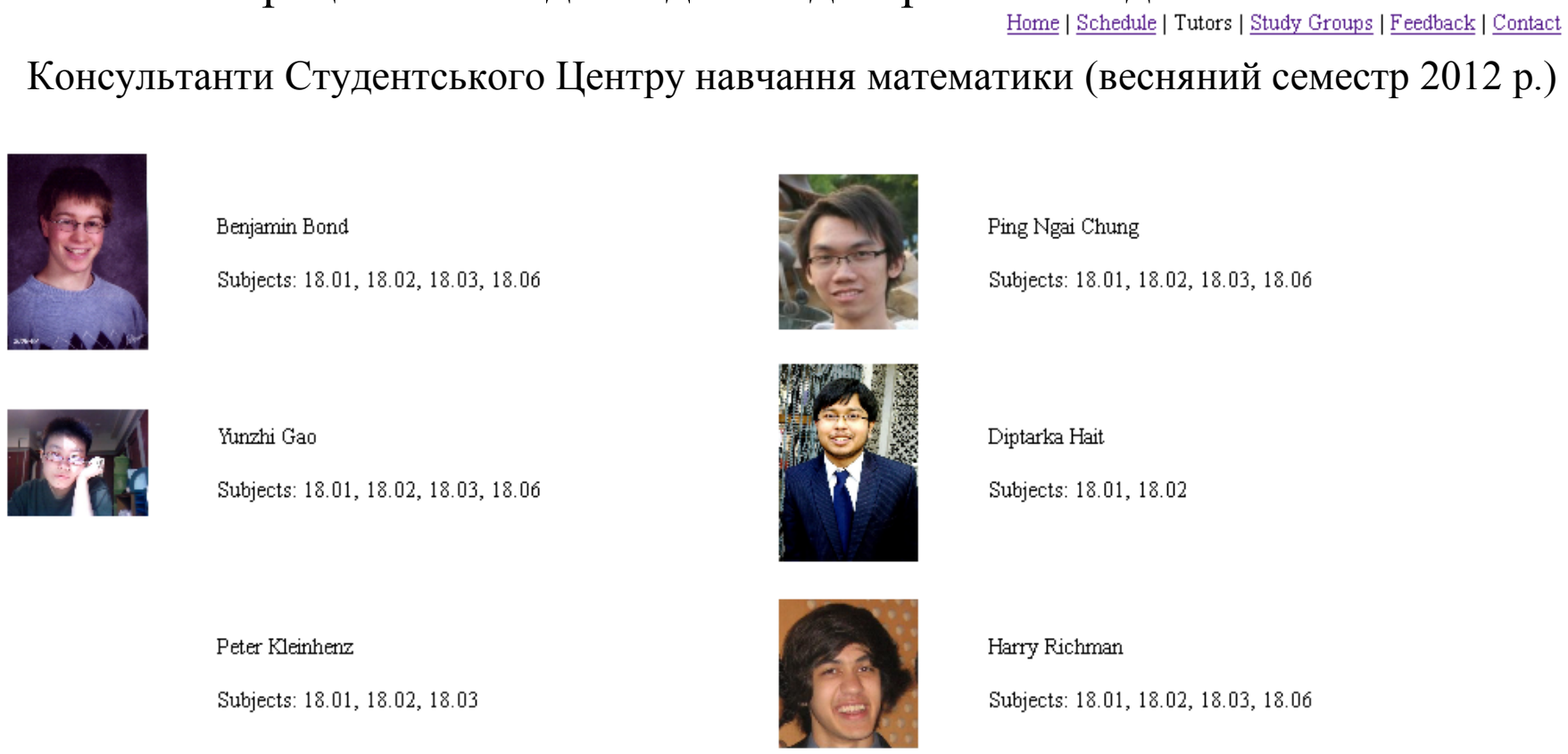

Рис. 2.4. Студентський Центр навчання математики (Math Learning Center) (http://math.mit.edu/learningcenter)

Крім організації консультативної допомоги в Центрі навчання математики, викладачами кафедри надано можливість студентам самостійно засвоювати навчальний матеріал, спираючись на допоміжні навчальні матеріали, приклади та розв’язання задач. Для кожного курсу як в онлайн перегляді, так і для завантаження файлів у форматі PDF доступно наступне:

– перелік тем курсу;

– викладення теоретичного матеріалу з наведеними розв’язаними прикладами;

– добірка прикладів з кожної теми курсу та відповіді до них;

– добірка прикладів для повторення та відповіді.

Зокрема на сторінці курсу «Числення однієї змінної» А. Маттука (Arthur Mattuck) і Д. Джерісона (David Jerison) можна обрати окремий розділ курсу для перегляду, або завантажити одразу всі файли (рис. 2.5) (скорочене посилання на динамічну сторінку у Internet Archive – http://goo.gl/CUnHRP).

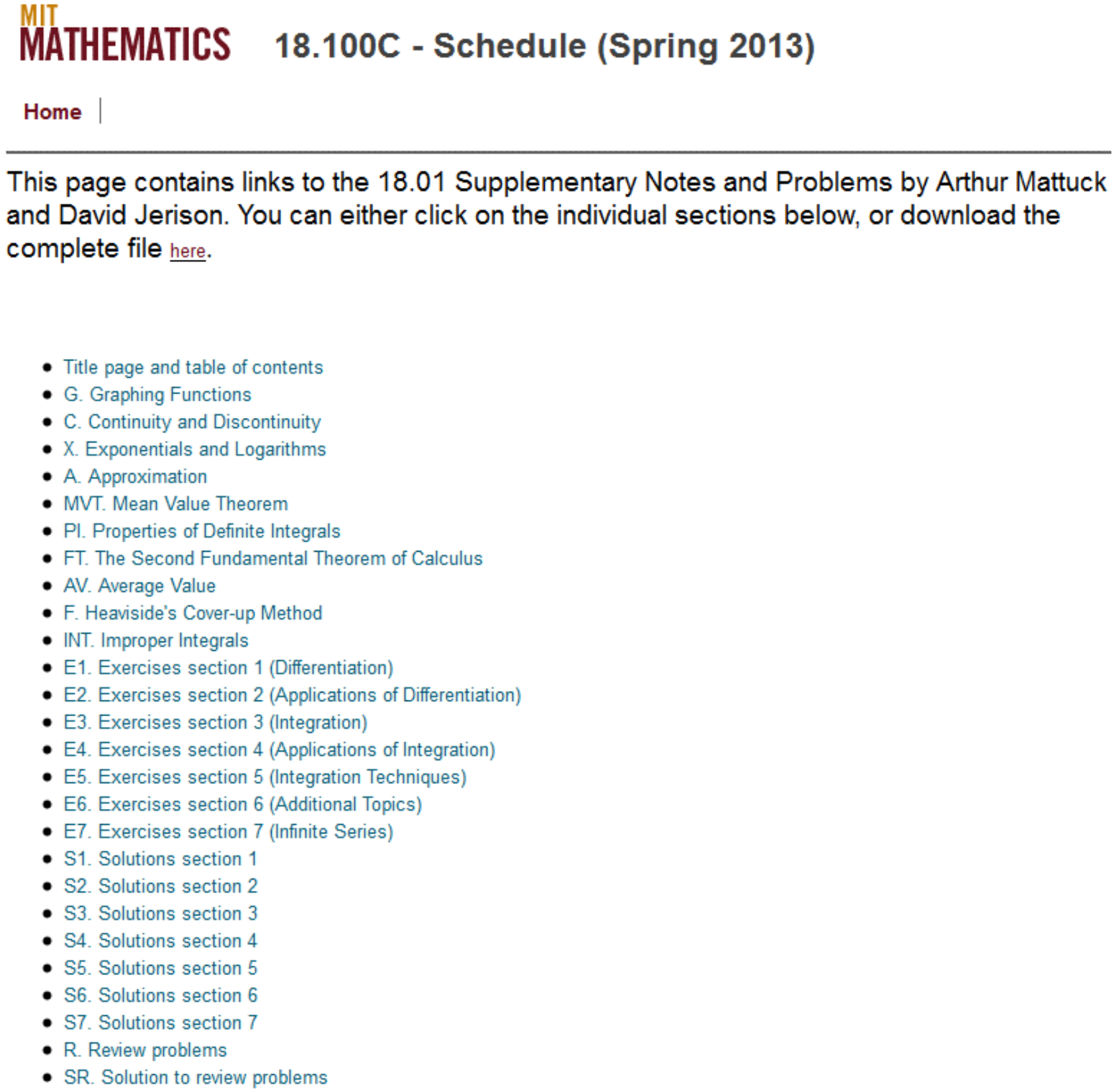
MIT MATHEMATICS 18.100C - Schedule (Spring 2013)

Home |

This page contains links to the 18.01 Supplementary Notes and Problems by Arthur Mattuck and David Jerison. You can either click on the individual sections below, or download the complete file here.

- Title page and table of contents
- G. Graphing Functions
- C. Continuity and Discontinuity
- X. Exponentials and Logarithms
- A. Approximation
- MVT. Mean Value Theorem
- PI. Properties of Definite Integrals
- FT. The Second Fundamental Theorem of Calculus
- AV. Average Value
- F. Heaviside's Cover-up Method
- INT. Improper Integrals
- E1. Exercises section 1 (Differentiation)
- E2. Exercises section 2 (Applications of Differentiation)
- E3. Exercises section 3 (Integration)
- E4. Exercises section 4 (Applications of Integration)
- E5. Exercises section 5 (Integration Techniques)
- E6. Exercises section 6 (Additional Topics)
- E7. Exercises section 7 (Infinite Series)
- S1. Solutions section 1
- S2. Solutions section 2
- S3. Solutions section 3
- S4. Solutions section 4
- S5. Solutions section 5
- S6. Solutions section 6
- S7. Solutions section 7
- R. Review problems
- SR. Solution to review problems

Рис. 2.5. Сторінка сайту з курсом «Числення однієї змінної» А. Маттука (Arthur Mattuck) і Д. Джерісона (David Jerison)

В якості форуму для спілкування учасників курсу, на якому будь-хто може розмістити питання чи коментарі з матеріалів курсу, домашніх завдань тощо, забезпечується використанням платформи Piazza, що виступає в ролі інтерактивної навчальної позааудиторної взаємодії.

*2.1.2 Платформа Piazza*, як зазначено на сайті (https://piazza.com/), є абсолютно безкоштовною, простою у використанні і швидкою в налаштуванні платформою (рис. 2.6) для інструкторів, що призначена ефективно управляти створеним на платформі Piazza класом. Студенти можуть задати питання і спільно редагувати відповіді на ці питання. Викладачі можуть також відповісти на запитання або підтримати відповіді студентів і редагувати їх, або видаляти будь-який розміщений зміст. Piazza призначена для імітації реального обговорення в аудиторії.

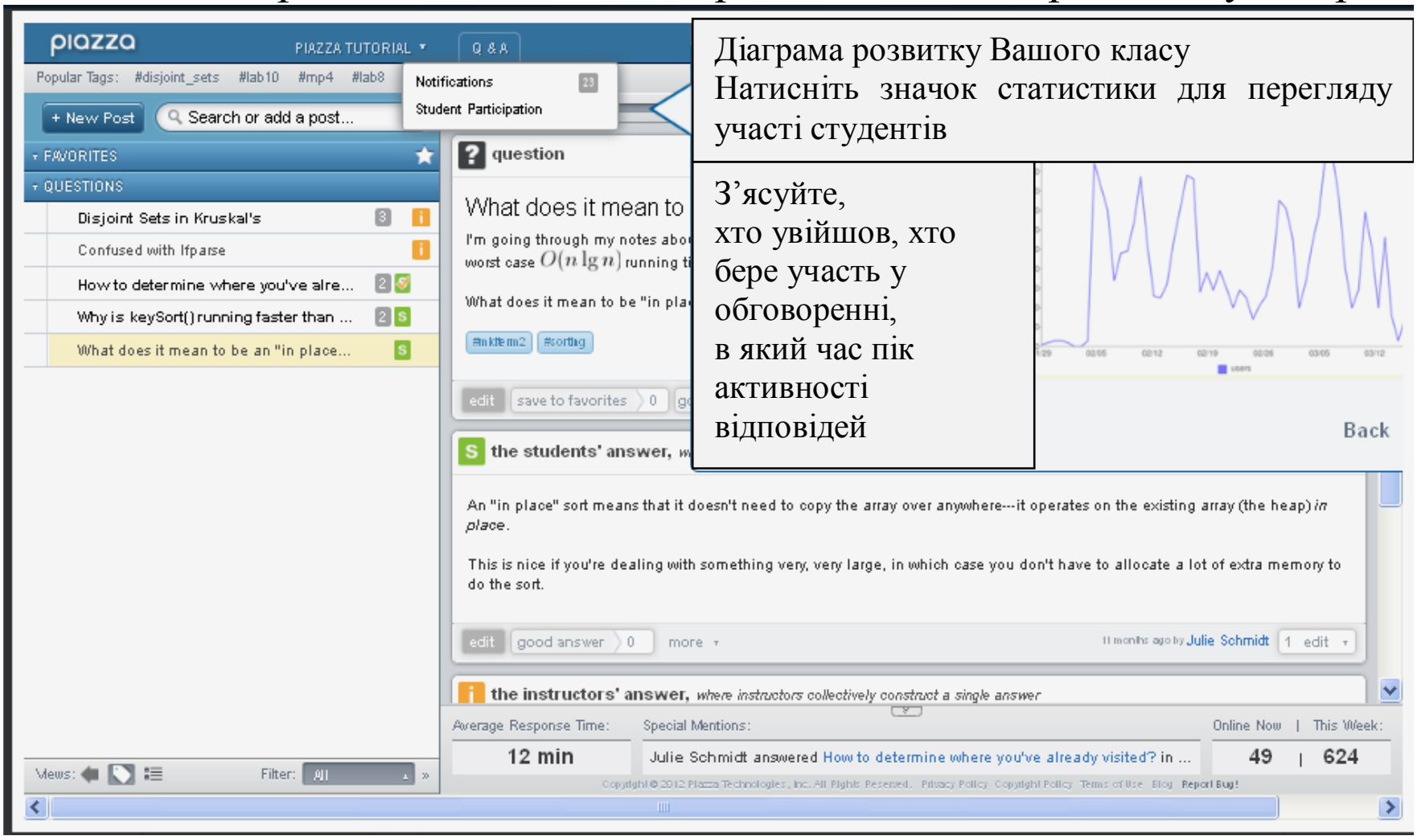

Рис. 2.6. Платформа Piazza для інтерактивної навчальної позааудиторної взаємодії (https://piazza.com/why-piazza-works)

Назва «Piazza» походить від італійського слова «площа» – загальної площі міста, де люди можуть збиратися разом, щоб поділитися знаннями та ідеями. Її призначення – стимулювати студентів до навчання і заощаджувати час викладача.

Використання платформи Piazza у процесі навчання створює умови для реалізації таких методів навчання, як навчання у групах, метод проектів, метод різнорівневого навчання, де студенти виступають у ролі викладачів-консультантів, оскільки бесіди в Piazza можна продовжувати і після закінчення занять в аудиторіях.

Використання платформи Piazza у процесі навчання створює умови для роботи студентам навіть анонімно. Задані питання аналізуються інструктором і викладаються на публічне обговорення у групі. Аналіз усіх отриманих відповідей та коментарів узагальнюється інструктором із зазначенням більш чіткої та правильної відповіді.

Питання та відповіді в Piazza редагуються спільними зусиллями. На кожне запитання можна дати лише одну відповідь студента і одну відповідь інструктора (рис. 2.7), після чого запитання розміщується для подальшої дискусії.

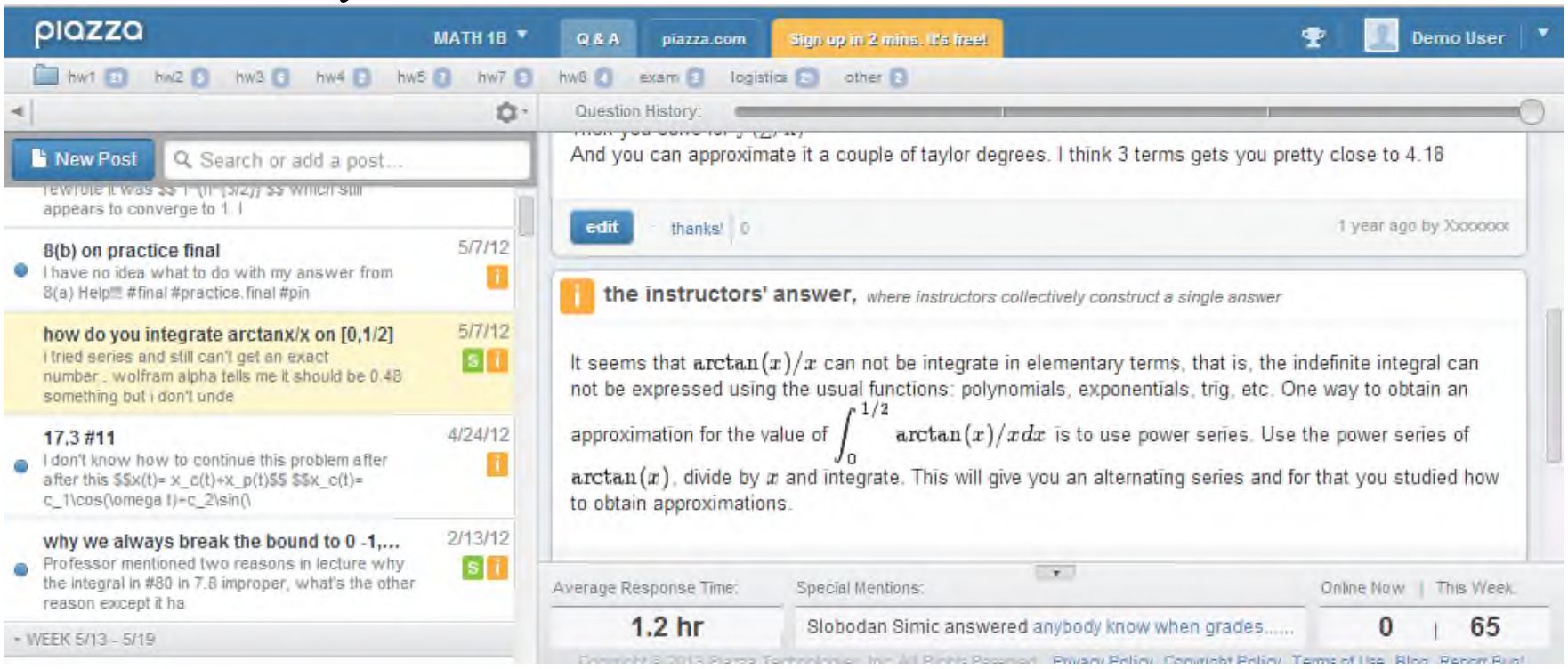

Рис. 2.7. Сторінка обговорення питання на платформі Piazza

Існує також мобільна версія платформи Piazza на базі додатків IOS і Android, що надає можливість віддаленим від аудиторії студентам брати участь у роботі групи, відсилаючи свої відповіді через повідомлення електронної пошти (рис. 2.8).

Платформа підтримує wiki-стиль Q&A (Questions and answers), що надає студенту можливість отримувати відповідь на запитання, не читаючи непотрібні довгі обговорення в форумі, прочитавши тільки одну відповідь, що має високу оцінку запитань і відповідей.

Відповіді студентів та викладачів з'являються швидко в зв'язку з швидкою роботою платформи Piazza. На відміну від Web-сайтів статичних класів, екран Piazza оновлюється, як тільки відбуваються зміни, так що можна побачити всі зміни в реальному часі. Студент може залишити свій браузер відкритим весь вечір і побачити, як відповідають на його запитання.

У викладача є повний редакторський контроль над змістом дискусій класу. Внески викладачів у класі виділені, так що студенти можуть легко визначити відповіді інструктора (рис. 2.9). Інструктор може заохочувати студентів задавати і давати відповіді на запитання, схвалюючи їх.

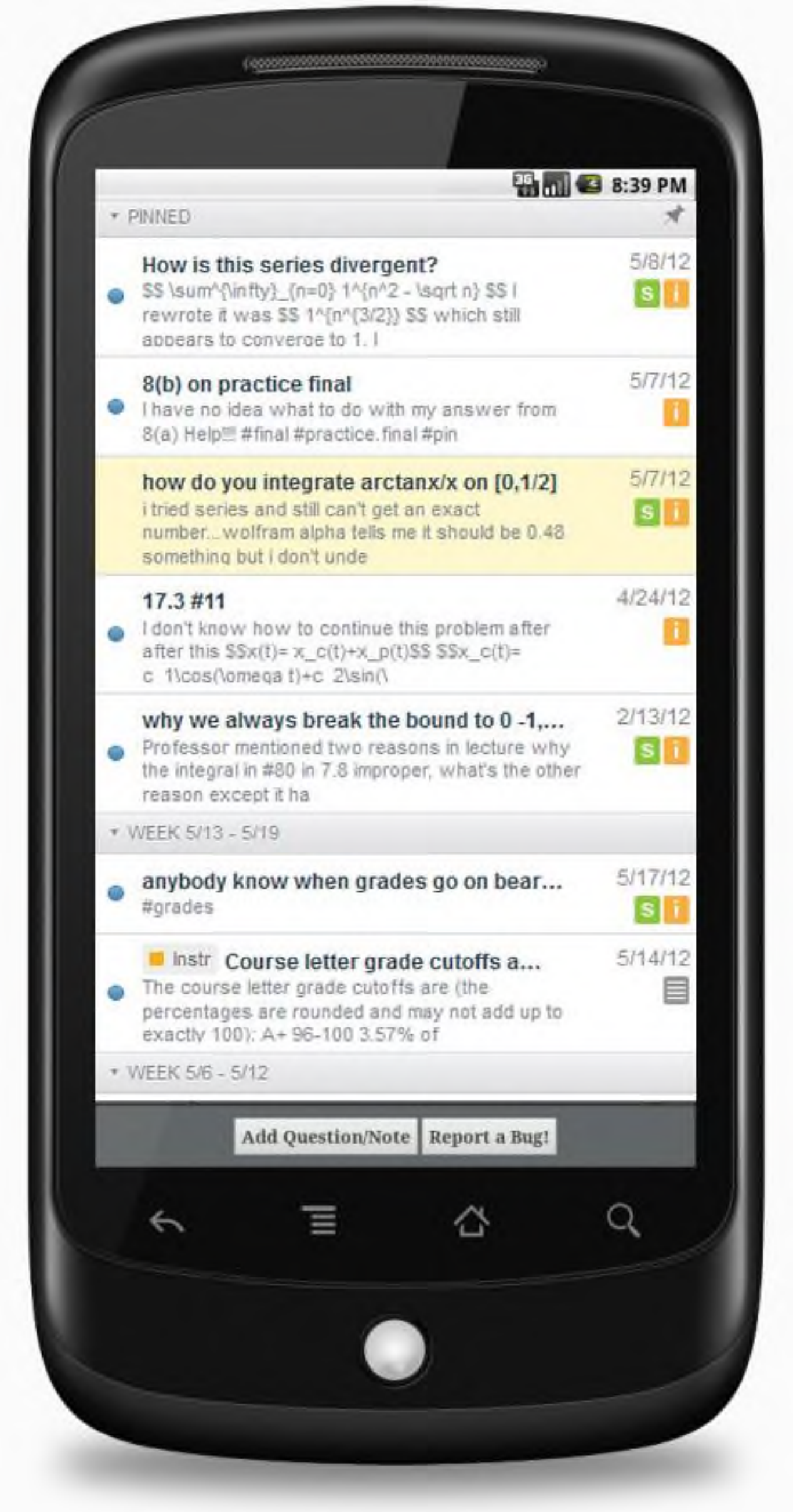

Рис. 2.8. Мобільна версія платформи Piazza

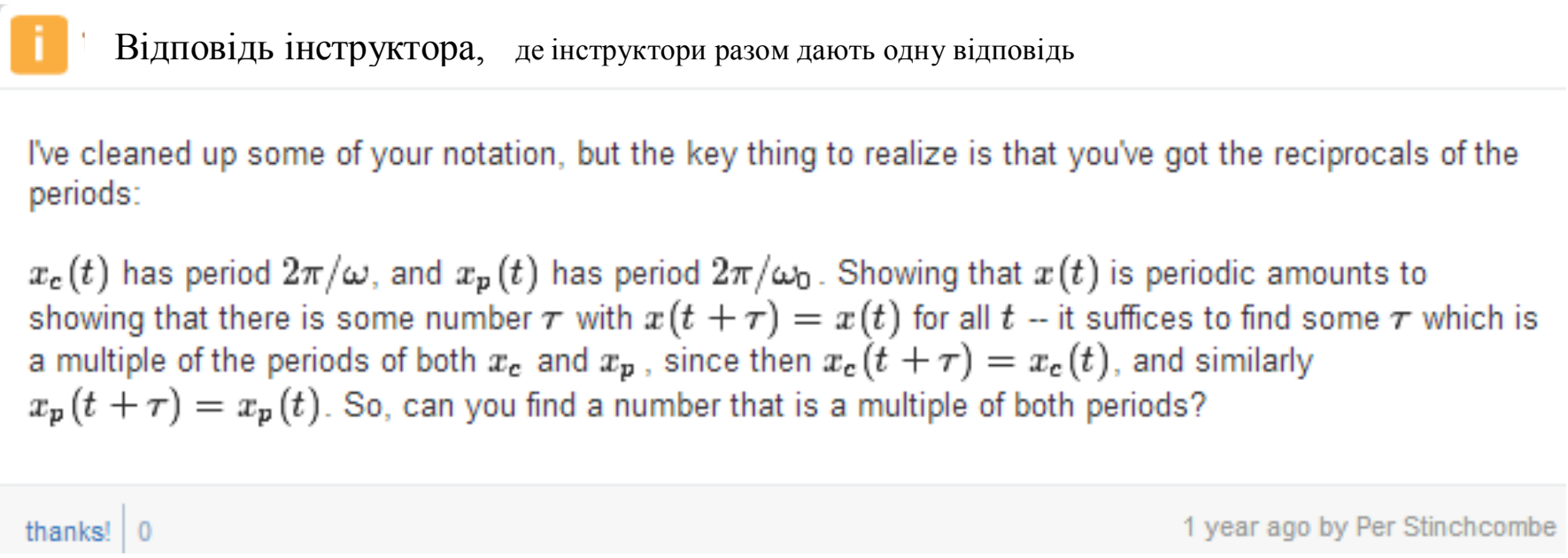

Рис. 2.9. Повідомлення інструктора на платформі Piazza

Сторінки курсу – це місце, де можна управляти оголошеннями, відомостями про навчальні програми та навчальні плани і ресурсами курсу (рис. 2.10). У розділі оголошень є можливість переглянути весь список електронних адрес курсу. Якщо є необхідність зробити повідомлення для всього класу, необхідно просто опублікувати повідомлення, позначивши його як «важливе оголошення», і вибрати операцію «відправити всім».

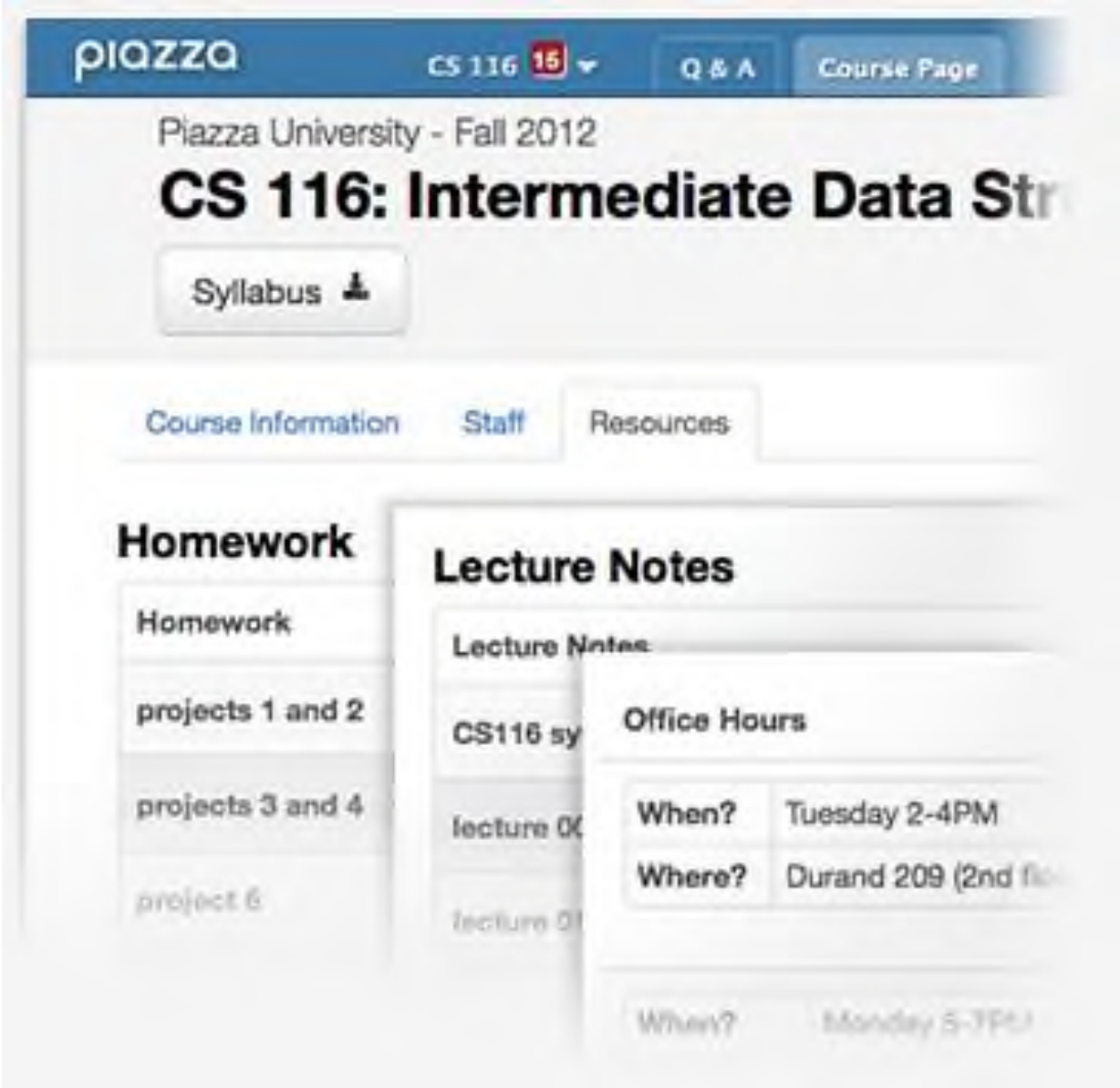

Рис. 2.10. Сторінки курсу Piazza, де можна управляти оголошеннями, відомостями про навчальні програми, навчальні плани і ресурси курсу

Розділ відомостей про навчальні програми – простий спосіб зберігати дані про курс та персональні контактні дані. Повідомивши студентам про місце зберігання цих даних, студенти можуть дізнаватися про розклад консультацій протягом тижня.

Студенти можуть відвідувати розділ «Ресурси курсу», щоб знайти такі матеріали: домашні завдання, конспекти лекцій, важливі посилання, що розміщуються всі в одному місці. Поставивши позначку про оновлення або внесення нових ресурсів, весь клас завжди знайде їх на відповідній сторінці.

Переглядаючи звіти про роботу студентів на платформі, можна визначити рейтинг роботи студентів групи, їх запитання, відповіді та редагування. В загальному звіті групи є дані про те, які питання найчастіше задаються і як вони узгоджуються з тим, що відбувається в групі (рис. 2.11).

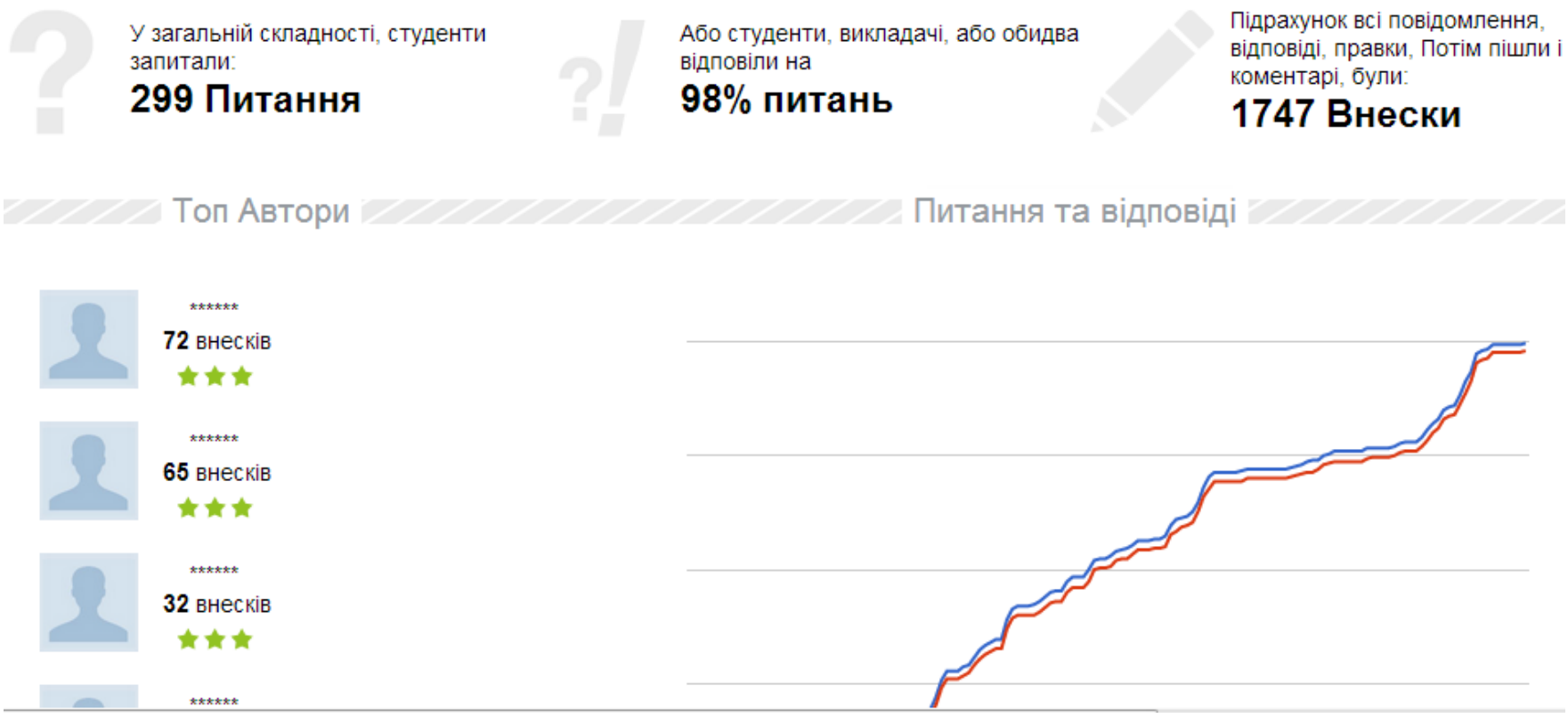

Рис. 2.11. На сторінці платформи Piazza можна сформувати звіти про роботу студентів

Використання платформи Piazza дає можливість проводити опитування серед студентів після заняття. Шляхом анонімних повідомлень інструктору можна отримати детальні коментарі про проведене аудиторне заняття.

Використовуючи опитування студентів, можна з’ясувати, які теми і завдання викликають труднощі, а які є простими (рис. 2.12); хто зі студентів потребує додаткових консультацій.

Контролюючи діяльність студентів за допомогою платформи Piazza, викладач має уявлення про проблемні теми, яким слід приділити більшу увагу [139].

Оцінювати роботу студентів у класі Piazza викладачам полегшують об’єднані засоби навчальної взаємодії (Learning Tools Interoperability – LTI). Додавши посилання на Piazza у LMS, що використовується, студентам не потрібні ім’я користувача і пароль для входу в систему.

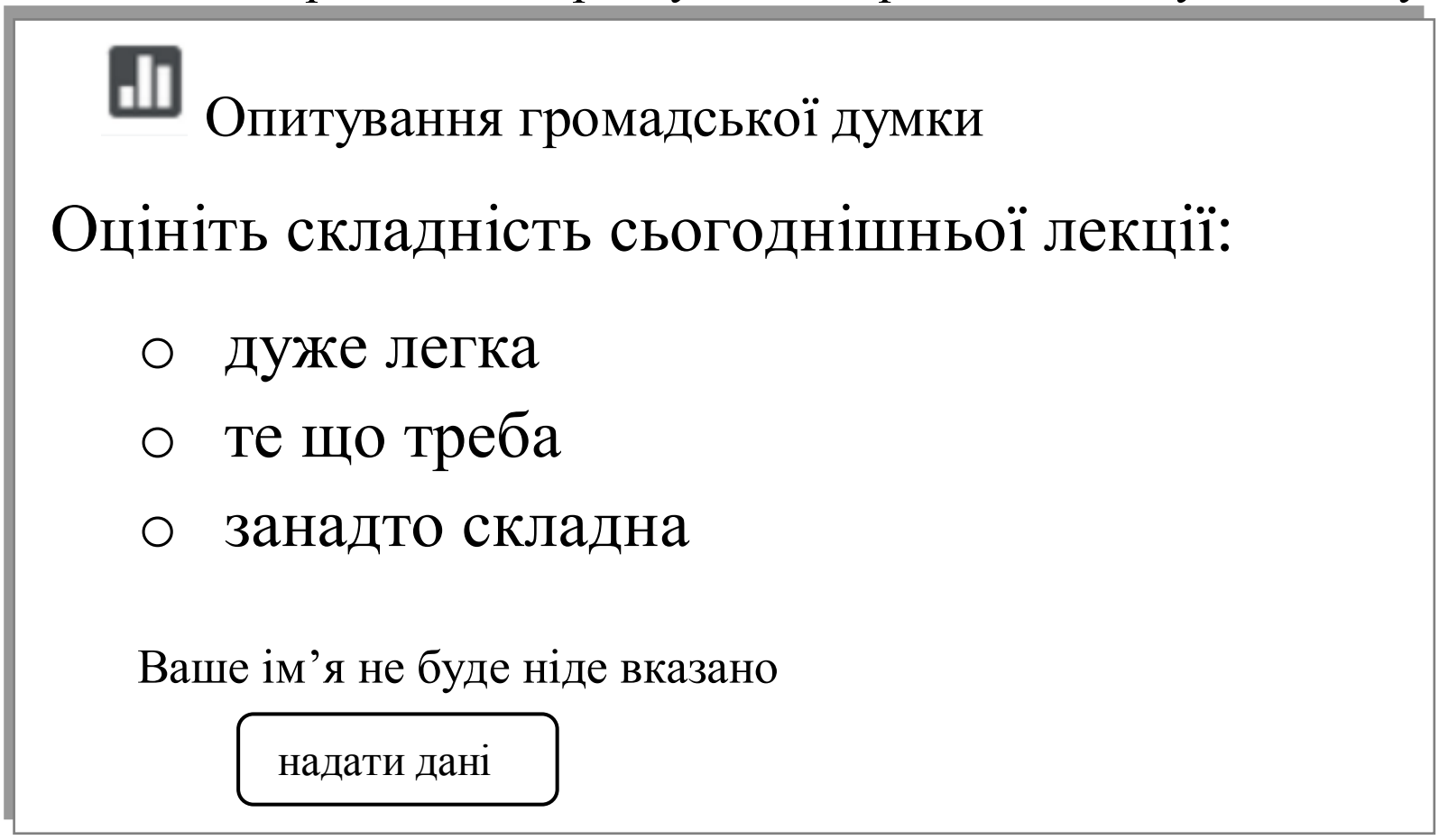

Рис. 2.12. Бланк для анонімного опитування студентів у Piazza про доступність викладання матеріалів на заняттях

Зручність використання Piazza з LMS полягає в тому, що:

– викладач може швидко активувати платформу Piazza для аудиторії;

– студенти та викладач автоматично реєструються з відповідними можливостями, встановленими при виконанні реєстрації;

– при натисканні на кнопку «Piazza» в LMS, користувачі автоматично заходять на платформу Piazza, тому немає необхідності запам’ятовувати окремі імена користувача та пароль.

До LMS, що підтримують платформу Piazza, відносять: Blackboard, Moodle, Sakai, Desire2Learn, Canvas, Jenzabar, Angel та інші.

Позитивними моментами використання платформи Piazza в організації процесу навчання з вищої математики є наступне [102]:

– платформа має надійний редактор формул LaTeX для простих математичних виразів (рис. 2.13);

– використання платформи Piazza надає викладачам можливість пояснити складний для розуміння математичний матеріал;

– надає студентам, які не зрозуміли лекцію, можливість повторного перегляду лекції з наступним обговоренням проблемних питань;

– у мобільній версії платформи студенти можуть робити фотографії їх письмових робіт, що усуває проблему введення формул та математичних виразів;

– платформа підтримує мультимедійні вкладення, що створює умови для обміну між викладачами і студентами діаграмами, зображеннями та відео.

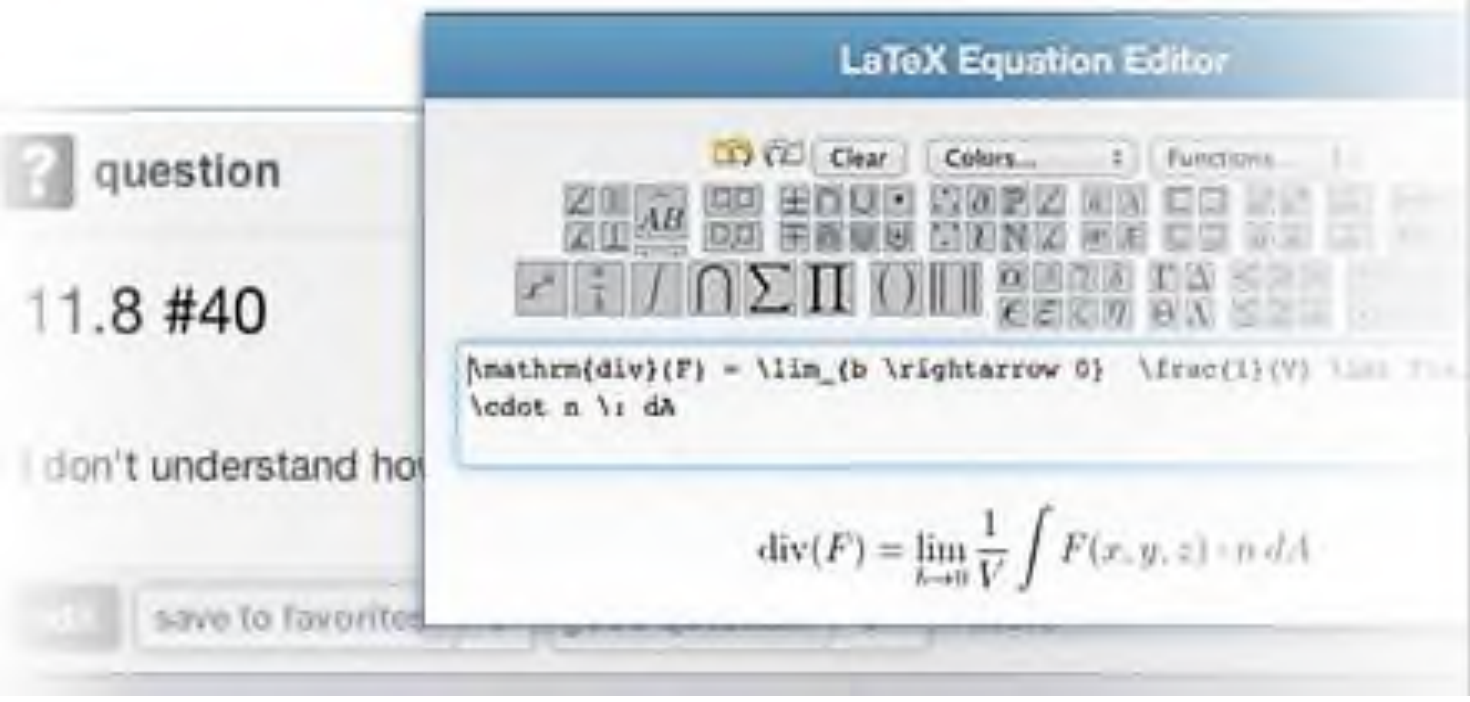

Рис. 2.13. Редактор формул LaTeX в Piazza

*2.1.3 Система управління навчанням Stellar* – це система управління курсами в МТІ, що може бути як дистанційною (рис. 2.14*а*), так і мобільною (рис. 2.14*б*).

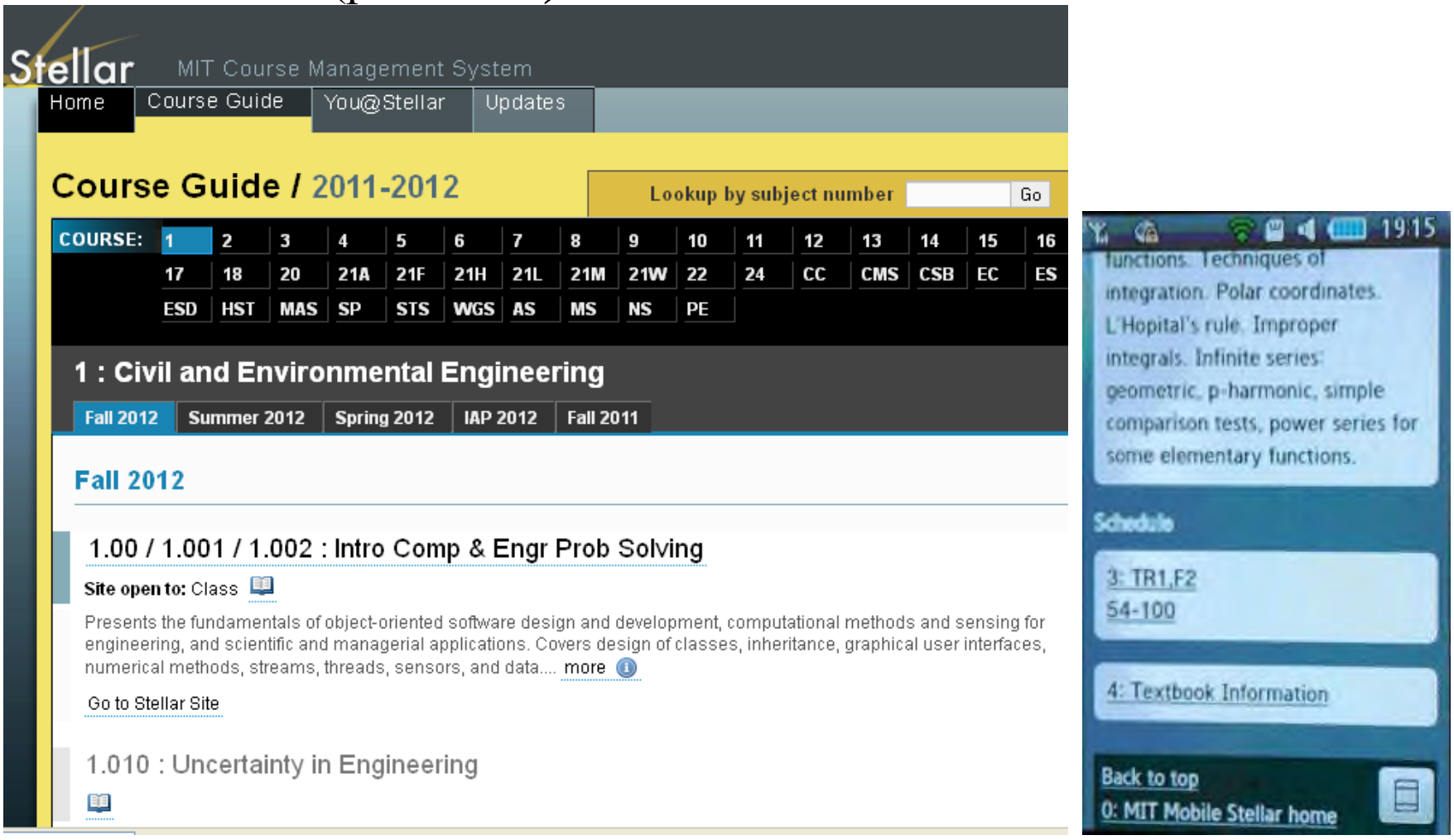

*а* *б*

Рис. 2.14. Дистанційна (*а*) та мобільна версія системи управління навчанням Stellar (http://stellar.mit.edu/mobile/atstellar)

Система управління навчанням Stellar надає студентам всі відомості про курси, які вони вивчають, у тому числі оголошення, навчальні програми та календар курсу, контактні дані викладачів, журнал оцінок тощо [75]. На сайтах OCW і Stellar викладачі вищої математики використовують велику кількість програм, у тому числі Web-додатків, для закріплення основних понять, які були представлені в аудиторії.

Метою розробки Stellar, як зазначено на офіційному сайті (http://stellar.mit.edu/about/), є інтеграція з інфраструктурою МТІ. Призначення Stellar полягає в наступному:

– зробити доступними навчальні відомості з кожного курсу, що пропонується в Массачусетському технологічному інституті;

– надавати можливість користувачам створювати інтуїтивно зрозумілий і багатофункціональний Web-сайт для класу або проекту без написання коду;

– пропонувати надійне онлайн управління і публікувати навчальні матеріали;

– надати користувачам можливість контролювати рівень доступу до свого сайту, що робить його як обмеженим, так і відкритим, за їх вибором;

– автоматизувати трудомісткі адміністративні задачі;

– забезпечити взаємодію зі студентами через оголошення, повідомлення електронної пошти і форуми;

– забезпечити взаємодію між користувачами шляхом інтеграції вікі;

– приймати рекомендації щодо Інтернет-управління процесом навчання;

– забезпечувати пошук і керівництво курсом;

– публікувати матеріали досліджень і опитувань.

Stellar розроблений спеціально для підтримки всіх курсів у МТІ, забезпечуючи основу для розміщення змісту курсів та інших матеріалів у мережі. Використовуючи програми, такі як Microsoft Word, PowerPoint або MATLAB, викладачі можуть створювати документи, а потім завантажувати їх на Web-сайт, де Web-сторінки генеруються автоматично. Викладач може легко управляти змістом курсу і доступом до нього; студенти мають можливість ініціювати і приєднуватися до дискусій.

Використання системи Stellar надає викладачам можливість контролювати темп розвитку змісту курсу. Вони можуть завантажити весь курс з усіма його вимогами і призначеннями на початку семестру, або можуть доповнювати або змінювати матеріали з курсу впродовж всього навчального періоду.

Матеріали курсів, що розміщені на Stellar, можна легко імпортувати на сайт OpenCourseWare МТІ для відкритого обміну змістом курсів МТІ серед викладачів і студентів для самоосвіти як по всій країні, так і у всьому світі [36].

Як правило, кожного тижня з курсу вищої математики у МТІ за розкладом 2-3 лекції та 2 семінари для закріплення матеріалу [30]. У таблиці Е.1 додатку Е наведено приклад розкладу занять на осінній семестр 2012-2013 н. р. з 18.014 – «Числення однієї змінної з теорією». На лекції вводяться і розглядаються нові математичні поняття з ілюстрацією прикладів, тоді як на практичному занятті поняття закріплюються прикладами. У випадку, якщо студенти мають проблеми із засвоєнням будь-якого матеріалу, вони можуть знайти викладачів у будь-який час робочого дня за допомогою сайту Stellar або сайту конкретного курсу [75].

Для більшості курсів MTI з вищої математики загальна оцінка курсу складається в основному з 2-3 проміжних екзаменів, підсумкового екзамену, і оцінок щотижневих домашніх завдань. При підготовці до екзаменів студентам рекомендується використовувати Інтернет-ресурси, щоб знайти зразки прикладів для опрацювання; переглядати матеріали, такі як екзамени за попередні роки або зразки онлайн-завдань. Щотижневі завдання, однак, спираються більше на онлайн-програми та математичні інструменти. Вони розміщені в Інтернет, на сайті Stellar чи на Web-сайті конкретного курсу, тому студенти можуть переглядати і опрацьовувати їх у зручному для них режимі. У той час як більшість завдань включає в себе традиційні обчислення (тобто розв’язання без використання засобів ІКТ), серед викладачів MTI існує тенденція розробки завдань з використанням програм, що надають студентам можливість розібратися і зрозуміти навчальний матеріал з вищої математики [75].

*2.1.4 Курси для самостійного вивчення.* Курси на *OpenCourseWare Scholar* розроблені спеціально для найбільшої аудиторії OCW: самоосвітян – тих, хто навчаються самостійно (independent learners). Курси OCW Scholar призначені для використання тими, хто має намір вивчати навчальний матеріал самостійно. Курси є більш повні, ніж подібні курси, розміщені на OCW і включають в себе новий, створений користувачами, зміст з раніше опублікованих курсів. Матеріали розташовані в логічній послідовності і часто включають в себе мультимедійні, такі як відео та модельні демонстрації. На даний час курси OCW Scholar є експериментом, тому адміністратори активно приймають відгуки про ефективність цих курсів.

Курси OCW Scholar є однією з цілого ряду інноваційних ініціатив, над якими працюють в MTI. Через OCW Scholar організатори експерименту сподіваються надати підтримку саме для тих людей, які навчаються самостійно. OCW спочатку був задуманий як набір ресурсів для педагогів з усього світу, що можна адаптувати для використання в процесі навчання. З часом виникла необхідність пропонувати матеріали курсу з метою полегшення самостійного навчання. OCW Scholar – перша спроба зробити це.

На початок 2012-2013 навчального року курси OCW Scholar є складовою OpenCourseWare. Вони є відкритими ліцензованими навчальними матеріалами, що, на відміну від курсів дистанційного навчання, не забезпечують двосторонній зв'язок з викладачами і призначені для підтримки самоосвіти. Вони відрізняються від інших курсів OCW тим, що навчальні матеріали добре структуровані, розташовані в логічній послідовності із зростанням ступені складності.

На сьогодні розроблено курси в галузях природничих наук, математики, техніки та економіки, оскільки ці курси, на думку розробників, мають найбільший попит серед студентів для організації самостійного навчання. До 2015 року передбачено в цілому розробити 15 курсів [44].

Незважаючи на відсутність взаємодії з викладачами MTI або студентами, ці курси включено в пілотний проект OpenStudy – навчальну соціальну мережу, де студенти задають питання, надають допомогу, і спілкуються з іншими студентами, що вивчають подібні розділи вищої математики (рис. 2.15).

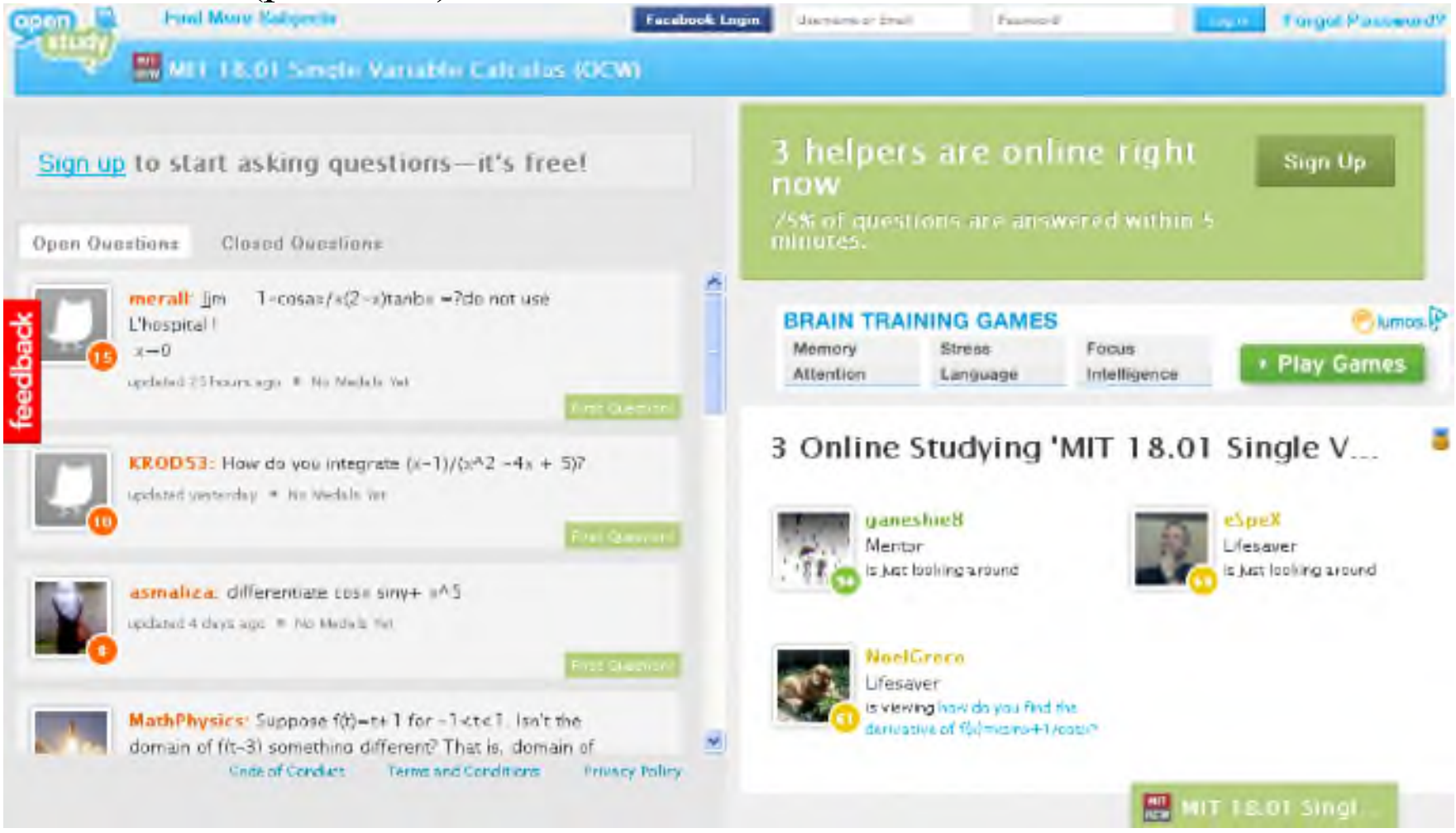

Рис. 2.15. Сторінка соціальної мережі Відкрите навчання (OpenStudy – http://openstudy.com/)

Як зазначено на офіційному сайті OpenStudy, його призначення в тому, щоб зробити світ єдиною великою навчальною групою, незалежною від закладів, місцерозташування користувачів та їх походження. Через OpenStudy користувачі курсів OCW можуть взаємодіяти один з одним, спільно працювати над завданнями і відповідати на запитання.

У зв'язку з тим, що матеріали курсів OCW Scholar структуровані для самостійного вивчення, з ними легко працювати в будь-якому темпі, зручному для студентів.

Не всі OCW курси опубліковані в форматі Scholar OCW, тому що ці курси є значно більш дорогими і трудомісткими у виробництві, а тому обмежують можливості їх створення. Виробництво OCW курсів фінансується за рахунок гранту від Stanton Foundation, що був створений Ф. Стентоном (Frank Stanton), який по праву вважається одним з найвидатніших керівників в історії електронних комунікацій [44].

Версія курсу Числення однієї змінної, розміщена у OCW Scholar, містить такі он-лайн-ресурси та засоби навчання:

1) відеолекції та відеоприклади розв’язання задач, доступні для завантаження та онлайн-перегляду (зокрема, на YouTube) (подібно до OCW);

2) статичні задачі з розв’язаннями (подібно до OCW);

3) Java-аплети («mathlets» – онлайн-програми математичного призначення), що демонструють ключові концепції курсу (рис. 2.16).

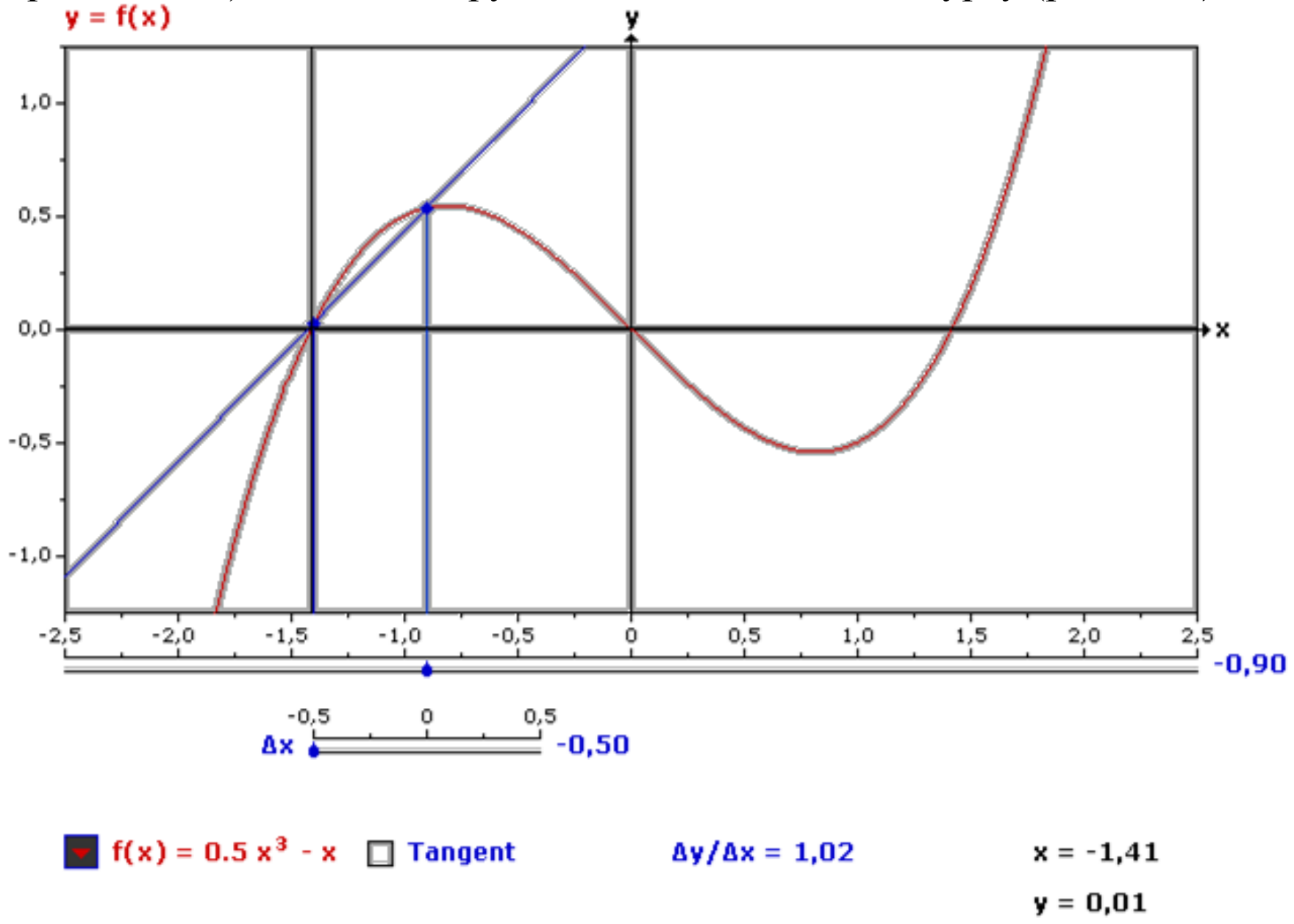

Рис. 2.16. Лекційна демонстрація до лекції «Поняття похідної функції однієї змінної»

На сайті кафедри вищої математики MTI студентам надано мінімальний набір додаткового програмного забезпечення, вказуючи його місце знаходження (табл. 2.2): текстові (Pine, elm, mail) та графічні (Thunderbird) поштові клієнти; текстові (links) та графічні (Firefox) Web-браузери; математичні пакети (MATLAB, Mathematica, Maple, R, MAGMA); текстові редактори (emacs, vim, nano, Xemacs, Gedit, OpenOffice); мережні засоби (VoIP – Skype, IM – Empathy, FTP – KFTP); наукові текстові процесори (LaTeX, Kile); компілятори (C/C++ – gcc, icc, Fortran – gfortran, ifort); системи відображення документів (PDF – Acrobat Reader, PostScript – evince, DVI – xdvi) [9].

*Таблиця* 2.2

**Додаткове програмне забезпечення курсів вищої математики в МТІ**

| **Програма** | **Відомості** | **Посилання на** |
|---|---|---|
| Pine | Поштовий клієнт для Unix-подібних систем. Програма для читання текстової електронної пошти з вікна терміналу типу Pine. | Інструкція з використання Pine |
| elm | Програма для читання текстової електронної пошти з терміналу типу elm. | Інструкція з використання Elm (текстовий файл) |
| mail | Програма для читання текстової електронної пошти з терміналу типу mail. | Керівництво з використання пошти |
| Thunderbird | Програма для графічної підтримки електронної пошти, що знаходиться у розділі Додатки → Інтернет (Applications → Internet). | Сторінка підтримки Thunderbird |
| Firefox | Сучасний графічний Web-браузер (перегляд web-сторінок), що знаходиться в розділі Додатки → Інтернет (Applications → Internet). | Сторінка підтримки Firefox |
| links | Перегляд web-сторінок на основі набраних посилань з терміналу типу links. | З терміналу набрати «man links» для доступу до сторінок довідника |
| MATLAB | Технічне обчислювальне програмне забезпечення, що знаходиться на сайті кафедри вищої математики → MATLAB (Math Department → MATLAB), або набрати в терміналі MATLAB. Для використання на | Швидкий запуск MATLAB |

| Програма | Відомості | Посилання на |
|---|---|---|
| | персональних машинах можна замовити копію з MTI (подано посилання з вказівкою «сертифікати вимагаються»). | |
| Mathematica | Технічне обчислювальне програмне забезпечення, що знаходиться на сайті кафедри вищої математики → Mathematica (Math Department → Mathematica), або в терміналі набрати mathematica. Для використання на персональних машинах викладачі та співробітники можуть отримати Mathematica з Help Desk у комп'ютерному відділі. Студенти можуть замовити копію з MTI (подано посилання з вказівкою «сертифікати вимагаються»). | Центр навчання Mathematica |
| Maple | Технічне обчислювальне програмне забезпечення, що знаходиться на сайті кафедри вищої математики → Maple (Math Department → Maple), або в терміналі набрати maple або xmaple для графічного інтерфейсу. Для використання на персональних машинах, можна отримати копію maple з MTI (подано посилання з вказівкою «сертифікати вимагаються»). | Maple Mastery (сертифікати вимагаються) |
| R | Статистичні обчислення і графіка. З терміналу набрати R. | Вступ до R |
| MAGMA | Програмне забезпечення комп'ютерної алгебри. З терміналу набрати magma. Це програмне забезпечення буде працювати тільки на серверах runge, hypatia, dwork, або archimedes. | Підручник MAGMA |
| emacs | Багатофункціональний розширюваний текстовий редактор. | Керівництво з Emacs |
| vim | Багаторежимний текстовий редактор, що добре налаштовується. | Документація до vim |
| nano | Простий консольний текстовий редактор для Unix та Unix-подібних | Довідник з nano. |

| Програма | Відомості | Посилання на |
|---|---|---|
| | операційних систем. | |
| xemacs | Графічна версія emacs. | |
| Gedit | Невеликий текстовий редактор, що потребує незначних потужностей; використовується для машин з GNOME; знаходиться в Стандартні → Текстовий редактор (Accessories → Text Editor). | Керівництво Gedit |
| OpenOffice.org | Пакет офісних програм для розробок презентацій, електронних таблиць, креслень. Знаходиться в розділі Додатки → Офіс (Applications → Office). | |
| Skype | Програмне забезпечення з закритим кодом, що забезпечує шифрований голосовий зв'язок і відеозв'язок через Інтернет. | Використання Skype |
| Empathy | Програмне забезпечення для миттєвого обміну повідомленнями. Сумісне з gchat, AIM, і майже зі всіма сервісами обміну миттєвими повідомленнями. Знаходиться в розділі Додатки → Інтернет (Applications → Internet). | Empath FAQ |
| KFTP | Графічний FTP клієнт. Знаходиться в розділі Додатки → Інтернет (Applications → Internet). | |
| LaTeX | Система підготовки математичних документів. | Вікіпідручник LaTeX |
| Kile | Текстовий редактор для LaTeX. | Довідник Kile |
| gcc | Компілятори GNU C/C++. | Використання компілятора GNU |
| icc | Компілятор Intel C. До початку компілювання треба конфігурувати середовище командою iccvars.sh. Для 32-бітної версії набрати icc32, для 64-бітної – icc64. | Приступаємо до роботи |
| gfortran | Компілятор GNU Fortran. | |
| ifort | Компілятор Intel Fortran. До початку компілювання треба конфігурувати | Приступаємо до роботи |

| Програма | Відомості | Посилання на |
|---|---|---|
| | середовище командою ifortvars.sh. Для 32-бітної версії набрати ifort32, для 64-бітної – ifort64. | |
| acroread | Для перегляду документів в PDF форматі. Також відомий як Acrobat Reader. | |
| evince | Простий PDF/PS переглядач документів. | |
| xdvi | Переглядач dvi-файлів. | |

За результатами проведеного дослідження найбільш популярними СКМ, що використовують у навчанні вищої математики у США є: Mathematica, MATLAB, Maple, GAUSS, Scilab, Mathcad, Maxima та Sage. У табл. 2.3 подано загальну характеристику систем комп'ютерної математики, що застосовуються в МТІ при навчанні вищої математики, спираючись на роботи [16; 20; 42; 70; 160; 213; 272; 325; 326].

*Таблиця* 2.3

**Загальна характеристика систем комп'ютерної математики, що використовують у навчанні вищої математики у США**

| Назва (остання версія) | Розробник | Основні характеристики | |
|---|---|---|---|
| | | ***переваги*** | ***недоліки*** |
| MATLAB R2013b (Version 8.2) (5 вересня 2013) | The MathWorks | – універсальна СКМ для здійснення швидких і точних чисельних розрахунків у різних предметних галузях;<br>– відкритість і розширюваність;<br>– підтримка 3D-графіки;<br>– сумісність з різними операційними платформами;<br>– підтримує роботу з базами даних | – відсутність у ядрі підтримки розв'язання нерівностей, діофантових рівнянь, рекурентних співвідношень;<br>– вимогливість до апаратних ресурсів інформаційної системи |
| Mathematica 9.0.1 (30 січня 2013) | Wolfram Research | – убудована підтримка паралельних обчислень;<br>– статистичний аналіз моделей; | – складність синтаксису;<br>– уявлення про дані як про сукупність окремих |

| Назва (остання версія) | Розробник | Основні характеристики | |
|---|---|---|---|
| | | *переваги* | *недоліки* |
| | | – унікальність 3D-графіки;<br>– сумісність з різними операційними платформами;<br>– висока швидкість виконання математичних операцій та обчислень;<br>– має розвинений графічний інтерфейс, що надає можливість працювати з багатьма документами;<br>– підтримує роботу з базами даних | виразів, що знижує продуктивність розв'язання складних задач |
| Maple 17 (13 березня 2013) | Waterloo Maple Inc. | – найкраще символьне ядро;<br>– висока точність обчислень;<br>– уведення (з 11 версії) математичних виразів у природній математичній нотації;<br>– структурованість документу;<br>– інтуїтивно зрозумілий інтерфейс;<br>– має розвинений графічний інтерфейс, що надає можливість працювати з багатьма документами;<br>– в останніх версіях є панель для розпізнавання символів, введених від руки, для швидкого пошуку потрібної команди або символу; | – незручність у роботі з великою кількістю числових даних;<br>– уявлення про дані як про сукупність окремих виразів, що знижує продуктивність розв'язання складних задач |

| **Назва (остання версія)** | **Розробник** | **Основні характеристики** | |
|---|---|---|---|
| | | ***переваги*** | ***недоліки*** |
| | | – взаємодія з CAD-системами, що надає можливість візуалізувати складні об'єкти, створювати креслення на підставі отриманих результатів обчислень та інше;<br>– підтримує роботу з базами даних | |
| R 3.0.2 (25 вересня 2013) | Р. Іака (Ross Ihaka), Р. Джентльмен (Robert Clifford Gentleman) | – розповсюджується безкоштовно;<br>– широкі можливості для проведення статистичних аналізів, включаючи лінійну і нелінійну регресію, класичні статистичні тести, аналіз часових рядів (серій), кластерний аналіз і інше;<br>– може використовуватися для матричних розрахунків;<br>– можливість підключати код, написаний на C, C++, або Fortran;<br>– містить засоби для візуалізації результатів обчислень (2-вимірні, 3-вимірні графіки, діаграми, гістограми, діаграми (схеми) Ганта тощо);<br>– функція Sweave, що надає можливість інтеграції і виконання коду R у документах, написаних за допомогою | – для роботи використовується командний інтерпретатор |

| **Назва (остання версія)** | **Розробник** | **Основні характеристики** | |
|---|---|---|---|
| | | ***переваги*** | ***недоліки*** |
| | | LaTeX з метою створення динамічних звітів;<br>– користувачі можуть розширювати функціонал за рахунок написання нових функцій;<br>– надають можливість створювати високоякісні графіки з різними атрибутами, включаючи також математичні формули і символи | |
| MAGMA V2.19-10 (28 серпня 2013) | Computational Algebra Group, School of Mathematics and Statistics, University of Sydney | – надає можливість розв'язувати завдань з алгебри, теорії чисел, геометрії і комбінаторики;<br>– працює на Unix-подібних і Linux операційних системах, а також Windows;<br>– має розвинений графічний інтерфейс, що надає можливість працювати з багатьма документами;<br>– підтримує роботу з базами даних | – відсутність у ядрі підтримки інтегрування, розв'язання нерівностей, диференціальних рівнянь, рекурентних співвідношень |

Таким чином, на сайтах OpenCourseWare і Stellar викладачі вищої математики використовують велику кількість програм, у тому числі Web-додатків, для закріплення основних понять, що були представлені в аудиторії.

При підготовці до іспитів студентам рекомендується використовувати Інтернет-ресурси, щоб знайти зразки прикладів для опрацювання; перегляд матеріалів, таких як екзамени за попередні роки або інші зразки онлайн завдань. Домашні та самостійні завдання базуються на онлайн програмах та математичних інструментах, що розміщені в Інтернеті, або у системі Stellar, або на Web-сайті конкретного курсу, так що студенти можуть переглядати і опрацьовувати їх у зручному для них режимі. У той час як більшість завдань може включати в себе традиційні обчислення (розв'язання задач без використання програмних засобів), існує тенденція серед викладачів розробляти свої власні завдання з використанням наявних програмних засобів. У цілому, доступні програми забезпечують альтернативний підхід для студентів, щоб мати можливість розібратися і зрозуміти матеріали з розділів вищої математики. І викладачі, і студенти погоджуються з тим, що використання таких програм, а також онлайн-інструментів, сприяє розвитку їх математичних компетентностей, а також їх розумінню матеріалу, необхідного для фундаментальних знань в їх подальшій кар'єрі [75].

Розглянуті ІКТ використовують не тільки для організації процесу навчання в МТІ, а й широко використовуються в інших технічних університетах США.

Так, наприклад, платформа для інтерактивного спілкування Piazza використовується в Каліфорнійському університеті в Берклі, Технологічному інституті Джорджії, Університеті штату Техас в Остіні, Стенфордському університеті, Корнельському університеті, Мічиганському університеті, Університеті Ватерлоо, Університеті штату Пенсільванія, Університеті штату Іллінойс в Урбана-Шампейн, Бостонському університеті, Університеті Карнегі-Меллона, Прінстонському університеті, Гарвардському університеті, Університеті Центральної Флориди та інші [103].

Уведення нових засобів навчання неминуче призводить до зміни методів та форм навчання. Із виникненням та впровадженням засобів ІКТ у процес навчання вищої математики, відбувається поєднання традиційних форм навчання із компютерно-орієнтованими (на основі систематичного, послідовного і логічного використання ІКТ у процесі навчання [317, 259-269]).

У моделі комбінованого навчання навчальний процес вищої математики містить основні традиційні форми організації навчального процесу: лекція, практичне заняття, самостійна робота студентів, контроль якості знань; крім того традиційні форми доповнюються формами організації дистанційного, мобільного навчання та навчання за допомогою Інтернет і мультимедіа: навчальні матеріали курсу, онлайн-спілкування, індивідуальні та групові онлайн-проекти, віртуальна класна кімната, аудіо- та відеолекції, анімація та симуляція, мобільні тренінги тощо.

На *лекційному* занятті відбувається виклад основного навчального матеріалу курсу. Пояснення, що супроводжується лекційними демонстраціями із використанням програмних засобів, сприяє кращому усвідомленню та засвоєнню розглянутих понять. Записи лекційних занять, розміщені в мережі Інтернет, надають студентам можливість допрацювати незрозумілий матеріал. У разі неможливості присутності викладача на лекційному занятті або винесенні матеріалу на самостійне опрацювання, традиційну лекцію доцільно замінити проведенням Web-конференції (вебінару, відеолекції), мультимедіа лекції.

*На практичному занятті* студенти закріплюють отриманні знання та набувають навички розв'язування практичних задач. Завдання практичного заняття мають бути підібрані із урахуванням напряму підготовки студентів, роблячи наголос на доцільності вивченого матеріалу у подальшій професійній діяльності студентів. Питання та задачі практичного заняття, приклади для самостійного розв'язання бажано довести до відома студентів заздалегідь, розмістивши їх у мережі.

Доцільно на практичному занятті використовувати засоби ІКТ для проведення експрес-тестування для визначення рівня засвоєного матеріалу студентами. Для проведення самоконтролю студентами необхідно надати чіткі інструкції по роботі з СКМ, що надасть їм можливість виконувати перевірку правильності відповідей у розв'язаних задачах і формуватиме у студента навички роботи з СКМ у подальшій професійній діяльності.

За Ю. В. Триусом, однією із форм організації навчальної діяльності, пов'язаної із набуттям студентами практичних навичок у відповідній галузі знань за допомогою засобів ІКТ, є *комп'ютерно-орієнтоване практичне заняття*: «будується на поєднанні традиційних і комп'ютерних форм навчання та контролю знань і орієнтовано на розв'язування задач, що забезпечують наступність між практичними, лабораторними і лекційними заняттями на основі внутрішніх і міждисциплінарних логічних зв'язків» [317, 261-262].

*Контроль якості знань* можна проводити в аудиторії та за її межами засобами Інтернет. Використання СКМ, тестових програм та тестових модулів у системах підтримки навчання забезпечує можливість студентам проведення самоконтролю, а викладачам автоматизувати процес навчання.

*Онлайн-спілкування* – це ефективний спосіб отримати своєчасну консультацію із проблемних питань. Таке спілкування здійснюється різними засобами ІКТ, серед яких можна виділити використання електронної пошти, платформ для інтерактивного спілкування, соціальних мереж для спілкування в групах. Використання Web-інструментів надає можливість викладачам підтримувати комунікацію зі студентами одним із наступних способів: у режимі реального часу, динамічно вирішуючи питання; із швидкою реакцією на питання за допомогою миттєвих повідомлень; засобами електронної пошти підтримувати затяжні діалоги (наприклад, протягом усього семестру чи навчального року ) і проводити віртуальні дискусії у класі за допомогою дошки обговорення або Web-конференцій.

Проведення своєчасної *консультації* є необхідною умовою для правильної організації самостійної роботи студентів. Використання засобів ІКТ значно розширює можливості викладача. Викладач може проводити індивідуальні та групові консультації за межами ВНЗ. Використовуючи відео записи занять, викладач може продемонструвати приклади розв'язання задач, що надасть можливість кращому засвоєнню навчального матеріалу студентами.

Для якісного процесу навчання необхідно обрати такі *засоби навчання*, що забезпечать найкращу інтеграцію аудиторного та позааудиторного навчання. Серед таких засобів можна виділити паперові та електронні книги та посібники; навчальні матеріали, що розміщені в мережі; навчальні аудіо- та відеоматеріали; системи підтримки математичної діяльності, тестові системи, тренажери для відпрацювання навичок та вмінь, електронні бібліотеки та довідники тощо.

Аналізуючи розглянуті форми організації навчання, можна зробити висновок, що в процесі навчання вищої математики студентів вищих технічних навчальних закладів відбувається гармонійне поєднання аудиторного навчання з віртуальним, очного навчання з дистанційним, індивідуального та групового навчання. При цьому викладач та студенти підтримують як синхронну, так і асинхронну комунікацію.

## 2.2 Динаміка розвитку теорії та методики використання інформаційно-комунікаційних технологій у навчанні вищої математики студентів інженерних спеціальностей у вищих навчальних закладах Сполучених Штатів Америки

Питанню інтеграції ІКТ у систему вищої освіти приділяється значна увага, воно розглядається на всіх рівнях: від департаменту системи освіти до викладачів ВНЗ. Головною метою інтеграції ІКТ у систему освіти є втілення у життя напряму освіти впродовж всього життя. Перед викладачами технічних ВНЗ США стоїть проблема підготовки інженера, що здатний працювати в умовах швидкого розвитку інформаційних та комунікаційних технологій.

Аналіз літератури з питання використання ІКТ в освіті надав можливість проаналізувати різні концептуальні схеми. На рис. 2.17 наведено приклад схеми використання ІКТ в освіті.

Педагогічна діяльність викладача повинна підкріплюватись ІКТ.

Однак педагогічна практика визначається не тільки такими якостями викладача, як їх академічна кваліфікація та компетентність у сфері ІКТ, але також факторами ВНЗ і факторами системного рівня. Якщо вважати, що на навчальні досягнення студентів впливають набуті ними навички та уміння, які використовуються у процесі навчання вищої математики, слід визнати, що результати (уявні чи справжні) впливають на подальші педагогічні рішення викладача [51].

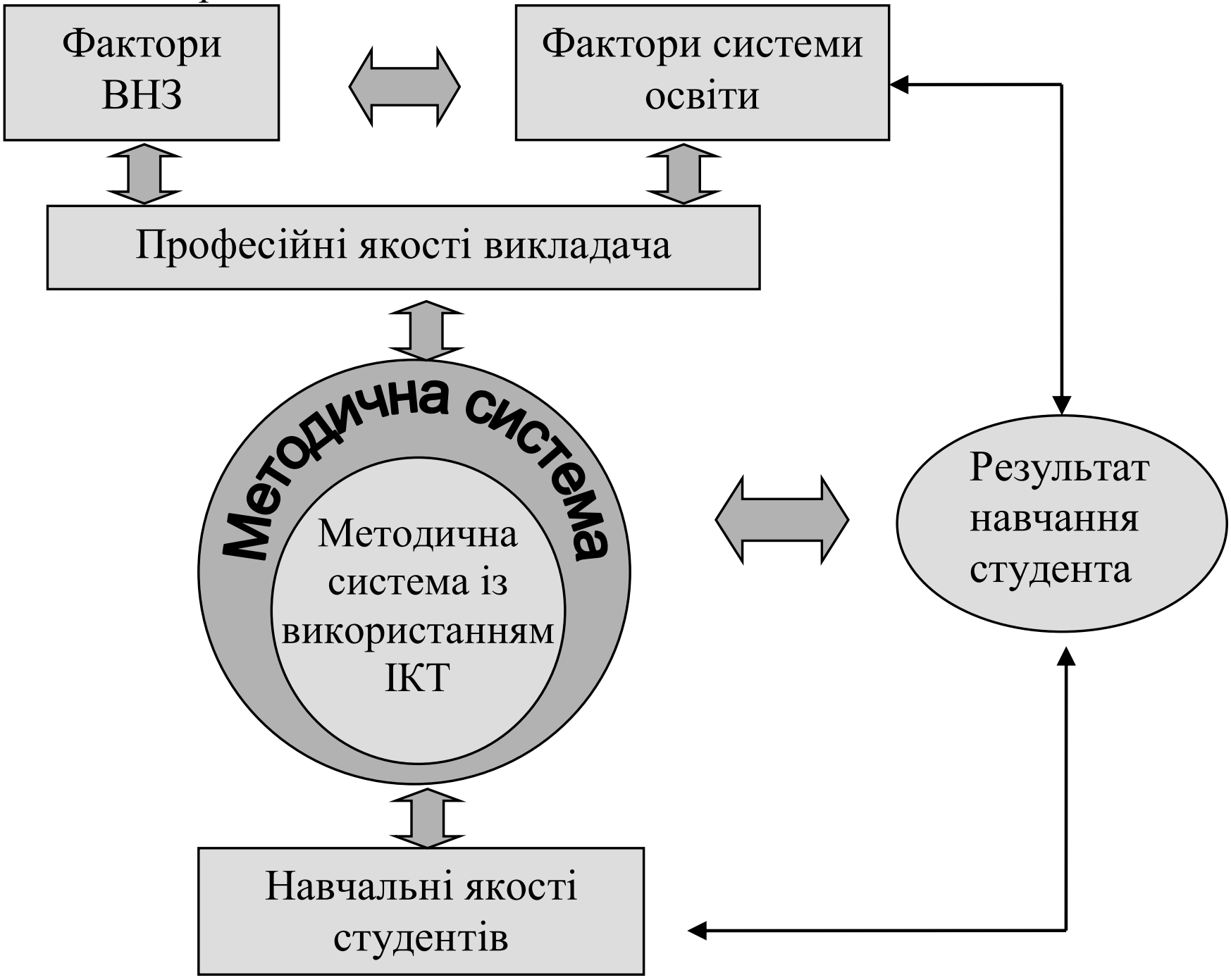

Рис. 2.17. Концептуальна схема використання ІКТ в освіті [51]

Для того, щоб інтеграція ІКТ у національну системи освіти стала ефективною, потрібно відповідне поєднання наступних політичних і практичних чинників [51]:

1) чіткі цілі та створення національної програми підтримки використання ІКТ в освіті;

2) допомога та стимулювання як державних, так і приватних навчальних закладів до придбання обладнання ІКТ (наприклад, шляхом цільового державного фінансування, включаючи кошти на технічне обслуговування; податкових знижок на обладнання ІКТ та програмне забезпечення для навчальних закладів; інвестицій або спонсорства досліджень з розвитку недорогого обладнання та програмного забезпечення ІКТ тощо);

3) пристосування навчальних програм до впровадження ІКТ, розвиток і придбання стандартних якісних електронних навчальних посібників та програмного забезпечення;

4) впровадження програм масової підготовки викладачів до використання ІКТ;

5) умотивованість викладачів та студентів організовувати процес навчання із використанням ІКТ;

6) адекватний рівень національного моніторингу та система оцінювання якості освіти, що надає можливість регулярно визначати результати та дієвість, а також заздалегідь виявляти недоліки з метою підвищення ефективності використання ІКТ в освіті.

Виданий Департаментом освіти США Національний план використання освітніх ІКТ у 2010 році – це є модель навчання, що базується на використанні ІКТ та включає в себе цілі і рекомендації в п’яти основних областях освітньої діяльності: навчання, оцінювання, викладацька діяльність, засоби і продуктивність [127].

Розглянемо, як у [127] інтерпретується ці основні області освітньої діяльності.

***Навчання***. Викладачі мають підготувати студентів до навчання впродовж всього життя за межами аудиторії, тому необхідно змінити зміст та засоби навчання для того, щоб відповідати тому, що людина повинна знати, як вона набуває знання, де і коли вона навчається, і змінити уявлення про те, хто повинен навчатися. В XXI столітті необхідно використовувати доступні ІКТ навчання для мотивації й натхнення студентів різного віку.

Складні і швидко змінні потреби світової економіки визначають необхідний зміст навчання та рівень знань студентів. Використання ІКТ надає можливість впливати на рівень знань студентів і розуміння ними навчального матеріалу.

На рис. 2.18 показана модель навчання, що базується на використанні ІКТ.

На відміну від традиційного навчання в аудиторії, де найчастіше один викладач передає один і той же навчальний матеріал всім студентам однаково, модель навчання із використанням ІКТ ставить студента у центр і дає йому можливість взяти під контроль своє індивідуальне навчання, забезпечуючи гнучкість у кількох вимірах. Основний набір стандартних знань, вмінь та навичок утворюють основу того, що всі студенти повинні вивчати, але, крім того, студенти та викладачі мають можливість вибору у навчанні: великі групи чи малі групи, діяльність у відповідності з індивідуальними цілями, потребами та інтересами.

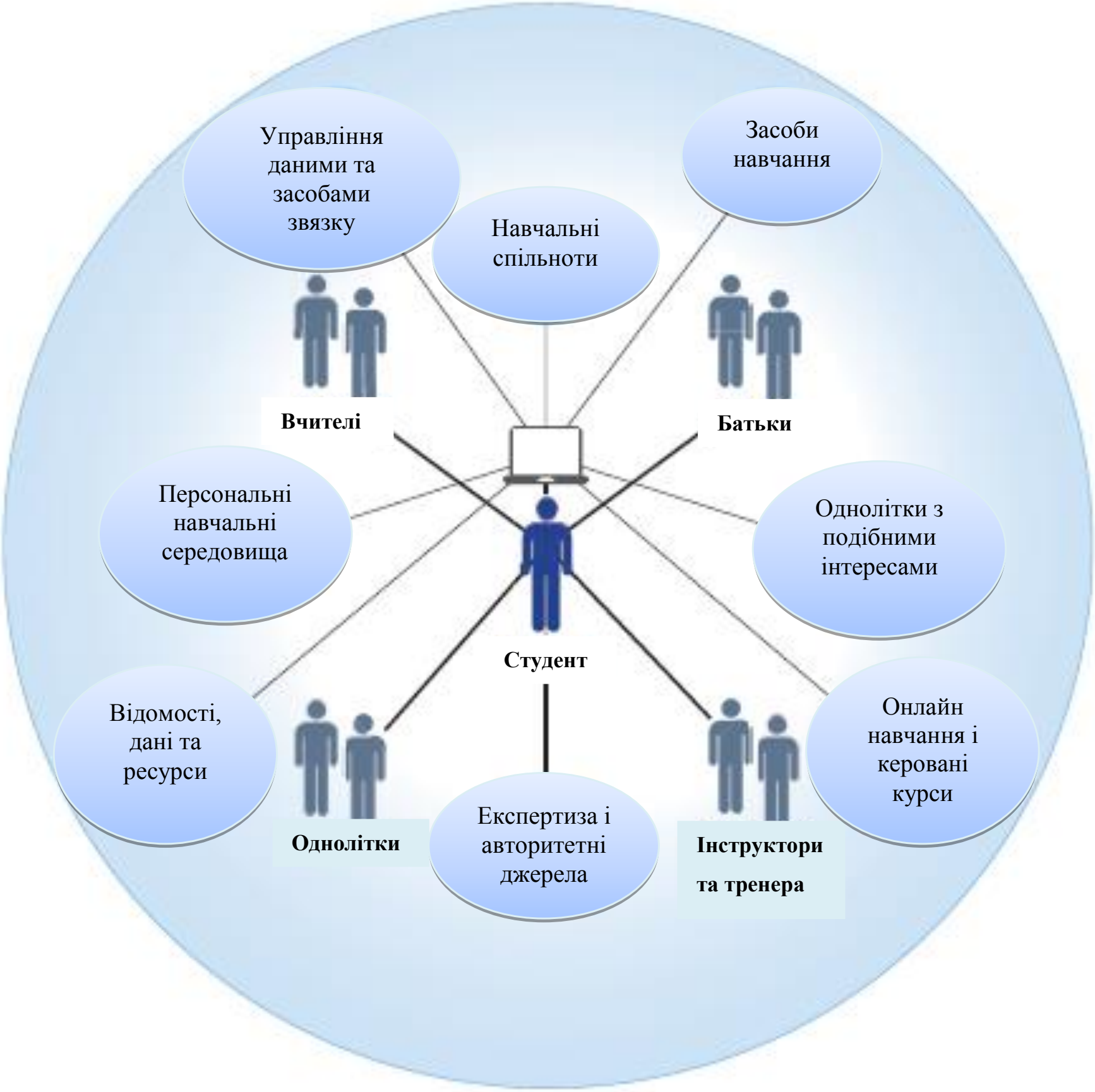

Рис. 2.18. Модель навчання із використанням ІКТ у США [127]

У цій моделі процес навчання підтримується ІКТ, що створюють умови для залучення навчальних середовищ та інструментів для ефективного розуміння і запам'ятовування навчального матеріалу. Залучення ІКТ до процесу навчання забезпечує доступ до більшого розмаїття набору навчальних ресурсів, ніж традиційне викладання в аудиторіях, при цьому до процесу організації та контролю за навчанням долучаються викладачі ВНЗ, батьки, експерти і наставники за межами аудиторії. Процес навчання стає індивідуалізованим та диференційованим, що створює умови для побудови персонального навчального середовища, шляхом обміну навчальними відомостями не тільки з викладачами, а й з іншими студентами.

Для кожної дисципліни, хоча й існують стандарти змісту навчання,

модель навчання із використанням ІКТ дає зрозуміти, яким чином можна проводити навчання. Серед всіх можливих варіантів будується власний проект навчання, що розв'язує проблеми реальної значимості. Добре продумані плани індивідуального навчання з дисциплін допомагають студентам отримати необхідні знання, а також підтримують розробку спеціалізованого адаптивного досвіду, що може бути застосований і в інших дисциплінах.

Використання ІКТ також дає студентам можливість контролювати та управляти власним навчанням. Студент, керуючи електронними навчальними портфоліо, може залишати постійні записи про навчання і допомагати студентам групи розвивати самосвідомість, встановлювати власні цілі навчання, висловлювати свої погляди, сильні і слабкі сторони, досягнення, нести за них відповідальність. Викладачі можуть використовувати ці записи для оцінки рівня розвитку досягнень студентів.

Згідно з Національним планом освітніх ІКТ Департаменту освіти США індивідуалізація, диференціація і персоналізація стали ключовими поняттями у сфері освіти [127].

*Індивідуалізація* розглядається як підхід, що визначає потрібний темп у навчанні різних студентів. При цьому навчальні цілі однакові для всіх студентів, але студенти можуть вивчати матеріал з різною швидкістю в залежності від їх потреб у навчанні.

*Диференціація* розглядається як підхід, що враховує особливості різних студентів. Цілі навчання однакові для всіх студентів, але методи навчання варіюються в залежності від уподобань кожного студента або потреб студентів.

*Персоналізація* розглядається як підхід, за якого вивчаються навчальні потреби студентів із урахуванням навчальних переваг та конкретних інтересів різних студентів. Персоналізація включає в себе диференціацію та індивідуалізацію.

Викладачі постійно мають визначати необхідний рівень знань та вмінь студентів. На сучасному етапі в навчанні, крім знань з конкретних дисциплін, студент має володіти критичним мисленням, умінням комплексно вирішувати проблеми, бути готовим до співпраці. Крім того, студент має відповідати таким категоріям: *інформаційна грамотність* (здатність ідентифікувати, знаходити, оцінювати та використовувати дані для різних цілей); *медіаграмотність* (здатність до використання і розуміння засобів масових відомостей, а також ефективного спілкування, використовуючи різні типи носіїв); *можливість оцінювати і використовувати ІКТ*, відповідно вести себе в соціально прийнятних Інтернет-спільнотах, а також розумітися в питанні навколишньої конфіденційності та безпеки. Все це вимагає базового розуміння самих ІКТ і здатності використовувати їх в повсякденному житті.

Навчаючи, викладачі мають враховувати те, що студенти не можуть вивчити все, що їм потрібно знати в житті, і економічна реальність така, що більшість людей будуть змінювати місце роботи протягом всього життя. Тому необхідно набувати та розвивати адаптивні навички навчання, що поєднують зміст знань із можливістю дізнатися щось нове.

В організації процесу навчання використовують Web-ресурси і технології, що створюють студентам умови для дослідження та підтримки співпраці і спілкування у роботі. Для студентів такі інструменти надають нові навчальні можливості і допомагають їм подолати реальні проблеми, розробити стратегії пошуку, оцінити рівень довіри до Web-сайтів і авторів, а також спілкуватися за допомогою мультимедіа. Так, при вивченні вищої математики використання динамічних графіків та статистичних програм полегшує сприйняття складних тем, роблячи їх більш доступними для студентів і допомагаючи студентам краще розуміти той навчальний матеріал, що має відношення до їх спеціальності.

ІКТ можуть бути використані для забезпечення більших можливостей у навчанні у поєднанні з традиційним методам навчання. Із використанням ІКТ можна подавати навчальні матеріали, вибираючи різні типи носіїв, та сприяти засвоєнню знань, вибираючи інструменти, до яких відносяться тематичні карти, хронології, що забезпечують візуальний зв’язок між наявними знаннями і новими ідеями.

Із використанням ІКТ розширюються засоби навчання студентів:

– забезпечується допомога студентам у процесі навчання;

– надаються інструменти для спілкування у процесі навчання (це можна зробити через Web-інтерфейс мультимедіа, мультимедійні презентації тощо);

– сприяють виникненню Інтернет-спільнот, де студенти можуть підтримувати один одного у дослідженнях та розвивати більш глибоке розуміння нових понять, обмінюватися ресурсами, працювати разом поза аудиторією і отримувати можливості проведення експертизи, керівництва та підтримки.

Для стимулювання викладачами мотивації до взаємодії викладачів і студентів із використанням ІКТ можна:

– підвищувати інтерес та увагу студентів;

– підтримувати зусилля та академічну мотивацію;

– розробляти позитивний імідж студента, який постійно навчається.

Оскільки людині впродовж всього життя доводиться навчатися, то ключовим фактором постійного і безперервного навчання є розуміння можливостей ІКТ. Використання ІКТ у навчанні надає студентам прямий доступ до навчальних матеріалів та надає можливість будувати свої знання організовано і доступно. Це дає можливість студентам взяти під контроль і персоналізувати їхнє навчання.

***Оцінювання.*** В системі освіти на всіх рівнях планується використовувати ІКТ для планування змісту навчального матеріалу, що є актуальним на момент навчання, і використовувати ці дані для безперервного вдосконалення навчальних програм.

Оцінювання, що проводиться сьогодні у ВНЗ, спрямоване на аналіз якості кінцевого результату процесу навчання. При цьому не відбувається оцінювання рівня розвитку мислення студента в процесі навчання, а це могло б допомогти їм навчатися краще.

У процесі організації оцінювання, на думку [127], необхідно вносити зміни, що включають в себе пошук нових та більш ефективних способів оцінювання. Необхідно проводити оцінювання в процесі навчання таким чином, щоб студенти підвищували свій якісний рівень знань шляхом мотивації зі сторони працедавців.

Президентом США Б. Обамою (Barack Hussein Obama II) у 2009 році було запропоновано губернаторам штатів і розробникам державних освітніх стандартів проводити оцінювання таким чином, щоб не вимірювати вміння студента виконувати стандартні завдання, а оцінювати володіння ними такими навичками, як розв’язування нестандартних задач, критичне мислення, підприємництво, творчість. Вимірювання цих складних навичок вимагає проектування і розробку таких моделей оцінювання, що показують весь спектр стандартних знань і вмінь студента. Когнітивні дослідження і теорії дають багато моделей і уявлень про те, як відбувається розуміння та осмислення студентами ключових понять навчальної програми.

Існує багато прикладів використання ІКТ для комплексного оцінювання знань студентів. Ці приклади ілюструють, як використання ІКТ змінило характер опитування студентів, воно залежить від характеру викладання та апробації теоретичного матеріалу. Впровадження ІКТ надає можливість представити навчальний матеріал, системи, моделі і дані різними способами, що раніше були недоступними. Із залученням ІКТ у процес навчання можна демонструвати динамічні моделі систем; оцінювати навчальні досягнення студентів, запропонувавши їм проводити експерименти із маніпулюванням параметрів, записом даних та графіків і описом їх результатів.

Ще однією перевагою використання ІКТ для оцінювання є те, що за їх допомогою можна оцінити навчальні досягнення студента в аудиторії та за її межами.

У рамках проекту «Національна оцінка освітніх досягнень» (The National Assessment of Educational Progress – NAEP) розроблено і представлено навчальні середовища, що надають можливість проводити оцінювання навчальних досягнень студентів при виконанні ними складних завдань і вирішенні проблемних ситуацій.

Використання ІКТ для проведення оцінювання сприяє поліпшенню якості навчання. На відміну від проведення підсумкового оцінювання, використання корекційного оцінювання (тобто оцінювання, що надає можливість студенту побачити та виправити свої помилки в процесі виконання запропонованих завдань, наприклад, тестування з фізики, запропоноване Дж. Р. Мінстрелом (J. R. Minstrell) [33] – http://www.diagnoser.com), може допомогти підвищити рівень знань студентів.

Під час аудиторних занять викладачі регулярно намагаються з'ясувати рівень знань студентів, проводячи фронтальне опитування. Але це надає можливість оцінити лише незначну кількість студентів, не оцінивши знання та розуміння навчального матеріалу всієї аудиторії студентів. Для вирішення такої проблеми вивчається можливість використання різних технологій на аудиторних заняттях в якості «інструменту» для оцінювання. Одним із прикладів є використання тестових програм, що пропонують декілька варіантів відповідей на питання, до складу яких включено як істинні, так і неправдиві відповіді. Студенти можуть отримати корисні відомості із запропонованих відповідей на подібні питання, якщо вони ретельно розроблені.

При навчанні студентів із використанням засобів мережі Інтернет існують різні варіанти використання доступних Інтернет-технологій для проведення оцінювання. Використовуючи онлайн програми, можна отримати детальні дані про рівень досягнень студентів, що не завжди можливо в рамках традиційних методів навчання. При виконанні завдань студентами програмно можна з'ясувати час, що витрачають студенти на виконання завдань, кількість спроб на розв'язання завдань, кількість підказок даних студенту, розподіл часу в різних частинах даного завдання. Прикладом такої системи є вільно розповсюджувана онлайн платформа ASSISTment (http://www.assistments.org) [7], що призначена для підтримки навчання з математики. Педагоги отримують докладні звіти про вміння студентів зі 100 математичних навичок, а також їх точність, швидкість виконання, звернення за допомогою, кількість спроб розв'язування. За допомогою системи ASSISTment можна оцінити якість знань кожного студента, визначити проблеми засвоєння матеріалу та продемонструвати аудиторії етапі правильного розв'язання завдання.

У моделі навчання, де студенти самі обирають доступні засоби навчання, оцінювання виступає в новій ролі – визначення рівня знань студента з метою розробки подальшого унікального плану навчання для конкретного студента. У використанні такого адаптивного оцінювання забезпечується диференціація навчання.

Однією з проблем, що виникає при використанні нових технологій для проведення оцінювання, є час і витрати на їх розробку, визначення достовірності і надійності оцінювання, його працездатності. Цю проблему можна вирішити також за допомогою використання ІКТ. Розроблені завдання можна випробовувати автоматично, розмістивши їх у навчальному Web-середовищі, де тисячі студентів можуть виконати завдання в режимі онлайн. Дані, зібрані таким чином, можна використати для з’ясування успішності проведення такого роду оцінювання навчальних досягнень студентів і вдосконалення завдань перед широкомасштабним використанням їх у процесі навчання.

Із використанням Інтернет-технологій для проведення оцінювання можна створити базу даних про рівень знань студента в різних галузях (подібно до портфоліо), що може в подальшому бути використана в професійній діяльності студентів.

У системі освіти в США на всіх рівнях застосовуються можливості Інтернет-технологій для вимірювання знань студентів, що надає можливість використовувати дані оцінки для безперервного вдосконалення процесу навчання.

Для проведення вдосконалення процесу навчання в плані розглядається необхідність наступних дій [127]:

1) на рівні держави та округів необхідно проектувати, розробляти і здійснювати оцінювання, що дає студентам, викладачам та іншим зацікавленим сторонам своєчасні та актуальні дані про навчальні досягнення студентів для підвищення рівня та навчальної практики студентів;

2) науковий потенціал викладачів освітніх установ, а також розробників Інтернет-технологій використовувати для поліпшення оцінювання в процесі навчання. Із використанням Інтернет-технологій можна проводити вимірювання ефективності навчання, забезпечуючи систему освіти можливостями проектування, розробки та перевірки нових і більш ефективних методів оцінювання;

3) проведення наукових досліджень для з’ясування того, як із використанням технологій, таких як моделювання, навчальні середовища, віртуальні світи, динамічні ігри та навчальні програми, можна заохочувати та підвищувати мотивацію студентів при оцінюванні складних навичок;

4) проведення наукових досліджень і розробок із проведення

об'єктивного оцінювання (без оцінювання сторонніх здібностей студента). Для того, щоб оцінювання було об'єктивним, воно повинно вимірювати потрібні якості та не повинно залежати від зовнішніх факторів;

5) перегляд практики, стратегії і правил забезпечення конфіденційності та захисту даних про одержані оцінки студентів, при одночасному забезпеченні моделі оцінювання, що включає в себе постійне збирання і обмін даними для безперервного вдосконалення процесу оцінювання.

***Викладацька діяльність.*** Викладачі можуть індивідуально або колективно підвищувати свій професійний рівень, використовуючи всі доступні технології. Вони можуть отримати доступ до даних, змісту, ресурсів, відомостей і передового досвіду навчання, що сприяє розширенню можливостей викладачів і надихає їх на забезпечення більш ефективного навчання студентів.

Багато викладачів працюють поодинці, не спілкуючись з колегами або викладачами з інших ВНЗ. Професійний розвиток зазвичай проводиться на короткому, фрагментарному і епізодичному семінарі, що пропонує мало можливостей для використання отриманих матеріалів на практиці. Основна аудиторна робота викладача на практиці зводиться до перевірки набутих знань студентами. Багато викладачів не мають часу та стимулу для постійного підвищення свого професійного рівня щороку.

Так само, як використання ІКТ може допомогти поліпшити процес навчання та оцінювання, використання ІКТ може допомогти краще підготуватися до ефективного викладання, підвищити професійний рівень. Використання ІКТ створює умови для переходу до нової моделі – зв'язаного навчання.

Під *зв'язаним навчанням* будемо розуміти організований двосторонній процес навчальної взаємодії викладача та студента, в результаті якого свій рівень знань та досягнень підвищують не тільки студенти, а й викладачі.

У зв'язаному навчанні викладачі мають отримувати повний доступ до даних про процес навчання студентів та аналітичні інструменти для обробки цих даних. Їм необхідно забезпечити комунікацію зі своїми студентами, доступ до даних, ресурсів і систем підтримки навчання, що дозволить їм створювати, управляти і оцінювати досягнення навчання студентів у позааудиторний час. Викладачі також можуть отримати доступ до ресурсів, що надають можливість покращити їх практику викладання, постійно підвищуючи свій професійний рівень (рис. 2.19). Як у моделі навчання студентів (рис. 2.18), описаної раніше, викладачі приймають участь в особистому навчанні, розвиваючи себе як високоваліфікованого педагога.

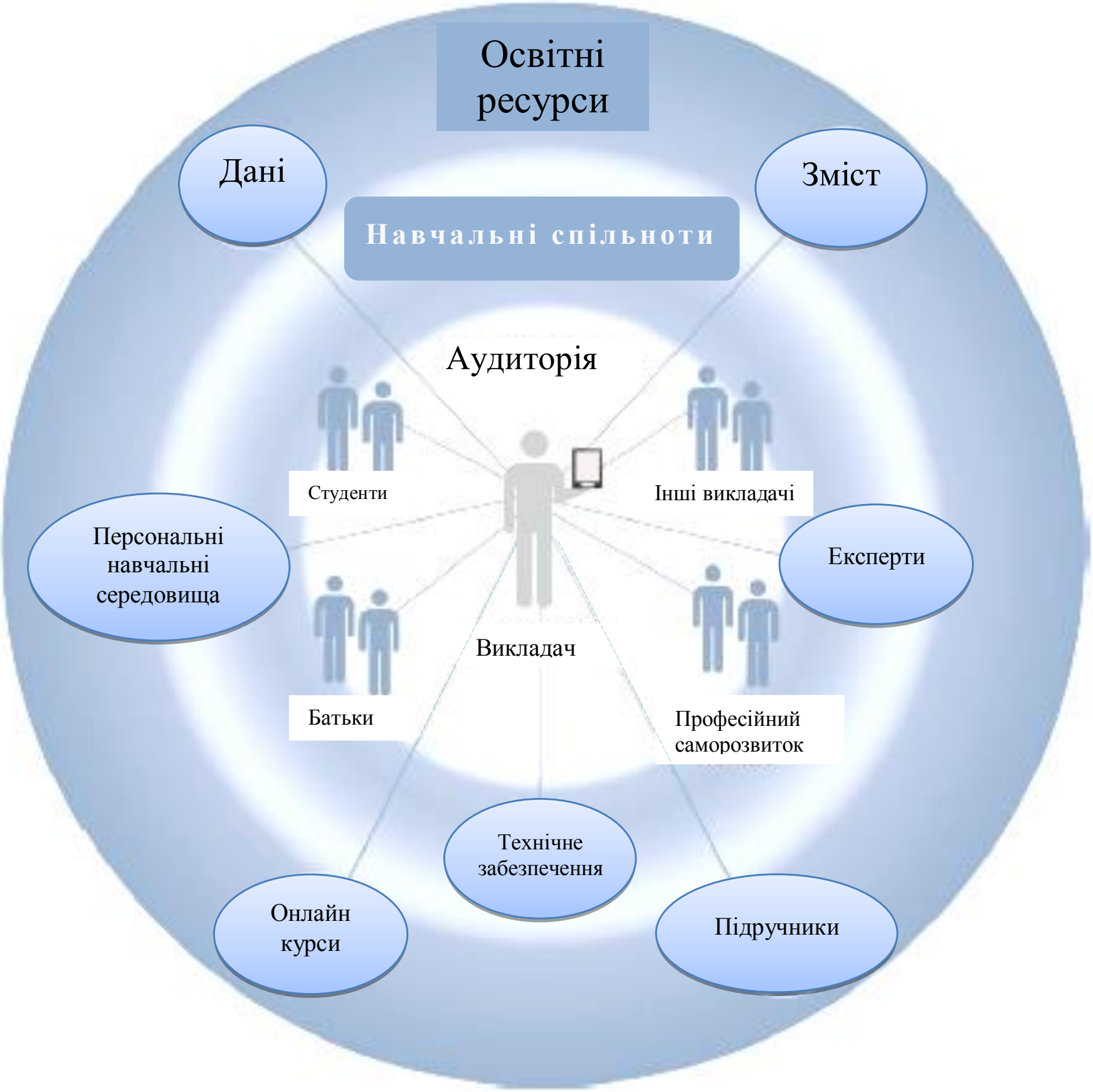

Рис. 2.19. Модель зв’язаного навчання у США [127]

У такому навчанні неефективні семінари для підвищення кваліфікації можна замінити на професійне навчання педагога, що є узгодженим і безперервним та може об’єднати матеріали найбільш ефективних семінарів, забезпечуючи розширення можливостей самоосвіти, оперативність, доступність та зручність притаманну онлайн-навчанню.

Оскільки середовище навчання постійно ускладнюється, у зв’язаному навчанні забезпечується підтримка викладачів в організації та управлінні навчальним процесом. Безкоштовні системи управління навчанням з відкритим доступом вже сьогодні широко використовуються у ВНЗ США. Такі інструменти надають можливість викладачам визначати зміст навчальних матеріалів, навчальних планів, завдань, дискусій, більше зосереджуючись на потребах студентів. Використання онлайн-середовищ надає можливість створювати віртуальні класи, у яких викладачі та студенти можуть взаємодіяти по-новому один з одним, управляючи змістом курсу. Такі навчальні середовища надають можливість залишати зразки домашніх, екзаменаційних робіт, проводити дискусії та розміщувати мультимедійні об'єкти, забезпечувати зворотній зв'язок зі студентами.

Інтернет-середовища сприяють підвищенню активності студентів у навчанні. Куратори груп можуть приєднатися до віртуальних класів для знаходження даних про навчальні досягнення студентів. Батьки або члени інших установ-партнерів також мають змогу входити в систему для віртуальної взаємодії з класом або доповнювати матеріали навчального середовища.

У зв'язаному навчанні пропонується широкий спектр можливостей забезпечення індивідуалізації навчання. Використовуючи всі доступні технології, викладачі отримують можливість спілкуватися зі студентами, дізнаватися про їх навчальні потреби та можливості, тим самим підвищуючи мотивацію студентів до навчання, можливості для творчості і самовираження. В Інтернет-спільнотах студенти отримують можливість спілкуватися з викладачами та однолітками зі всього світу, а викладачі можуть заохочувати студентів до навчальної діяльності.

У зв'язку з тим, що онлайн-навчання стає все більш важливим у системі вищої освіти, це створює необхідність в викладачах, які є кваліфікованими в онлайн-навчанні та обізнані у передовому досвіді викладання із застосуванням онлайн-технологій. Вирішальне значення для розв'язання цієї потреби, забезпечуючи при цьому ефективне навчання, є впровадження відповідних стандартів для онлайн-курсів та навчання.

Досягнення поставлених цілей навчання можливо за рахунок забезпечення підтримки особистого та групового професійного розвитку викладачів із використанням онлайн-технологій, що надає можливість отримати їм доступ до навчальних матеріалів, змісту, ресурсів, досвіду навчання і підвищує ефективність навчального процесу. Для досягнення цієї мети необхідні такі дії [127]:

1) розширення можливостей викладачів у доступі до даних про технології, що використовуються для збереження, зміни та обробки навчальних матеріалів у довільний зручний для них час. Використання Інтернет-технологій надає можливості викладачам використовувати ресурси і підтримувати спілкування серед викладачів різних ВНЗ у межах держави та за її межами, не витрачаючи на це багато зусиль та часу.

Сьогодні викладачі повинні мати доступ до всіх можливих технологічних ресурсів, що забезпечить більш привабливі та ефективні умови для навчання студентів;

2) використання соціальних мереж і платформ для створення співтовариств викладачів, що забезпечують професійний ріст педагогів протягом всієї їх кар'єри, доступу до засобів навчання та ресурсів, що роблять професійне навчання своєчасним і актуальним, а також передової діяльності викладачів, що постійно вдосконалюється і розвивається. В Інтернет-спільнотах повинна надаватися можливість викладачам проходити онлайн-курси, а також надаватися доступ до платформ та засобів проектування та розробки матеріалів для підвищення кваліфікації;

3) використання Інтернет-технологій для підтримки онлайн-навчання студентів, підвищуючи ефективність навчання та створюючи умови для організації комбінованого навчання. Студентам, для самореалізації, потрібні навчальні курси з різних дисциплін, що викладаються на високому рівні та відповідають сьогоденню. В онлайн навчанні необхідно передбачити використання кращих навчальних практик, доступних всім студентам;

4) забезпечення якісної підготовки викладачів з використання технологій в їх навчальній практиці, підвищувати їх грамотність в області цифрових технологій і надання їм можливості створювати завдання для студентів, які покращують навчання, оцінювання та навчальну практику. Використання технологій може допомогти не лише залучити та мотивувати студентів до навчання, а також може допомогти в процесі підготовки і безперервного навчання педагогів. Використання технологій також повинно бути невід'ємною складовою методики навчання курсів, а не розглядатися як набуття окремого навику;

5) підготовка кваліфікованих фахівців з онлайн-навчання. Оскільки онлайн-навчання набуває все більшого значення в системі освіти, необхідно вивчати досвід комбінованого та Інтернет-навчання, що є більш персоналізованими і втілюють кращі практики для залучення всіх студентів. Це створює необхідність підготовки викладачів, які є кваліфікованими в створенні навчальних курсів із залученням новітніх технологій. Вирішальне значення для задоволення цієї потреби відіграють відповідні стандарти для онлайн-курсів.

***Засоби.*** Всі студенти та викладачі мають отримати доступ до усіх можливих засобів навчання в довільний зручний для них час і місці.

Використання Інтернет-технологій надає можливість відійти від традиційного навчання студентів лише в аудиторіях, лабораторіях, бібліотеках, перемістивши робочі місця додому в зручний для них час, де вони мають доступ до мережі Інтернет. Засоби навчання покликані підтримувати навчання суспільства впродовж всього життя.

До *засобів навчання* автори плану віднесли освітні ресурси, стратегії і стійкі моделі для постійного поліпшення зв'язку, сервери, програмне забезпечення, системи управління та адміністрування. Питання забезпечення засобами навчання вимагає участі та співпраці викладачів усіх дисциплін і типів установ по всім рівням освіти.

Використання засобів навчання надає можливість розкрити нові способи збирання і поширення знань, заснованих на мультимедійних нерухомих і рухомих зображеннях, текстах, аудіододатках, що працюють на різних пристроях. Це звільняє навчання від монотонної моделі передавання навчальних матеріалів (з книги або викладача до студента) і сприяє підвищенню навчальної активності студентів.

На більш оперативному рівні засоби навчання об'єднують і надають можливість отримати доступ до даних з різних джерел, забезпечуючи при цьому належний рівень безпеки та конфіденційності. До засобів навчання відносять також комп'ютерну техніку, дані і мережу, інформаційні ресурси, сумісне програмне забезпечення, проміжні послуги та інструменти, а також пристрої, що забезпечують підтримку зв'язку у міждисциплінарних групах викладачів, відповідальних за використання засобів, розвиток, технічне обслуговування і управління навчанням, змінюючи підходи до викладання та навчання.

Використання відкритих освітніх ресурсів відіграє важливу роль в процесі навчання, в проведенні науково-дослідної роботи. Відкриті освітні ресурси розглядають як навчальні або наукові ресурси, що є загальнодоступними засобами Інтернет і відкриті всім користувачам і студентам. До відкритих освітніх ресурсів відносять тести, тренажери, презентації, електронні бібліотеки підручників, ігри, дистанційні курси, навчальні відео- та аудіоматеріали, а також інші електронні навчальні засоби.

Для забезпечення процесу навчання потрібними засобами навчання, необхідно почати перехід до наступного покоління архітектури обчислювальної системи. ВНЗ необхідно розглянути варіанти скорочення кількості серверів. Віртуалізація надає можливість на одному сервері працювати з декількома додатками безпечно і надійно, так що можна зменшити кількість серверів у мережах, скорочуючи витрати і зробивши мережу менш складною, простішою в управлінні.

Перехід до хмарних обчислень включає в себе перехід від закупівель і технічного обслуговування серверів у центрах обробки даних місцевих покупців програмного забезпечення в хмарі. Використання хмарних обчислень стає більш бажаним і можливим у зв'язку з тим, що так відбувається наближення до розвитку Інтернет-технологій і забезпечуються потреби у більш потужних платформах з меншою затратою ресурсів. Використання хмарних обчислень може забезпечити підтримку академічних та адміністративних послуг, необхідних для навчання та освіти. Це надає можливість студентам і викладачам отримувати доступ до ресурсів навчання, використовуючи різні пристрої.

***Продуктивність.*** У системі освіти на всіх рівнях проводиться перебудова процесу і структури навчання з метою використання ІКТ для покращення результатів навчання та більш ефективного використання часу, грошей і персоналу. Використання ІКТ для планування, управління, контролю і звітності надасть можливість отримати точне і повне уявлення про рівень знань студентів, фінансові показники в системі освіти на всіх рівнях, що має важливе значення для підвищення її продуктивності.

Для організації перебудови в системі освіти необхідно забезпечити можливість проведення таких дій:

1) розробити та прийняти загальне тлумачення поняття продуктивності праці в сфері освіти і більш актуальних і значущих показників результатів, разом з поліпшенням стратегій і використанням ІКТ для управління витратами, зокрема для закупівель. Безперервне вдосконалення та підвищення продуктивності праці не можливо без визначення й вимірювання його поточного стану. Необхідно визначитися з бажаними результатами навчання, забезпечення можливості встановлювати відношення результатів навчання до витрат на нього, що може бути відстежене з плином часу;

2) переосмислити основні чинники в системі освіти, які перешкоджають використанню ІКТ для покращення навчання, починаючи з поточної практики організації роботи студентів і викладачів. Для реалізації всіх можливостей використання ІКТ для поліпшення продуктивності праці студентів, необхідно усунути перешкоди, що стоять на шляху впровадження ІКТ. У системі освіти необхідно виявити і переосмислити основні її положення. Деякі з них включають в себе вимірювання якості освіти кількістю аудиторного часу, структурою окремих навчальних дисциплін, організацією навчання в групах з однаковою кількістю студентів, а також використанням факультативних занять;

3) розробити інформативну діагностику використання ІКТ в освіті по всій території США. Поточні дані про використання освітніх та інформаційних технологій в системі освіти складаються з відомостей про здійснені закупки і кількість комп'ютерів, підключених до Інтернету. Дуже мало даних про те, як ІКТ фактично використовуються для підтримки викладання, навчання та оцінювання. Отримавши дані про те, як і коли використовуються ІКТ, можна визначити ефективність використання, що забезпечить покращення результатів і продуктивності системи освіти;

4) розробити технологію використання ІКТ у системі освіти таким чином, щоб отримати випускників шкіл, готових до вступу у ВНЗ, і висококваліфікованих спеціалістів, готових до професійної діяльності. Співпраця адміністрації шкіл та адміністрації вищих навчальних закладів у напряму використання ІКТ в освіті має вирішальне значення для вирішення цієї проблеми.

У розглянутій моделі навчання, що базується на використанні ІКТ, припускається, що в системі освіти будуть розроблені і прийняті дії по забезпеченню вільного доступу до всіх можливих технологій навчання, що забезпечить ефективність процесу навчання, оцінювання, викладацької діяльності, використання засобів навчання. Крім того, використання технологій значно спрощує процес управління в освіті, підвищуючи її продуктивність.

Використання у процесі навчання інтегрованих систем, що забезпечують можливість доступу до навчальних матеріалів, важливе значення для здійснення індивідуалізованого, диференційованого та персоналізованого навчання. Інтегрована система повинна надавати можливість: отримувати доступ до необхідних освітніх ресурсів; налаштовувати ресурси відповідно до рівня знань та вмінь студентів; вибирати необхідні налаштування для забезпечення підтримки навчання студентів та можливості його доступу до навчальних матеріалів.

Використання ІКТ має підтримувати:

– систематичний аналіз рівня знань, вмінь та навичок студентів (у тому числі набуття вмінь проводити дослідження, розмірковування, проектування та комунікації), визначені стандартами освіти та види критеріїв, необхідних для визначення набутих навичок;

– визначення критеріїв оцінювання завдань і ситуацій, які б забезпечили необхідні вміння, знання та навички студентам;

– управління процесом оцінювання;

– розробка та застосування правил та статистичних моделей для одержання надійних висновків про набуті вміння, знання та навички студентів за підсумками оцінювання завдань.

Перед системою освіти стоїть завдання визначити і затвердити принципи побудови дієвої та ефективної системи навчання в режимі онлайн і комбінованого навчання, що підтримують розвиток знань студентів всередині і за межами академічних груп. Існує багато підтверджень того факту, що навчання може бути прискорене та вдосконалене із використанням онлайн навчання, перебудовою навчальних програм, а також шляхом надання можливості комунікації для поліпшення процесу навчання.

Одним із перспективних підходів до організації процесу навчання вищої математики є модель інтеграції технологій навчання: традиційного та дистанційного, мобільного та навчання за допомогою Інтернет і мультимедіа. Інтеграція аудиторної та позааудиторної роботи в процесі навчання можлива за рахунок використання педагогічних технологій та сучасних ІКТ, зокрема, засобів дистанційного, мобільного навчання та навчання за допомогою Інтернет і мультимедіа. Для того, щоб процес інтеграції був ефективним, викладачу необхідно управляти, регулювати та контролювати навчальну діяльність студентів [272, 89].

Традиційна форма навчання у вищий технічній школі реалізується за денною формою навчання та передбачає ознайомлення студентів із навчальним матеріалом, з найважливішими проблемами курсу, що потребують пояснення викладача; дискусії, роботу в групах (діяльність, пов'язана з безпосереднім контактом на різних рівнях); контрольні роботи (деякі види проміжного тестування рівня сформованості тієї чи іншої навички краще проводити за технологіями дистанційного навчання); захист проектів (при очному захисті бажано подавати необхідні матеріали на сайті). Основними формами за традиційної моделі навчання у ВНЗ є лекції, практичні та семінарські заняття, лабораторні роботи, заняття з контролю та перевірки знань.

Дистанційна форма навчання полягає в самостійному оволодінні поданого навчального матеріалу, дослідницькій діяльності з використання ресурсів Інтернет; виконання додаткових завдань, що сприяють засвоєнню навчального матеріалу; тестів, лабораторних та практичних робіт; спільного виконання завдань творчого характеру; дистанційних консультацій викладача та інше засобами ІКТ [272].

Саме використання комбінованого навчання демонструють викладачі при навчанні вищої математики студентів інженерних спеціальностей у США. Впровадження комбінованого навчання в процес навчання вищої математики надає викладачам можливість використовувати великий набір засобів для підтримки процесу навчання, що стає неперервним (рис. 2.20). Організація комбінованого навчання – це кропітка робота викладачів та фахівців, які відповідають за постійне оновлення засобів, що підтримують процес навчання [3].

Викладачі розглядають свою роботу також як навчання, що спрямоване на задоволення унікальних потреб як своїх, так і студентів. Одним із головних завдань викладача є спрямування своїх знань на підтримку процесу навчання та консультацій студентів. Викладання спрямоване на численні можливості навчання, готуючи студентів до майбутньої професії [96].

Розглянемо реалізацію технології комбінованого навчання шляхом виокремлення окремих компонент моделі використання ІКТ у навчанні вищої математики студентів інженерних спеціальностей США.

Як зазначає Н. І. Євтушенко [171], побудова моделі навчання повинна відбуватися з урахуванням таких *етапів*:

1) входження в процес і вибір методологічних засад для моделювання, якісний опис предмета дослідження;

2) обґрунтування завдань моделювання;

3) конструювання моделі з уточненням залежності між основними елементами досліджуваного явища, визначенням параметрів об’єкта та критеріїв оцінки змін цих параметрів, вибір методик діагностики;

4) дослідження валідності моделі;

5) використання моделі у педагогічному експерименті;

6) змістова інтерпретація результатів моделювання.

Метод моделювання має бути теоретично і практично структурованим, відтворюючи дійсність у схематизованій і наочній формі. Розробка моделі використання ІКТ у навчанні вищої математики майбутніх інженерів в університетах США потребує реалізації *принципів*, відповідних до [143; 171]:

– відповідність загальній стратегічній меті фундаментальної підготовки інженера в технічному ВНЗ та відкритість щодо інформаційно-освітнього простору;– реалізація особистісно-орієнтованого навчання в умовах використання ІКТ у процесі навчання;

– цілісність, зумовлена інтегральними властивостями моделі, послідовність, неперервність, наступність структурно-функціональних елементів моделювання;

– достовірність моделювання (модель описує або проектує сегмент педагогічної реальності, очікуваний результат від її використання має, як правило, експериментальне підтвердження);

– ймовірність очікуваного ефекту (результат практичної реалізації розробленої моделі має імовірнісний характер і підкоряється законам статистики).

Проведений аналіз джерел [22; 51; 95; 131; 133; 219; 235; 294; 330] надав можливість виокремити дидактичні моделі використання ІКТ у процесі навчання вищої математики студентів інженерних спеціальностей США на різних етапах розвитку теорії та методики використання ІКТ у навчанні вищої математики студентів інженерних спеціальностей.

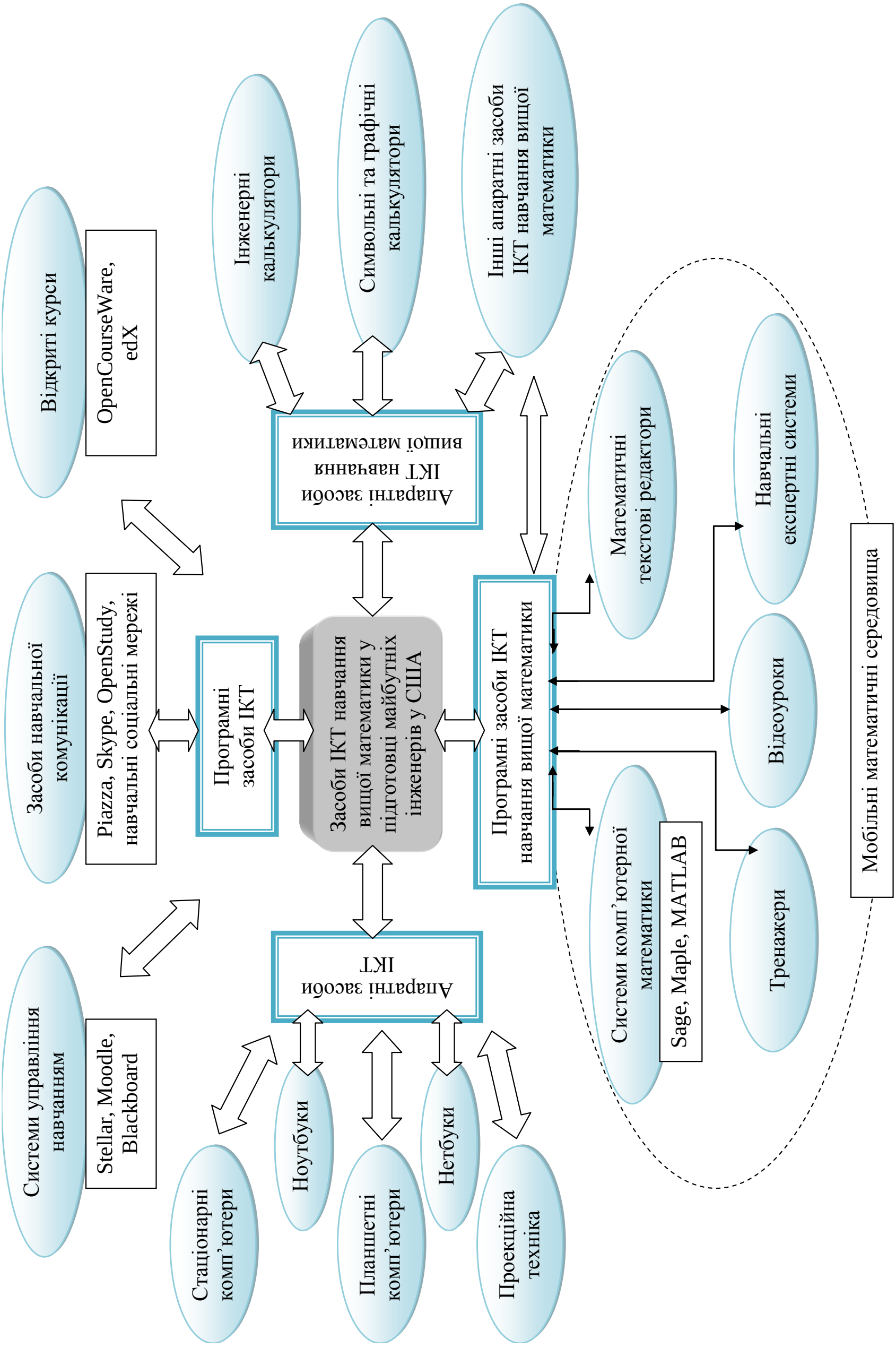

Рис. 2.20. Засоби ІКТ навчання вищої математики в США

Оскільки на першому та другому етапах розвитку теорії та методики використання ІКТ у навчанні вищої математики студентів інженерних спеціальностей у США використання ІКТ не носило системного характеру, тому виділимо дидактичні моделі використання ІКТ у навчанні вищої математики студентів інженерних спеціальностей США на таких етапах:

– *третьому етапі* – 1981–1989 рр. – пов'язаному із поширенням персональних комп'ютерів;

– *четвертому етапі* – 1989–1997 рр. – пов'язаному із створенням World Wide Web та використанням технологій Web 1.0;

– *п'ятому етапі* – 1997–2003 рр. – пов'язаному із появою та розробкою систем управління навчанням;

– *шостому етапі* – з 2003 р. по теперішній час – пов'язаному із перенесенням у Web-середовище засобів підтримки математичної діяльності та становленням і розвитком хмарних технологій навчання.

Запропоновані моделі визначають динаміку процесу використання засобів ІКТ у навчанні вищої математики. Моделі поєднують у собі описання типів засобів ІКТ за напрямом їх використання, провідних форм організації навчання вищої математики, а також специфіку впливу засобів ІКТ на процес навчання вищої математики.

Розглядаючи **третій етап** розвитку теорії та методики використання ІКТ у навчанні вищої математики студентів інженерних спеціальностей у США (пов'язаному із поширенням персональних комп'ютерів) (рис. 2.21), можна виділити такі *напрями використання засобів*: засоби підтримки математичної діяльності, засоби для подання навчальних відомостей та засоби для контролю знань та умінь студентів.

Використання персональних комп'ютерів у процесі **навчання** вищої математики студентів інженерних спеціальностей технічних ВНЗ США:

– сприяло індивідуалізації навчального процесу;

– надало можливість викладачам представляти навчальні матеріали студентам в електронному вигляді, створювати презентації лекцій;

– надало можливість організувати взаємодію викладача зі студентами таким чином, щоб досягти максимального засвоєння та усвідомлення навчального матеріалу студентами. На цьому етапі на заняттях з вищої математики студенти починають створювати комп'ютерні програми для проведення досліджень з різних тем математичних дисциплін [133];

– організація та підтримка роботи із застосуванням комп'ютера сприяє покращенню процесу взаємодії між викладачем та студентом.

Застосування комп'ютерів у процесі викладання вищої математики надає можливість виділити такі *переваги*:

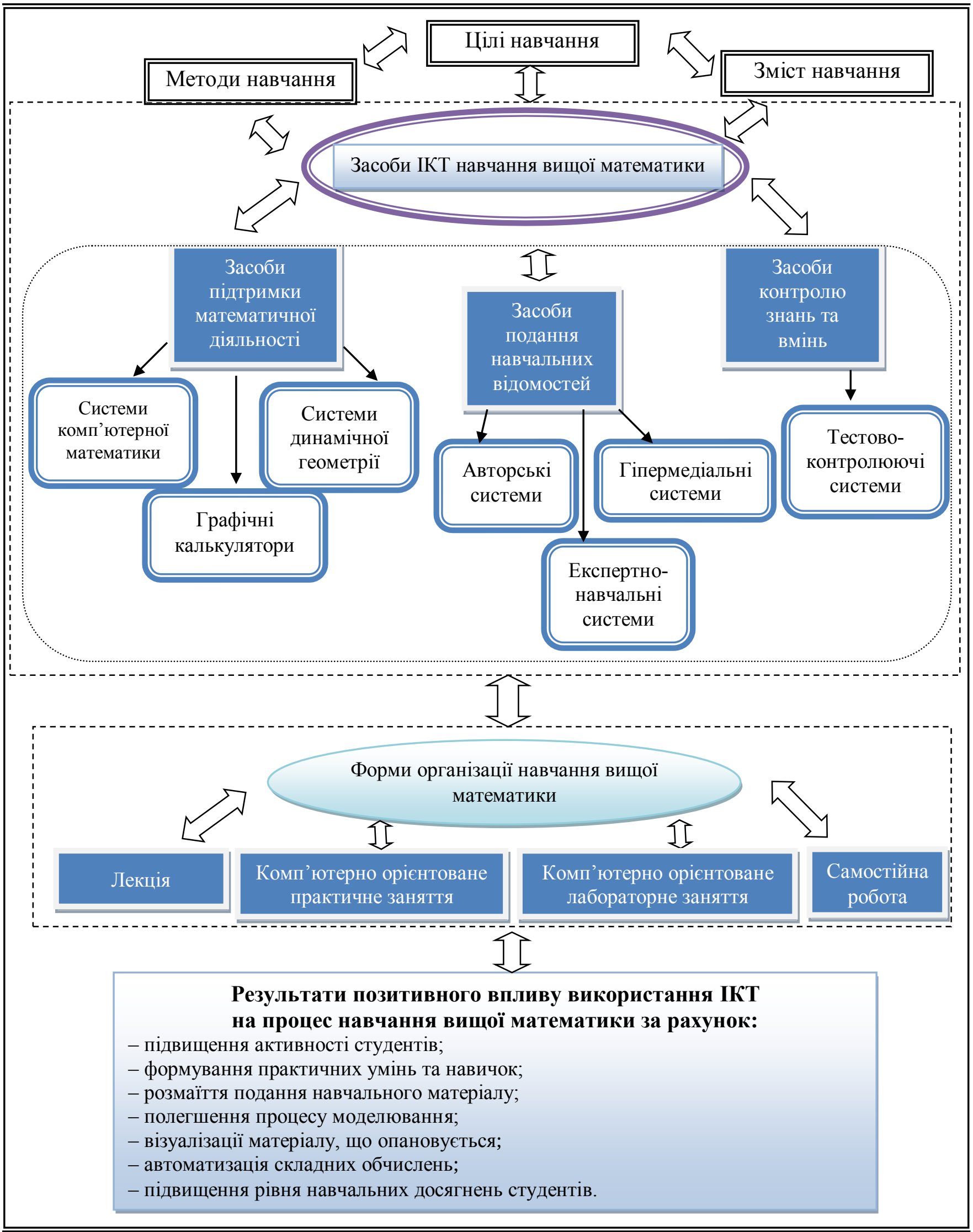

Рис. 2.21. Дидактична модель використання ІКТ при вивченні вищої математики в США на третьому етапі (1981–1989 рр.: поширення персональних комп'ютерів)

– підвищення активності студентів;
– формування практичних умінь та навичок;
– розмаїття подання навчального матеріалу;

– полегшення процесу моделювання;

– візуалізації матеріалу, що опановується;

– автоматизація складних обчислень;

– підвищення рівня навчальних досягнень студентів.

Використання комп’ютера у процесі навчання вищої математики сприяє вирішенню *часової проблеми* організації **контролю** знань студентів. Викладач здатний визначати рівень знань та умінь студентів, розуміння ними матеріалу, використовуючи програмні засоби (зокрема системи тестування для автоматизації).

Виділено **форми** організації навчання вищої математики із використанням засобів ІКТ: лекційне заняття, комп’ютерно орієнтоване практичне заняття, комп’ютерно орієнтоване лабораторне заняття, самостійна робота студентів.

Розглядаючи процес використання персональних комп’ютерів у навчанні вищої математики, можна стверджувати, що їх використання здійснює позитивний вплив на процес навчання вищої математики за рахунок:

– наявних переваг використання комп’ютера на заняттях, що значно скорочує час на розв’язання задач, полегшує процес моделювання;

– ознайомлення студентів з роллю та місцем вищої математики в наукових та прикладних дослідженнях, використовуючи доступні комп’ютерні програми;

– отримувати первісні навички математичного дослідження, оцінювати отримані результати, розвивати математичне мислення та підвищувати загальний рівень математичної культури студентів за допомогою програмних засобів;

– можливості самоосвіти студентів за рахунок самостійного аналізу наданої викладачем навчальної літератури з математики, поданої в електронному вигляді;

– розуміти місце отриманих знань та навичок для подальшого їх використання у навчанні та в майбутній професійній діяльності.

Аналіз експериментальних досліджень США [86] показав, що використання комп’ютера у процесі навчання математики надає можливість зменшити час на вивчення тем розділів вищої математики та підвищити активність студентів на заняттях.

*Напрями використання* засобів ІКТ на **<u>четвертому етапі</u>** (пов’язаному із створенням World Wide Web та використанням технологій Web 1.0) (рис. 2.22) розподіляються на *навчальні, контролюючі та управлінські*.

У процесі **навчання** вищої математики із використанням засобів ІКТ на четвертому етапі:

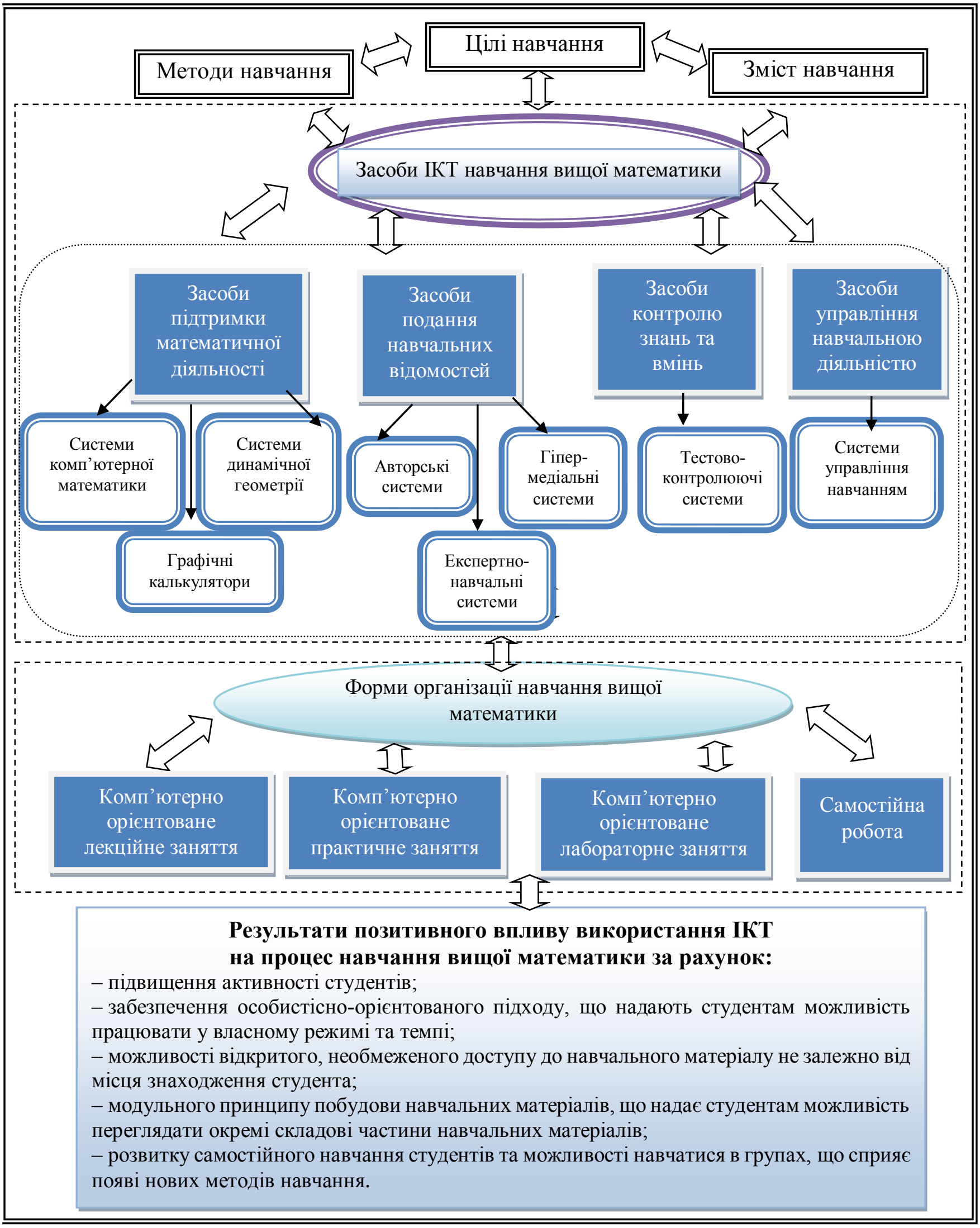

Рис. 2.22. Дидактична модель використання ІКТ при вивченні вищої математики в США на четвертому етапі (1989–1997 рр.: створення World Wide Web та використання технологій Web 1.0)

– викладачі мають можливість викладати навчальні матеріали курсу в мережу Інтернет і створювати власні сайти з поданням електронного наочного матеріалу для занять (схем, таблиць, рисунків тощо);

– студенти отримують доступ до глобальних ресурсів мережі;

– студенти отримують можливість опрацьовувати надані викладачем навчальні матеріали в зручний для себе час, зберігаючи їх на власному комп'ютері. Асинхронна природа створених викладачами онлайн-курсів надає студентам можливість самостійно вибирати час і місце для виконання завдань, надає можливість працювати у своєму власному темпі.

Використання засобів ІКТ надає можливість вирішити *часову та просторову проблеми* організації **контролю** знань студентів. Викладач здатний визначати рівень знань і вмінь студентів, розуміння ними матеріалу в зручний для нього час, використовуючи доступні засоби комунікацій (телефон, миттєві повідомлення, електронну пошту, програмні засоби тощо).

У системі **управління** процесом навчання:

– організовано контроль за доступом як до навчального матеріалу, так і до персональних даних студентів;

– організовано управління змістом, ресурсами, навчальною діяльністю студентів.

Основними призначеннями управління процесу навчання є: організація навчальної діяльності студентами, адміністрування процесу навчання, контроль за виконанням навчального плану.

Серед **форм організації** навчання вищої математики із використанням засобів ІКТ є: комп'ютерно орієнтоване лекційне заняття, комп'ютерно орієнтоване практичне заняття, комп'ютерно орієнтоване лабораторне заняття, самостійна робота студентів.

Із метою організації самостійної роботи студентів та підтримки комунікації викладачів та студентів використовується модель з інтеграцією ІКТ та технологій дистанційного навчання з технологіями традиційного навчання. В наслідок такої інтеграції викладач має можливість не тільки організовувати самостійну роботу студентів за допомогою забезпечення їх електронними навчальними матеріалами, а й керувати та оцінювати результати їх роботи на відстані, у зручний час, у зручному місці.

Використання ІКТ у процесі самостійної роботи студентів надає їм можливість дистанційно, через Інтернет, ознайомитися з навчальним матеріалом, який подається у вигляді різнотипних інформаційних ресурсів (текст, відео, анімація, презентація); отримувати консультаційну підтримку з курсу; одержувати електронні завдання для самоконтролю та підготовки до поточного контролю; дізнатися відомості про організацію навчального курсу.

Викладач має можливість: створювати власні авторські курси,

наповнювати їх необхідними навчальними матеріалами, давати консультації на відстані, надсилати повідомлення студентам, розподіляти завдання, налаштовувати різноманітні ресурси навчального курсу тощо. Доступ до ресурсів курсу – відкритий.

Електронні навчальні матеріали, розміщені в мережі, використовуються студентами для організації самостійної роботи, підготовки до виконання домашніх та екзаменаційних робіт.

Організація та підтримка процесу навчання вищої математики із застосуванням ІКТ надає можливість:

– створювати електронні освітні ресурси;

– розширювати доступ до цих ресурсів як студентам, так і викладачам;

– створювати організаційну та технологічну базу для впровадження технологій дистанційного навчання;

– підтримувати процес спілкування викладачів і студентів, студентів між собою.

Розглядаючи процес використання засобів ІКТ у навчанні вищої математики на четвертому етапі, можна стверджувати, що їх використання здійснює позитивний вплив на організацію навчальної діяльності студентів за рахунок:

– підвищення активності студентів;

– забезпечення особистісно-орієнтованого підходу, що надають студентам можливість працювати у власному режимі та темпі;

– можливості відкритого, необмеженого доступу до навчального матеріалу не залежно від місця знаходження студента;

– модульного принципу побудови навчальних матеріалів, що надає студентам можливість переглядати окремі складові частини навчальних матеріалів;

– розвитку самостійного навчання студентів та можливості навчатися в групах, що сприяє застосуванню методів інтерактивного навчання: кейс-методу, методу портфоліо, методу навчання у групах зі змінним складом, навчання у співробітництві.

Наукові дослідження того часу надали можливість зробити висновок: системне використання ІКТ у процесі організації навчання математики сприяє появі можливості застосування інтерактивних методів навчання, що робить студента активним учасником процесу навчання, знання якого формуються у процесі дослідження [123].

*Напрями використання* засобів ІКТ на **п’ятому етапі** (пов’язаному із появою та розробкою систем управління навчанням) (рис. 2.23) розподіляються, як і на попередньому етапі, на *навчальні, контролюючі та управлінські.*

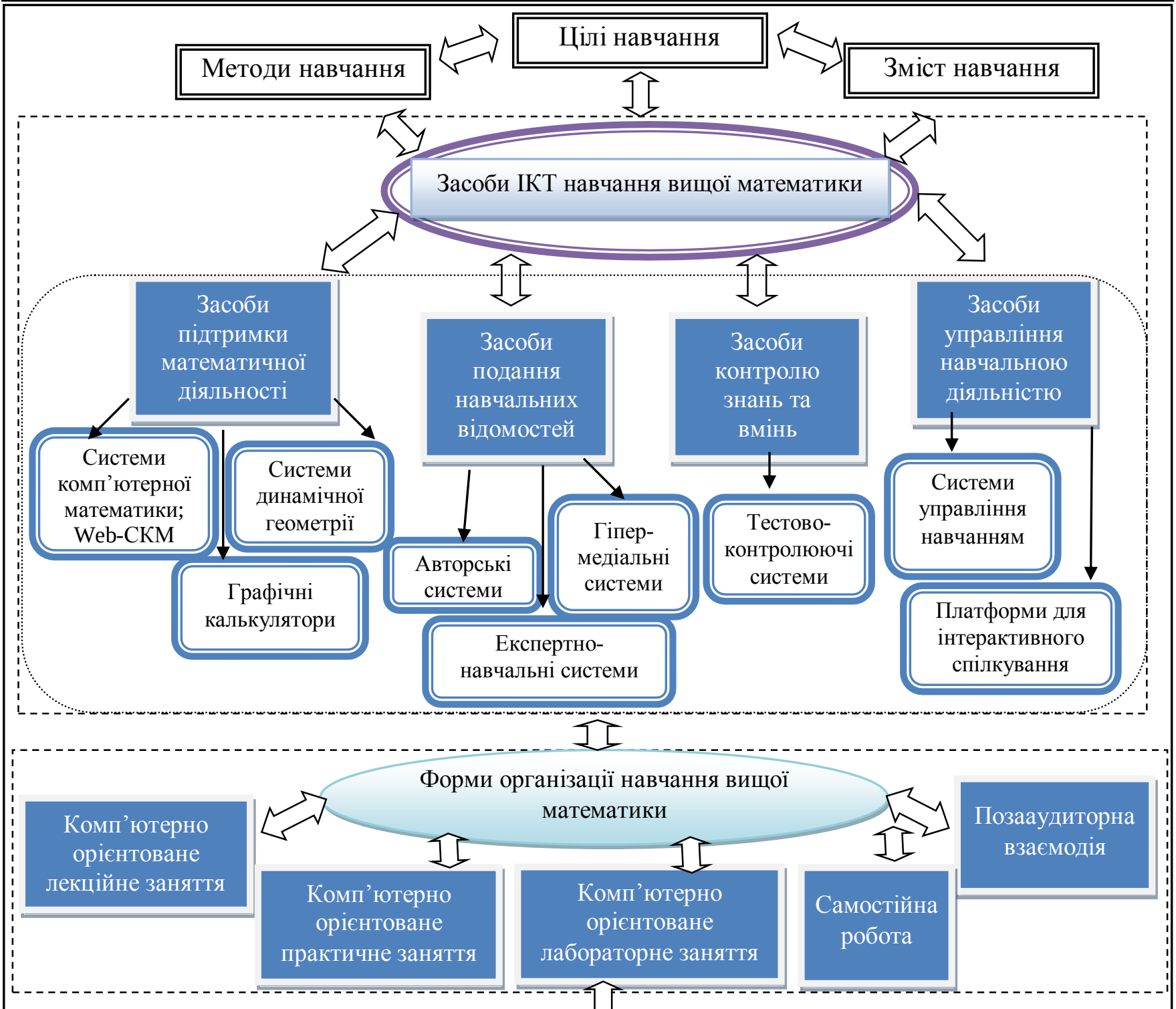

**Результати позитивного впливу використання ІКТ на процес навчання вищої математики за рахунок:**

– студент може самостійно визначати темп своєї роботи, з урахуванням власних можливостей і здібностей;
– студент має можливості необмеженого доступу до навчального матеріалу незалежно від місцезнаходження;
– різноманітність способів надання навчального матеріалу (відео, текст, інтерактивне спілкування);
– регулювання інтенсивності навчання на різних етапах навчального процесу;
– доступу до навчальних матеріалів попередніх семестрів;
– візуалізації матеріалу, що опановується;
– модульний принцип побудови навчальних матеріалів;
– організація самостійного навчання та можливість навчатися в онлайн-групах;
– здатність підтримувати комунікацію у процесі навчання як з викладачами, так і з студентами, що вивчають ті ж самі теми та розділи, проводити консультації, незалежно від місцезнаходження.

Рис. 2.23. Дидактична модель використання ІКТ при вивченні вищої математики в США на п’ятому етапі (1997–2003 рр.: поява та розробка систем управління навчанням)

У процесі **навчання** вищої математики із використанням засобів ІКТ п’ятого періоду організовано двосторонню взаємодію викладача зі студентами, діяльність викладача спрямована на досягнення максимального засвоєння та усвідомлення навчального матеріалу і подальшого застосування отриманих знань, умінь та навичок на практичних заняттях та в їх подальшій професійній діяльності. Використання засобів ІКТ надають можливість підтримувати індивідуальну, групову і колективну форми роботи студентів. До засобів цього етапу можна додати динамічні Java-аплети та відеолекції, що сприяють візуалізації навчального матеріалу.

Використання засобів ІКТ сприяють вирішенню *часової та просторової проблеми* організації **контролю** знань студентів. За допомогою апаратних та програмних засобів (телефон, миттєві повідомлення, електронна пошта, форум, чат, блог тощо) викладач у зручний для себе час здатен визначати рівень знань і вмінь студентів, розуміння ними теоретичного матеріалу та установити зворотній зв’язок із ними.

У системі **управління навчанням**: забезпечується цілеспрямований вплив на об’єкти навчальної діяльності, тобто на обсяг навчального матеріалу, що має самостійну логічну структуру та зміст; організовано **контроль** за доступом як до навчального матеріалу, так і до персональних даних студентів; організовано управління змістом, ресурсами, навчальною діяльністю студентів.

На цьому етапі виділяємо такі **форми організації** навчання вищої математики із використанням засобів ІКТ: комп’ютерно орієнтоване лекційне заняття, комп’ютерно орієнтоване практичне заняття, комп’ютерно орієнтоване лабораторне заняття, самостійна робота студентів та позааудиторна взаємодія викладача та студента.

Використання ІКТ у процесі навчання вищої математики надає студентам такі можливості:

– опрацьовувати навчальні матеріали, що подаються у вигляді різнотипних інформаційних ресурсів (текст, відео, анімація, презентація);

– у зручний час отримувати від викладача консультаційну підтримку з курсу;

– одержувати завдання в електронному вигляді для самоконтролю та підготовки до поточного контролю;

– дізнаватися про організацію навчального курсу.

Викладач має змогу наповнювати необхідними навчальними матеріалами курси і давати консультації на відстані, надсилати повідомлення студентам, розподіляти завдання, вести електронний журнал обліку активності студентів у навчанні, налаштовувати різноманітні ресурси навчального курсу тощо. Доступ до ресурсів курсу – відкритий.

Електронні навчальні матеріали, розміщені в мережі,

використовуються студентами для організації індивідуальної роботи, підготовки до виконання домашніх та екзаменаційних робіт.

Розглядаючи процес використання засобів ІКТ у навчанні вищої математики на п'ятому етапі, можна стверджувати, що їх використання здійснює позитивний вплив на: якість навчання; ступінь індивідуалізації, диференціалізації та персоналізації навчання; доступність матеріалу; самостійність студентів; мотивацію студентів.

Такий позитивний вплив при вивченні вищої математики з використанням ІКТ п'ятого етапу відбувається за рахунок того, що:

– студент може самостійно визначати темп своєї роботи, з урахуванням власних можливостей і здібностей;

– студент має можливості необмеженого доступу до навчального матеріалу не залежно від місця знаходження студента;

– різноманітності методів засвоєння навчального матеріалу (відео, текст, інтерактивне спілкування);

– регулювання інтенсивності навчання на різних етапах навчального процесу;

– доступу до навчальних матеріалів попередніх семестрів;

– візуалізації матеріалу, що вивчається;

– модульного принципу побудови навчальних матеріалів, що надає можливість переглядати окремі складові курсу;

– організації самостійного навчання та можливості навчатися в онлайн-групах;

– здатність підтримувати комунікацію у процесі навчання як з викладачами так і зі студентами, що вивчають ті ж самі теми та розділи, проводити консультації, не залежно від місця знаходження.

На **шостому етапі** (пов'язаному із перенесенням у Web-середовище засобів підтримки математичної діяльності та становленням і розвитком хмарних технологій навчання) (рис. 2.24) *напрями використання* розподіляються, як і на попередніх двох етапах, на *навчальні, контролюючі та управлінські.*

У процесі **навчання** вищої математики із використанням Web-орієнтованих засобів ІКТ організовано співпрацю викладача зі студентами, діяльність викладача спрямована на допомогу студентам в оволодінні навчального матеріалу і подальшого застосування отриманих знань, умінь та навичок в їх подальшій професійній діяльності.

Використання Web-орієнтованих засобів ІКТ надає можливість підтримувати індивідуальну, групову і колективну форми роботи студентів поза межами аудиторії. До засобів цього етапу можна додати платформи для проведення вебінарів, хмарні сховища даних, Web-орієнтовані засобів підтримки математичної діяльності.

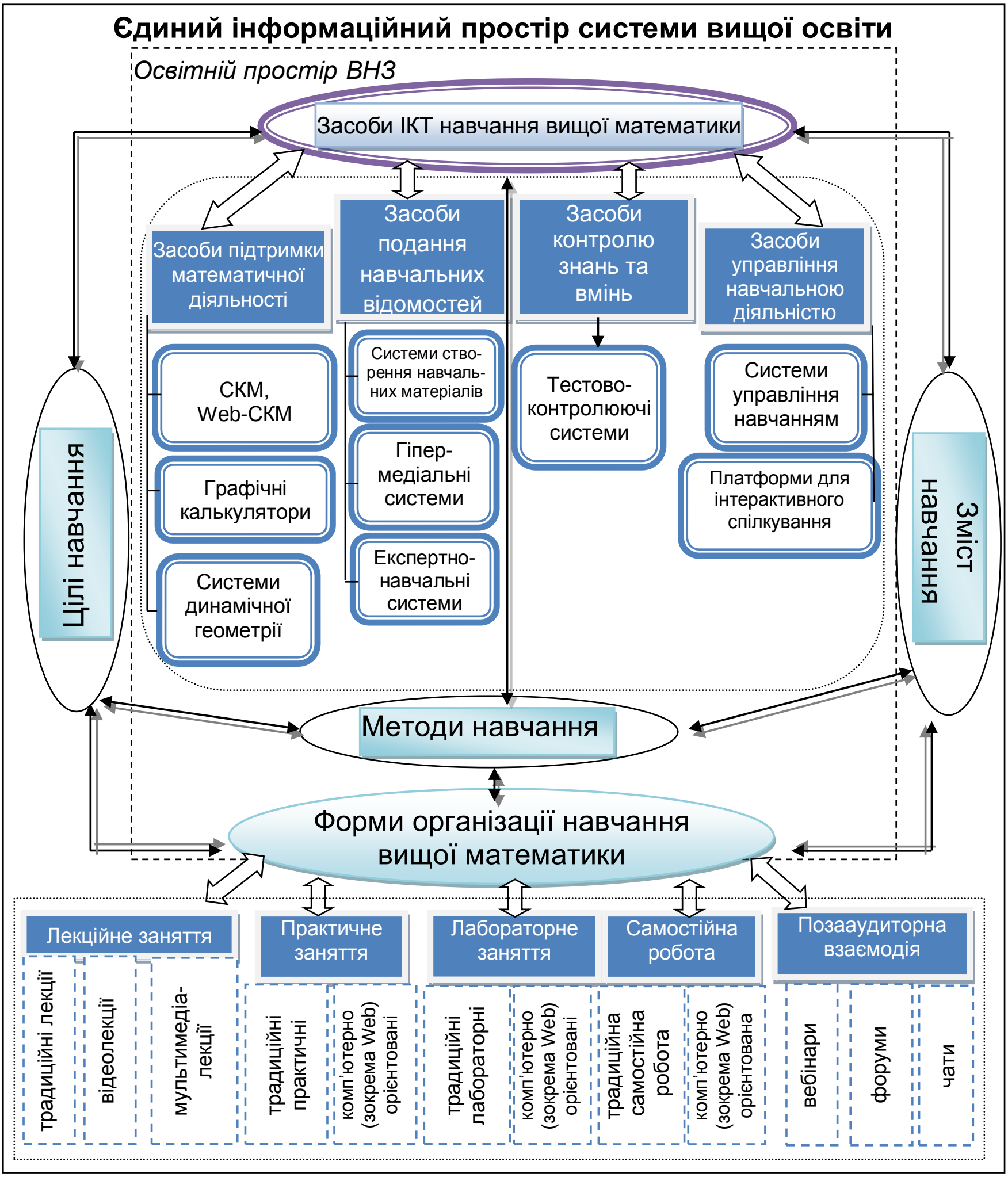

Рис. 2.24. Дидактична модель використання ІКТ при вивченні вищої математики в США на шостому етапі (з 2003 р. по теперішній час: перенесення у Web-середовище засобів підтримки математичної діяльності та становлення і розвиток хмарних технологій навчання)

Використання Web-орієнтованих засобів ІКТ сприяє вирішенню *часової та просторової проблеми* організації **контролю** знань студентів. За допомогою апаратних та програмних засобів викладач у зручний для себе час здатен провести контроль рівня знань і вмінь студентів, розуміння ними теоретичного матеріалу та установити зворотній зв’язок із ними.

У процесі **управління** забезпечується цілеспрямований вплив на об'єкти навчальної діяльності, тобто на обсяг навчального матеріалу, що має самостійну логічну структуру та зміст, і надає можливість оперувати цим матеріалом у процесі навчальної діяльності.

У системі управління навчанням організовано контроль за доступом як до навчального матеріалу, так і до персональних даних студентів; організовано управління змістом, ресурсами, навчальною діяльністю студентів. Організоване управління забезпечує підвищення ефективності функціонування процесу навчання. Основними призначеннями управління є координація студентами, адміністрування процесом навчання, контроль за виконанням навчального плану.

На цьому етапі виділяємо такі **форми організації** навчання вищої математики із використанням ІКТ: відео лекції та мультимедіа-лекції, комп'ютерно орієнтовані та Web-орієнтовані практичні, лабораторні заняття та самостійна робота студентів, позааудиторна взаємодія викладача зі студентами, що включає вебінари, форуми, чати тощо.

Використання Web-орієнтованих засобів ІКТ у процесі навчання вищої математики надає викладачам можливості здійснювати ефективну підтримка процесу навчання; надавати доступ до навчальних матеріалів всім учасникам навчального процесу; здійснювати систематичне управління навчального процесу.

Викладач має змогу наповнювати необхідними навчальними матеріалами курси і давати консультації на відстані, надсилати повідомлення студентам, розподіляти завдання, вести електронний журнал обліку активності студентів у навчанні, налаштовувати різноманітні ресурси навчального курсу тощо. Доступ до ресурсів курсу – відкритий.

Для забезпечення студентів електронними навчальними матеріалами, організації та керування самостійною роботою студентів, підтримкою комунікаційних можливостей використовується інтеграція традиційних методів навчання з методами дистанційного, мобільного навчання та навчання за допомогою Інтернет і мультимедіа.

Використання Web-орієнтованих засобів ІКТ у процесі навчання вищої математики надає студентам можливості виділити базові знання, необхідні для їх професійної діяльності; здійснювати інтеграцію аудиторної та позааудиторної роботи; здійснювати перебудову процесу навчання, надаючи йому принципів мобільності.

Електронні навчальні матеріали, розміщені в мережі, використовуються студентами для організації індивідуальної роботи, підготовки до виконання домашніх та екзаменаційних робіт. Організація та підтримка роботи із застосуванням Web-орієнтованих ІКТ надає можливість активізувати використання наявних і створювати нові освітні ресурси; розширити доступ до цих ресурсів студентам та викладачам; створити організаційну та технологічну базу для впровадження у процес навчання технологій дистанційного, мобільного навчання та навчання із використанням Інтернет та мультимедіа; покращити процес взаємодії між викладачем та студентом.

Розглядаючи процес використання засобів ІКТ у навчанні вищої математики шостого етапу, можна стверджувати, що їх використання здійснює позитивний вплив на якість навчання, ступінь індивідуалізації навчання, доступність матеріалу, самостійність студентів, мотивацію студентів, мобільність подання та передавання навчальних матеріалів, здатність підтримувати комунікацію у процесі навчання.

У статті Д. С. Пердью [100] виділено позитивні зміни, що викликані впровадженням Web-орієнтованих технологій у процес навчання вищої математики: по-перше, використання засобів ІКТ навчання надає можливість приділяти більше часу обговоренню навчального матеріалу, зменшуючи обсяг матеріалу, не зрозумілий для студентів; по-друге, студенти отримують можливість стати активними учасниками власного процесу навчання – кожному студенту подано весь навчальний матеріал та завдання з курсу, що переносить відповідальність за виконання поставлених вимог на студента. Використання Web-технологій у процесі навчання вищої математики надає можливість досягти кращого засвоєння навчального матеріалу студентами.

Запропоновані моделі об’єднують багато факторів, що мають важливе значення, таких як напрями використання засобів ІКТ, способи їх використання та напрями позитивного впливу використання технологій на процес навчання. Аналіз динаміки розвитку методик використання ІКТ у навчанні вищої математики студентів інженерних спеціальностей у вищих навчальних закладах Сполучених Штатів Америки показав, що висновок Ю. В. Триуса про те, що «використання ІКТ суттєво впливає на всі компоненти методичної системи навчання (цілі, зміст, методи, засоби і форми організації навчання)» [317, 269], є універсальним.

Поява нових ІКТ призводить до змін у теорії та методиці навчання вищої математики, що обумовлює виникнення нових цілей, засобів, форм, методів організації процесу навчання та доповненням змісту навчання, що має позитивний вплив на процес навчання вищої математики, значно розширюючи можливості студентів (табл. 2.4). Динаміка розвитку методик використання ІКТ у навчанні вищої математики студентів інженерних спеціальностей у ВНЗ США мала діалектичний характер: на кожному етапі її розвитку зміна засобів ІКТ супроводжувалась набуттям нового змісту способу доступу до навчальних ресурсів: *термінального – віддаленого – автономного – мережного – онлайн – хмарного*.

*Таблиця* 2.4

**Загальна характеристика змін у методиці використання ІКТ при вивченні вищої математики в США на різних етапах**

| Етапи / Напрями змін | Перший | Другий | Третій | Четвертий | П'ятий | Шостий |
|---|---|---|---|---|---|---|
| **Форми організації навчання** | традиційні | | | | | |
| | | комп'ютерно орієнтоване практичне заняття | | | | |
| | | | комп'ютерно орієнтоване лабораторне заняття | | | |
| | | | | комп'ютерно орієнтоване лекційне заняття | | |
| | | | | | позааудиторна взаємодія, комбіноване навчання | |
| | | | | | | Web-орієнтовані, масові відкриті дистанційні курси |
| **Методи навчання** | традиційні, використання програмованого навчання | | | | | |
| | | | дослідницьке навчання, зокрема – метод проектів | | | |
| | | | | кейс-метод, метод портфоліо, метод навчання у групах зі змінним складом, навчання у співробітництві | | |
| | | | | | комбіновані методи навчання | |

| Етапи / Напрями змін | Перший | Другий | Третій | Четвертий | П'ятий | Шостий |
|---|---|---|---|---|---|---|
| **Вплив на зміст** | зміст підручників доповнювався фрагментами програм (переважно мовою FORTRAN). | зміст підручників доповнювався алгоритмами обчислень | зміст підручників доповнювався розділами, пов'язаними з комп'ютеризацією усіх видів інженерної діяльності | зміст підручників переноситься у гіпертекстове Web-середовище | зміст підручників у Web-середовищі доповнюється мультимедійними демонстраціями та тестовими системами | зміст підручників у Web-середовищі адаптується до потреб студента |
| **Провідні засоби навчання** | комп'ютери PDP-8 та KENBAK-1, мови програмування FORTRAN, BASIC, LOGO | комп'ютери Apple II та TRS-80, мережна операційна система UNIX, мова програмування BASIC, СКМ | персональні комп'ютери, мови програмування BASIC, Pascal, FORTRAN, Algol, СКМ | технології Web 1.0 – статичні та динамічні HTML, JavaScript, Java; апаратні засоби (КПК, ноутбуки), СКМ | системи управління навчанням; СКМ | відкриті курси, СКМ та Web-СКМ; системи управління навчанням, засоби навчальної комунікації |

**Висновки до розділу 2**

1. Провідним напрямом організації процесу комп'ютерно орієнтованого навчання вищої математики у вищій технічний школі США є інтеграція традиційного навчання і навчання за допомогою Інтернет та мультимедіа, що сприяє гармонійному поєднанню теоретичної та практичної складових процесу навчання.

2. На сучасному етапі розвитку вищої інженерної школи США провідними засобами підтримки навчальної діяльності з вищої математики майбутніх інженерів є Web-орієнтовані ІКТ загального призначення (системи управління навчанням, системи розміщення відкритих навчальних матеріалів, засоби комунікації та спільної роботи) та спеціального призначення (системи комп'ютерної математики, лекційні демонстрації, динамічні навчальні матеріали).

3. Поява нових ІКТ призводить до змін у теорії та методиці навчання вищої математики, що обумовлює виникнення нових цілей, засобів, форм, методів організації процесу навчання та доповненням змісту навчання, що має позитивний вплив на процес навчання вищої математики, значно розширюючи можливості студентів. Динаміка розвитку методик використання ІКТ у навчанні вищої математики студентів інженерних спеціальностей у ВНЗ США мала діалектичний характер: на кожному етапі її розвитку зміна засобів ІКТ супроводжувалась набуттям нового способу доступу до навчальних ресурсів: *термінального – віддаленого – автономного – мережного – онлайн – хмарного*.

4. На першому етапі розвитку теорії та методики використання ІКТ у навчанні вищої математики студентів інженерних спеціальностей у ВНЗ США термінальний спосіб доступу до навчальних програм та даних надавав можливість одночасної роботи групи студентів у наперед визначений термін у заданому місці. Зміст підручників з вищої математики доповнювався фрагментами програм (переважно мовою FORTRAN) та міг реалізовувати програмоване навчання. Програмні засоби навчання були орієнтовані на певний клас комп'ютерів. Апаратне забезпечення було нестандартизоване, мобільність засобів ІКТ була низькою та забезпечувалась лише мобільними навчальними класами. Форми організації та методи навчання вищої математики залишались традиційними, підтримка самостійної роботи студентів засобами ІКТ – відсутня або обмежена доступом до них. Системний підхід до організації комп'ютерно орієнтованого навчання вищої математики був відсутнім.

5. На другому етапі розвитку теорії та методики використання ІКТ у навчанні вищої математики студентів інженерних спеціальностей у ВНЗ США віддалений спосіб доступу до навчальних програм та даних розширював можливості групової та автономної роботи студентів зі спеціально обладнаної аудиторії ВНЗ до університетського кампусу та за його межами. Зміст підручників з вищої математики доповнювався алгоритмами обчислень. Програмні засоби навчання набули мобільності за рахунок використання мобільних мов програмування. Нестандартизованість апаратного забезпечення компенсувалось появою мобільних операційних систем. Виникає нова форма організації навчання – комп'ютерно орієнтоване практичне заняття з використанням систем програмування та систем комп'ютерної математики. Самостійна робота студентів з вищої математики підтримується спеціалізованими засобами ІКТ – інженерними калькуляторами. З'являється системний підхід до організації комп'ютерно орієнтованого навчання обчислювальної математики.

6. На третьому етапі розвитку теорії та методики використання ІКТ у навчанні вищої математики студентів інженерних спеціальностей у ВНЗ США автономний спосіб доступу до навчальних програм та даних розширив можливості самостійної роботи студентів з вищої математики. Поява нового засобу навчання – персональних комп'ютерів – сприяла підвищенню активності студентів на заняттях з вищої математики, формуванню практичних умінь та навичок, розмаїттю подання навчального матеріалу, полегшенню процесу моделювання, візуалізації навчального матеріалу, автоматизації складних обчислень. У змісті навчання вищої математики з'явились розділи, пов'язані з комп'ютеризацією усіх видів інженерної діяльності. З'являються підручники з вищої математики, орієнтовані на спільне використання систем комп'ютерної математики та систем динамічної геометрії доповнювався алгоритмами обчислень. Виникла нова форма організації навчання з вищої математики – комп'ютерно орієнтована лабораторна робота. Цілі навчання вищої математики доповнились набуттям умінь застосовувати персональні комп'ютери для розв'язання прикладних задач та проведення навчальних досліджень. У навчанні вищої математики набувають поширення методи дослідницького навчання, зокрема – метод проектів. З'являється системний підхід до організації комп'ютерно орієнтованого навчання вищої математики.

7. На четвертому етапі розвитку теорії та методики використання ІКТ у навчанні вищої математики студентів інженерних спеціальностей у ВНЗ США мережний спосіб доступу до навчальних програм та даних сприяв інтеграції персональних обчислювальних ресурсів студентів та ВНЗ у навчально-інформаційне середовище на основі Web-технологій. Зміст підручників з вищої математики переноситься у гіпертекстове Web-середовище. У програмних засобах навчання використовуються технології Web 1.0 – статичний та динамічний HTML, JavaScript, Java. Програмне забезпечення мобільних апаратних засобів ІКТ (КПК, ноутбуки тощо) підтримує окремі складові навчальної діяльності з вищої математики. Використання Web-технологій сприяло забезпеченню особистісно-орієнтованого підходу до навчання; можливості відкритого, необмеженого доступу до навчального матеріалу незалежно від місця знаходження студента; модульного принципу побудови навчальних матеріалів, що надає студентам можливість переглядати окремі складові частини навчальних матеріалів; розвитку самостійного навчання студентів, що сприяло появі нових методів навчання: кейс-метод, метод портфоліо, метод навчання у групах зі змінним складом, навчання у співробітництві.

8. На п'ятому етапі розвитку теорії та методики використання ІКТ у навчанні вищої математики студентів інженерних спеціальностей у ВНЗ США онлайн доступ до електронних освітніх ресурсів зумовив появу і розвиток систем управління навчанням. Відбувається інтеграція подання навчальних матеріалів, планування процесу навчання та контролю навчальних досягнень студентів з вищої математики у єдиному Web-середовищі, доступному онлайн через різноманітні Інтернет зорієнтовані засоби ІКТ. Використання систем управління навчанням створює умови для програмно-методичної підтримки індивідуальної, групової та колективної форм організації навчання студентів, надаючи їм можливість самостійно визначати темп своєї роботи з урахуванням власних можливостей і здібностей, урізноманітнюючи способи подання навчального матеріалу на різних типах пристроїв, зокрема мобільних, регулювання інтенсивності навчального навантаження на різних етапах навчального процесу, доступ до навчальних матеріалів попередніх семестрів та курсів, організації самостійного навчання та можливість навчатися в онлайн-групах зі змінним складом. Здатність підтримувати навчальну комунікацію як з викладачами, так і з студентами за межами аудиторії у системі управління навчанням сприяє виникненню комбінованих форм організації та методів навчання.

9. На шостому етапі розвитку теорії та методики використання ІКТ у навчанні вищої математики студентів інженерних спеціальностей у ВНЗ США використання хмарних технологій сприяє досягненню студентами та викладачами високого рівня мобільності. Зміст підручників з вищої математики разом з іншими електронними освітніми ресурсами (зокрема, засоби підтримки математичної діяльності) переносяться у Web-середовище, завдяки чому суттєво розширюється спектр засобів ІКТ, що можуть бути використані для навчання вищої математики. Виникає нова форма організації навчання – масові відкриті дистанційні курси.

Провідними формами організації навчання стають Web-орієнтовані лекційні, практичні та лабораторні заняття. Ураховуючи незавершеність даного етапу, він містить у собі багато рис попереднього.

# РОЗДІЛ 3
# ПЕРСПЕКТИВИ ВИКОРИСТАННЯ ДОСВІДУ СПОЛУЧЕНИХ ШТАТІВ АМЕРИКИ ЗАСТОСУВАННЯ ІНФОРМАЦІЙНО-КОМУНІКАЦІЙНИХ ТЕХНОЛОГІЙ У НАВЧАННІ ВИЩОЇ МАТЕМАТИКИ СТУДЕНТІВ ІНЖЕНЕРНИХ СПЕЦІАЛЬНОСТЕЙ В УКРАЇНІ

## 3.1 Інформаційно-комунікаційні технології навчання вищої математики студентів інженерних спеціальностей університетів України

Однією з особливостей нашого часу в Україні, як і в багатьох інших країнах світу, є рух до формування інформаційного суспільства. Велику роль в цьому процесі відіграє просування ІКТ в освітню сферу, що зумовлює необхідність постійного підвищення ефективності використання новітніх ІКТ у навчальному процесі, впровадження ІКТ та підтримки їх на сучасному рівні для управління освітньою галуззю, включаючи кожен навчальний заклад, своєчасного оновлення змісту освіти та підвищення якості підготовки фахівців з ІКТ [295].

За останні десять років в Україні зроблено чимало позитивних кроків щодо впровадження і ефективного використання у вищій освіті нових освітніх технологій, що спираються на ІКТ. У першу чергу, це стосується нормативно-правового забезпечення даного напряму. Так, за цей час було прийнято 4 Укази, 14 Законів України, 2 Постанови Верховної Ради України, 16 Постанов Кабінету Міністрів України, 13 Наказів Міністерства освіти і науки України, які в тій чи іншій мірі регламентують розвиток ІКТ у вищій освіті [254].

Важливим завданням перспективної системи освіти є підготовка фахівців для професійної інженерної діяльності в інформаційній сфері суспільства та педагогічних кадрів, здатних розробляти і застосовувати інформаційні технології навчання [156]. Поточний стан, проблеми, перспективи розвитку використання ІКТ у вищій освіті в Україні розглянуто у додатку Ж.

До ІКТ, що можуть бути використані у процесі навчання математики можна віднести:

– мережні технології, що використовують локальні мережі та глобальну мережу Інтернет (електронні методичні рекомендації, платформи дистанційного навчання, що забезпечуютьпідтримку інтерактивного зв’язку зі студентами, зокрема, в онлайн режимі);

– технології, що зорієнтовані на локальні комп’ютери (навчальні програми, комп’ютерні моделі реальних процесів, демонстраційні програми, електронні задачники, тестові системи);

– мобільні технології, що надають студенту та викладачу високу ступінь свободи (рис. 3.1).

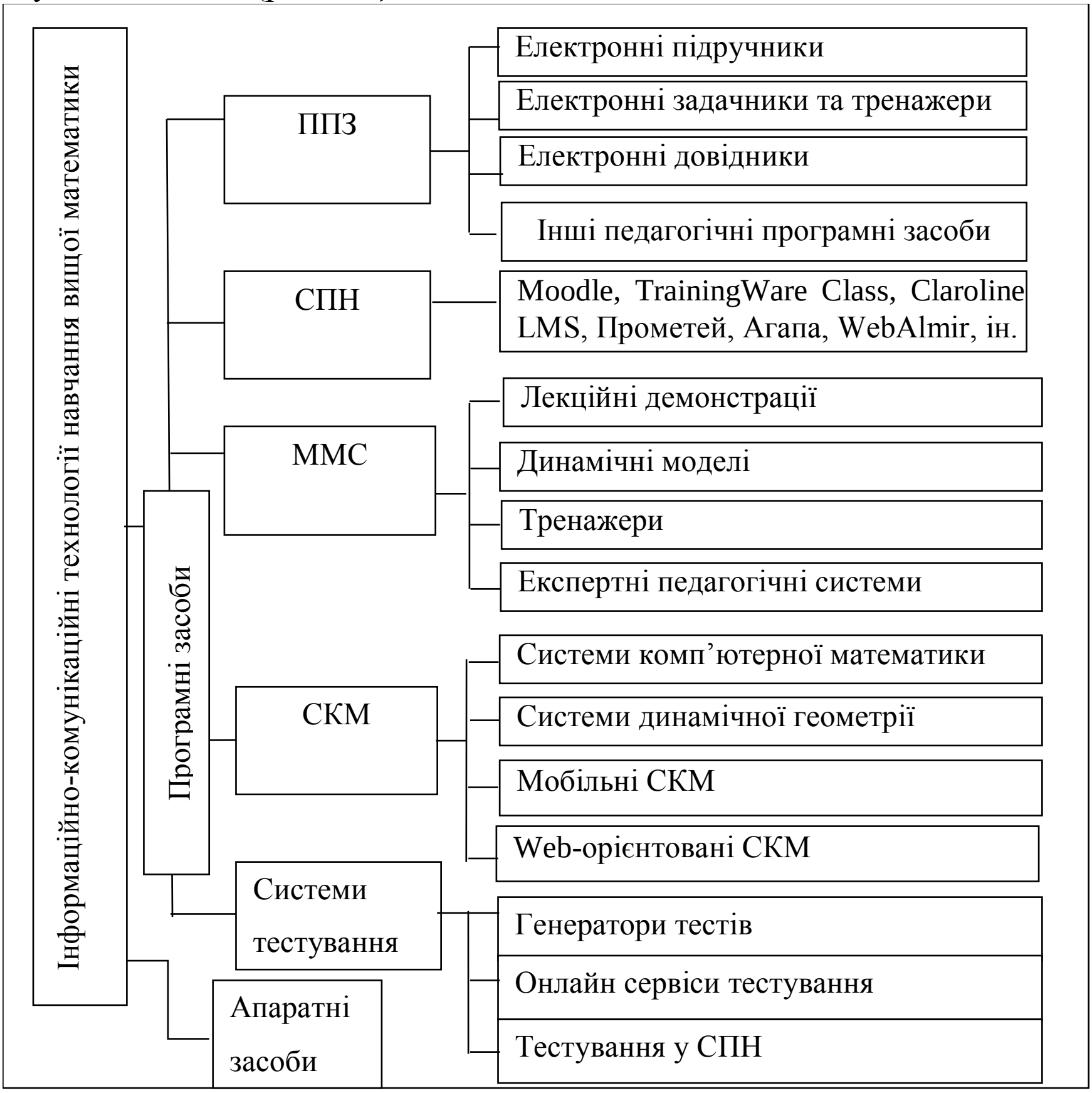

Рис. 3.1. ІКТ навчання вищої математики

*3.1.1 Електронні освітні ресурси* (ЕОР) – навчальні, наукові, інформаційні, довідкові матеріали та засоби, розроблені в електронній формі та представлені на носіях будь-якого типу або розміщені у комп’ютерних мережах, які відтворюються за допомогою електронних цифрових технічних засобів і необхідні для ефективної організації навчально-виховного процесу, в частині, що стосується його наповнення якісними навчально-методичними матеріалами [259].

До основних видів ЕОР належать [259]:

– *електронний документ* – документ, відомості в якому подані у формі електронних даних і для використання якого потрібні технічні

засоби;

– *електронне видання* – електронний документ, який пройшов редакційно-видавниче опрацювання, має вихідні відомості й призначений для розповсюдження в незмінному вигляді;

– *електронний аналог друкованого видання* – електронне видання, що в основному відтворює відповідне друковане видання, зберігаючи розташування на сторінці тексту, ілюстрацій, посилань, приміток тощо;

– *електронні дидактичні демонстраційні матеріали* – електронні матеріали (презентації, схеми, відео- й аудіозаписи тощо), призначені для супроводу навчально-виховного процесу;

– *інформаційна система* – організаційно впорядкована сукупність документів (масивів документів) та інформаційних технологій, в тому числі з використанням технічних засобів, що реалізують інформаційні процеси та призначені для зберігання, обробки, пошуку, розповсюдження, передачі та надання відомостей;

– *депозитарій електронних ресурсів* – інформаційна система, що забезпечує зосередження в одному місці сучасних ЕОР з можливістю надання доступу до них через технічні засоби, у тому числі в інформаційних мережах (як локальних, так і глобальних);

– *комп'ютерний тест* – стандартизовані завдання, представлені в електронній формі, призначені для вхідного, проміжного і підсумкового контролю рівня навчальних досягнень, а також самоконтролю та/або такі, що забезпечують вимірювання психофізіологічних і особистісних характеристик випробовуваного, обробка результатів яких здійснюється за допомогою відповідних програм;

– *електронний словник* – електронне довідкове видання упорядкованого переліку мовних одиниць (слів, словосполучень, фраз, термінів, імен, знаків), доповнених відповідними довідковими даними;

– *електронний довідник* – електронне довідкове видання прикладного характеру, в якому назви статей розташовані за абеткою або в систематичному порядку;

– *електронна бібліотека цифрових об'єктів* – набір ЕОР різних форматів, в якому передбачено можливості для їх автоматизованого створення, пошуку і використання;

– *електронний навчальний посібник* – навчальне електронне видання, використання якого доповнює або частково замінює підручник;

– *електронний підручник* – електронне навчальне видання з систематизованим викладом дисципліни (її розділу, частини), що відповідає навчальній програмі;

– *електронні методичні матеріали* – електронне навчальне або виробничо-практичне видання роз'яснень з певної теми, розділу або

питання навчальної дисципліни з викладом методики виконання окремих завдань, певного виду робіт;

– *курс дистанційного навчання* – інформаційна система, яка є достатньою для навчання окремим навчальним дисциплінам за допомогою опосередкованої взаємодії віддалених один від одного учасників навчального процесу у спеціалізованому середовищі, яке функціонує на базі сучасних психолого-педагогічних технологій та ІКТ;

– *електронний лабораторний практикум* – інформаційна система, що є демонстраційною моделлю природних і штучних об’єктів, процесів та їх властивостей із застосуванням засобів комп’ютерної візуалізації.

Опишемо ІКТ навчання вищої математики.

**Програмно-педагогічні засоби (ППЗ)** – це комплекси прикладних програм, що призначені для організації та підтримки процесу навчання із використанням комп’ютера. ППЗ призначені для подання навчальних відомостей, вони надають можливість організувати індивідуальний підхід до кожного студента за допомогою налагодженого зворотного зв’язку користувача з програмою [237].

Використання ППЗ у процесі навчання вищої математики надають можливість студентам самостійно проводити дослідження, опрацьовувати отримані дані та інтерпретувати одержані результати, що сприяє більш глибинному сприйняттю абстрактного навчального матеріалу і надає йому професійної спрямованості [148].

До основних ППЗ можна віднести електронні програмно-методичні комплекси, електронні підручники, електронні довідники, електронні задачники та тренажери, але усі вони повинні бути розроблені з урахуванням наступних дидактичних засад [282]:

– *інтегрованість*: одну й ту саму наочність можна використовувати з різним цільовим призначенням;

– *конструктивність* забезпечується аналізом комп’ютерних зображень реальних предметів, перенесенням їх властивостей на відповідні їм моделі, де увага приділяється поелементному їх створенню, внаслідок чого студент самостійно формулює означення нових понять;

– *інтерактивність* забезпечується використанням сучасних форм організації проведення занять (лекція з ілюстраціями, групова, парна, індивідуальна робота, семінарське заняття тощо), підтримка активних методів навчання (проведення посильних навчальних досліджень, моделювання і конструювання професійних задач;

– *візуалізація* забезпечується розробленими комп’ютерними динамічними моделями.

*Електронні програмно-методичні комплекси* (ЕПМК) – це комп’ютерні засоби організації самонавчання студентів (наодинці чи під керівництвом викладача) у процесі їх самостійної навчально-пізнавальної діяльності з урахуванням майбутнього фаху, етапу навчання, робочої програми дисципліни, визначених форм, видів, методів і підходів до навчання та стратегій управління траєкторією учіння студентів [215]. ЕПМК містять навчальні відомості, конструктори занять, словник термінів і понять, історичну довідку, різноманітні таблиці.

*Електронний підручник* використовують для самостійного вивчення студентом теоретичного матеріалу курсу. Електронні підручники містять всі види навчальної діяльності, що спрямовані на підтримку як аудиторної так і самостійної роботи студента.

Одним із електронних підручників, що використовують у процесі навчання вищої математики студентів інженерних спеціальностей Донецького національного університету, є розроблений О. Г. Євсєєвою, де навчальні відомості подаються невеликими частинами в структурованому вигляді. Кожний розділ семантичного конспекту містить посилання на підтему тематичного компоненту, якій він відповідає. Кожне висловлювання містить посилання на інші висловлювання, від яких воно залежить, визначається і з яких виходить. Після кожної порції даних наводяться завдання, спрямовані їх на засвоєння й одночасне виконання математичних предметних дій. Завдання представлені тестами різних типів: закриті, відкриті, завданнями на відповідність, на встановлення правильної послідовності. Окрім цього, система завдань містить завдання евристичного характеру, спрямовані на розвиток у студентів логічного мислення і виконання математичних предметних дій теоретичного характеру.

Особливістю подання матеріалу є використання О. Г. Євсєєвою процедури орієнтування, що складається із загального орієнтування (визначення, що треба робити і що для цього треба знати) і загального орієнтування (визначення, які дії необхідно виконати і за допомогою чого), що сприяє освоєнню математичних дій. При цьому для кожного типу задач, що розв’язуються, пропонується складати так звану схему орієнтування [170].

*Електронний довідник* надає можливість студенту у будь-який час оперативно отримати довідкові відомості. Електронні довідники, на думку Л. М. Наконечної [237; 238], відносяться до інформаційно-довідкових джерел, що забезпечують загальну інформаційну підтримку. Такі електронні ресурси використовують при розв’язку творчих навчальних задач, в тому числі тих, що виходять за рамки навчальної програми. Довідкові джерела наділені основними дидактичними якостями: автоматичністю та відкритістю змісту, можливістю копіювання окремих частин матеріалу, що подається в будь-яких поєднаннях.

*Електронна бібліотека* (англ. Digital library) надає можливість зберігати і використовувати різнорідні колекції електронних документів (текст, графіка, аудіо, відео тощо) завдяки глобальним мережам передачі даних в зручному, для кінцевого користувача, вигляді [185].

О. М. Спірін [306] зазначає, що електронні бібліотеки відіграють особливу роль у розширенні доступу до відомостей і даних, що забезпечують подання інформаційних ресурсів в електронному вигляді, віддалений доступ до них з використанням ІКТ. Вони значно підвищують рівень надання користувачам бібліотечних послуг, а саме:

– сприяють вільному доступу до наявних електронних інформаційних ресурсів в мережі Інтернет, насамперед до бібліотек і періодичних видань, а також до зарубіжних електронних ресурсів;

– забезпечують якісно новий рівень задоволення інформаційних потреб науковців завдяки використанню ІКТ.

*Комп'ютерні моделі, конструктори і тренажери* надають можливість закріпити отримані на лекції чи практичному занятті знання; відпрацювати розв'язки типових задач, що надають можливість наочно пов'язати теоретичні знання з конкретними професійними проблемами, на вирішення яких вони можуть бути спрямовані.

Одним з найперших вітчизняних засобів візуалізації математичної задачі та її розв'язку, що робить діалог учня (студента) та викладача більш доступним та евристичним, є педагогічний програмний засіб **GRAN**, розробка якого розпочалася у 1989 році авторським колективом під керівництвом М. І. Жалдака.

До складу ППЗ GRAN входять педагогічні програмні засоби: GRAN1, що призначена для графічного аналізу функцій; GRAN-2D – для комп'ютерної підтримки навчання планіметрії; GRAN-3D – для комп'ютерної підтримки навчання стереометрії.

За допомогою <u>GRAN1</u> (GRaphic ANalysis) можна розв'язувати такі класи задач [175]: побудова графіків функцій, обчислення значення виразів; графічне розв'язування рівнянь та систем рівнянь; графічне розв'язування нерівностей та систем нерівностей; відшукання найбільших та найменших значень функції на заданій множині точок; побудова січних та дотичних до графіків; обчислення визначених інтегралів, обчислення площ довільних фігур, обчислення довжини дуги кривої, обчислення об'ємів та площ поверхонь обертання; елементи статистичного аналізу експериментальних даних.

<u>GRAN-2D</u> (GRaphic ANalysis 2-Dimension) відноситься до розряду програм динамічної геометрії та призначений для графічного аналізу геометричних об'єктів на площині. Використання пакету GRAN-2D надає можливість: створювати динамічні моделі геометричних фігур та їхніх комбінацій аналогічно класичним побудовам за допомогою циркуля та лінійки, а також використовуючи елементи аналітичної геометрії (систему координат, рівняння прямих і кіл, алгебраїчні залежності між частинами побудови, графіки функцій тощо); проводити вимірювання геометричних величин, досліджувати геометричні місця точок; аналізувати динамічні вирази, висувати припущення, встановлювати закономірності; будувати графічні зображення, використовуючи коментарі, кнопки, підказки та гіперпосилання; експортувати рисунки у графічні формати для вбудовування їх у інші додатки і для створення геометричних ілюстрацій тощо.

GRAN-3D (GRaphic ANalysis 3-Dimension) для графічного аналізу тривимірних об'єктів. Використання пакету GRAN-3D надає можливість: створювати та перетворювати моделі базових просторових об'єктів; виконувати перерізи многогранників площинами; обчислювати об'єми та площі поверхонь многогранників і тіл обертання; вимірювати відстані та кути [174].

Складовою програмно-методичного комплекту GRAN є також посібник для вчителів, у якому наведена значна кількість математичних прикладів, що унаочнюють графічні зображення задач і вправ для самостійного виконання, питання для самоконтролю. Ці завдання можна використовувати також для вивчення деяких розділів вищої математики.

**Програмне середовище (ПС) «Системи лінійних рівнянь»** розроблено авторським колективом під керівництвом О. В. Співаковського [277]. Основним призначенням ПС «Системи лінійних рівнянь» є комп'ютерна підтримка практичних занять і лабораторних робіт з алгебри при вивченні теми «Системи лінійних рівнянь», а також при розв'язанні арифметичних, фізичних та інших задач, в яких математична модель є системою лінійних рівнянь. У процесі такого роду діяльності студент використовує теоретичні знання, придбані на попередніх стадіях навчання, для рішення практичних задач, що надає можливість вирішити задачі формування необхідних вмінь та навичок з даної теми [304].

До складу ПС «Системи лінійних задач» входять такі модулі: теоретичний матеріал курсу, завдання для контрольних та самостійних робіт, завдання для тестового контролю (рис. 3.2).

Використання ППЗ «Системи лінійних задач» під час вивчення вищої математики надає можливість обчислювати визначники, зводити матриці до трикутного виду, знаходити власні значення вектора, розв'язувати системи лінійних рівнянь графічним способом, способом додавання, способом підстановки, а також текстових задач, математичною моделлю яких є системи лінійних рівнянь.

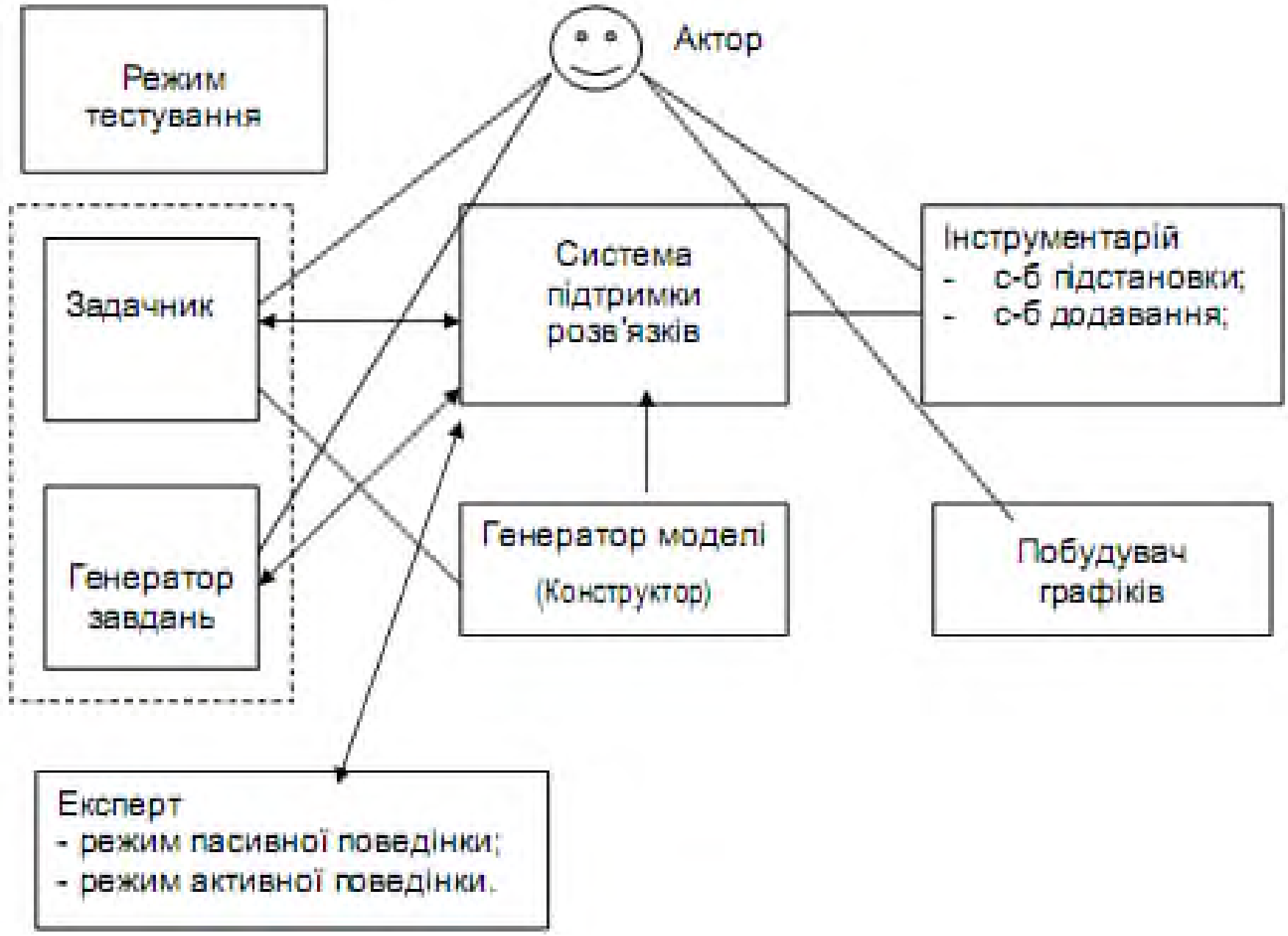

Рис. 3.2. Структура ППЗ «Системи лінійних рівнянь»

**Електронний навчально-методичний посібник** «Вища математика», створений у Донбаській державній машинобудівній академії під керівництвом К. В. Власенко [157], може бути використаний як електронний дистанційний курс.

Цей програмно-методичний комплекс забезпечує можливість самостійно засвоїти навчальний модуль та з дидактичної точки зору містить декілька структур [157]:

1. *Лінійні програми*, що допомагають засвоїти декларативні та процедурні знання, необхідні для майбутніх інженерів-машинобудівників.

2. *Розгалужені програми*, під час роботи з якими вибір правильних відповідей вимагає від студентів більших розумових здібностей, ніж повторення відомих теоретичних відомостей.

3. *Циклічні програми* ґрунтуються на використанні принципу зворотного зв’язку, пов’язаного з контролем засвоєння навчальних відомостей та, за необхідності, наступною корекцією цього процесу.

4. *Узагальнені ієрархічно-інформаційні програми*, що надають можливість створення більш інерційного зворотного зв’язку,

забезпечують контроль і корекцію засвоєння навчального матеріалу в межах декількох занять як із спільною так і різною тематикою.

**Інструментальні програмні засоби**, такі як Xtremum, XtremumND, Extremum, Nonline, Asimplex, Master of Logic, AlgoMachines, розроблені під керівництвом Ю. В. Триуса [321] і призначені для розв'язування задач оптимізації.

*Програмний засіб ASimplex* створено для розв'язування задач лінійного програмування, який має зручний інтерфейс, режим встановлення параметрів середовища та систему допомоги. Користувачеві надається можливість вводити математичну модель задачі, редагувати і зберігати її, автоматично одержувати канонічну форму задачі та двоїсту до неї задачу, візуально спостерігати за процесом розв'язування задачі лінійного програмування за звичайним симплекс-методом, методом штучного базису і двоїстим симплекс-методом, які обираються автоматично в залежності від математичної моделі, а також сформувати у вигляді HTML-файлу повний протокол з відомостями про всі поточні обчислення, вибір елементів розв'язання, перевірку поточних планів на оптимальність і побудову поточних симплекс-таблиць [321].

*Програмний засіб Extremum* доцільно використовувати при викладанні курсу «Математичні методи оптимізації», оскільки, на думку розробників, вона дає можливість поглибити знання студентів з питань, що стосуються математичного моделювання і обчислювального експерименту, дослідження ефективності методів розв'язування математичних задач за допомогою комп'ютера, аналізу та інтерпретації отриманих результатів, підвищити інформаційну культуру студентів [321].

Програмний продукт *AlgoMachines* надає можливість користувачу працювати в двох основних режимах: інструментальному та режимі контролю. В інструментальному режимі програма надає можливість користувачу конструювати алгоритми в алгоритмічних системах Тьюрінга, Поста та Маркова для розв'язування різноманітних задач і виконувати їх як автоматично, так і покроково. При цьому користувач звільняється від виконання громіздких і рутинних операцій, що на певному етапі навчання є несуттєвими. Користувач також може одержувати необхідні теоретичні відомості через контекстну систему допомоги. Програмний продукт в інструментальному режимі надає можливість: створити новий алгоритм в одній з алгоритмічних систем; відкрити вже існуючий алгоритм та відредагувати його; виконати алгоритм при будь-яких вхідних даних, заданих у певному алфавіті; проаналізувати його в режимі покрокового виконання; зберегти алгоритм; змінити мову інтерфейсу програми під час роботи (передбачено три мови: українську, російську та англійську); подавати алгоритми різними способами, зокрема алгоритми Маркова можна подати у вигляді граф-схеми, скороченої граф-схеми чи блок-схеми, а програму машини Тьюрінга у вигляді таблиці чи послідовності команд [321].

Для комп'ютерної підтримки вивчення основ математичної логіки пропонується авторська інструментально-контролююча програма *Master of Logic*, що призначена для розвитку умінь і навичок розв'язування основних класів задач алгебри висловлень та організації проміжного та підсумкового контролю.

Використання Master of Logic дозволяє: будувати таблицю істинності; визначати, до якого класу формул належить задана формула; встановлювати рівносильність двох формул алгебри висловлень; для формул алгебри висловлень що задаються в аналітичному вигляді або таблично, будувати досконалі кон'юнктивні або диз'юнктивні нормальні форми; встановлювати, чи є задана формула алгебри висловлень логічним наслідком скінченої кількості посилок; розв'язувати систему булевих рівнянь; створювати і працювати з файлами, що містять задачі вказаних типів і можуть використовуватися для проведення занять, зокрема практичних; створювати завдання для проведення самостійних і контрольних робіт з відповідних тем або для комплексного контролю знань [321].

*3.1.2 Системи підтримки навчання* (СПН) являють собою платформи дистанційного навчання чи мобільні модулі для навчання вищої математики, що можуть бути інтегровані в систему дистанційного навчання.

У системі вищої освіти України для організації процесу навчання вищої математики використовують вільно поширюванні платформи для дистанційного навчання. До більш поширювальних СПН в Україні відносять: Moodle, TrainingWare Class, Claroline LMS, Прометей, Агапа, WebAlmir, деякі з них розглянуто в додатку К.

На базі Херсонського державного університету розроблено проект **«Херсонський віртуальний університет 2.0»**, що призначений для здійснення взаємодії за допомогою мережі Інтернет між фасилітатором та слухачами дистанційного курсу у процесі розподіленого навчального процесу, під час якого викладач може обмінюватися зі студентами мультимедійною у реальному часі. Даний проект передбачає також спілкування учасників дистанційного курсу, прослуховування аудіо файлів, трансляцію та перегляд відеороликів, необхідних для проведення відео конференцій.

Проект складається із клієнтської частини, створеної на основі Flash-технології та серверної частини, що представляє віддалений Web-сервіс, створений на основі технології NET. Зв'язок між клієнтською частиною додатка і серверною частиною реалізовано на основі технології Flash Remoting MX [216].

Важливим аспектом проекту є те, що клієнтська частина є платформонезалежною і в електронній дошці можна працювати з будь-якої операційної системи, достатньо лише мати на клієнтському комп'ютері браузер і Flash-плеєр.

*3.1.3 Мобільне математичне середовище* (ММС) – відкрите модульне мережне мобільне інформаційно-обчислювальне програмне забезпечення, що надає користувачу (викладачу, студенту) можливість мобільного доступу до інформаційних ресурсів математичного і навчального призначення, створюючи умови для організації повного циклу навчання (зберігання та подання навчальних матеріалів; проведення навчальних математичних досліджень; підтримка індивідуальної та колективної роботи; оцінювання навчальних досягнень тощо) та інтеграції аудиторної і позааудиторної роботи у безперервний процес навчання [299].

У процесі розробки ММС з вищої математики особливу увагу слід приділити вибору математичного пакету, що складає ядро ММС (рис. 3.3). Головними критеріями вибору математичного пакету для побудови ММС з вищої математики є: *розширюваність* (система повинна надавати можливість користувачеві доповнювати її для задоволення професійних потреб); *наявність Web, WAP-інтерфейсів, XML-RPC, SOAP* та інших Web-сервісів (для забезпечення мобільного доступу); *крос-платформеність* (мобільність програмного забезпечення); можливість створення програм із стандартними елементами управління – лекційні демонстрації, динамічні моделі, тренажери та навчальні експертні системи; можливість інтегрувати у себе різноманітне програмне забезпечебння для навчання математики; підтримка технологій Wiki; можливість локалізації та вільне поширення.

*3.1.4 Системи комп'ютерної математики* – це сукупність методів і засобів, що забезпечують максимально комфортну й швидку підготовку алгоритмів і програм для розв'язування математичних завдань будь-якої складності з високим ступенем візуалізації усіх етапів розв'язування [168].

Ю. В. Триус виділяє основні об'єктивні та суб'єктивні причини низького рівня використання СКМ при вивченні математичних дисциплін [320]:

– до *об'єктивних* причин відносить: недостатній рівень забезпечення сучасною комп'ютерною технікою математичних кафедр для

регулярного її використання в навчанні математичних дисциплін; відсутність коштів у ВНЗ на придбання ліцензованого програмного забезпечення (навіть студентських версій, які коштують значно дешевше, ніж комерційні та академічні версії); відсутність коштів у ВНЗ і викладачів на придбання навчальної, методичної і довідкової літератури з СКМ;

– до *суб'єктивних*: недостатню обізнаність викладачів з можливостей використання СКМ, особливо тих, що вільно розповсюджуються, їх роль в математичних дослідженнях і математичній освіті; певний консерватизм викладачів у підходах до навчання математичних дисциплін; недостатній рівень інформаційної культури викладачів математичних дисциплін і студентів некомп'ютерних спеціальностей.

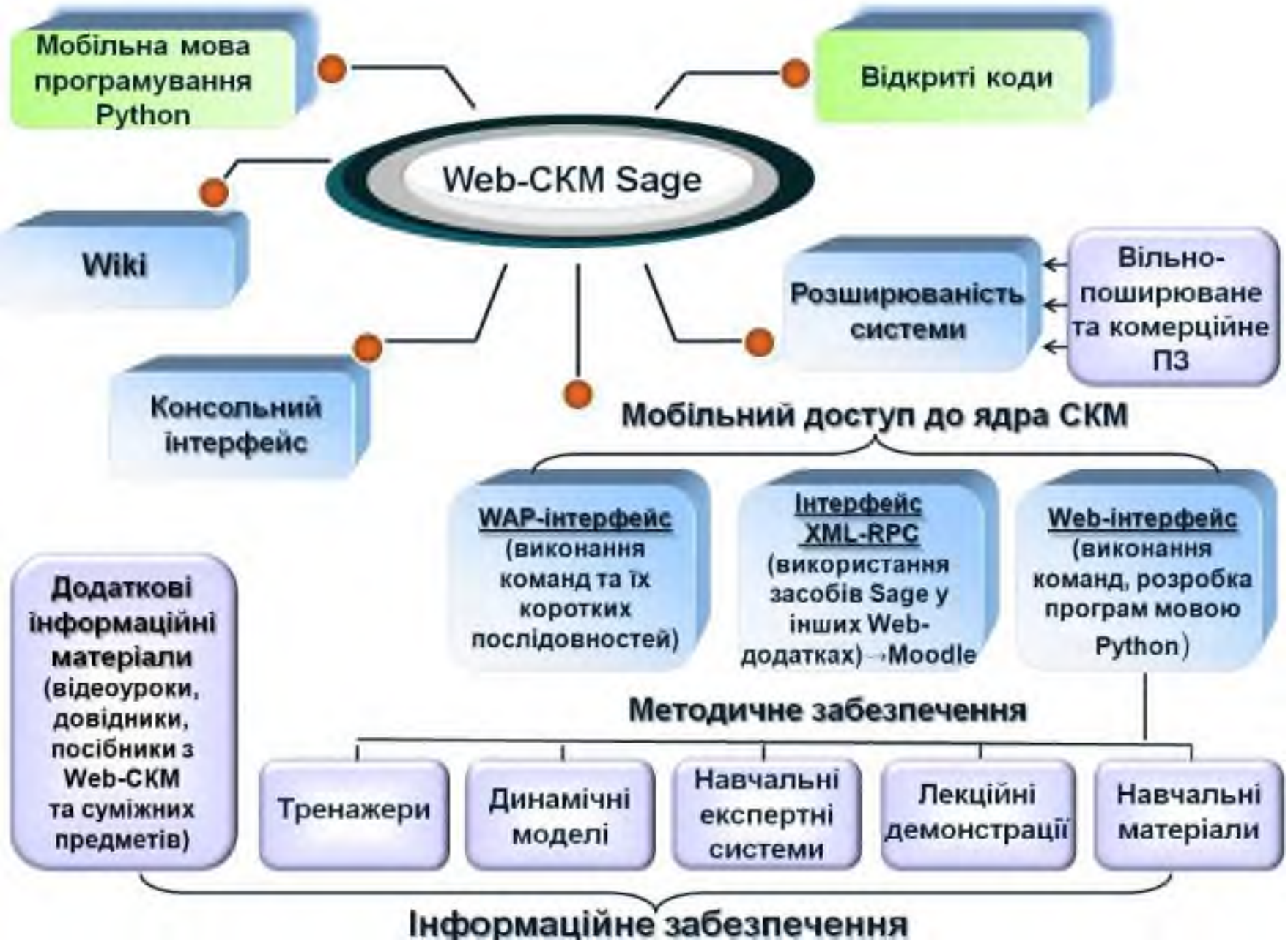

Рис. 3.3. Архітектура мобільного математичного середовища на основі Web-СКМ Sage [299]

У вищих технічних навчальних закладах України у процесі навчання вищої математики найчастіше використовуються такі СКМ: Mathematica, Mathcad, Maple, Derive, SmathStudio, GeoGebra, MathPiper, Sage, DG, що є як комерційними, так і вільно поширювальними.

*3.1.5 Системи тестування.* Забезпечення підготовки фахівців із заданою якістю кваліфікації можливо лише у спеціально організованій педагогічній системі, яка містить прогресивні педагогічні технології навчання, вивчення, організації та управління навчальним процесом. Ці технології передбачають необхідність розробки не тільки специфічної системи засобів і способів їх побудови та впровадження, але й контроль результатів застосування, тобто контроль якості навчання, що сприяє підготовці фахівців на усіх етапах та ступенях освіти. Роль такої системи в сучасних умовах виконує система *тестування,* що є способом одержання відомостей про певний об’єкт і його характеристики шляхом випробовувань [151].

Комп’ютерне тестування є одним із сучасних засобів вимірювання навчальних досягнень. Основними перевагами тестової форми контролю є об’єктивність, економічність, точність, технологічність перевірки виконання завдань. Щодо застосування тестів у математичних дисциплінах, використання комп’ютера надає значні переваги у створенні малюнків, об’ємних зображень, ілюстрацій до питань [279].

І. А. Сверчевською [279] була розроблена програма комп’ютерного тестування «Universal Test System», що призначена для здійснення контролю за процесом засвоєння знань та результатами навчальних досягнень та проведення подальшої корекції знань при вивчені розділу стереометрії «Геометричні тіла». Тестова програма «Universal Test System» складається з чотирьох частин: «Многогранники», «Тіла обертання», «Об’єми» та «Повторення».

Для проведення тестування існують онлайн сервіси, що надають можливість створювати, редагувати і проводити тестування. Для створення тестів у більшості сервісів необхідна реєстрація. Ці сервіси з російськомовним інтерфейсом. У них міститься велика колекція різноманітних тестів, що можуть бути використані у практичній діяльності. Серед онлайн сервісів *Майстер-тест* (www.master-test.net/uk), *БанкТестов.РУ* (www.banktestov.ru), *Твой тест* (www.make-test.ru), *Tests Online* (tests-online.ru), *Тестируем Все* (testing-all.ru), *Анкетёр* (www.anketer.ru/) та інші.

Для організації тестування існує можливість створювати тести у середовищах дистанційного навчання, зокрема Moodle. Тести в Moodle створюються з різними типами питань, використовуються тестові завдання закритого та відкритого типів, допускаються завдання на відповідність, передбачається коротка тестова відповідь, а також числова або обчислювана. Запитання зберігаються в базі даних тестів, є можливість генерації тесту зі списку бази даних питань, а також повторне проходження тесту (декілька спроб тестування), можуть бути використані повторно в цьому курсі або в іншому [151].

Зроблений аналіз надає можливість стверджувати, що на сьогодні впровадження ІКТ у навчання вищої математики студентів інженерних спеціальностей в Україні має велику кількість розробок та впроваджень. Існує досить велика кількість інформаційних технологій для застосування в процесі навчання. Проте ефективне застосування сучасних ІКТ у навчальному процесі можливе лише за умови обґрунтованого та гармонійного інтегрування відповідних технологій у даний процес, забезпечуючи нові можливості і викладачам, і студентам. «Широкий спектр аналітичних, обчислювальних і графічних операцій, що підтримується в сучасних математичних пакетах, роблять їх одними з основних інструментів у професійній діяльності математика-аналітика, інженера, економіста-кібернетика тощо. Тому їх використання у навчальному процесі ВНЗ при вивченні математичних дисциплін надає можливість підвищувати рівень професійної підготовки студентів, рівень їхньої математичної та інформаційної культури» [320].

Проведений аналіз доступних ІКТ навчання вищої математики в Україні показав, що не всі розглянуті впровадження реалізують у повній мірі системний підхід до використання ІКТ, що властивий системі вищої освіти США, де в процес навчання вищої математики студентів інженерних спеціальностей впроваджено математичні пакети, системи підтримки навчання, Інтернет-технології, мультимедійні програмні засоби та інше програмне забезпечення, спрямоване саме на інженерні спеціальності. Крім того, ІКТ навчання вищої математики, що використовуються в США, мають такі визначальні характеристики: тривалий процес розробки та тестування засобів навчання; урахування тенденцій розвитку світових систем (які продовжують бути флагманами в світовому масштабі); реалізація планів по забезпеченню концепції STEM – інтеграції науки, технології, інженерії та математики; прикладна спрямованість навчання вищої математики.

У зв'язку з цим існує об'єктивна необхідність вдосконалення методик використання ІКТ у навчанні вищої математики студентів інженерних спеціальностей ВНЗ України з урахуванням позитивного досвіду США.

## 3.2 Рекомендації щодо використання досвіду Сполучених Штатів Америки у навчанні вищої математики студентів інженерних спеціальностей в Україні

Сьогодні до фахівця інженера роботодавцями поставлено серйозні вимоги, що передбачають співпрацю з творчою особистістю, здатною легко пристосовуватися та добре орієнтуватися в сучасних ринкових умовах. У зв'язку з цим перед ВНЗ постає завдання підготувати спеціалістів, які здатні самостійно приймати рішення, вміти відслідковувати появу нових відомостей, необхідних для їх саморозвитку, вміти оцінювати свою потребу у нових знаннях та навичках, критично підходити до своєї самоосвіти та постійно вдосконалюватися. Однин із засобів, що може допомагати викладачам у вирішенні цього питання є впровадження в навчальний процес ІКТ навчання. Використання ІКТ у фундаментальній підготовці майбутніх спеціалістів у технічних ВНЗ потребує розробки та впровадження змін у методики навчання всіх дисциплін. Це обумовлене тим, що викладач перестає бути єдиним джерелом отримання знань для студента. Сьогодні студенту надано багато можливостей для самоосвіти, оскільки необхідні матеріали можна знайти в мережі Інтернет. Запам'ятовування та відтворення, властиві традиційному навчанню, сьогодні стають не актуальними. Перед студентами виникають нові вимоги: вміння співставляти, аналізувати, оцінювати, прогнозувати та планувати.

Впровадження і поширення освітніх технологій, що використовують ІКТ у процесі навчання, сприяє підвищенню якості вищої освіти. Освітні технології в процесі розвитку, стають більш зручними для використання, вони швидкими темпами проникають у всі дисципліни, оскільки все більше викладачів визнає потребу в підготовці фахівців, обізнаних у призначенні та можливостях використання ІКТ для професійної діяльності.

Зміни у методиці підготовки майбутніх інженерів мають торкнутися як фундаментальних, так і фахових дисциплін. Використання засобів ІКТ у навчанні посилює роль методів активного навчання.

Вивчення вищої математики є основою для фахових дисциплін студентів інженерних спеціальностей (рис. 3.5), тому використовувати ІКТ необхідно починати саме при вивчені вищої математики, основною метою якої є навчити складати математичні моделі процесів і конструкцій, пов'язаних з подальшою діяльністю фахівців, вивчати такі моделі, інтерпретувати відповідно здобуті результати. До знань, вмінь та навичок, якими мають оволодіти студенти після вивчення дисципліни «Вища математика» відносять: знання елементів лінійної алгебри та аналітичної геометрії, вміти використовувати елементи математичного аналізу, диференціального та інтегрального числення, розв'язувати диференціальні рівняння, знати теорію ймовірностей та статистичний аналіз.

З метою вивчення стану впровадження засобів ІКТ у процес навчання вищої математики студентів технічних ВНЗ, нами було проведено анкетування (текст анкети було розроблено Ю. В. Триусом [317] та допрацьовано нами, додаток Л).

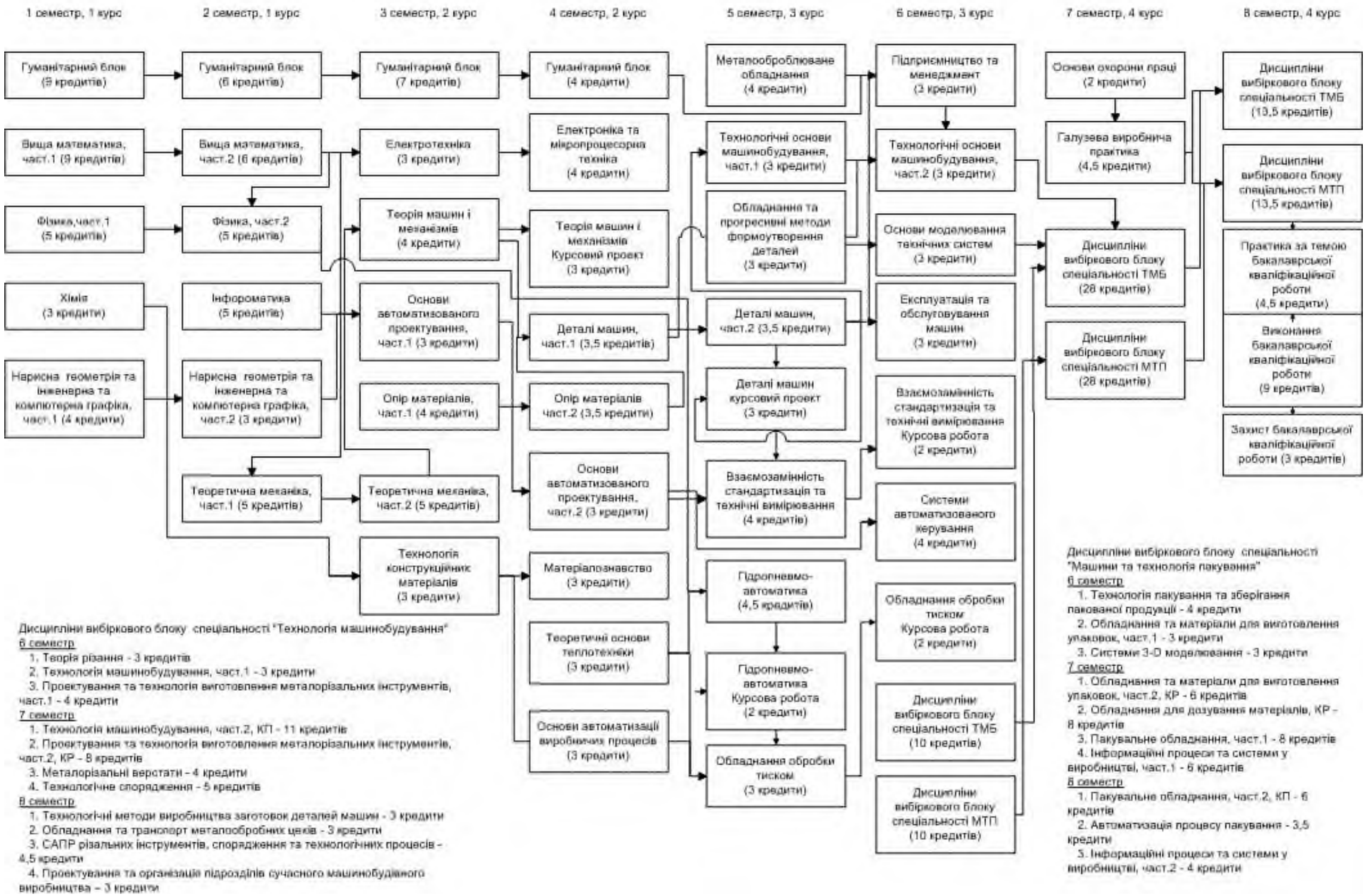

Рис. 3.5. Структурно-логічна схема навчання бакалаврату «Інженерна механіка» [239]

Дане анкетування, проходило в два етапи: опитування викладачів математичних дисциплін за допомогою електронного анкетування в документах Google та опитування викладачів математичних дисциплін за допомогою анкетування, проведеного очно у паперовому вигляді.

В електронному анкетуванні прийняло участь 71 викладач (https://docs.google.com/spreadsheet/gform?key=0ArDq7crTpuk7dFJ3dE9U SjEzeTA0OHlEdTJZbUtJbkE&gridId=0-chart). Спираючись на результати електронного опитування (рис. 3.6, рис. 3.7), можна стверджувати, що опитані викладачі, які користуються електронною поштою, мають навички використання Інтернет, активно впроваджують засоби ІКТ у процес навчання. Серед основних напрямів впровадження засобів ІКТ у навчання вищої математики виділено: створення методичних та дидактичних матеріалів з дисципліни, зокрема мультимедійних (77 %); як джерело відомостей через Internet (70 %); як інструмент для розв’язування задач (48 %); для активізації самостійної роботи студентів (39 %); для вимірювання навчальних досягнень студентів (комп’ютерне тестування, автоматизований контроль) (35 %); як засіб дистанційного навчання математичних дисциплін (23 %).

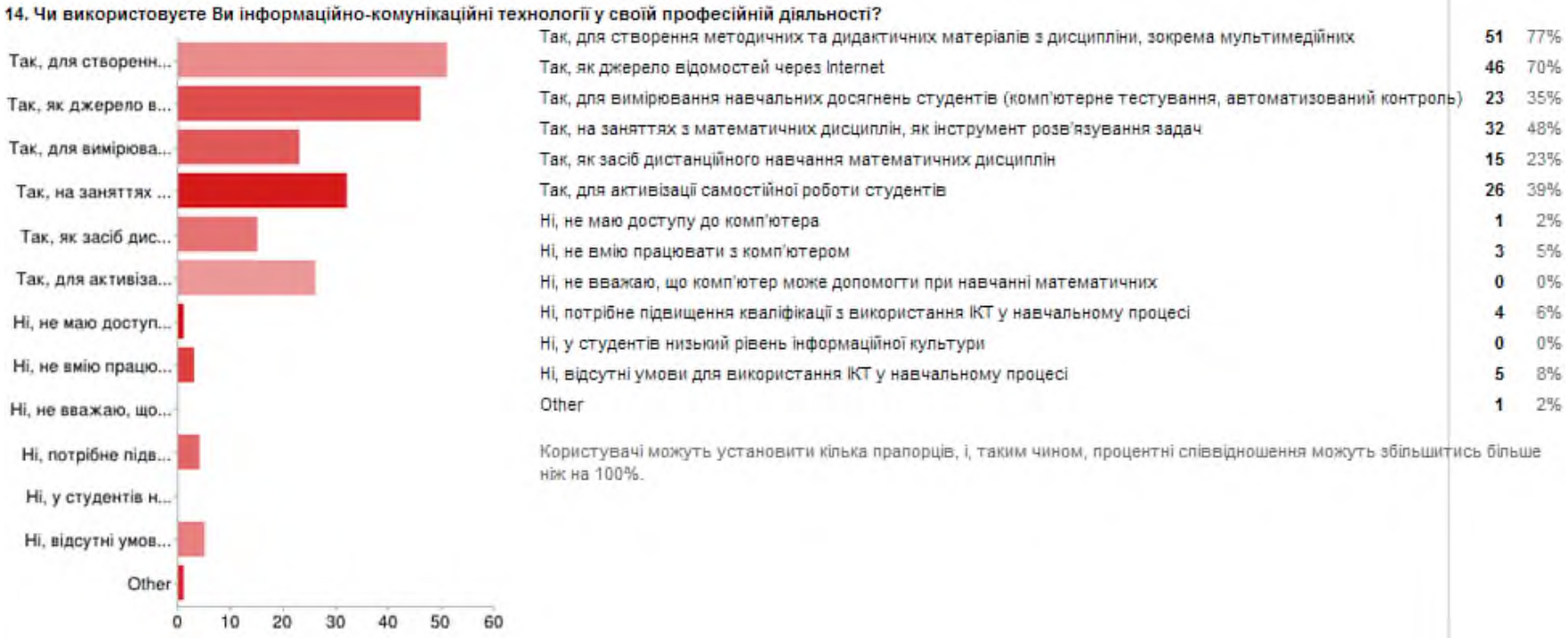

Рис. 3.6. Результати анкетування в документах Google

У паперовому вигляді анкету проходили викладачі кафедри інженерної математики ДВНЗ «Криворізький національний університет» та викладачі математики з технікумів, які проходили стажування в ДВНЗ «Криворізький національний університет». Всього в опитуванні прийняло участь 37 викладачів. Як показали результати паперового опитування, при викладанні вищої математики викладачі здебільшого не використовують засоби ІКТ. Серед основних причин було виділено: не вмію працювати з комп’ютером (16 %); потрібне підвищення кваліфікації з використання ІКТ у навчальному процесі (65 %). Спираючись на результати останнього анкетування, можна стверджувати, що існує необхідність у підготовці рекомендацій та проведенні спецкурсів з теорії та методики використання ІКТ у процесі навчання вищої математики.

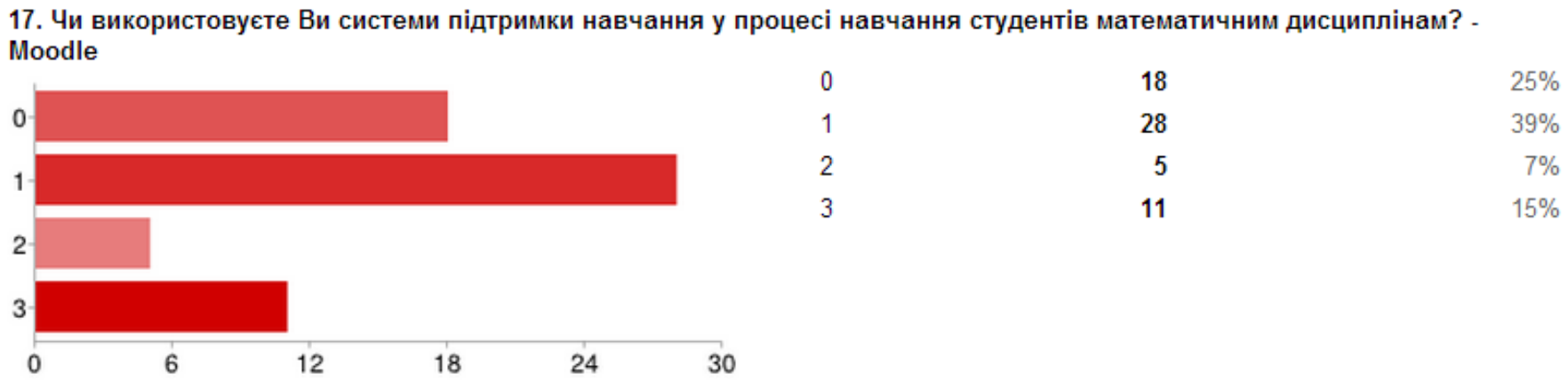

(0 – невідомі, 1 – знаю, але не застосовую в роботі, 2 – знаю і застосовую для своєї роботи, 3 – знаю і застосовую при навчанні математичних дисциплін)

Рис. 3.7. Результати анкетування в документах Google

Аналіз досліджень у галузі використання ІКТ у процесі навчання вищої математики показав, що на сьогоднішній день рівень використання ІКТ викладачами вищої математики у ВНЗ України є досить високим, але не має системного підходу, як зазначалося раніше. Широке використання ІКТ у ВНЗ викладачами та студентами обмежене у зв'язку з відсутністю навичок їх використання, мотивації та інтересів до використання ІКТ для полегшення їх роботи. Хоча, на нашу думку, використання ІКТ викладачами вищої математики у технічних університетах України сприятиме підвищенню рівня професійної підготовки студентів і поліпшенню процесу викладання і навчання, створюючи покоління інженерів, конкурентоспроможне на ринку праці.

Використання ІКТ у процесі навчання надає викладачам вищої математики можливість урізноманітнити лекційні та практичні заняття, проводити демонстрації навчальних матеріалів, організовувати самостійну роботу студентів, підвищувати їх активність та мотивацію. Крім того, використання різних технологій надає можливість викладачам економити час і активізувати увагу студентів під час аудиторних занять.

Проведений аналіз теорії та практики використання засобів ІКТ навчання вищої математики студентів інженерних спеціальностей надав можливість обґрунтувати *фактори, що зумовлюють доцільність поширення досвіду США* у вітчизняному освітньому просторі: перші позиції у рейтингу країн за наявними розробленими та впровадженими засобами ІКТ у процес навчання; демократично-гуманістичні орієнтири освітньої політики щодо реалізації процесу інформатизації вищої освіти, що надає можливість забезпечити індивідуалізацію, диференціацію і персоналізацію процесу навчання; професійна спрямованість заходів щодо інформатизації системи технічної освіти в ВНЗ; тривалий етап розвитку та впровадження ІКТ, що обумовлює їх надійність при використанні в процесі навчання; реалізація системного підходу до використання ІКТ у процесі навчання вищої математики студентів інженерних спеціальностей.

Ґрунтуючись на досвіді використання ІКТ у підготовці майбутніх інженерів у США, слід зазначити, що системне використання ІКТ надає можливість організувати сам процес навчання таким чином, щоб відбувався перехід:

– від навчання під керівництвом викладача до особистісно-орієнтованого навчання;

– від виконання прямих вказівок викладача до інтерактивного навчання;

– від передавання навчальних відомостей до самостійного здобування знань;

– від послідовного повідомлення теоретичних навчальних матеріалів до дослідницького навчання;

– від навчання строго за програмою до навчання за освітніми потребами тих, хто навчається;

– від друкованих навчальних матеріалів до матеріалів в електронному вигляді;

– від традиційного навчання до комбінованого;

– від навчання в аудиторіях до навчання в Інтернет-спільнотах;

– від навчання в університеті до навчання протягом всього життя.

На підставі проведеного педагогічного дослідження розвитку ІКТ навчання вищої математики студентів інженерних спеціальностей у США визначено основні *напрями можливого застосування* позитивного досвіду США у вітчизняному освітньому просторі, а саме:

– впровадження індивідуально-диференційованого підходу до організації навчання студентів через створення належних організаційно-педагогічних умов: наявності у програмі елективних курсів, навчально-методичного забезпечення і налагодженої системи консультування студентів із використанням Інтернет-ресурсів для інтерактивної позааудиторної взаємодії викладачів та студентів;

– урізноманітнення способів подання навчального матеріалу (використання наочних засобів навчання; систем комп'ютерної математики, математичних онлайн ресурсів; навчальних електронних курсів; персональних педагогічних Web-ресурсів);

– урізноманітнення способів перевірки знань студентів із залученням ІКТ – проведення онлайн тестування з використанням тестових систем;

– забезпечення розвитку процесу інформатизації інженерної освіти через інформатизацію змісту вищої математики для забезпечення гнучкості, рефлективності, доступності інженерної освіти;

– визначення напрямів розвитку ІКТ-компетентності викладачів вищої математики, розвиваючи здатність педагогів використовувати ІКТ на всіх етапах процесу навчання.

Розроблено структурно-функціональну схему використання ІКТ у навчанні вищої математики студентів інженерних спеціальностей України (рис. 3.8).

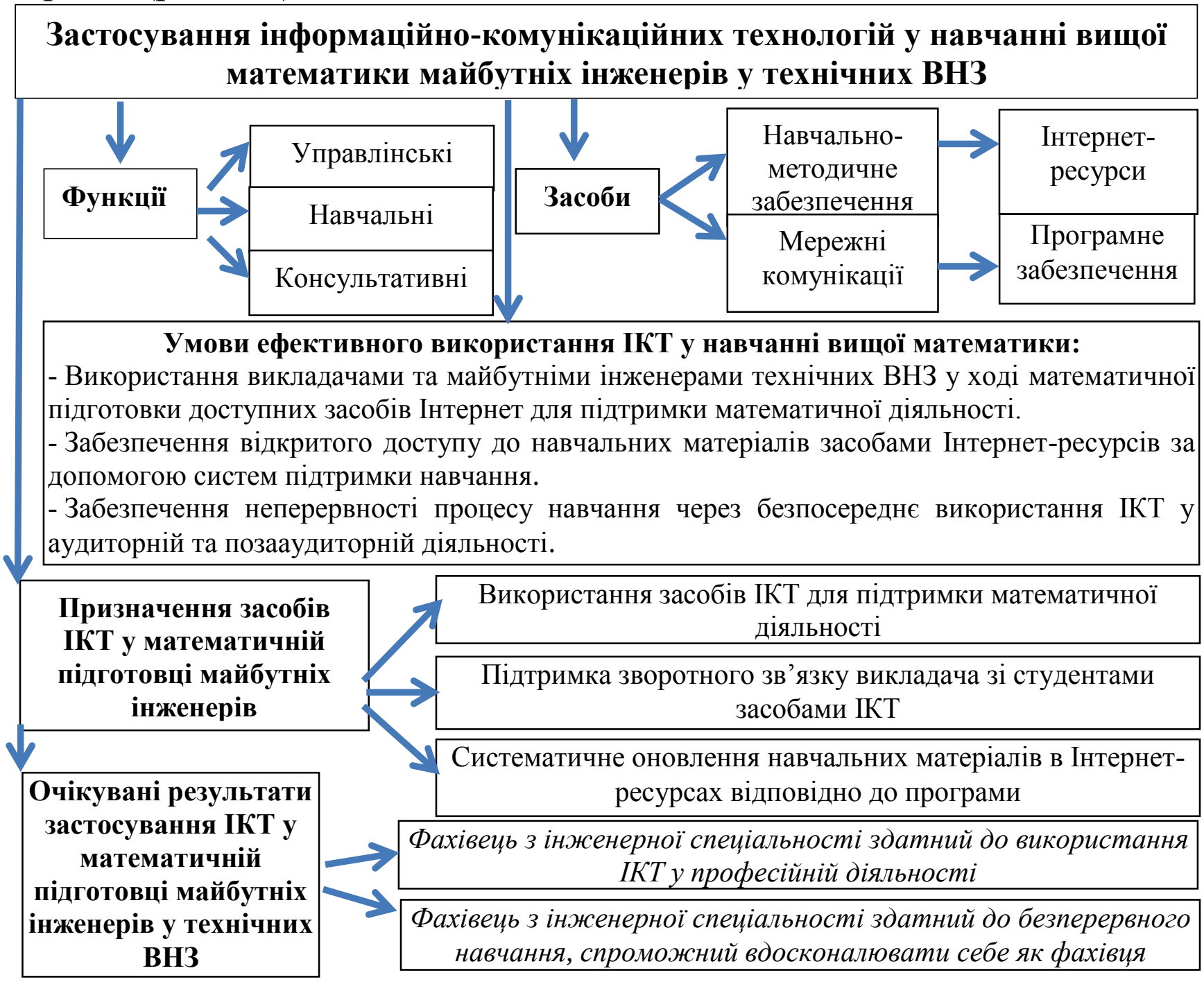

Рис. 3.8. Структурно-функціональна схема використання ІКТ у навчанні вищої математики студентів інженерних спеціальностей

На основі досвіду США виділено основну мету запропонованої схеми – використання ІКТ у навчанні вищої математики студентів інженерних спеціальностей у технічних ВНЗ України.

Застосування ІКТ навчання у математичній підготовці майбутніх інженерів супроводжується управлінськими, навчальними та консультативними функціями з урахуванням навчально-методичних, організаційних і дидактичних вимог, що впливають на результат процесу навчання. Тому функції застосування ІКТ навчання у математичній підготовці майбутніх інженерів розподілені на управлінські, навчальні та консультативні.

*Навчальна функція* застосування ІКТ у процесі навчання вищої математики студентів інженерних спеціальностей технічних ВНЗ України полягає у наданні всіх доступних ІКТ для забезпечення високого рівня засвоєння знань, умінь і навичок з дисципліни. З цією метою пропонується використовувати засоби ІКТ для візуалізації навчального матеріалу (робота над проектами та презентації із використанням сервісів Google, онлайн навчання у WiZiQ); системи комп'ютерної математики, базу знань Wolfram|Alpha, математичні онлайн ресурси; системи підтримки навчання для розробки навчальних електронних курсів; Інтернет-ресурси для інтерактивної позааудиторної взаємодії викладачів та студентів (платформа Piazza, Skype, Google+); персональні педагогічні Web-ресурси – авторські сайти викладачів. Доцільно використовувати Інтернет-технології для надання відомостей щодо навчальних планів, програм, графіків роботи студентів, доступу до навчального матеріалу, до персональних даних студентів; проводити адміністрування роботи студентів, проводити контроль рівня їх знань та вмінь із використанням тестових систем (сервіси Майстер-тест, Банк Тестов.RU, Твой тест, система електронного тестування Tests Online.

*Управлінська* діяльність викладача із використанням засобів ІКТ включає планування, організацію, стимулювання, поточний контроль, регулювання діяльності студентів, аналіз результатів їх навчальних досягнень. Використання засобів ІКТ для управління навчальної діяльності студентів допомагає студентам розвивати вміння аналізувати, корегувати та вдосконалювати власний процес навчання. Маючи відомості про стан успішності студента, про його досягнення та помилки, викладач може правильно скорегувати роботу студента у процесі навчання та роботу з організації самого процесу навчання та методики навчання того чи іншого розділу курсу.

*Консультативна* функція використання ІКТ у процесі навчання вищої математики студентів інженерних спеціальностей ВНЗ України сприяє кращому засвоєнню навчального матеріалу. Студенти мають усвідомлювати необхідність використання отриманих знань з вищої математики в процесі професійної діяльності. Тому для забезпечення високої якості навчання викладач може залучати студентів до участі в тематичних чатах та форумах у системах підтримки навчання або соціальних мережах, використовувати електронну пошту для особистого листування. Викладач має змогу наповнювати необхідними матеріалами навчальні курси і проводити консультації на відстані в зручний для себе і студентів час, розподіляти завдання, вести електронний журнал обліку активності студентів у навчанні, налаштовувати різноманітні ресурси навчального курсу тощо. Використовуючи соціальні мережі, викладач має змогу створювати web-сторінки, що містять різноманітні навчальні матеріали.

Виокремимо *умови ефективного використання ІКТ* у навчанні вищої математики студентів інженерних спеціальностей:

– використання викладачами та майбутніми інженерами технічних ВНЗ у ході математичної підготовки доступних засобів Інтернет для підтримки математичної діяльності;

– забезпечення відкритого доступу до навчальних матеріалів засобами Інтернет-ресурсів за допомогою систем підтримки навчання;

– забезпечення неперервності процесу навчання через безпосереднє використання ІКТ у аудиторній та позааудиторній діяльності.

Отже, використання ІКТ у навчанні вищої математики надає можливість регулювати навчальну діяльність студентів, розвивати їх навчальні інтересів. Впровадження в процес навчання ІКТ і надання студентам інженерних спеціальностей та викладачам необхідних методичних рекомендацій щодо використання ІКТ у математичній підготовці створює умови для підвищення якості підготовки майбутніх інженерів у технічних ВНЗ України.

При застосуванні ІКТ у математичній підготовці майбутніх інженерів необхідно враховувати *професійну спрямованість* застосування запропонованих ІКТ, а навчальні матеріали в Інтернет-ресурсах мають *систематично оновлюватися* відповідно до програми навчання та спеціальності студентів.

Засоби ІКТ у математичній підготовці майбутніх інженерів доцільно використовувати за таким призначенням:

– використання засобів ІКТ для підтримки математичної діяльності студентів;

– підтримка зворотного зв’язку викладача зі студентами засобами ІКТ;

– систематичне оновлення навчальних матеріалів в Інтернет-ресурсах відповідно до програми.

С. Тренхолм (Sven Trenholm) [75] зазначає, що більшість університетів по всьому світу в даний час інтегрує віртуальні навчальні середовища (virtual learning environments – VLE) в свої програми вищої освіти. Ці Web-інструменти можуть бути використані для розробки альтернативних та доповнення традиційних форм навчання і надають можливість зробити доступними інструкцію для студентів, які обмежені в часі або місці. Це, в свою чергу, розглядається як потенційний і ефективний засіб для вирішення проблеми у деяких розвинених країнах, пов’язаної з зростаючим попитом на вищу освіту. Такі платформи навчання за допомогою Інтернет і мультимедіа в даний час забезпечують, наприклад, студентам можливість зручно отримати доступ до всіх або частини необхідних їм навчальним матеріалам, проходити тестування, доступ до комплекту домашнього завдання, брати участь у різних індивідуальних та спільних навчальних заходах, задаючи питання викладачам або спільній групі студентів та вирішувати навчальні проблеми.

Виділимо основні чинники, що необхідно враховувати для успішної розробки методики навчання та розвитку навчання з вищої математики з використанням ІКТ:

– *основний навчальний матеріал курсу* – є джерелом відомостей для студентів у процесі навчання студентів. Він повинен бути розрахований на самостійне навчання. Він також повинен надавати студентам можливість зрозуміти усі суттєві аспекти курсу. Бажано, щоб всі ресурси були розроблені і написані однією і тією ж людиною (або проектною групою). У зв'язку з характерними властивостями математичних дисциплін, де студентам необхідно багато зрозуміти і вивчити під час процесу навчання, необхідно забезпечити доступ до основних матеріалів курсу в електронному вигляді, зручному для друку. Ці основні навчальні матеріали можуть (і повинні) бути доповнені додатковими навчальними матеріалами і математичними ресурсами, такими як аплети, статті, симулятори;

– *роль викладачів як Інтернет-інструкторів*: існує необхідність у забезпеченні керівництва і підтримки процесу навчання студентів, а також безперервного та своєчасного зворотного зв'язку під час всього процесу навчання. Це керівництво має здійснюватися на основі надісланих повідомлень (наприклад, на початку кожного тижня) з чіткими інструкціями про те, який зміст і вид діяльності повинен бути виконаний в найближчий час. У той час як студенти проходять курс, підтримка повинна бути забезпечена швидкою реакцією на запитання студентів у форумах, електронною поштою, засобами синхронного зв'язку. Такий зворотний зв'язок повинен бути своєчасним, в той же або наступний день, коли було задано запитання. За необхідності, має бути забезпечена координація між різними викладачами одного і того ж курсу для забезпечення однорідності навчального процесу;

– *ефективне використання програмного забезпечення навчання математики*: важливо, щоб студенти розуміли, що математичні дисципліни забезпечують практичні знання і навички, необхідні для професійної діяльності. Такий підхід буде цінуватися студентами і сприятиме підвищенню рівня мотивації. Наявне математичне програмне забезпечення має подібні можливості, тому мова йде не стільки про те, яке спеціальне програмне забезпечення використовувати, скільки про те, як ефективно інтегрувати його у навчальний процес. А яке програмне забезпечення використовувати, залежить від умов, уподобань та знань викладача, що мають постійно розширюватися;

– *доступність віртуального навчального середовища*: найважливішим аспектом будь-якого навчального середовища є його доступність, студенти та викладачі повинні відчувати себе комфортно, використовуючи дане навчальне середовище, а всі основні його параметри повинні бути інтуїтивно зрозумілими. Для курсів з вищої математики, як і для будь-якої іншої дисципліни, в навчальному просторі має бути забезпечена можливість розміщувати відзначені викладачем завдання – офіційні повідомлення від викладачів до студентів, забезпечена можливість надсилати нотатки та вести дискусії та обговорення, що стосуються змісту курсу. Інші бажані функції навчального середовища – наявність засобів забезпечення математичної діяльності і функція моніторингу навчальної діяльності студентів.

У зв'язку з тим, що розвиток ІКТ відбувається постійно і безперервно, особливо стрімко в останні роки, необхідно враховувати це в методиці використання ІКТ у процесі навчання.

Комбінований підхід до вибору технологій для навчання студентів надає можливість для створення зрозумілого, доступного, особистісно-орієнтованого процесу навчання, в якому студенти зорієнтовані на співпрацю між собою та з викладачем. В структуру комп'ютерно орієнтованого навчально-методичного комплексу (КОНМК) викладача вищої математики необхідно включити новітні ІКТ. На рис. 3.9 запропоновано структуру КОНМК з вищої математики для студентів інженерних спеціальностей. До КОНМК включено служби та сервіси мережі Інтернет, що допомагає організувати самостійну роботу студентів.

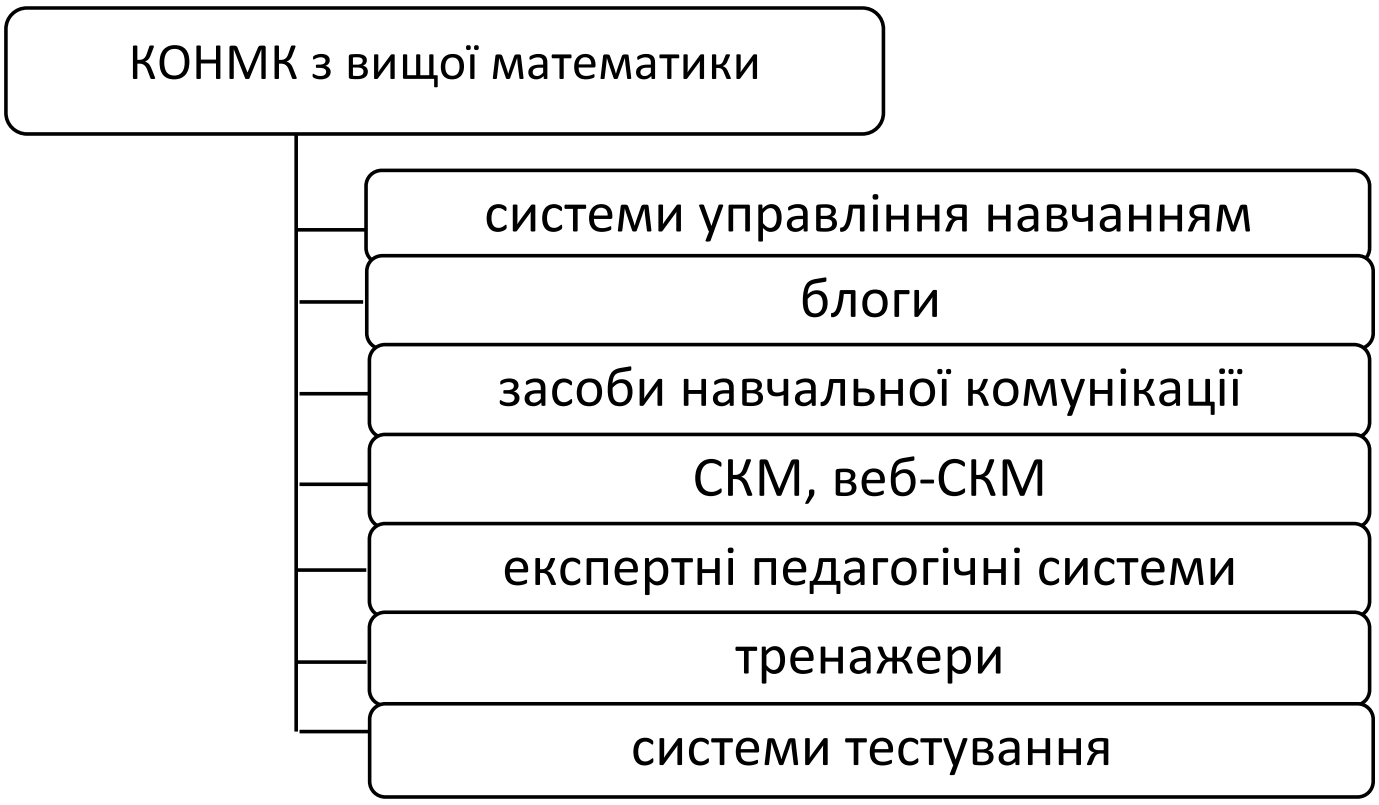

Рис. 3.9. Структура КОНМК з вищої математики

Опишемо запропоновані компоненти КОМНК.

Застосування системи управління навчанням у математичній підготовці майбутніх інженерів надає можливість покращити процес управління. Вибір системи управління навчання все частіше обумовлюється ціновою політикою. У ВНЗ України викладачі використовують найчастіше вільно поширюванні платформи для організації навчання, що вже зазначалося раніше.

Важливу роль у організації навчального процесу, на наш погляд, відіграють авторські сайти викладачів, СКМ та Web-СКМ, експертні педагогічні системи, тренажери, системи тестування.

На авторських сайтах викладач має змогу надавати відомості про навчальний курс, розміщувати методичні матеріали, проводити інтерактивне спілкування зі студентами.

На сьогодні існує велика кількість конструкторів сайтів, що не вимагають знань з програмування і надають можливість користувачу швидко і легко створювати власні сайти. Одним з таких конструкторів є сайти Google (sites.google.com). Використовуючи конструктор сайтів Google, можна легко створювати і оновлювати власні сайти. За допомогою конструктору сайтів Google можна інтегрувати та агрегувати різноманітні дані в одному місці, включаючи відео, слайдшоу, календарі, презентації, вкладення, додатки, текст, надаючи можливість перегляду і редагування цих даних невеликій групі осіб, організації або всім відвідувачам сайту.

Серед доступних можливостей роботи з сайтами Google є: налаштування сайту відповідно до власних потреб, створення підсторінок для упорядкування змісту сайту, вибір типів сторінок (Web-сторінка, оголошення, картотека), централізоване зберігання Web-вмісту та автономних файлів, можливість управління правами доступу до сайту, можливість пошуку змісту на сайтах Google з використанням пошукових технологій Google.

Нами було створено сайт для підтримки навчання студентів вищій математиці (https://sites.google.com/site/visamatematikavidkianovskoie/home/golovn) (рис. 3.10). Детальний опис створення сайту подано в додатку М.

До сайту включено такі сторінки: новини, журнали груп, лекції, література, матеріали, питання та завдання, практичні заняття, самостійна робота, тестування, форум та карта сайту. На сторінках «Лекції», «Практичні заняття» та «Самостійна робота» створено підсторінки, що містять навчальні матеріали, доступні для завантаження (рис. 3.11).

Сторінка «Журнали груп» містить посилання на сторінку з таблицями Google Диску. В Таблицях Google надається можливість використовувати ті ж формули, що підтримуються в більшості поширених програм для роботи з електронними таблицями. За допомогою цих формул можна створювати функції для обробки даних про навчальні досягнення студентів.

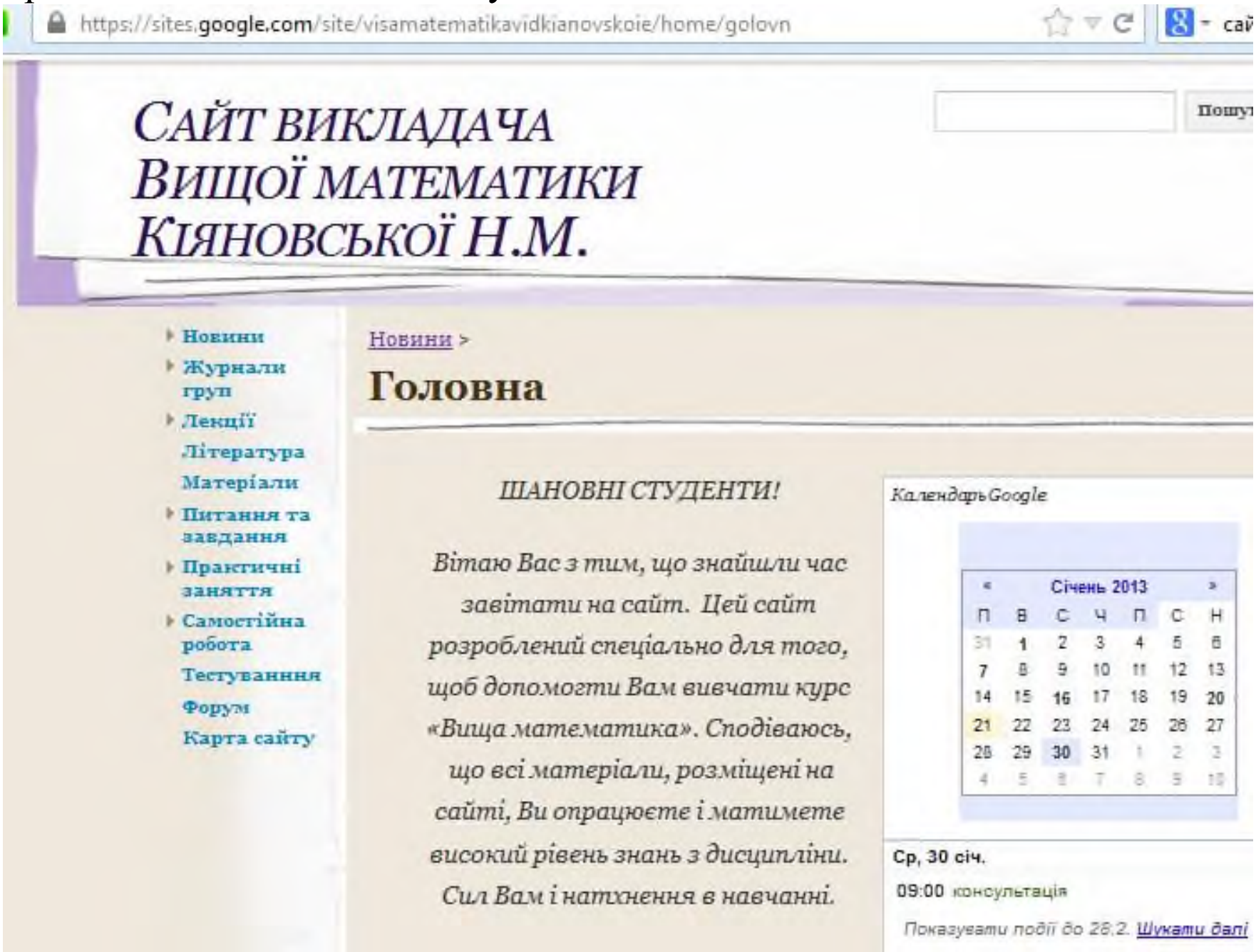

Рис. 3.10. Фрагмент головного вікна сайту для підтримки навчання студентів з вищої математики

Сторінка «Матеріали» містить допоміжні навчальні матеріали для підготовки до лекційних та практичних занять, до підсумкового контролю знань, вмінь та навичок, а також посилання на сайти, що можуть стати в нагоді студентам при вивченні вищої математики (рис. 3.12).

При проведення складних обчислень проміжного характеру, на розв'язування яких витрачається велика кількість часу, доцільно *використовувати СКМ.* Існування великої кількості наявних математичних пакетів надає можливість викладачу вибрати зручний, доступний та зрозумілий для нього ресурс, враховуючи зазначені переваги та недоліки цього програмного забезпечення.

Використовуючи в процесі навчання вищої математики СКМ, необхідно враховувати наступне:

– надавати студентам відомості про різні доступні математичні

пакети із наданням методичних матеріалів щодо їх використання, забезпечуючи можливість студентів самостійно вибирати зручний для них продукт;

– пропонувати студентам задачі, що мають бути розв'язані із використанням СКМ;

– використовувати СКМ на заняттях для розв'язування проміжних обчислень, що займають велику кількість часу;

– пропонувати студентам завдання з вищої математики професійної спрямованості, що мають бути розв'язані із використанням СКМ;

– поступове ускладнення та розширення застосування СКМ в процесі навчання;

– розвивати у студентів мотивацію до самостійного оволодіння вмінням застосовувати СКМ в процесі навчання;

– використовувати СКМ для демонстрації властивостей понять, що вивчаються та є складними для розуміння;

– використовувати СКМ для візуалізації процесу навчання математики.

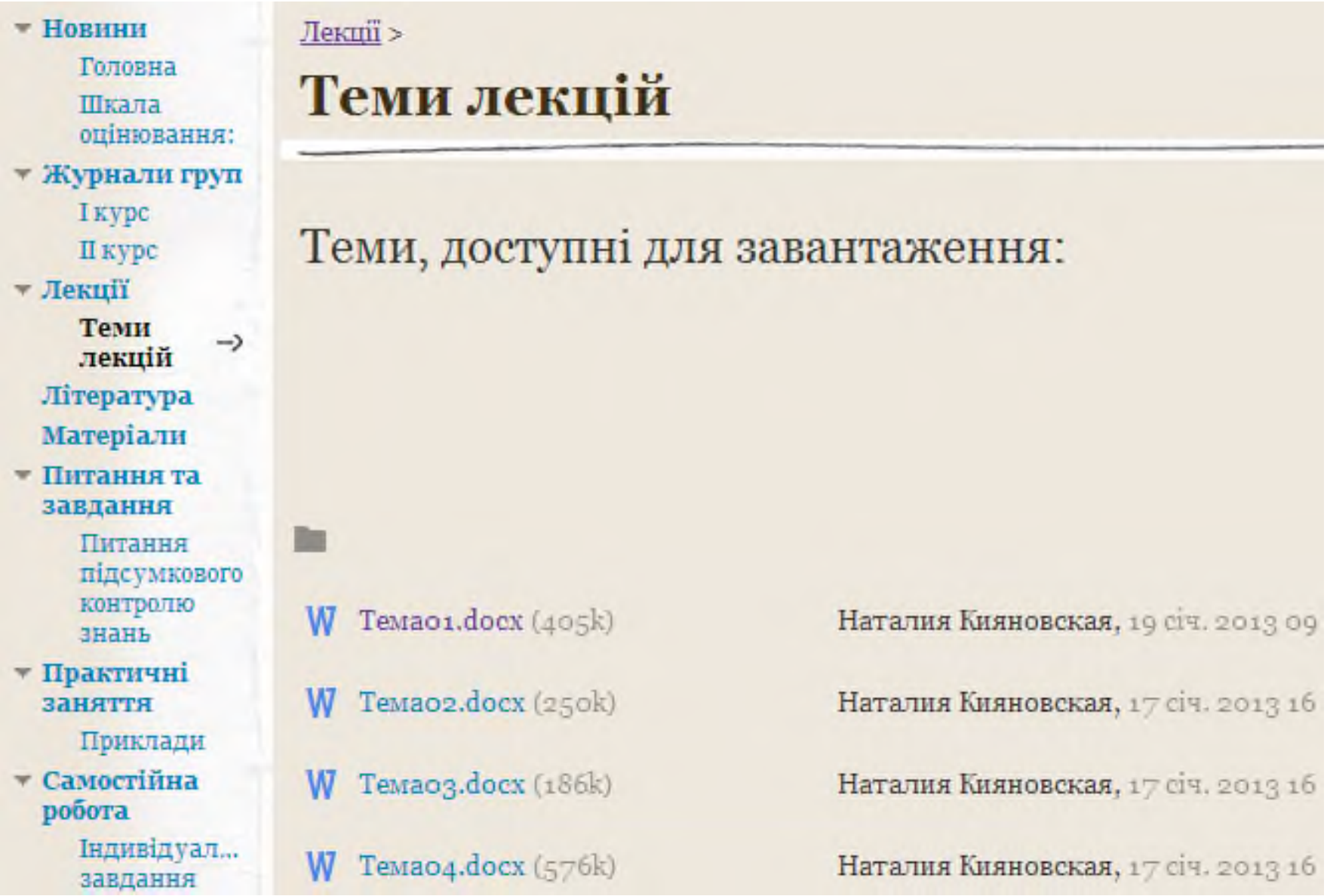

Рис. 3.11. Фрагмент сторінки сайту «Теми лекцій»

Існує велика кількість методичних розробок з впровадження СКМ в процес навчання. Так, наприклад, використання математичного пакету MATLAB для розв'язування прикладних задач студентами механіко-математичного факультету розглянуто в посібнику Б. П. Довгого [155], розв'язування задач із використанням Mathcad розглянуто в посібнику В. О. Доровського [167], символьні обчислення в системі Maple в посібнику А. В. Кузьміна [220], крім того, на освітньому сайті www.exponenta.ru містяться керівництва з використання MATLAB, Mathematica, Mathcad, Maple, Statistica та з інших пакетів.

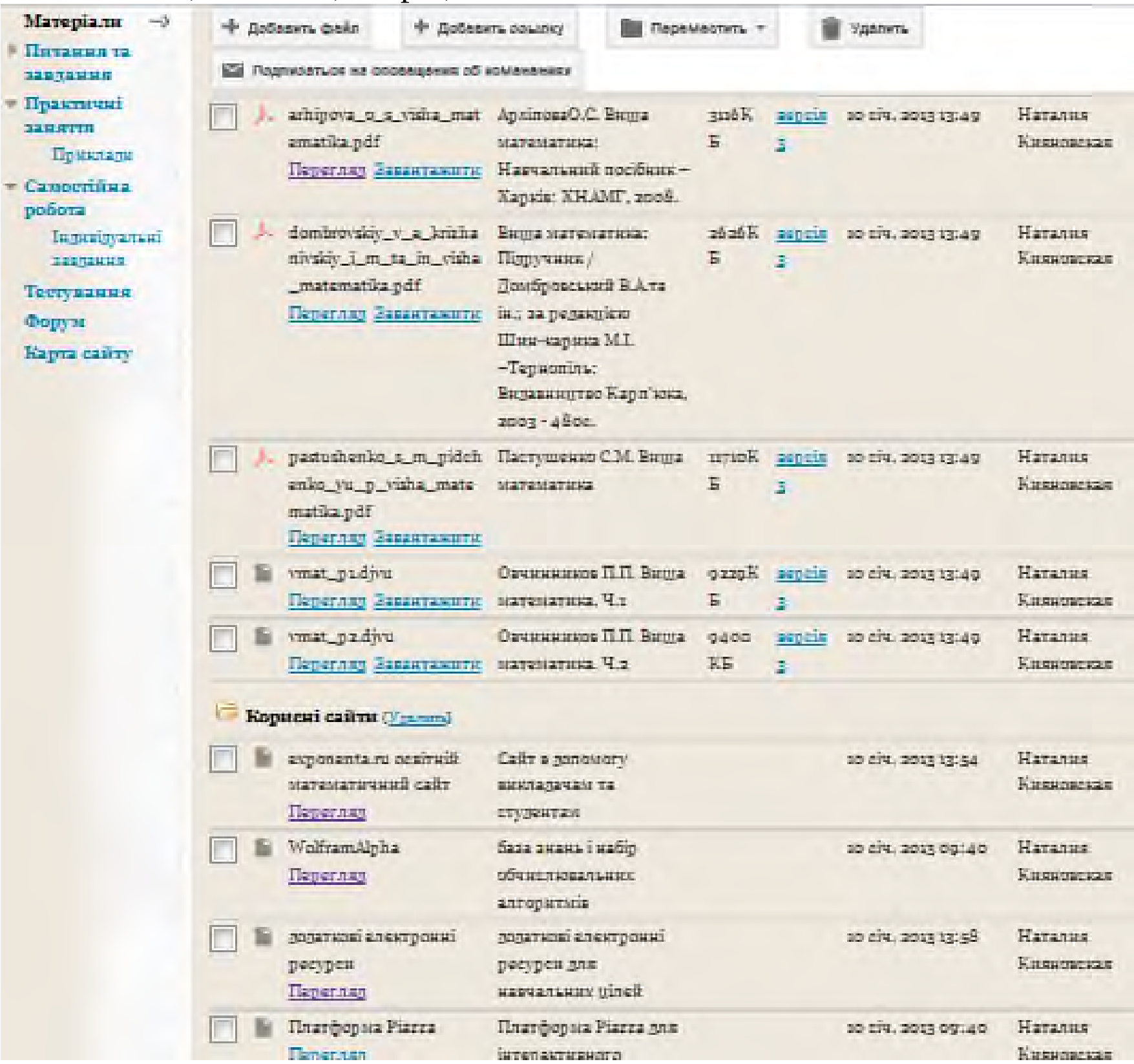

Рис. 3.12. Фрагмент сторінки «Матеріали» створеного сайту

При вивченні тем «Циліндричні поверхні», «Конічні поверхні», «Поверхні обертання», «Поверхні другого порядку», зображення поверхонь доцільно будувати в одній із СКМ і на лекційному занятті проводити демонстрації одержаних зображень поверхонь. Так у системі Maple для тривимірних побудов використовується функція plot3d: plot3d(expr, x=a..b, y=c..d), plot3d(f, a..b, c..d), де параметри **expr** – вираз відносно **x** та **y**; **f** – процедура або оператор; **a**, **b** – дійсні константи, процедури, або вирази в **y**; **c**, **d** – дійсні константи, процедури, або вирази в **x**; **x**, **y**, **t** – імена.

За допомогою функції plot3d можна побудувати зображення складних поверхонь (рис. 3.13), приклади написання цієї функції подано в довідковій системі Maple.

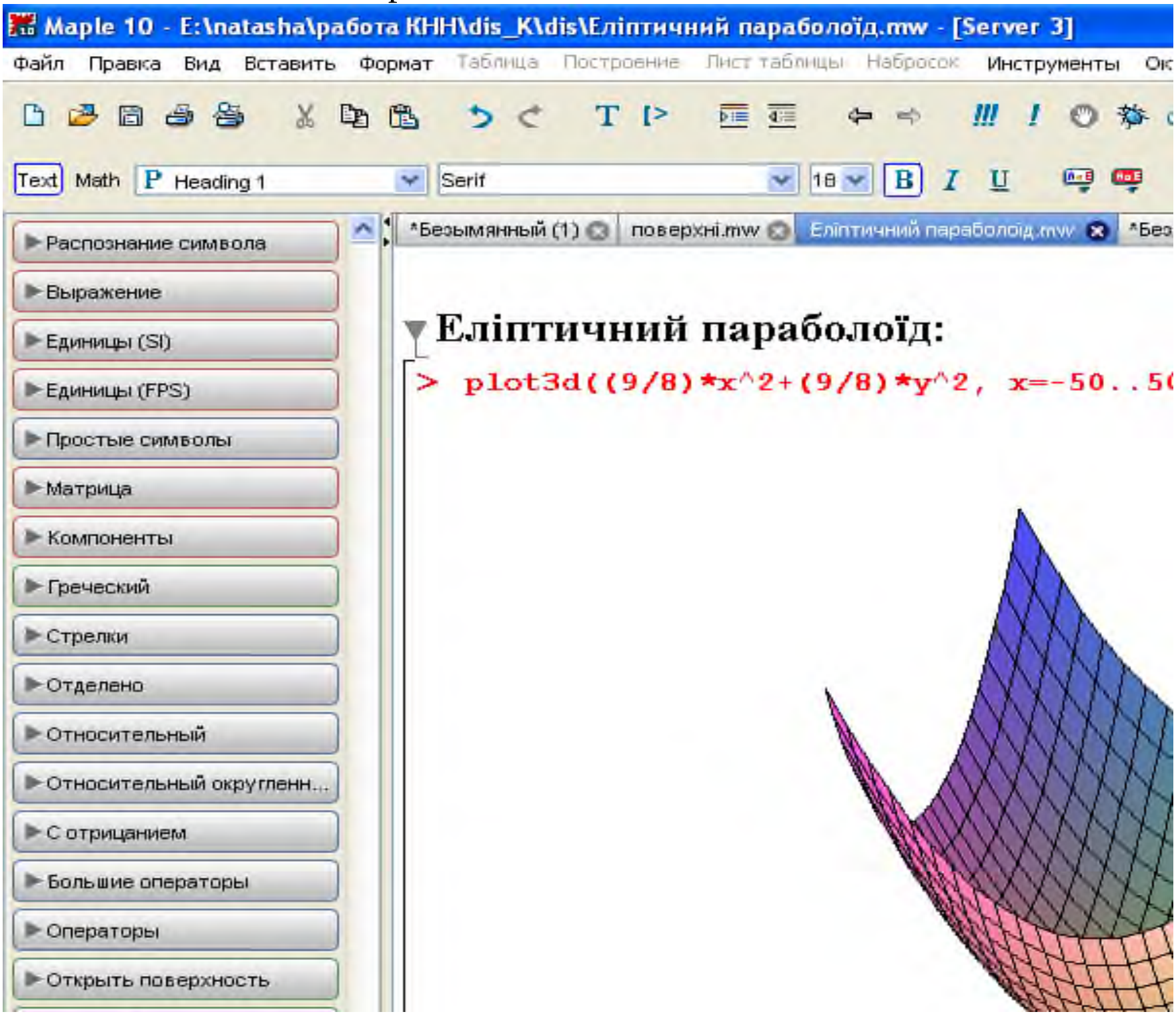

Рис. 3.13. Побудова еліптичного параболоїда в Maple

Крім СКМ існує велика кількість різноманітних онлайн-ресурсів для підтримки математичної діяльності, що є простими та зручними у використанні.

Використання математичних пакетів в процес навчання вищої математики студентів інженерних спеціальностей ВНЗ України надає можливість розширити знання студентів в області новітніх програмних засобів з математики та вміння застосовувати ці програмні засоби при розв'язуванні задач як з вищої математики, так і задач з дисциплін фахової направленості впродовж всього навчання у ВНЗ.

Потужними засобами для осмислення навчального матеріалу, виділення основних його компонентів та систематизації знань студентів є *експертні педагогічні системи, тренажери, системи тестування*.

Основне призначення тренажерів, на думку [226] – це осмислення і закріплення теоретичного матеріалу, контроль знань з теми дослідження. У системі вищої освіти розроблено достатню кількість тренажерів з вищої математики: *тренажер «Системи лінійних рівнянь»*, створений О. В. Тутовою [322]; *евристичний тренажер «Gauss»,* створений Т. С. Максимовою [157]; також Т. С. Максимовою створено *тренажери «Continuity and Graphics»* та *«Limit»*, за допомогою яких студенти повторюють і систематизують такі поняття: область визначення функції, область значень, монотонність, нулі функції, неперервність функції, графік функції однієї змінної та оволодівають вміннями по обчисленню границь функцій і послідовностей.

*Експертні педагогічні системи* – спеціалізовані комп'ютерні програми для обробки і структурування аморфних (неповних, невизначених, ненадійних, складних, заплутаних, суперечливих) даних та формування висновків у галузях, де такі дані циркулюють і використовується, тобто в педагогіці, психології, соціології, економіці, екології, податковій справі, географії, політиці, юриспруденції тощо [248].

Розроблення експертних систем спрямоване на використання комп'ютерів для опрацювання даних у тих галузях науки і техніки, де традиційні методи моделювання малопридатні. Основу експертних систем складає база знань про предметну область, що накопичується в процесі побудови й експлуатації експертної системи [145].

Виділяють декілька різновидів експертних педагогічних систем: експертні системи для підтримки навчального діалогу, експертні системи навчання мов або системи перекладу, експертні системи навчання предметних або штучних мов, експертні системи класифікації, проблемно-орієнтовані експертні системи, експертні системи доведення теорем [312].

У процесі вивчення вищої математики доцільно використовувати експертну оболонку eXpertise2Go, що є вільно поширюваним Web-орієнтованим програмним засобом (Web-HEC). Зазначена експертна оболонка надає можливість генерувати базу знань у форматі e2gRuleEngine за допомогою інструменту для створення та перевірки таблиць розв'язків e2gRuleWriter.

Як зазначалось раніше, існує велика кількість засобів для створення та проведення онлайн тестування, розроблено чимало начальних експертних систем, тренажерів, використання яких може допомогти студентам підготуватися якісно до залікових модульних робіт з вищої математики.

Проходячи тематичне онлайн тестування, студенти мають можливість: систематизувати отримані знання з теми, виділити основний матеріал теми, з’ясувати питання, що викликають ускладнення, здійснити самоконтроль, пройти тестування в зручний час.

Наприклад, на сайті Mathtest.ru (рис. 3.14), розробниками якого є викладачі кафедри математики національно-дослідного Іркутського державного технічного університету, запропоновано готові тести зі шкільної математики та з курсу вищої математики.

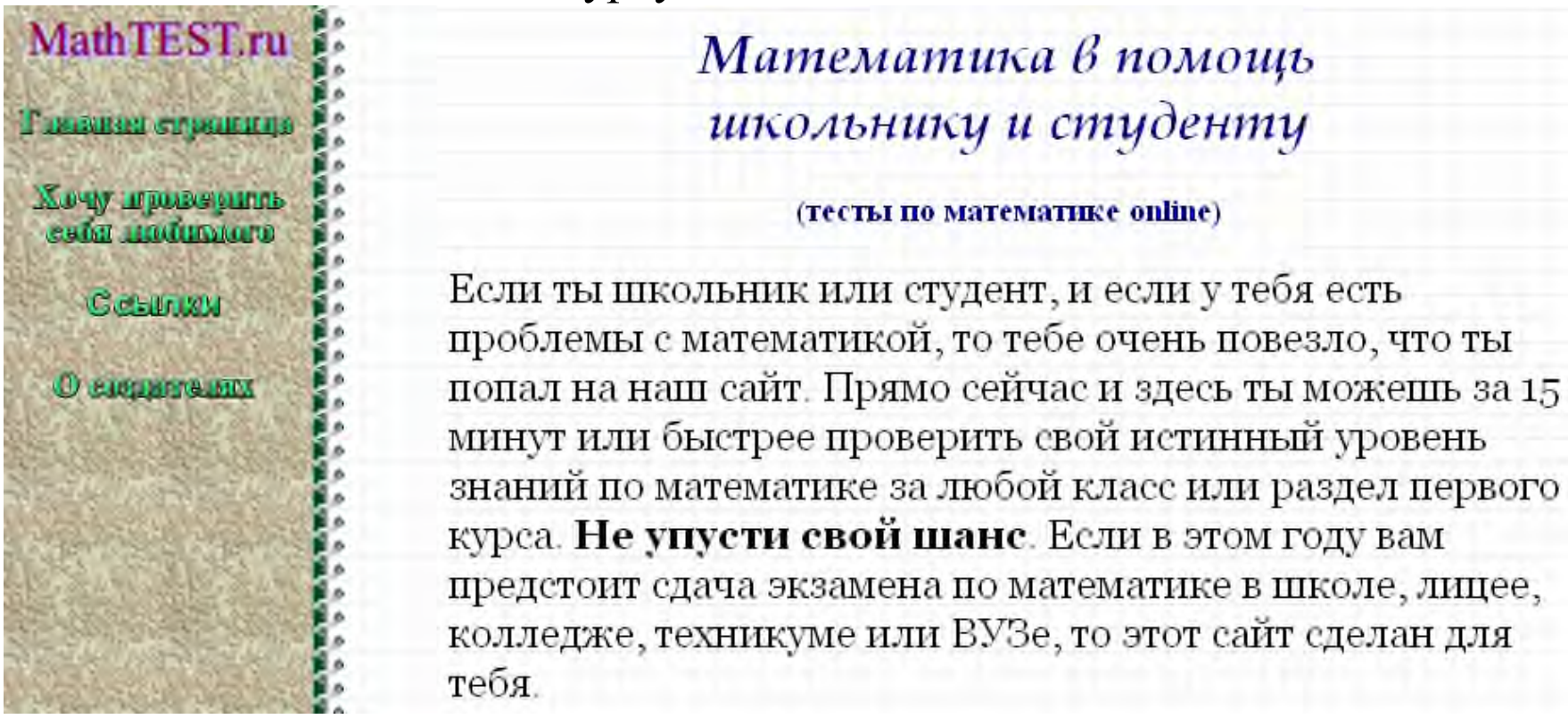

Рис. 3.14. Сайт Mathtest.ru для проведення тестування з математики

Перевагою запропонованого тестування є те, що в ньому вказуються неправильні відповіді студента і питання задаються до тих пір, поки студент не дасть правильну відповідь на питання (рис. 3.15). Таким чином, проходячи тестування, студент має змогу вчити теоретичний матеріал або корегувати свій рівень знань.

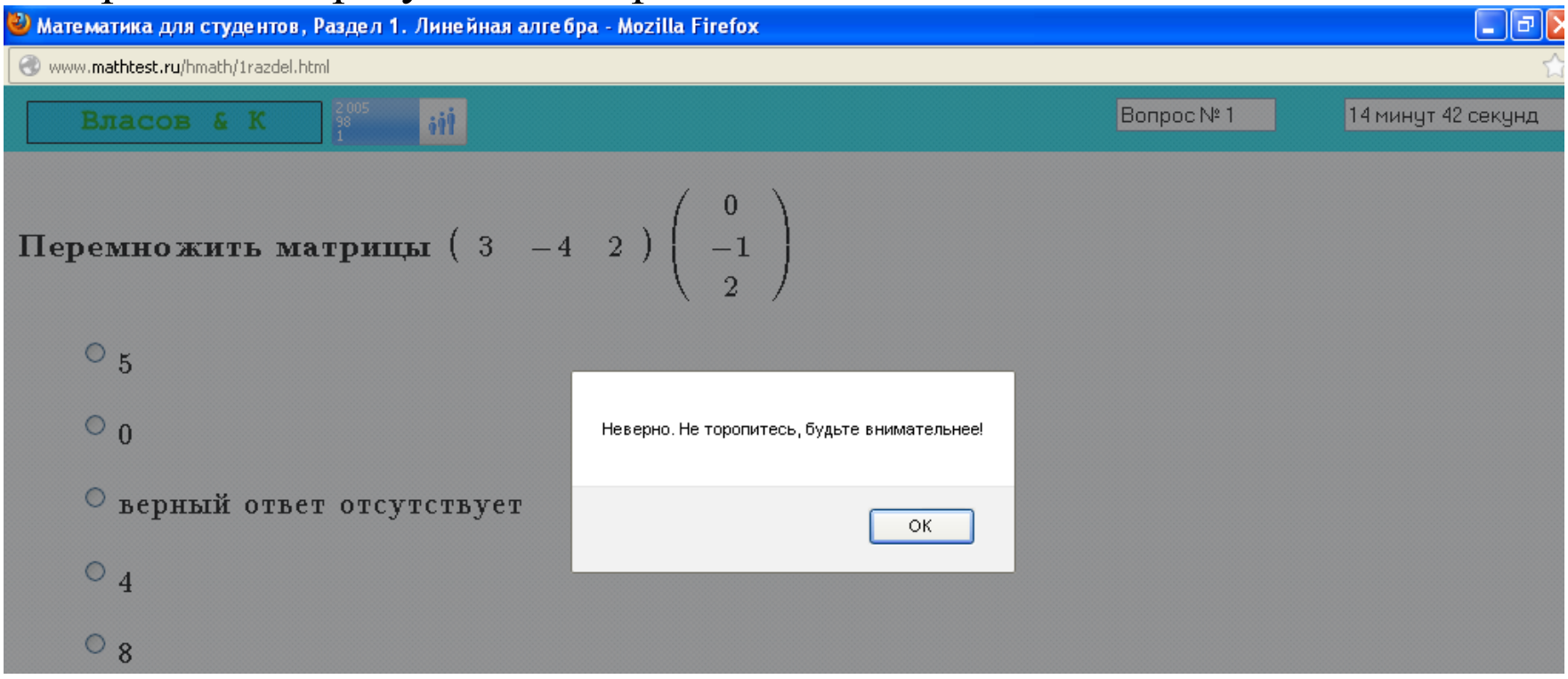

Рис. 3.15. Коментарі до відповідей на сайті Mathtest.ru

Крім використання власного сайту, електронної пошти, пошукових систем, освітніх порталів, вікі, блогів, корисним, на наш погляд, буде впровадження в процес навчання для підтримки позааудиторної навчальної комунікації *платформи Piazza*, що є розповсюдженою в багатьох ВНЗ США.

Використання платформи Piazza надає можливість викладачам та студентам продовжувати роботу поза аудиторією, підтримуючи інтерактивне спілкування зі студентами. Розміщуючи питання на Piazza, студенти можуть отримати відповідь на них як від викладача, так і від одногрупників. Імітуючи реальне обговорення в аудиторії, використання Piazza допомагає підвищити активність студентів у навчанні і заощаджує час викладача та студентів, оскільки всі проблеми вирішуються колективно та оперативно.

Створення класу та реєстрація користувачів в платформі Piazza відбувається англійською мовою, але проведення спілкування, надання всіх даних та навчальних матеріалів можна як українською мовою, так і іншою, зручною користувачам. Для створення свого класу в Piazza необхідно створити або вибрати з існуючих навчальний заклад (рис. 3.16 а).

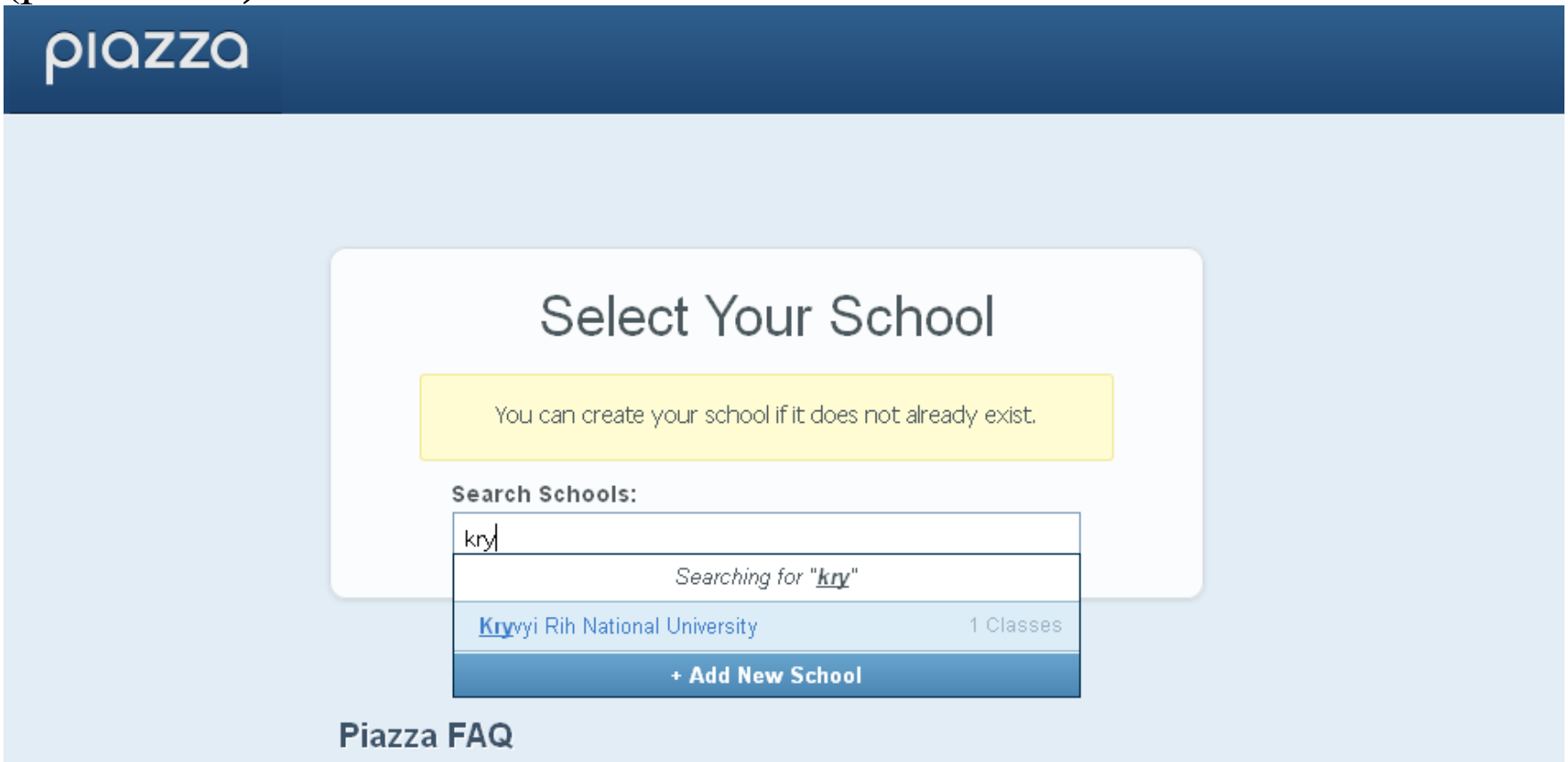

Рис. 3.16 а) Вибір навчального закладу в Piazza

На наступному кроці необхідно створити або вибрати з існуючих клас в Piazza (рис 3.16 б). Створивши клас, викладачу пропонується внести данні про новий клас: назва класу, номер класу, кількість студентів у класі. (рис 3.16 в). Далі пропонується зазначити роль користувача класу, який реєструється (як студент або як викладач), зазначивши електронну адресу (рис. 3.16 г). На електронну пошту приходить повідомлення для активації класу Piazza, після проходження активації можна починати наповнювати сторінку класу навчальними матеріалами, додати для завантаження студентами навчальний план (рис. 3.16 д). Викладачам надається можливість вносити данні про час та місце проведення аудиторних консультацій в навчальному закладі (рис. 3.16 є).

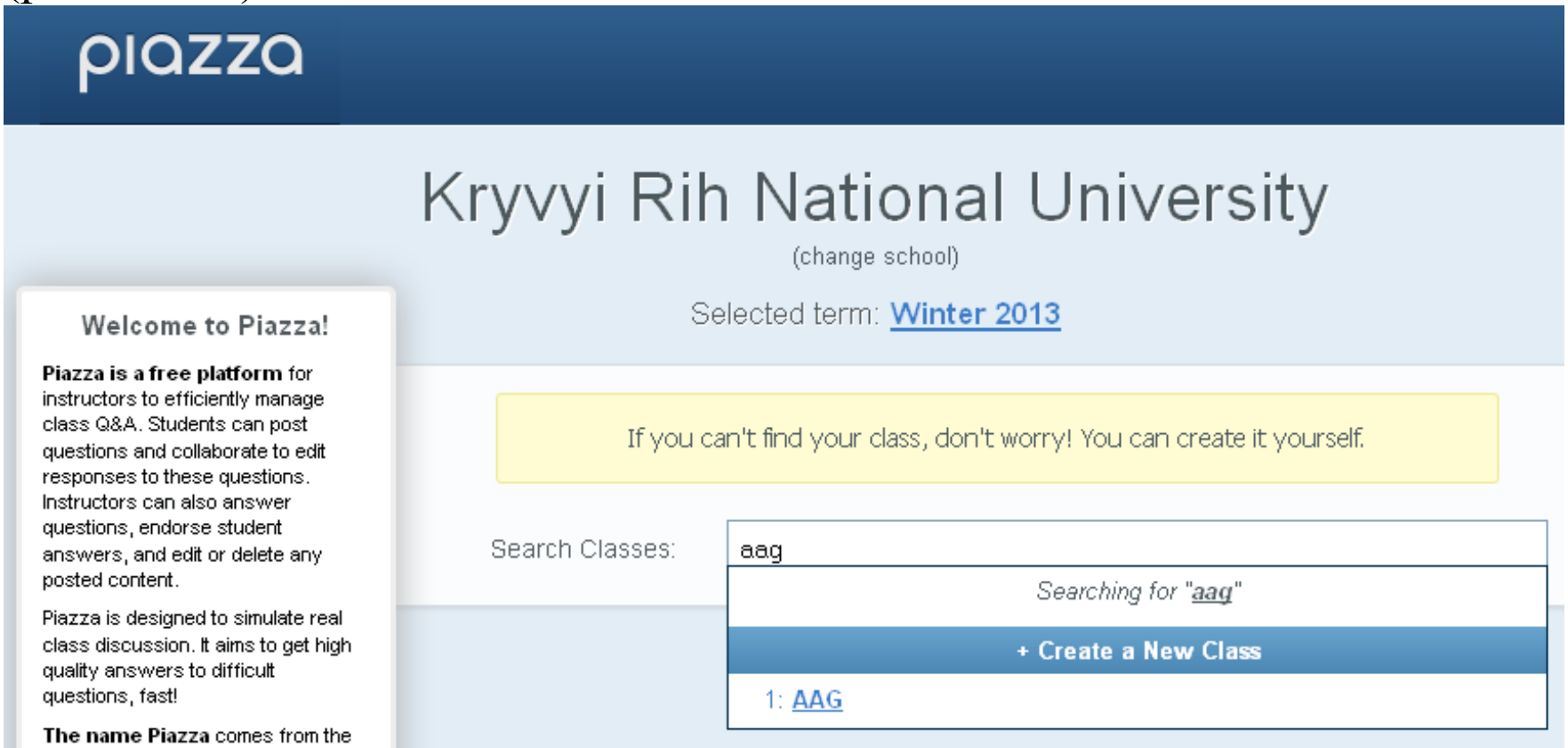

Рис. 3.16 б) Створення або вибір з існуючих класу в Piazza

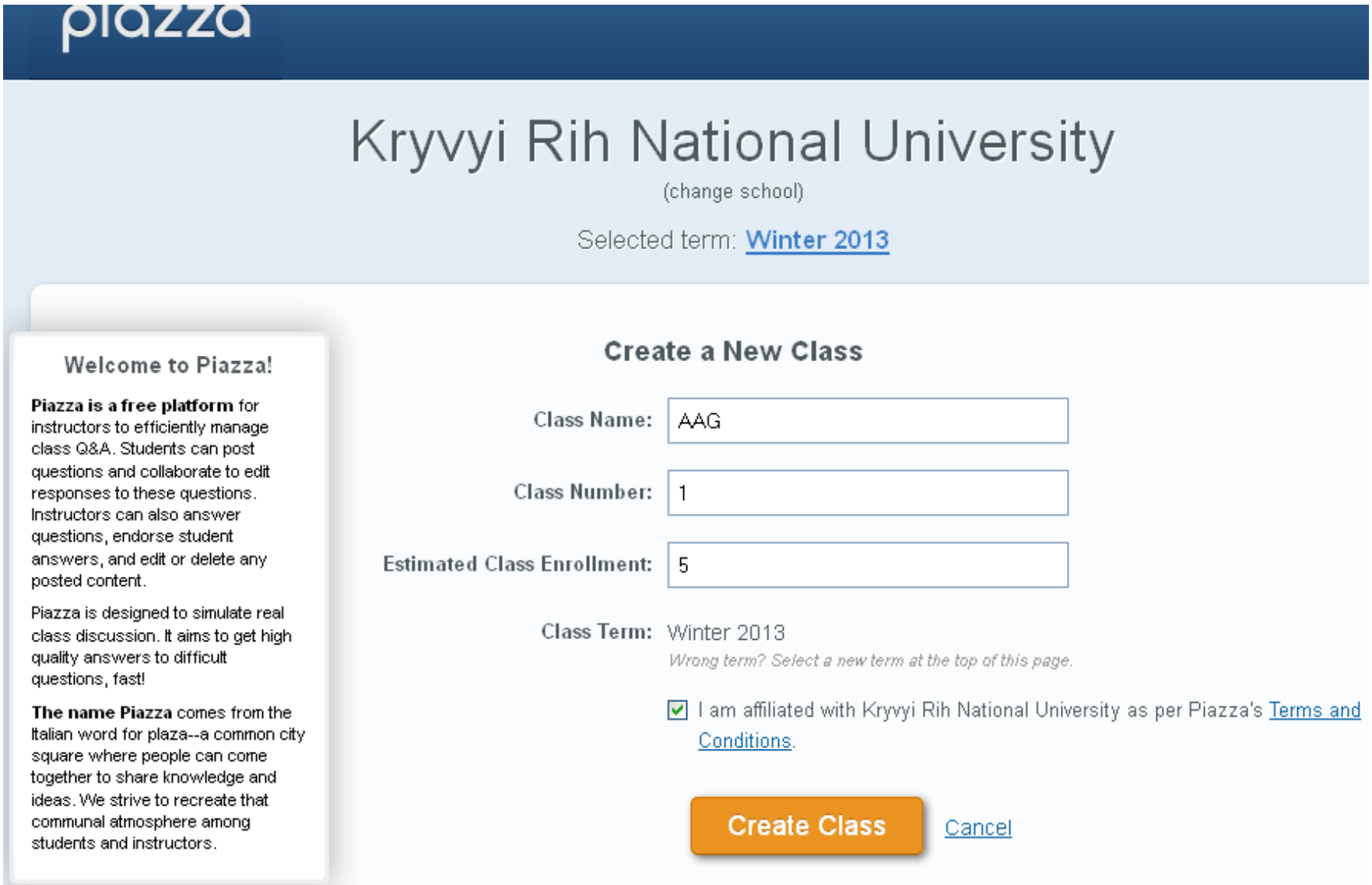

Рис. 3.16 в) Внесення даних про новий клас в Piazza

На сторінці ресурсів викладач може зробити посилання або завантажити файли з домашнім завданням для студентів, прикладом розв’язання домашнього завдання, лекційні матеріали та додаткові матеріали (рис. 3.17).

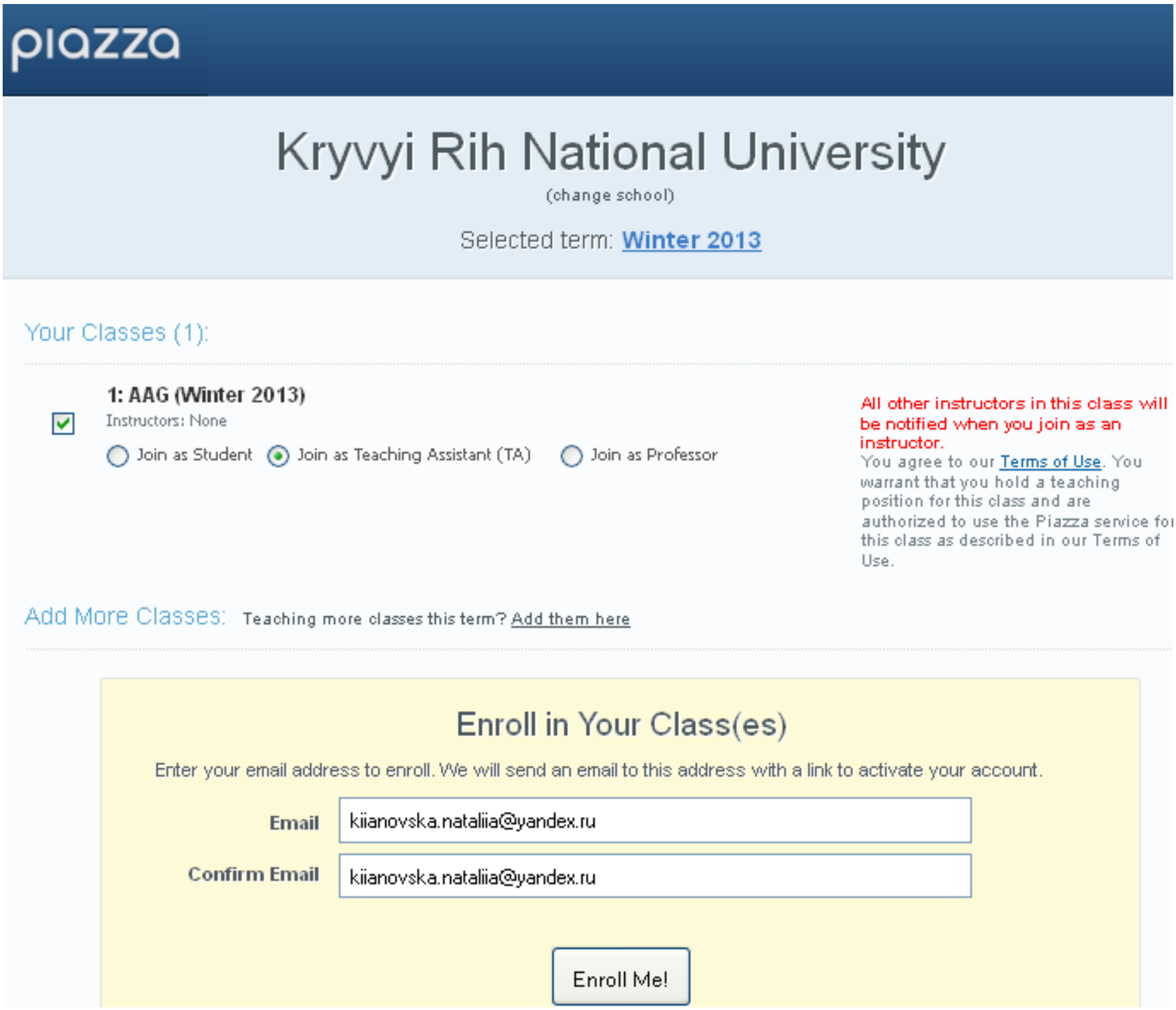

Рис. 3.16 г) Визначення ролі користувача в класі Piazza

Таким чином, використання платформи Piazza надає широкі можливості як викладачам так і студентам, активізуючи пізнавальну діяльність студента при вивченні курсу «Вища математика». Внутрішня мотивація студентів із використанням платформи Piazza зростає за рахунок організації нової, незвичної форми спілкування, при якій студенти почувають підтримку колективу та викладача в організації навчання.

Досягнення якісно нового рівня у підготовці фахівців з вищою освітою неможливе без забезпечення розвитку вищої школи на основі нових прогресивних концепцій, запровадження сучасних педагогічних та інформаційних технологій, науково-методичних розробок, відходу від засад авторитарної педагогіки і застарілих технологій навчання. В дослідженнях Ю. В. Триуса [317] показано, що серед педагогічних інновацій, використання яких може забезпечити підвищення якості вищої математичної освіти, сприяти пізнавальній активності студентів і набуттю ними комунікативних навичок й умінь, тобто вмінь працювати в різноманітних групах, виконуючи різні соціальні ролі (лідера, виконавця, посередника та ін.), формування вмінь самостійно конструювати свої знання та орієнтуватися в інформаційному просторі, є продуктивне навчання, навчання в співпраці, метод проектів [318]. Аналіз методів навчання вищої математики студентів інженерних спеціальностей розглянуто у додатку Н.

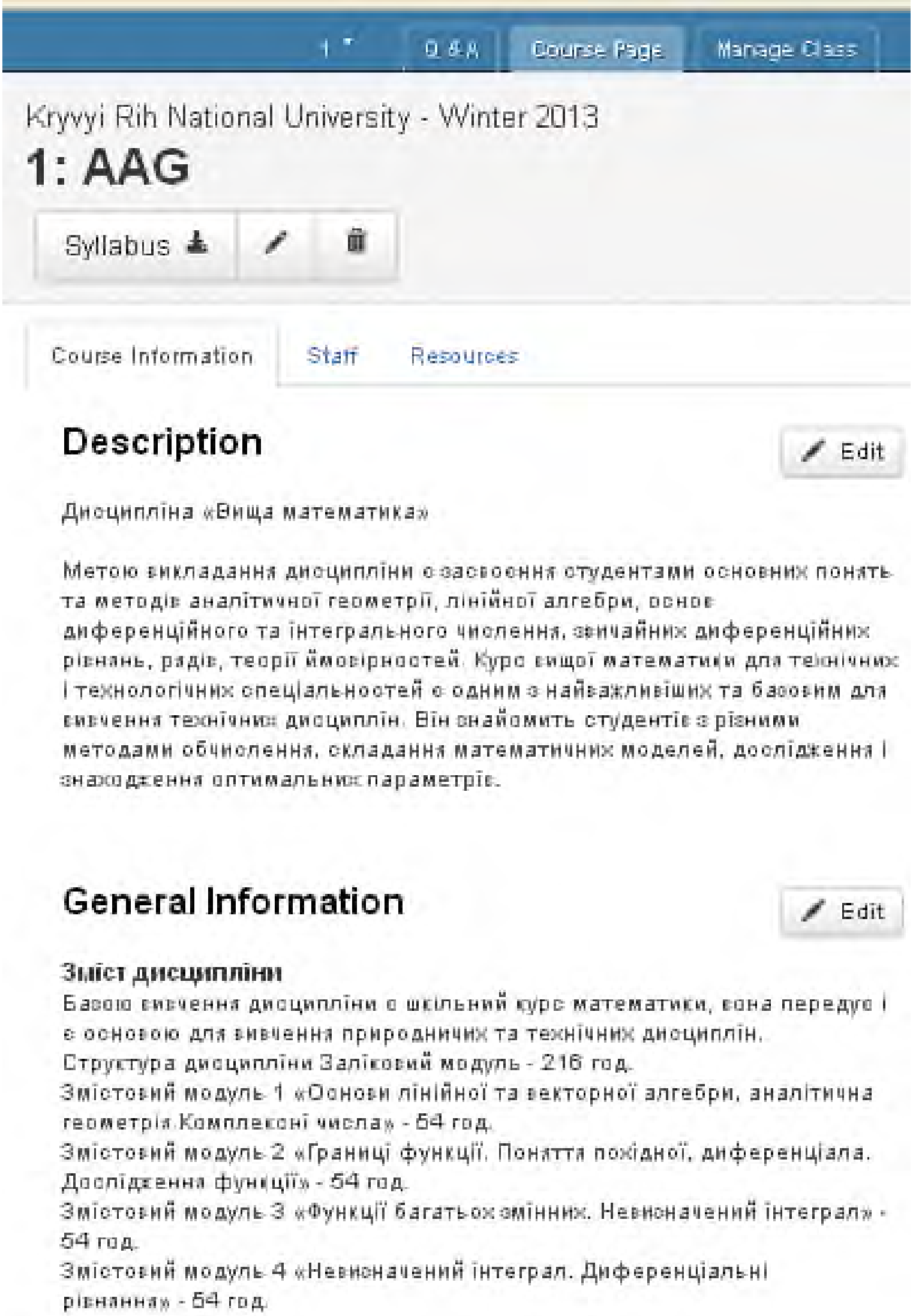

Рис. 3.16 д) Дані про навчальний курс в Piazza

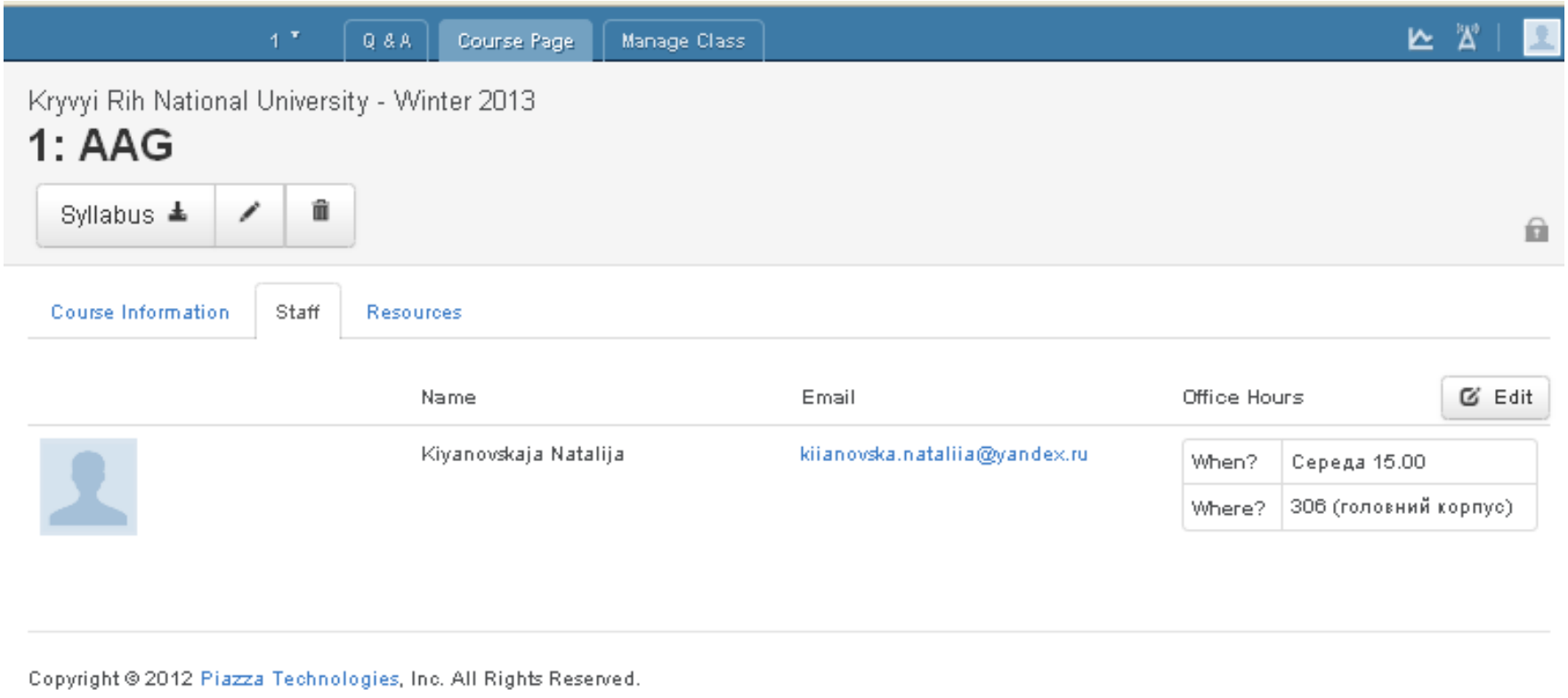

Рис. 3.16 є) Дані про графік аудиторних консультацій викладача

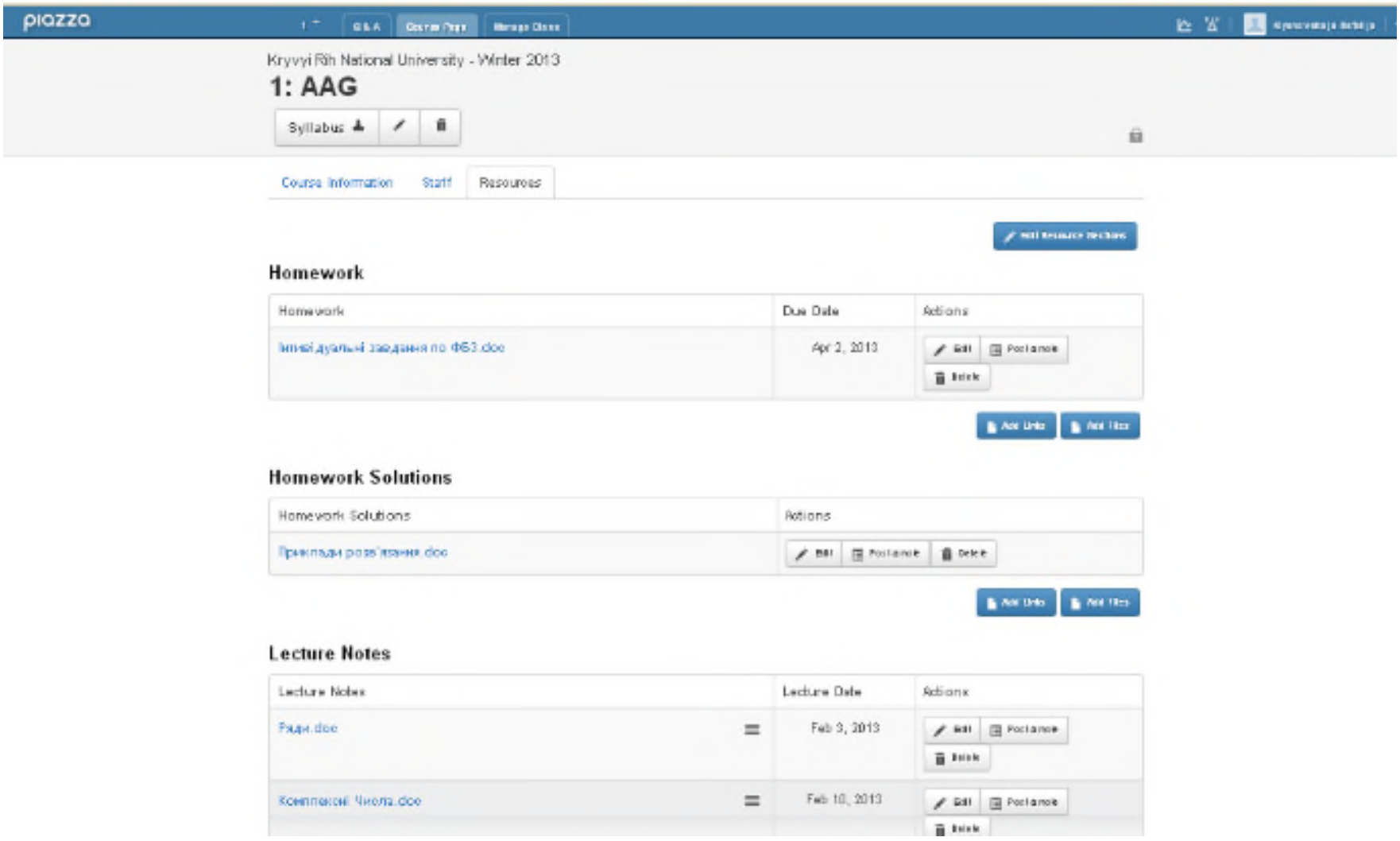

Рис. 3.17. Навчальні матеріали з курсу

Якщо кожна із зазначених вище інноваційних педагогічних технологій навчання, інтегруючись з ІКТ, займе своє місце в навчально-виховному процесі ВНЗ, поступово витісняючи методи і форми пасивного навчання, то згодом вдасться виробити досить ефективні підходи до організації навчального процесу у вищих навчальних закладах [318].

## 3.3 Інформаційно-комунікаційна компетентність викладачів вищої математики технічних вищих навчальних закладів в Україні

Педагоги у всіх країнах світу дуже добре усвідомлюють переваги, що надає методично обґрунтоване використання сучасних інформаційних і комунікаційних технологій у сфері освіти. Використання ІКТ допомагає вирішувати проблеми в тих галузях, де суттєве значення мають знання і комунікація. Сюди входять: вдосконалення процесів учіння / навчання, підвищення рівня навчальних досягнень студентів та їх навчальної мотивації, покращення взаємодії викладачів зі студентами, спілкування в мережі і виконання спільних проектів, вдосконалення організації та управління освітою та навчанням [309]. У зв'язку з цим перед викладачами постає завдання бути обізнаними в останніх досягненнях комп'ютерно орієнтованих технологій, розвивати свою ІКТ-компетентність.

Базуючись на дослідженнях О. М. Спіріна [305] та О. В. Овчарук [243], під *ІКТ-компетентністю* будемо розуміти підтверджені знання, вміння, ставлення та здатність особистості автономно і відповідально використовувати на практиці ІКТ для задоволення власних індивідуальних потреб і розв'язування суспільно значущих, зокрема професійних, задач у певній предметній галузі або виді діяльності.

Інформаційно-комунікаційна компетентність включає свідоме та критичне застосування технологій інформаційного суспільства для роботи, навчання, відпочинку та спілкування. Вона побудована на застосуванні базових інформаційно-комунікаційних навичок: використання засобів ІКТ для доступу, накопичення, вироблення, подання та обміну даними і відомостями та для спілкування, участі в спільнотах через мережу Інтернет [245, 46].

Основні знання, вміння та ставлення, що відносяться до цієї компетентності [245, 47]:

– інформаційно-комунікаційна компетентність вимагає свідомого *розуміння* та *знання* природи, ролі та можливостей технологій інформаційного суспільства в особистісному та соціальному житті, навчанні та роботі. Це включає використання комп'ютерних технологій (як, наприклад, текстових редакторів, електронних таблиць, баз даних, масивів даних локального та хмарного зберігання), розуміння можливостей та потенціальних ризиків Інтернету та спілкування через електронні медіа для роботи, навчання, відпочинку, обміну даними і відомостями та колаборативного мережного спілкування, навчання та досліджень;

– особистості повинні також *усвідомлювати*, як технології інформаційного суспільства можуть підтримувати творчість та

інноваційність, бути обізнаними про валідність та відповідність даних і відомостей, що на етичних та правових принципах є доступними та залучають до їх використання;

– *уміння* передбачають здатність знаходити, збирати та опрацьовувати дані, відомості і повідомлення та використовувати її систематично та критично, відповідно до реального та віртуального середовища. Особистості повинні володіти вмінням використовувати засоби для розробки, подання та усвідомлення комплексу даних та здатністю до доступу, пошуку та використання сервісів мережі Інтернет;

– також особистості повинні бути *здатними* використовувати засоби ІКТ для підтримки критичного мислення та відповідного ставлення до доступних даних і відомостей та відповідально використовувати медіа. Ця компетентність передбачає здатність входження до соціальних, культурних, професійних спільнот та мереж.

Особистості також повинні бути здатними використовувати ІКТ для підтримки не лише критичного мислення, а й творчості та інновацій [245, 47-48]:

1) ІКТ-бачення: розуміння та усвідомлення ролі та значення ІКТ для роботи та навчання впродовж життя;

2) ІКТ-культура: спосіб розуміння, конструювання, світоглядного бачення цифрових технологій для життя та діяльності в інформаційному суспільстві;

3) ІКТ-знання: фактичні та теоретичні знання, що відображають галузь ІКТ для навчання та практичної діяльності;

4) ІКТ-практика: практика застосування знань, умінь, навичок у галузі ІКТ для особистих та суспільних професійних та навчальних цілей;

5) ІКТ-удосконалення: здатність удосконалювати, розвивати, генерувати нове у сфері ІКТ та засобами ІКТ для навчання, професійної діяльності, особистого розвитку;

6) ІКТ-громадянськість: підтверджена якість особистості демонструвати свідоме ставлення через дію, пов'язану із застосуванням ІКТ для відповідальної соціальної взаємодії та поведінки.

Доцільно, щоб ці характеристики були максимально можливо відображені наскрізно і на всіх рівнях ІКТ-компетентностей у процесі їх набуття.

Існує безліч причин, які заважають студентам і викладачам в повній мірі використовувати можливості, що з'являються із використанням ІКТ. Це і брак коштів на закупівлю обладнання, і обмежений доступ в Інтернет, і відсутність цифрових освітніх ресурсів на рідній мові. Але головна причина в тому, що викладачі не завжди знають, як ефективно використовувати ІКТ [309].

Під *ІКТ-компетентністю викладача* будемо розуміти професійно значуще особистісне утворення – здатність педагогічно виважено та методично обґрунтовано виокремлювати, добирати, досліджувати, проектувати та використовувати ІКТ з метою усебічного забезпечення процесу навчання.

В рекомендаціях ЮНЕСКО щодо інформаційно-комунікаційної компетентності викладачів виділено такі напрями [309, 11]:

– *розуміння ролі ІКТ в освіті:* педагоги повинні бути знайомі з освітньою політикою і вміти пояснити на професійній мові, чому їх педагогічні практики відповідають освітній політиці і як її реалізують;

– *навчальна програма і оцінювання:* педагоги повинні відмінно знати освітні стандарти та вимоги щодо оцінювання навчальних досягнень зі свого навчального предмету. Крім того, педагоги повинні бути здатні включити використання засобів ІКТ у свою навчальну програму;

– *педагогічні практики*: педагоги повинні знати передовий досвід теорії та методики використання ІКТ у навчальній роботі та для подання навчального матеріалу;

– *технічні і програмні засоби ІКТ:* педагоги повинні знати базові прийоми роботи з технічними і програмними засобами; програмні засоби, що підвищують продуктивність роботи; Web-браузери; комунікаційні програмні засоби; засоби презентаційної графіки; програми для вирішення завдань управління;

– *організація і управління освітнім процесом:* педагоги повинні вміти використовувати засоби ІКТ для роботи з усім класом, у малих групах, а також для індивідуальної роботи. Вони повинні надавати всім студентам рівний доступ до використання ІКТ;

– *професійний розвиток:* педагоги повинні володіти навичками роботи з ІКТ і знати Web-ресурси, щоб отримувати додаткові навчально-методичні матеріали, необхідні для їх професійного розвитку.

Рівень розвитку інформаційно-комунікаційної складової педагогічної компетентності викладача визначається [152]:

– мотивацією: прагненням до якісного виконання інформаційної діяльності, ефективного використання ІКТ, самовдосконалення тощо;

– інформаційно-науковими знаннями з оперуванням поняттями, актами, властивостями, закономірностями, методами, алгоритмами та ін.;

– уміннями і навичками досвіду роботи з інформаційними джерелами;

– інформаційним світоглядом, інноваційним мисленням та інтуїцією, мобільністю, ціннісними орієнтаціями щодо доцільності використання ІКТ.

Спираючись на запропоноване О. М. Спіріним у дослідженні [245]

опис рівнів ІКТ-компетентностей, пропонуємо три *рівні ІКТ-компетентності викладача*:

*I рівень, базовий.* Систематично використовувати стандартні засоби ІКТ для підтримки навчання. Самостійно добирати засоби ІКТ для реалізації цілей навчання. Правильно добирати і використовувати ІКТ для розв'язування основних навчальних задач.

*II рівень, поглиблений.* Проводити проектування процесу навчання із використанням ІКТ. Створювати предметно орієнтоване навчальне середовище, сприяти розвитку персональних навчальних середовищ студентів. Застосувати ІКТ для комбінування форм організації, методів та засобів навчання. Уміти розв'язувати професійні задачі підвищеної складності з використанням ІКТ, адаптувати засоби ІКТ для розв'язування основних професійних задач, зокрема бути здатним проектувати, конструювати і вносити інновації до елементів наявних ІКТ навчання.

*III рівень, дослідницький.* Засвоїти та демонструвати повне володіння методикою використання ІКТ у предметній галузі. Досліджувати, добирати та проектувати засоби ІКТ організації навчального процесу. Зробити оригінальний вклад у розвиток теорії та методики використання ІКТ у процесі навчання, розробляти інноваційні ІКТ навчання.

Запропоновані рівні було визначено, виходячи з вимог державної цільової програми впровадження у навчально-виховний процес загальноосвітніх навчальних закладів ІКТ на період до 2015 року «Сто відсотків».

У визначенні компетентностей викладача вищої математики доцільно скористатися результатами дослідження С. А. Ракова, де вказується на необхідність формування [266; 86]:

– процедурної компетентності як умінь розв'язувати типові математичні та інформатичні задачі;

– логічної компетентності як володіння дедуктивним методом доведення та спростування тверджень;

– технологічної компетентності як умінь застосовувати у професійній діяльності засобів ІКТ;

– дослідницької компетентності як володіння методами дослідження соціально та індивідуально значущих задач математичними методами;

– методологічної компетентності як умінь оцінювати доцільність використання математичних методів для розв'язування індивідуально і суспільно значущих задач.

Викладач з вищої математики технічного ВНЗ має володіти такими **компетентностями** [152; 315]:

– *інформаційна* – здатність викладача до проведення критичного

аналізу джерел даних, пошуку необхідних ресурсів, синтезу, узагальненню та структуруванню опрацьованих відомостей;

– *технічна* – здатність та готовність викладача до ефективного використання та опанування апаратних та програмних засобів ІКТ;

– *технологічна* – здатність та готовність викладача до інформаційно-технологічної діяльності, а саме: постановка цілей створення електронних освітніх ресурсів, використанню існуючої або розробки нової технології для створення електронного освітнього ресурсу, тестуванню створеного продукту на відповідність до певних вимог тощо;

– *педагогічна* – здатність та готовність викладача до педагогічного проектування, змістового наповнення та використання електронних освітніх ресурсів у власній професійній діяльності;

– *мережна та телекомунікаційна* – здатність до опанування основними принципами побудови і використання локальних мереж та глобальної мережі Інтернет;

– *дослідницька* – здатність проводити дослідження доступними засобами ІКТ;

– *в питаннях інформаційної безпеки* – здатність запобігти можливим інформаційним атакам у комп’ютерних системах, володіти знаннями з принципів захисту даних, вміти проводити апаратні та програмні методи захисту даних.

Відповідно до означених компетентностей викладачі з вищої математики технічного ВНЗ повинні знати про: засоби ІКТ для навчання вищої математики у інженерній освіті; можливості та призначення засобів ІКТ; техніку безпечного користування засобами ІКТ; структуру мережі Інтернет та її значення для освіти.

У викладачів вищої математики технічного ВНЗ повинні бути сформовані такі основні **вміння** використання ІКТ:

– працювати з електронною поштою;

– працювати з Web-браузерами;

– використовувати математичні пакети для проведення обчислень та моделювання;

– застосовувати мережні засоби для підтримки спілкування;

– працювати з науковими текстовими процесорами (зокрема LaTeX, Kile);

– використовувати системи відображення документів;

– користуватися програмами автоматизації роботи з даними;

– працювати з периферійним комп’ютерним обладнанням (принтер, сканер, модем, Web-камера тощо);

– проектувати та створювати нові засоби навчання математики.

Ураховуючи, що ефективною та дієвою формою розвитку ІКТ-

компетентності викладачів вищої математики технічних ВНЗ є підвищення кваліфікації, було спроектовано навчальний спецкурс «Інформаційно-комунікаційні технології навчання вищої математики студентів інженерних спеціальностей», що розрахований на використання його для підвищення кваліфікацій викладачів математичних дисциплін технічних ВНЗ та підготовки магістрів математики.

*Метою спецкурсу* є розвиток інформаційно-комунікаційної компетентності викладачів вищої математики технічних ВНЗ, зокрема: навчити створювати текстові документи, таблиці, малюнки, діаграми, презентації, використовуючи Інтернет-технології, локальні мережі, бази даних; здійснювати анкетування, діагностування, тестування, пошук даних; розробляти власні електронні освітні ресурси та поєднувати готові.

Серед *завдань* спецкурсу виділено:

– сприяти формуванню стійкої педагогічної позиції щодо необхідності використання засобів ІКТ у власній професійній діяльності;

– систематизувати технічні, педагогічні, інформаційні, технологічні, мережні та телекомунікаційні, дослідницькі знання та вміння викладачів, знання викладачів з питань інформаційної безпеки для створення та використання електронних освітніх ресурсів в професійній діяльності;

– удосконалити навички роботи з апаратними та програмними засобами ІКТ;

– розвинути та удосконалити інформаційно-технологічні вміння щодо створення та використання електронних освітніх ресурсів у професійній діяльності.

У контексті даного спецкурсу слухачі повинні знати про можливості створення електронних освітніх ресурсів та використання створених або готових ресурсів в організаційно-методичній діяльності або навчальному процесі.

Очікуваним результатом спецкурсу є: ефективне використання інноваційних методів, форм організації навчання та засобів ІКТ у навчальному процесі з вищої математики; впровадження нових прийомів роботи у навчальний процес; розвиток ІКТ-компетентності кожного викладача.

Курс призначений для підвищення ефективності освітнього процесу та навчальних досягнень з вищої математики у вищих технічних закладах. Зміст курсу може бути адаптований для будь-якої категорії слухачів очної або очно-дистанційної форми організації навчання. На вивчення спецкурсу пропонується виділити 54 годин / 1,5 кредити ECTS.

Робочу програму спецкурсу «Інформаційно-комунікаційні технології навчання вищої математики студентів інженерних спеціальностей» подано у додатку П.

Рекомендації щодо вибору методів, форм та засобів навчання розділів вищої математики студентів інженерних спеціальностей у технічних ВНЗ розглянуто у додатку Р.

Використання ІКТ у навчанні вищої математики надає такі переваги:

– гнучкість і зручність у поданні навчального матеріалу;

– мобільність, доступність до навчального матеріалу в довільний час і будь-якому місці, і доступність до навчання взагалі;

– навчання є диференційованим, надається можливість організації індивідуальної роботи студента з викладачем за власним графіком навчання;

– адаптованість процесу навчання до студентів, до їх потреб і можливостей;

– скорочення розриву між теорією і практикою за допомогою програмного забезпечення, що надає можливість у вищій математиці моделювати і вирішувати реальні проблеми та приклади;

– забезпечення безперервного процесу навчання: студенти отримують своєчасні відомості про навчальний курс, про свої навчальні досягнення в ході курсу;

– прогресивна система оцінювання навчальних досягнень студентів поліпшує академічні результати студентів, а також сприяє мотивації їх навчання;

– розвиток таких рис особистості студента як самостійність, відповідальність, мобільність;

– можливість дистанційної роботи зі студентами як індивідуально, так і колективно;

– можливість співпраці викладача та студента, викладач спонукає студента до самостійної роботи, («тим, хто веде за руку», оскільки педагог в Давній Греції παιδαγωγός – «той, хто веде дитину»), а не управляє процесом навчання.

Однак існують, звичайно, і недоліки при використанні ІКТ у процесі навчання вищої математики. До їх числа можна віднести:

– технічні засоби не завжди працюють безперебійно, тому можливі збої в роботі ІКТ, що може привести до зниження ефективності навчання;

– необхідність інвестицій у розробку та впровадження проектів, подібних до OCW;

– необхідність наявності спеціальних знань у професорсько-викладацького складу ВНЗ та час на опанування нових технологій;

– необхідність у забезпеченні комп'ютерами та вільним доступом до Інтернет всіх учасників навчального процесу: викладачів та студентів;

– витрачання часу на розробку та впровадження ІКТ у процес навчання.

Та не зважаючи на наявні недоліки, переваги, які отримують учасники процесу навчання, набагато вагоміші і тому питанню використання ІКТ у процесі навчання вищої математики мають присвячуватися семінари, методичні засідання кафедр, конференції тощо з метою просування цієї політики навчання, що набула широкого розмаху в усьому світі.

Хоча переваги, що надаються при використанні ІКТ у процесі навчання такі значимі, ІКТ ні в якому разі не можуть замінити викладача. Роль викладача залишається провідною і лише від нього залежить якість отриманих знань студентів: від уміння зацікавити та організувати студентів, від акторської майстерності проводити заняття, якості підібраних матеріалів. Тому, як вже зазначалось раніше, доцільно інтегрувати технології традиційного навчання з технологіями навчання із використанням ІКТ.

## Висновки до розділу 3

1. Використання ІКТ у навчанні вищої математики сприяє якісному та своєчасному поданню навчальних матеріалів та відомостей про процес навчання. Використання різноманітних онлайн ресурсів освітніх мереж підвищує рейтинг ВНЗ, що може відігравати вагому роль для абітурієнтів при виборі місця навчання.

2. Використання засобів ІКТ (зокрема, онлайн) у математичній підготовці майбутніх інженерів в Україні сприятиме поглибленню розуміння матеріалу з фундаментальних основ інженерії, активізації навчальної діяльності з вищої математики, надаючи: *процесу навчання* вищої математики – властивостей мобільності, безперервності та адаптивності, *викладачам* – нових можливостей із комбінування форм організації та методів навчання вищої математики, *студентам* – вільний доступ до навчальних матеріалів, мобільну навчальну підтримку та варіативність процесу навчання вищої математики.

3. Використання ІКТ у математичній підготовці майбутніх інженерів технічних ВНЗ надає змогу підготувати фахівця з інженерної спеціальності, який володіє ІКТ, що мають професійну спрямованість, здатний до безперервного навчання, спроможній вдосконалювати себе як фахівця.

4. На сьогоднішній день існує велика кількість розробок та впроваджень різних засобів ІКТ у навчання вищої математики студентів інженерних спеціальностей. Але всі розглянуті впровадження не реалізують системний підхід до використання ІКТ, як це було розглянуто на прикладі системи вищої освіти США, де в процес навчання вищої математики студентів інженерних спеціальностей впроваджено математичні пакети, системи підтримки навчання, Інтернет-технології, мультимедійні програмні засоби, офісне програмне забезпечення, зі спрямованістю саме на інженерні спеціальності. Натомість, в Україні існує тенденція допрацьовувати та адаптувати засоби, що не були заздалегідь призначені для навчання студентів інженерних спеціальностей. Крім професійної спрямованості , розглянуті засоби мають давню історію розвитку та впровадження, що обумовлює їх надійність при використанні в процесі навчання.

5. Запропонований спецкурс «Інформаційно-комунікаційні технології навчання вищої математики студентів інженерних спеціальностей» враховує особливості досвіду США використання ІКТ у навчанні вищої математики студентів інженерних спеціальностей та розрахований на використання його для підвищення кваліфікацій викладачів математичних дисциплін технічних ВНЗ та підготовки магістрів математики. Вивчення спецкурсу сприятиме розвитку інформаційно-комунікаційної компетентності викладачів вищої математики технічних ВНЗ та магістрів математики.

## ВИСНОВКИ

1. Незважаючи на недержавну систему акредитації, відсутність державних галузевих стандартів та традиційне різноманіття пропонованих математичних курсів (як обов'язкових, так і факультативних), навчання вищої математики майбутніх інженерів у США здійснюється за схожими навчальними програмами. Системи підготовки інженерів у ВНЗ США та у ВНЗ України мають такі спільні риси: високий рівень математизації та комп'ютеризації загальноінженерних та спеціальних дисциплін; навчання фундаментальних дисциплін, зокрема вищої математики, відбувається переважно на молодших курсах, загальнопрофесійних дисциплін – на середніх курсах та спеціальних професійних – на старших; зміст математичної підготовки є професійно орієнтованим та диференційованим за рівнями початкової підготовки студентів; у навчанні вищої математики широко використовуються засоби ІКТ.

2. Вивчення історико-педагогічної джерельної бази надало можливість відобразити еволюцію та виокремити шість етапів розвитку теорії та методики використання ІКТ у навчанні вищої математики студентів інженерних спеціальностей у США: 1) 1965–1973 рр. – пов'язаний із появою достатньої кількості комп'ютерних засобів різного рівня, оснащених мовами високого рівня та специфікою апаратного забезпечення ІКТ (використання мейнфреймів з обмеженим мережним доступом); 2) 1973–1981 рр. – пов'язаний із поширенням в університетах США мережної операційної системи UNIX, використанням міні- та мікрокомп'ютерних систем; 3) 1981–1989 рр. – пов'язаний із поширенням персональних комп'ютерів; 4) 1989–1997 рр. – пов'язаний із створенням World Wide Web та використанням технологій Web 1.0; 5) 1997–2003 рр. – пов'язаний із появою та розробкою систем управління навчанням; 6) з 2003 р. по теперішній час – пов'язаний із перенесенням у Web-середовище засобів підтримки математичної діяльності та становленням і розвитком хмарних технологій навчання.

Поява нового типу апаратних чи програмних засобів ІКТ впливає на процес організації навчання вищої математики і на сучасному етапі створює умови для реалізації Web-орієнтованого навчання вищої математики, що стає провідним напрямом як у США, так і в Україні.

3. На сучасному етапі розвитку вищої інженерної школи США провідними засобами підтримки навчальної діяльності з вищої математики майбутніх інженерів є Web-орієнтовані ІКТ загального призначення (системи управління навчанням, системи розміщення відкритих навчальних матеріалів, засоби комунікації та спільної роботи) та спеціального призначення (системи комп'ютерної математики, лекційні демонстрації, динамічні навчальні матеріали). Виокремлено такі особливості сучасного етапу розвитку теорії та методики використання ІКТ у навчанні вищої математики студентів інженерних спеціальностей у США: перенесення математичної діяльності викладачів та студентів у мережне середовище; застосування засобів хмарних технологій для підтримки навчальної діяльності; становлення Web-орієнтованих методичних систем навчання вищої математики; розвиток масових відкритих дистанційних курсів.

4. Поява нових ІКТ призводить до змін у теорії та методиці навчання вищої математики, що обумовлює виникнення нових цілей, засобів, форм, методів організації процесу навчання та доповненням змісту навчання. Динаміка розвитку теорії й методики використання ІКТ у навчанні вищої математики студентів інженерних спеціальностей у ВНЗ США мала діалектичний характер: на кожному етапі її розвитку зміна засобів ІКТ супроводжувалась набуттям нового способу доступу до навчальних ресурсів.

5. Обґрунтовано фактори, що зумовлюють доцільність поширення досвіду США у вітчизняному освітньому просторі: перші позиції у рейтингу країн за наявними розробленими та впровадженими засобами ІКТ у процес навчання; демократично-гуманістичні орієнтири освітньої політики щодо реалізації процесу інформатизації вищої освіти, що надають можливість забезпечити індивідуалізацію, диференціацію і персоналізацію процесу навчання; професійна спрямованість заходів щодо інформатизації системи технічної освіти в ВНЗ; тривалий етап розвитку та впровадження ІКТ, що обумовлює їх надійність при використанні в процесі навчання; реалізація системного підходу до використання ІКТ у процесі навчання вищої математики студентів інженерних спеціальностей.

Визначено передумови поширення позитивного досвіду США у вітчизняні освітні практики: актуалізація процесів підвищення якості освіти на всіх рівнях; потреби подальшого розвитку теорії й практики використання ІКТ у навчанні вищої математики студентів інженерних спеціальностей; стрімкий розвиток системи інженерної освіти в Україні, що супроводжується зростанням вимог до гнучкості, комфортності, доступності її організації.

На основі обґрунтованої структурно-функціональної схеми використання ІКТ у навчанні вищої математики розроблено рекомендації для викладачів вищої математики технічних ВНЗ України, що включають опис структури комп'ютерно орієнтованого навчально-методичного комплексу викладача вищої математики та передбачають включення інноваційних технологій навчання, впровадження комп'ютерно-орієнтованих форм організації навчання вищої математики, використання викладачами та майбутніми інженерами технічних ВНЗ у ході математичної підготовки доступних засобів Інтернет для підтримки математичної діяльності. Розроблені рекомендації щодо застосування досвіду США з використання ІКТ у навчанні вищої математики студентів інженерних спеціальностей в Україні узагальнено у спецкурсі «ІКТ навчання вищої математики студентів інженерних спеціальностей» для підвищення кваліфікації викладачів математичних дисциплін.

Виконане дослідження не вичерпує всіх аспектів досліджуваної проблеми. Продовження наукового пошуку за даною проблематикою доцільно у таких напрямах:

– удосконалення змісту, форм і методів підготовки студентів інженерних спеціальностей до використання ІКТ у їх професійній діяльності;

– вивчення, розробка, обґрунтування та впровадження технологій, спрямованих на розвиток ІКТ-компетентності майбутніх інженерів.

## СПИСОК ВИКОРИСТАНИХ ДЖЕРЕЛ

# ДОДАТКИ

*А Порівняльна таблиця розподілу кредитів на вивчення вищої математики*

*Таблиця* А.1

**Навчальні курси, що вивчають студенти інженерних спеціальностей у США**

| № з/п | Назва спеціальності | Навчальні курси у кредитах | | | | |
|---|---|---|---|---|---|---|
| | | Елементарне числення | Числення однієї змінної | Числення багатьох змінних | Calculus III (диференціальні рівняння) | Лінійна алгебра |
| *Massachusetts Institute of Technology* [79] | | | | | | |
| 1. | Aeronautics and astronautics | | 12 | 12 | 12 | |
| 2. | Civil and Environmental Engineering | | 12 | 12 | 12 | |
| 3. | Mechanical Engineering | | 12 | 12 | 12 | |
| 4. | Materials Science and Engineering | | 12 | 12 | 12 | |
| 5. | Nuclear Science and Engineering | | 12 | 12 | 12 | |
| *Purdue University at West Lafayette* [106] | | | | | | |
| 6. | Aeronautical and Astronautical Engineering | | 4 | 8 | 6 | 3 |
| 7. | Computational Science and Engineering | | 4 | 8 | 3 | 3 |
| 8. | Computer Graphics and Visualization | | 4 | 8 | | 3 |
| 9. | Database and Information Systems | | 4 | 8 | | 3 |
| 10. | Foundations of Computer Science | | 4 | 8 | | 3 |
| 11. | Machine Intelligence | | 4 | 8 | | 3 |
| 12. | Software Engineering | | 4 | 8 | | 3 |
| 13. | Systems Programming | | 4 | 8 | | 3 |
| *New York Institute of Technology* [129] | | | | | | |
| 14. | Computer Science | | 4 | 4 | | 3 |
| 15. | Computer Engineering | | 4 | 4 | 7 | 3 |

| № з/п | Назва спеціальності | Навчальні курси у кредитах | | | | |
|---|---|---|---|---|---|---|
| | | Елементарне числення | Числення однієї змінної | Числення багатьох змінних | Calculus III (диференціальні рівняння) | Лінійна алгебра |
| 16. | Mechanical Engineering | | 4 | 4 | 7 | |
| *West Virginia University* [138] | | | | | | |
| 17. | Electrical Engineering | | 4 | 8 | 3 | 3 |
| 18. | Aerospace Engineering | | 4 | 8 | 4 | |
| 19. | Mechanical Engineering | | 4 | 8 | 4 | |
| *California State Polytechnic University, Pomona* [19] | | | | | | |
| 20. | Mechanical Engineering | | 4 | 10 | 4 | 4 |
| 21. | Industrial and Manufacturing Engineering | | 4 | 10 | 4 | 4 |
| 22. | Electrical and Computer Engineering | | 4 | 10 | 4 | 4 |
| Georgia Southern University [48] | | | | | | |
| 23. | Civil Engineering Technology | 4 | 4 | 4 | | |
| 24. | Electrical Engineering Technology | 4 | 4 | 4 | | |
| 25. | Mechanical Engineering Technology | 4 | 4 | 4 | | |
| *University of Pittsburgh at Johnstown* [104] | | | | | | |
| 26. | Civil Engineering Technology | | 4 | 4 | 8 | |
| 27. | Computer Engineering Technology | | 4 | 4 | 8 | |
| 28. | Mechanical Engineering Technology | | 4 | 4 | 8 | |
| *Howard University* [23] | | | | | | |
| 29. | Electrical Engineering | | 4 | 4 | 8 | 3 |
| 30. | Computer Engineering | | 4 | 4 | 4 | 3 |
| 31. | Civil & Environmental Engineering | | 4 | 4 | 8 | |
| 32. | Mechanical Engineering | | 4 | 4 | 8 | |
| 33. | Chemical Engineering | | 4 | 4 | 8 | |
| *Illinois Institute of Technology* [61] | | | | | | |
| 34. | Computer Engineering | | 4 | 8 | 4 | |
| 35. | Electrical Engineering | | 4 | 8 | 4 | 3 |

| № з/п | Назва спеціальності | Навчальні курси у кредитах | | | | |
|---|---|---|---|---|---|---|
| | | Елементарне числення | Числення однієї змінної | Числення багатьох змінних | Calculus III (диференціальні рівняння) | Лінійна алгебра |
| 36. | Mechanical Engineering | | 4 | 8 | 4 | |
| 37. | Architectural Engineering | | 4 | 8 | 4 | |
| *Lawrence Technological University* [73] | | | | | | |
| 38. | Civil Engineering | | 4 | 4 | 7 | |
| 39. | Mechanical Engineering | | 4 | 4 | 7 | |
| 40. | Computer Engineering | | 4 | 4 | 7 | |
| *Michigan Technological University* [22] | | | | | | |
| 41. | Civil Engineering | | 4/5 | 4 | 6 | 2 |
| 42. | Mechanical Engineering | | 4/5 | 4 | 6/7 | 2/3 |
| 43. | Computer Engineering | | 4/5 | 4 | 5 | 2 |
| 44. | Electrical Engineering | | 4/5 | 4 | 6 | 2 |
| *Oregon Institute of Technology* [24] | | | | | | |
| 45. | Civil Engineering | | 8 | 4 | 4 | |
| 46. | Mechanical Engineering | | 8 | 4 | 4 | 3 |
| 47. | Computer Engineering | 8 | 8 | 4 | | |
| 48. | Electrical Engineering | | 8 | 4 | 4 | 3 |
| *Polytechnic Institute of New York University* [76] | | | | | | |
| 49. | Civil Engineering | | 6 | 4 | 2 | |
| 50. | Mechanical Engineering | | 4 | 8 | 2 | 2 |
| 51. | Computer Engineering | | 4 | 6 | 2 | 2 |
| 52. | Electrical Engineering | | 4 | 8 | 2 | 2 |

*Б Аналіз категоріального поля дослідження*

Навчання у вищий школі – це цілеспрямований спільний процес діяльності викладача і студентів, спрямований на досягнення дидактичних, розвивальних, виховних цілей та професійного розвитку і становлення студентської молоді. Для управління активністю студентів; демонстрації явищ і процесів, що вивчаються і взаємозв’язків між ними; для концентрації уваги на основному змісті матеріалу на заняттях доцільно використовувати різноманітні засоби навчання.

*Засоби навчання* – будь-які засоби, прилади, обладнання та устаткування, що використовуються для передачі відомостей в процесі навчання [169, 313–314].

Засоби навчання формують матеріальну та інформаційну складові навчального середовища, впливають на діяльність суб’єктів навчання і організацію процесу навчання, створюють умови для забезпечення можливості набуття знань та формування навчальних умінь з метою досягнення наперед визначених цілей навчання. З іншого боку, засоби навчання можуть використовуватися за різними формами організації та методами навчання, вибір яких зумовлений частинними методиками навчання.

На кожному етапі розвитку педагогічної науки адекватно розвиваються і засоби навчання, що акумулюють та відтворюють науково-технічні, психолого-педагогічні та методичні досягнення свого часу. Еволюція засобів навчання визначається потребами педагогічної практики, а їх розвиток спрямовується на задоволення цих потреб. Засоби навчання у процесі навчання виступають як засоби пізнання. Вплив наукових і технічних досягнень людства на зміст, структуру і організацію процесу навчання опосередковується і має матеріальний вираз у засобах навчання як знаряддях навчальної діяльності.

Технічний прогрес зумовив появу принципово нових засобів навчання, що сприяють формуванню навчальних середовищ на базі ІКТ. Суттєвою ознакою таких середовищ є те, що їх можна використовувати як частину внутрішнього навчального середовища ВНЗ (замкненого), так і як складову відкритого навчального середовища. Якщо в замкненому середовищі у процесі навчання можна не застосовувати засоби ІКТ, то у відкритому навчальному середовищі навчання відбувається тільки із впровадженням цих засобів у навчальний процес. Одні й ті ж самі засоби ІКТ можуть бути використані при навчанні різних навчальних дисциплін за умови відповідної зміни методик навчання.

З розвитком технологій у всіх галузях виробництва, зв’язку, спілкування освітяни також почали вивчати технологічні особливості педагогічного процесу та здійснювати його на основі технологічного підходу. Особливо актуальною є реалізація технологічного підходу для вітчизняної освіти на сучасному етапі її розвитку, що характеризується створенням Європейського простору вищої освіти, що відзначається високим рівнем технологізації, використанням освітньо-технологічних інновацій, переходом до «суспільства знань» [228].

Поняття «*технологія*» походить з грецької τέχνη – мистецтво, майстерність, ремесло + λόγος – знання, вчення, тобто «знання про майстерність», але на сьогодні дане поняття тлумачать по різному.

На думку авторів [222, 855], *технологія* – це сукупність виробничих операцій, методів, процесів у певній галузі виробництва, способів, що використовуються у певній справі; або це сукупність знань про них, послідовність окремих виробничих операцій у процесі виробництва; навчальний предмет, у якому викладаються ці знання, відомості.

Як зазначає В. Т. Бусел [154, 1448], *технологією* є сукупність способів обробки чи переробки матеріалів, даних, виготовлення виробів, проведення різних виробничих операцій, надання послуг тощо.

За Н. В. Морзе [234, 10], *технологія* – це сукупність методів, засобів і реалізації людьми конкретного складного процесу шляхом поділу його на систему послідовних взаємопов'язаних процедур і операцій, що виконуються більш або менш однозначно і мають на меті досягнення високої ефективності певного виду діяльності.

Поняття «технології в освіті» вперше було використане в науковій літературі та педагогічній періодиці (насамперед, у зарубіжній) наприкінці XIX – початку XX ст. [169, 906]. Розвиток науки й техніки забезпечив появу та використання в навчально-виховному процесі технічних засобів: кінопроектора, епідіаскопа, мікроскопа, фотообладнання тощо.

Так, С. А. Лисенко [228], вивчаючи питання технологічного підходу до навчально-виховної діяльності, зазначає, що зміна термінів відбувалася у такій послідовності: «технологія в освіті», «технологія виховання», «технологія навчання», «технологія освіти», «педагогічна технологія», що відображено й у змінах змісту цих понять та виникненні різних підходів до розуміння їх сутності.

Розглянемо кілька основних тлумачень поняття «педагогічна технологія» (табл. Б.1).

Г. К. Селевко [281] основними структурними елементами педагогічної технології вбачає наступне: концепція; змістовний компонент, що включає цілі та зміст навчання та виховання; процесуальний (технологічний) компонент, що включає організацію навчально-виховного процесу, методи й форми діяльності студентів, методи та форми діяльності педагога, управлінську діяльність педагога в навчально-виховному процесі та діагностику навчально-виховного процесу.

*Таблиця* Б.1

**Тлумачення поняття «педагогічна технологія»**

| **Автор**<br>***Назва поняття*** | **Тлумачення поняття** |
|---|---|
| П. І. Підкасистий<br>*педагогічна технологія* | Новий напрямок у педагогічній науці, що займається конструюванням оптимальних навчальних систем, проектуванням навчальних процесів [247]. |
| В. П. Беспалько<br>*педагогічна технологія* | Змістовна техніка реалізації навчального процесу, в яку входять процеси навчання, організація навчання та засоби навчання [144, 7]. |
| Б. Т. Ліхачов<br>*педагогічна технологія* | Сукупність психолого-педагогічних настановлень, що визначають спеціальний набір та компонування форм, методів, способів, прийомів навчання, виховних засобів; педагогічна технологія є організаційно-методичним інструментарієм педагогічного процесу [229, 166–167]. |
| Г. К. Селевко<br>*педагогічна технологія* | Система функціонування всіх компонентів педагогічного процесу, що побудована на науковій основі, запрограмована в часі та просторі і веде до намічених результатів [281, 31]. |
| В. І. Загвязинський<br>*педагогічна технологія* | Область знань, що охоплює сферу практичних взаємодій викладача і студента в будь-яких видах діяльності, організованих на основі чіткого визначення мети, систематизації, алгоритмізації прийомів навчання [179]. |
| М. В. Кларін<br>*педагогічна технологія* | Системна сукупність і порядок функціонування всіх особистісних, інструментальних і методологічних засобів, що використовуються для досягнення педагогічних цілей [212]. |
| Ю. В. Триус<br>*педагогічна технологія* | Система способів, прийомів педагогічних дій і засобів, що охоплюють цілісний навчально-виховний процес від визначення його мети до очікуваних результатів і які цілеспрямовано, систематично й послідовно впроваджуються в педагогічну практику з метою підвищення якості освіти [316]. |

Аналіз наукової та навчальної літератури засвідчує, що не лише категоріальна, а й більш широка наукова база теорії педагогічної технології розроблена досить ґрунтовно: виділено ознаки та критерії, основні якості та аспекти, здійснено класифікацію, змодельовано та описано конкретні оптимальні технології навчання й виховання. Проте всі вони, на думку С. А. Лисенка [228], підтверджують тлумачення, наведені у табл. Б.1, розглядають педагогічні технології з позиції педагога як їх творця, організатора, керівника.

Ю. В. Буган [150] під *технологією навчання* розуміє певний порядок, логічність і послідовність викладу змісту навчання відповідно до поставленої мети; певною мірою алгоритмізація спільної діяльності викладача та студентів у процесі навчання, узгодженість їхніх дій та взаємовідносин.

Дж. Дж. Пір (Joseph J. Pear) визначає такі властивості технології навчання [99]:

1) навчання є спостережуваним процесом, а його результати – вимірювані;

2) процес досягнення бажаних результатів навчання визначається у термінах чітко виділених операцій;

3) ці операції визначаються на основі фактичних даних.

Автор зазначає, що комп’ютерно орієнтовані технології навчання найбільш повно реалізують ці властивості. Дійсно:

1) спостереження за процесом навчання опосередковується, а його результати – фіксуються, опрацьовуються та узагальнюються засобами ІКТ;

2) процес навчання алгоритмізується та реалізується у комп’ютерній програмі;

3) дії з управління процесом навчання обираються відповідно до результатів навчання.

Таким чином, доцільним є вибір засобів ІКТ для реалізації технології навчання.

Поняття інформаційної технології з’явилося з виникненням інформаційного суспільства, основою соціальної динаміки в якому є не традиційні, матеріальні, а інформаційні ресурси – знання, наука, організаційні чинники, інтелектуальні здібності людей, їх ініціатива і творчість. Вперше поняття і перспективи розвитку інформаційних технологій докладно проаналізував академік В. М. Глушков, який визначив *інформаційну технологію* як людино-машинну технологію збирання, опрацювання та передавання даних. До інформаційних технологій відносять усі види технологій, що використовуються для створення, збереження, обміну і використання даних в усіх можливих формах [189, 10].

М. І. Жалдак [172, 2] вводить поняття *інформаційної системи*, яка на його думку є системою збору, збереження, переробки, передачі та подання відомостей, що базуються на електронній техніці та засобах телекомунікацій.

Під *інформаційною технологією* М. І. Жалдак [177] розуміє сукупність методів та технічних засобів, що використовують для збирання, створення, організації, зберігання, опрацювання, передавання, подання і використання відомостей, розширюючи знання людей і розвиваючи їх можливості в управлінні технічними і соціальними процесами.

За Н. В. Морзе [234, 12] *інформаційна технологія* – це сукупність методів, засобів, прийомів, що забезпечують пошук, збирання, зберігання, опрацювання, подання, передавання відомостей між людьми на основі електронних засобів, комп'ютерної техніки та зв'язку.

Н. В. Макарова [182] під *інформаційною технологією* розуміє процес, що використовує сукупність засобів і методів збору, обробки і передачі даних (первинних відомостей) для отримання даних нової якості про стан об'єкта, процесу або явища.

У тлумачному словнику з інформаційно-педагогічних технологій говориться, що *інформаційна технологія* – це сукупність засобів і методів, за допомогою яких здійснюється процес обробки даних [218].

Аналіз тлумачень поняття «нові інформаційні технології» [164; 182; 184; 225] надав можливість сформулювати узагальнене тлумачення даного поняття. На нашу дамку, *нові інформаційні технології* – це сукупність засобів і методів з отримання, опрацювання, зберігання та передавання даних з використанням електронної техніки некваліфікованим користувачем.

Щодо загальної користі впровадження *комп'ютерних («нових») інформаційних технологій* у педагогічний процес С. П. Новіков [242, 32] зазначає, що вони допоможуть підвищити якість підготовки студентів, підготовки фахівця, який володіє сучасним науковим світоглядом і досвідом емоційно-ціннісних відносин до світу знань; використання засобів НІТ допомагає вирішенню різноманітних психолого-педагогічних проблем, у тому числі формування умінь і навичок здійснення експериментально-дослідної діяльності, вибору змісту навчання, а також можливостей застосування засобів НІТ в якості засобу навчальної, науково-дослідної та управлінської діяльності.

Тлумачний словник з інформаційно-педагогічних технологій визначає *інформаційно-комунікаційні технології* як сукупність методів, виробничих процесів і програмно-технічних засобів, інтегрованих з метою збору, обробки, зберігання, поширення, відображення й використання даних в інтересах її користувачів [218].

М. І. Жалдак [165] визначає *інформаційно-комунікаційні технології* як сукупність методів, засобів і прийомів, використовуваних для збирання, систематизації, зберігання, опрацювання, передавання, подання різних повідомлень і даних.

Н. В. Морзе [234, 12] визначає *інформаційно-комунікаційні технології* як інформаційні технології на базі персональних комп'ютерів, комп'ютерних мереж і засобів зв'язку, для яких характерна наявність зручного для роботи користувача середовища.

Федеральним агентством з технічного регулювання та метрології Російської Федерації було розроблено державний стандарт «Інформаційно-комунікаційні технології в освіті. Терміни та тлумачення», згідно якого *інформаційно-комунікаційна технології* визначають як інформаційні процеси та методи роботи з даними, що здійснюються за допомогою засобів обчислювальної техніки та засобів комунікацій [183].

Ю. В. Триус [316] визначає *інноваційні інформаційно-комунікаційні технології навчання* як оригінальні технології (методи, засоби, способи) створення, передавання і збереження навчальних матеріалів, інших інформаційних ресурсів освітнього призначення, а також організації і супроводу навчального процесу (традиційного, за допомогою Інтернет і мультимедіа, дистанційного, мобільного) за допомогою телекомунікаційного зв'язку та комп'ютерних систем і мереж, що цілеспрямовано, систематично й послідовно впроваджуються в педагогічну практику з метою підвищення якості освіти.

*В Елементи навчання за допомогою Інтернет і мультимедіа, його переваги, недоліки та основні проблеми*

*Елементи системи навчання за допомогою Інтернет і мультимедіа,* що є *спільними з дистанційним* [294, 91-92]:

– змістові об'єкти: навчальний матеріал поділений на модулі, що містять об'єкти різної природи – текст, графіку, зображення, аудіо, анімацію, відео тощо. Як правило, вони зберігаються в базі даних і доступні в залежності від потреб суб'єктів навчання. Результатом є індивідуалізація навчання – студенти отримують лише те, що їм потрібно, засвоюючи знання у бажаному темпі;

– спільноти: студенти можуть створювати Інтернет-спільноти для взаємодопомоги та обміну повідомленнями;

– експертна онлайн-допомога: викладачі або експерти (інструктори з курсу) доступні в мережі для проведення консультацій, відповіді на питання, організації обговорення;

– можливості для співпраці: за допомогою відповідного програмного забезпечення можна організувати онлайн-конференції, спільну роботу над проектом студентів, географічно віддалених один від одного;

– мультимедіа: сучасні аудіо- та відеотехнології подання навчальних матеріалів з метою стимулювання прагнення студентів до набуття знань та підвищення ефективності навчання.

До *переваг навчання*, що проводиться з використанням технологій навчання за допомогою Інтернет і мультимедіа відносять [37]:

1. Персоніфікація. Слухач навчання, що проводиться з використанням технологій навчання за допомогою Інтернет і мультимедіа, може самостійно: визначити швидкість вивчення навчального матеріалу; визначити, коли він хоче проходити навчання; визначити які саме розділи навчального матеріалу і в якій послідовності йому необхідно вивчити.

2. Можливість проходження навчання без відриву від виробництва.

3. Можливість комбінування навчального контенту для формування різноманітних навчальних програм, адаптованих під конкретного учня.

4. Можливість отримати набагато більше відомостей, необхідних для оцінки знань, навичок і умінь, отриманих в результаті проведеного навчання. У тому числі: час витрачається на питання, кількість спроб, питання або завдання, які викликали найбільші труднощі, тощо. Наявність таких відомостей надає можливість набагато гнучкіше управляти проведеним навчанням.

5. Вартість. Незважаючи на необхідність високих початкових інвестицій, навчання, яке проводиться з використанням технологій навчання за допомогою Інтернет і мультимедіа, виявляється значно дешевшим порівняно з традиційним очним навчанням.

6. Використання широкого діапазону різноманітних засобів навчання. Всі ці кошти можуть бути використані і при проведенні традиційного очного навчання, але частіше всього цього не відбувається, а навчання за допомогою Інтернет і мультимедіа вимагає обов’язкового їх використання. В результаті цього навчання, яке проводиться з використанням технологій навчання за допомогою Інтернет і мультимедіа, виявляється найчастіше більш ефективним в порівнянні з традиційним очним навчанням.

7. Можливість його використання для навчання осіб з особливими потребами.

8. Надання доступу до якісного навчання особам, які з тих чи інших причин не мають можливості навчатися за традиційною очною формою. Наприклад, в місці їх проживання немає якісного навчального закладу.

9. Побудова ефективної системи управління навчанням на основі можливості збору значно більших відомостей про проходження навчання студентом у порівнянні з традиційним очним навчанням.

До *недоліків навчання*, що проводиться з використанням технологій навчання за допомогою Інтернет і мультимедіа, слід віднести [37]:

1. Складність внесення оперативних змін, у випадку якщо навчання вже розпочалося.

2. Необхідність формування додаткової мотивації у студентів, порівняно з іншими формами навчання.

3. Необхідність високих інвестицій у розбудову середовища навчання за допомогою Інтернет і мультимедіа.

4. Висока залежність від технічної інфраструктури: збій в інфраструктурі може призвести до зниження ефективності чи взагалі зриву навчання.

5. Відсутність достатньої кількості фахівців у сфері технологій навчання за допомогою Інтернет і мультимедіа.

6. Високі інвестиції на конструювання якісного змісту навчання.

Швидкий розвиток принципово нового напряму в освітній сфері неминуче привів до появи великої кількості проблем. Швидкість подальшого розвитку технологій навчання за допомогою Інтернет і мультимедіа багато в чому залежить від того, наскільки успішно будуть вирішені існуючі на сьогодні проблеми. Можна виділити такі *основні проблеми* в сфері технологій навчання за допомогою Інтернет і мультимедіа [37]:

– проблема визначення еквівалентності дистанційних курсів і визнання дистанційної освіти наряду з традиційною очною освітою;

– мовна проблема при імпорті (експорті) освіти. Електронні курси,

розроблені на одній мові, вимагають значних інвестицій для їх перекладу на іншу мову, включаючи необхідність урахування соціальних, культурологічних та інших особливостей регіону, де буде проводитися навчання з використанням технологій навчання за допомогою Інтернет і мультимедіа;

– нерівномірний розвиток інформаційних технологій, особливо, в частині каналів передавання даних. Недостатня пропускна спроможність каналів передавання даних обмежує можливість застосування засобів навчання за допомогою Інтернет і мультимедіа;

– висока вартість розробки і підтримки в актуальному стані електронних навчальних курсів;

– різниця в часі в разі проведення навчання на великих територіях. Особливо актуальним це стає при використанні засобів навчання за допомогою Інтернет і мультимедіа, що функціонують в режимі реального часу.

*Г Розвиток та впровадження інформаційно-комунікаційних засобів у процес навчання*

Педагогічна психологія вищої школи розглядає процес навчання у комплексі інформаційно-навчальної, розвивальної і виховної функцій.

На думку М. І. Жалдака [176], широке використання сучасних ІКТ у навчальному процесі дає можливість розкрити значний потенціал усіх дисциплін, завдяки формуванню наукового світогляду, розвитку аналітичного і творчого мислення, суспільної свідомості і свідомого ставлення до навколишнього світу. При цьому в основу інформатизації навчального процесу слід покласти створення нових комп'ютерно-орієнтованих методичних систем навчання на принципах поступового і неантагоністичного вбудовування ІКТ у діючі дидактичні системи, гармонійного поєднання традиційних та комп'ютерно-орієнтованих технологій навчання, не заперечування і відкидання здобутків педагогічної науки минулого, а, навпаки, їх удосконалення і посилення за рахунок використання досягнень у розвитку комп'ютерної техніки і засобів зв'язку.

Використання ІКТ докорінно змінює роль і місце викладача та студента в процесі навчання, сприяє реалізації індивідуального підходу в навчанні. У такій моделі викладач перестає бути основним джерелом знань, значне місце починають займати сучасні технології навчання [147].

Використання ІКТ, на думку М. І. Жалдака [178, 72], надає можливість значно підвищити ефективність осмислення і засвоєння повідомлень і даних, що циркулюють в навчально-виховному процесі, за рахунок їх своєчасності, корисності, доцільного дозування, доступності (зрозумілості), педагогічно доцільної надлишковості, оперативного використання джерел навчального матеріалу, адаптації темпу подання навчального матеріалу до швидкості його осмислення і засвоєння, врахування індивідуальних особливостей студентів, ефективного поєднання індивідуальної і колективної навчально-пізнавальної діяльності, методів і засобів навчання, організаційних форм навчального процесу.

На думку Ю. І. Машбиця [232], використання ІКТ у навчальному процесі спричиняє суттєві зміни в методах навчання. Ефективність методів навчання при цьому підвищується завдяки тому, що:

– використання ІКТ та засобів навчання надає широкі зображувальні можливості в розкритті способу вивчення об'єкта, у наочному поданні прийомів аналізу умови завдання, контролю за власними діями тощо;

– значно розширюється коло навчальних завдань, зокрема, професійного змісту;

– використання мобільних ІКТ та засобів навчання створює умови

для надання масового характеру індивідуальному навчанню;

– використання мобільних ІКТ та засобів навчання надає можливість моделювати спільну діяльність викладача і студента на всіх етапах вивчення предмета.

Для більш повного та об'ємного бачення проблем і перспектив впровадження ІКТ в систему освіти, необхідно розглянути передісторію цього питання, звернути увагу на всі зміни в даному напряму аж до теперішнього часу. У більшості країн за останні десятиліття проблема використання ІКТ в сфері освіти стала одним із пріоритетних напрямів діяльності урядових органів. Задача розширення використання ІКТ в освітніх установах вирішується в рамках взятого курсу на комп'ютеризацію суспільства в цілому. І хоча на цьому шляху приймається чимало важливих і необхідних заходів (національні концепції, стратегії і програми), все ж є ряд негативних факторів, що суттєво уповільнюють досягнення високих результатів в області використання ІКТ в освіті. До таких факторів слід віднести [265]:

– фінансову нестабільність країн, що приводить до неповного або безсистемного фінансування проектів з упровадження ІКТ в освіту;

– низький рівень досвіду міжнародного співробітництва.

Зазначені фактори є взаємопов'язаними і, відповідно, вирішення одного спричинить позитивну динаміку і в іншому. Таким чином, розширення міжнародного співробітництва, одночасно з привнесенням успішного досвіду, може зняти частину фінансових зобов'язань з держави, що в свою чергу сприятиме більш динамічним темпам комп'ютеризації.

Розглянемо періодизацію розвитку та впровадження ІКТ у процес навчання, запропоновану С. А. Лисенком [228]. Критеріями даної періодизації є суттєві зміни в технологічному проектуванні, розробках технічних засобів, практиці їх застосування, а також динаміку відображення цих змін у педагогічній теорії та практиці.

Отже, *перший період* (кінець XIX ст. – 30-ті рр. XX ст.) характеризується впровадженням у педагогічну практику технічних засобів навчання, а в теорію – понять «технологія в освіті» та «педагогічна технологія». Останнє з'являється в 1930-ті рр. у працях радянських психологів і педагогів, які стояли на засадах педології (О. П. Нечаєв, О. О. Ухтомський, С. Т. Шацький). Термін «технологія виховання» вперше використав А. С. Макаренко.

Перші програми аудіовізуального навчання, створені в цей період у США, заклали підґрунтя технологічної революції в освіті та спонукали розробки теоретичних засад педагогічної технології як науки. У професійній підготовці педагогічних кадрів з'явився технологічна складова, що активно опрацьовували зарубіжні центри фахової підготовки.

*Другий період* (40-ві рр. – середина 50-х рр. XX ст.) пов'язаний із упровадженням в освітні процеси розвинутих країн різноманітних технічних засобів презентації аудіовізуальних даних. У 1940-х рр. виник термін «технологія навчання», який тлумачили як використання продукту інженерно-технічної думки в навчанні (магнітофонів, передавачів, телевізорів, проекторів). Перший освітньо-технологічний проект – інтегрована програма аудіовізуальної освіти Л. К. Ларсона (L. C. Larson) 1946-1947 н. р. в університеті штату Індіана, США [71].

*Третій період* (друга половина 50-х – 70-ті рр. XX ст.) – динамічне наростання, розширення та поглиблення процесів технологізації в педагогічній практиці, що стало підґрунтям для подальшого розвитку теоретичних засад на основі введення поняття «технологічний підхід до навчання» [169, 906]. Формується два аспекти: теоретичний, що ототожнюється з поняттям «технологія освіти» або науковий опис педагогічного процесу (мета, форми, методи, засоби, прогнозований результат); та практичний, що асоціюється з поняттям «технологія навчання», яке, у свою чергу, сприймається як використання технічних засобів навчання.

С. А. Лисенко виділяє такі особливості цього етапу [228]:

– розвиток матеріально-технічної бази освіти, насамперед, аудіовізуальної: у 1950–1960-их рр. – засоби зворотного зв'язку, електронні класи, навчальні машини, лінгафонні кабінети, тренажери; у 1970-х рр. – відеомагнітофон, кадропроектор, електронна дошка, синхронізатори звуку й зображення, кабінети програмованого навчання тощо;

– виникнення технології програмованого навчання (автор – Б. Ф. Скіннер, 1954 р.);

– використання в навчальних цілях мови програмування Logo (Массачусетський технологічний інститут США, С. Пейперт (S. Papert) та інші, 1968 р.);

– створення масового персонального комп'ютера Apple (автори – С. Джобс (S. Jobs), С. Возняк (S. Wozniak), 1976 р.);

– розширення теоретичної бази: поява таких наукових галузей і навчальних дисциплін, як інформатика, системний аналіз, теорія телекомунікацій, педагогічна кваліметрія тощо. Виникає поняття «технологія освіти» як науково-педагогічний опис сукупності засобів і методів педагогічного процесу;

– виникнення поняття «педагогічна технологія», формування тенденції розуміння цього терміна як педагогічної системи – цілісного,

взаємоузгодженого педагогічного процесу, який шляхом використання технологічних засобів та урахування принципу оптимізації підвищує його ефективність. Суттєвий вплив на розвиток технології як науки має нова галузь – теорія систем;

– розширення методичної бази – розробка методів аудіовізуального навчання, використання інших технічних засобів;

– активізація фахової підготовки педагогів-технологів.

*Четвертий період* (80-ті рр. XX ст.) характеризується комплексністю у використанні педагогічних технологій: створюються комп'ютерні й дисплейні класи, лабораторії; розробляються, збагачуються та стають більш якісними педагогічні програмні засоби та їх комплекси.

*П'ятий період* (з 1990 р. по теперішній час) – пов'язаний з використанням нових інформаційних технологій в освіті: персональні комп'ютери нових поколінь, комп'ютерні мережі, глобальна мережа Інтернет, мультимедійні засоби, комп'ютеризовані засоби комунікацій тощо; розробляються конкретні педагогічні технології на основі таких процесів, як інформатизація (задоволення інформаційних потреб) [222, 368], медіатизація (удосконалення засобів збирання, збереження та поширення даних) [324, 183], комп'ютеризація (удосконалення та впровадження в практику засобів пошуку й обробки даних) [222, 412], інтелектуалізація (розвиток розумових, пізнавальних, творчих здібностей людини) [222, 364].

Як зазначає М. І. Жалдак [177], широке впровадження інноваційних інформаційних технологій в навчальний процес породжує ряд проблем, що стосуються змісту, методів, організаційних форм і засобів навчання, інтеграції навчальних предметів і фундаменталізації знань, підготовки і удосконалення кваліфікації педагогічних кадрів, удосконалення системи неперервної освіти, зокрема системи самоосвіти і самовдосконалення викладачів, що забезпечує оволодіння ними основами сучасної інформаційної культури.

*Д LMS, що застосовуються у вищій освіті в США*

Розглянемо деякі LMS, що використовують для організації та підтримки процесу навчання в технічних ВНЗ США та опишемо їх характеристики.

<u>*Moodle*</u> (http://moodle.org) – це система управління курсами (CMS) з відкритим вихідним кодом, також відома як система управління навчанням (LMS) або віртуальне навчальне середовище (VLE).

Проект Moodle є найбільш поширеною системою для підтримки навчання у світі і має за мету забезпечення безперервного процесу навчання не тільки для студентів дистанційної форми організації навчання, але й для підтримки традиційного навчання за моделлю комбінованого навчання.

До особливостей системи Moodle відносять наступне:

1) процес навчання за допомогою системи Moodle можна організовувати для досить великої кількості студентів (>100000);

2) використання Moodle у процесі навчання сприяє організації самостійної роботи студентів та можливості організовувати та підтримувати процес навчання за моделлю комбінованого навчання;

3) система є досить популярною серед навчальних закладів усього світу;

4) система має модульну структуру, постійно оновлюється та має можливість інтегрувати в себе різноманітні додатки.

До можливостей, що надає Moodle, відносять [87]:

1. Якісну підтримку всіх основних інструментів, типових для середовищ дистанційного навчання:

– можливість публікувати ресурси в будь-яких форматах і управляти доступом до них;

– потужна і гнучка система тестування із створеним банком завдань;

– зручна система форумів та розсилок;

– можливість створення індивідуальних завдань в різних форматах (текст, файл, кілька файлів, завдання поза сайту, повідомлення в форумі, запис в глосарії, заповнена анкета, відкрите питання в тестах тощо);

– глосарії з підтримкою автопосилань на тлумачення терміну в навчальних матеріалах курсу, включаючи форуми;

– опитування.

2. Додаткові інструменти, що підвищують зручність і якість процесу навчання:

– гнучка система оцінок, з налаштовуваною шкалою і можливістю задання правил виведення проміжних і підсумкових оцінок;

– анкети для збору даних (з можливістю підтвердження і оцінювання викладачем);

– модуль «лекція» для створення сценаріїв адаптивного навчання (технологічна педагогічна система форм і методів, що сприяє ефективному індивідуальному навчанню);

– модуль «семінар» для розвитку критичного й аналітичного мислення студентів.

3. Можливість встановлення модулів і плагінів сторонніх розробників, що розширюють функціональність системи.

4. Налаштовувані шаблони оформлення.

5. Гнучка ієрархічна система управління повноваженнями користувачів на основі ролей (адміністратор, розробник курсу, викладач, асистент, студент, гість та інше).

6. Можливість інтеграції з зовнішніми базами даних за списками користувачів та підписками на курси.

7. Детальне протоколювання всіх дій в системі.

8. Можливість зберігати і відновлювати курси з файлу.

9. Можливість оновлення версії із збереженням усіх даних в системі.

10. Більше 140 типів модулів і плагінів, що надають можливість розширювати всі аспекти функціональності Moodle без модифікації коду ядра і без втрати можливості простого переходу на майбутні версії.

11. Підтримка концепції соціально-педагогічного конструкціонізму (взаємодія, навчальна діяльність, критичне осмислення).

12. Підтримка стандарту SCORM.

*Sakai* (http://www.sakaiproject.org/) – навчальне середовище (Collaboration and Learning Environment – CLE), що являє собою гнучкий корпоративний додаток, який підтримує викладання, навчання і наукове співробітництво в цілому або частково в онлайновому середовищі. Sakai CLE поширюється безкоштовно, має відкритий вихідний код, пропонує максимальну гнучкість і надає можливість уникати ризику зростаючих витрат на ліцензію.

Навчальне середовище Sakai містить такі складові:

– адміністрування системи (система оцінювання студента, система підтримки двостороннього зв’язку: дискусії, чати, система тестування);

– управління курсом: інструмент для дослідження та групових проектів. Для підтримки цієї функції в Sakai містяться інструменти, що надають можливість змінити налаштування усіх інструментів на основі розподілених ролей;

– wiki – Web-сайт або інша гіпертекстова збірка документів, що надає користувачам можливість змінювати самостійно вміст сторінок через браузер, використовуючи спрощену і зручнішу, порівняно з HTML, wiki-розмітку тексту;

– розсилки та архівування новин, анонсів статей, зображень, аудіо та

відео матеріалів.

До основних інструментів навчального середовища Sakai відносять [69]:

– інструменти налаштування сайту;

– дошку оголошень, на якій інформують учасників процесу навчання про поточні питання;

– DropBox, що надає можливість викладачам і студентам обмінюватися документами в особистій папці, яка створена для кожного учасника;

– архів електронної пошти зберігає всі повідомлення, відправлені на адресу електронної пошти сайту;

– навчальні ресурси, що можуть бути доступні тільки дозволеним користувачам курсу;

– інструмент для створення презентацій курсу;

– чат та форум для організації обговорення поставленої проблеми;

– Message Center – інструмент комунікації, що надає учасникам сайту можливість спілкуватися за допомогою внутрішньої пошти курсу;

– інструмент для проведення опитування студентів, що надає користувачам можливість голосувати за правильну відповідь;

– інструменти для організації спільних проектів;

– інструменти для створення портфоліо.

*Desire2Learn* (http://www.desire2learn.com/) є навчальним середовищем, створеним для підтримки процесу навчання, що містить набір необхідних інструментів для створення якісного курсу:

– інструмент налаштування середовища, що містить макети системи із яких можна обрати зручний для викладача вигляд системи: колір та фон сторінки, шрифт, панель навігації;

– панель навігації знаходиться у верхній частині кожної сторінки і надає студентам можливість переходити від одного курсу до іншого та переходити до новин і повідомлень;

– інструмент побудови курсу – завдяки інтуїтивно зрозумілому і доступному інтерфейсу можна будувати індивідуальний навчальний курс;

– інструмент вбудовування медіа ресурсу – аудіо та відеоматеріали можуть бути легко вбудовані в курс, що надає можливість створювати аудіо та відеолекції;

– календар надає адміністраторам, викладачам і студентам можливість управляти своїм часом у навчальному середовищі, викладачам переглядати час проходження завдань студентами, студентам створювати індивідуальну траєкторію навчання;

– інструмент ClassList містить дані про студентів курсу, їх доступ до

навчальних матеріалів, листування за допомогою електронної пошти чи миттєвих повідомлень;

– інструмент навчальних досягнень студентів, що надає викладачу можливість визначати навчальні досягнення кожного студенту курсу, порівнювати його з іншими студентами курсу і повідомляти студента про його успіхи та невдачі;

– інструмент поступового вивчення матеріалу автоматично керує процесом навчання студента: студент не може переходити до наступного блоку навчального матеріалу, поки не дасть правильні відповіді за попередньою темою.

*WebCT* – система дистанційного навчання, що являє собою інструментальне середовище для створення навчальних курсів і призначена для організації та супроводу процесу навчання в мережі Інтернет.

Система дистанційного навчання WebCT вперше була протестована у 1997 р.

WebCT притаманні такі функціональні можливості:

– серверна частина розрахована на значну кількість студентів і має розширений інструментарій для формування змісту курсу: глосарій, пошукова система, предметний покажчик, база даних графіків та рисунків;

– різноманітні інструменти для створення якісних дистанційних курсів: шаблони для створення курсу, пошукова система за заданими шаблонами, гіпертекстовий словник термінів курсу, бібліотека мультимедійних файлів;

– тестова система для самотестування та тестування викладачем студентів курсу, архів результатів тестування з візуалізацією помилок;

– система моніторингу знань, що складається з моніторингу поточної успішності студентів та подання їх робіт на «електронній дошці оголошень» курсу;

– різноманітні комунікаційні засоби: текстовий діалог в рамках інструменту «Завдання», електронну пошту, форуми, чат.

Серед позитивних рис системи слід зазначити наступне: інструментальна база системи WebCT ніяким чином не обмежує викладача у виборі інструментів для формування авторських дистанційних курсів. Розроблені курси можуть бути як простими послідовно структурованими курсами, так і динамічними; використання системи у процесі навчання забезпечує поступову адаптацію студента до середовища навчання і може бути використане для організації процесу навчання за моделлю комбінованого навчання.

Зараз пакет придбаний компанією BlackBoard [329].

*Blackboard* (www.blackboard.com) – платформа для дистанційного навчання, що орієнтована на підтримку процесу навчання в асинхронному режими і складається з трьох частин:

– частина 1 містить елементи для організації роботи студентів курсу в режимі онлайн: розклад занять, відомості про курс і оголошення, електронні підручники, лекції, завдання, тести, рейтинги студентів курсу; убудований інструментарій для комунікацій і спільної роботи включає дискусії, електронну пошту і текстовий чат. Студенти мають можливість робити замітки у записнику (online notebook); перевіряти розклад групи за календарем; здавати домашні завдання, використовуючи зошит завдань (Digital Drop Box) і виконувати наукову роботу, використовуючи інтегровані академічні web ресурси (Academic Web Resources);

– частина 2 додає до платформи модулі, що можуть інтегрувати раніше створені курси, комунікаційні й адміністративні сервіси через Інтернет доступ;

– частина 3 додає до платформи адміністративний пакет, який надає можливість управляти навчальним процесом усього університету.

Рішення сумісне з Microsoft .NET. Надає можливість розміщувати матеріали курсів у форматах Microsoft Office, Adobe Acrobat PDF, HTML, різні формати графіки, аудіо та відео, а також анімаційні ролики (Flash, Shockwave, Authorware) [137].

*Pearson LearningStudio* (http://www.pearsonlearningsolutions.com/). Офіційний сайт Pearson LearningStudio зазначає, що Pearson LearningStudio реалізує навчання в режимі онлайн без обмежень. Технологічна платформа оптимізована для розширення програм навчання. Сайт має понад 45 мільйонів використаних хвилин на день і більше 9 мільйонів активних користувачів. Pearson LearningStudio має багатий досвід в обслуговуванні великомасштабних програм навчання в режимі онлайн.

Pearson LearningStudio встановлює стандарти для забезпечення безпеки і надійності, дозволяючи клієнтам не турбуватися про програмне забезпечення, обслуговування, питання працездатності, масштабованості, безпеки, аварійного відновлення і забезпечення безперервності бізнесу.

Економія Pearson LearningStudio, притаманна хмарним платформам, виключає витрати на апаратну установку, оновлення ліцензії програмного забезпечення, системи управління та необхідності зберігати оновлення.

Pearson LearningStudio надає клієнтам всю технічну підтримку, необхідну для живлення успішних онлайн середовищі навчання.

*Е Приклад розкладу занять на осінній семестр 2012-2013 н. р. у МТІ з курсу 18.014 – «Числення однієї змінної з теорією»*

*Таблиця* Е.1

**Розклад занять на осінній семестр 2012-2013 н. р. у МТІ з 18.014 – «Числення однієї змінної з теорією»**

| Дата | Тема | Основна література | Додаткова література |
|---|---|---|---|
| 05.09 | Послідовності і функцій, натуральні числа та математична індукція | I 2, Конспект | I 3.6, I 4.1-6, 1.2 |
| 06.09 | Операції над натуральними числами, принцип упорядкованості для цілих чисел | Конспект | |
| 07.09 | Раціональні числа | Конспект | |
| 11.09 | Дійсні числа | Конспект | |
| 13.09 | Основна теорема Дедекінда та повноти | Конспект, I 3,8-9 | I 3.10-15 |
| 14.09 | Функції дійсних чисел, аксіоми площини | 1.2-11 | |
| 18.09 | Інтеграл Рімана I | 1.12-16 | Необов'язково: гл. 2 |
| 20.09 | Інтеграл Рімана II | 1.17-23, примітка C OCW | Необов'язково: гл. 2 |
| 21.09 | Канікули | немає занять | |
| 25.09 | Інтеграл Рімана III | 1.24-27 | 2.1-4 |
| 27.09 | Границя та неперервність I | 3.1-6, 7.15 | |
| 28.09 | Презентації | | |
| 02.10 | Границя та неперервність II | 3,7-11 | |
| 04.10 | Границя та неперервність III | 3.12-16 | Примітки Г від OCW |
| 05.10 | Границя та неперервність IV | 3.17-20 | Примітки Н від OCW |
| 09.10 | День Колумба | немає занять | |
| 11.10 | Диференціювання | 4.1-9 | |
| 12.10 | Основна теорема диференціального числення | 5.1-6 | Примітки К від OCW |
| 16.10 | Екзамен I | все перераховане | |

| Дата | Тема | Основна література | Додаткова література |
|---|---|---|---|
| | | вище | |
| 18.10 | Теорема про середнє значення | 4.10-15 | |
| 19.10 | Презентації | | |
| 23.10 | Методи і застосування диференціювання | 4.16-21 | |
| 25.10 | Логарифм і експонента | Нотатки M від OCW, | 6,1-5,9, 6.1-19 |
| 26.10 | Тригонометричні функції | Нотатки L з OCW | 2,5-8 |
| 30.10 | Обернені функції та інше | 6.20-22 | |
| 01.11 | Методи інтегрування | 6,23-26,5.7-11 | |
| 02.11 | Презентації | | |
| 06.11 | Поліном Тейлора | 7,1-6 | Примітки O від OCW |
| 08.11 | Правило Лопіталя | 7.12-17, | Примітки P від OCW |
| 09.11 | Послідовності та ряди I | 10.1-9 | |
| 13.11 | Послідовності та ряди II | 10.11-16 | Нотатки Q від OCW |
| 15.11 | Послідовності та ряди III | 10.17-20,23 | |
| 16.11 | Презентації | | |
| 20.11 | Рівномірна збіжність I | 11.1-5 | |
| 22.11 | День подяки | немає занять | |
| 23.11 | День подяки | немає занять | |
| 27.11 | Екзамен II | Все з моменту екзамен I | |
| 29.11 | Рівномірна збіжність II | | |
| 30.11 | Презентації | | |
| 04.12 | Степеневі ряди, ряд Тейлора | 11.6-11 | Примітки R і S від OCW |
| 06.12 | Ряди Фур’є | | Примітки T від OCW |
| 07.12 | Лишки та огляд | | |
| 11.12 | Огляд | | |

| Date | Topic | Reading | Supplementary reading |
|---|---|---|---|
| Sep 5 | Sets and functions, natural numbers, and mathematical induction | I 2, Lecture notes | I 3.6, I 4.1-6, 1.2 |
| Sep 6 | Operations of natural numbers, the well-ordering principle, the integers | Lecture notes | |
| Sep 7 | Rational numbers | Lecture notes | |
| Sep 11 | Real numbers | Lecture notes | |
| Sep 13 | Dedekind's main theorem and completeness | Lecture notes, I 3.8-9 | I 3.10-15 |
| Sep 14 | Functions on the real numbers, area axioms | 1.2-11 | |
| Sep 18 | The Riemann integral I | 1.12-16 | Optional: Ch 2 |
| Sep 20 | The Riemann integral II | 1.17-23, Notes C from OCW | Optional: Ch 2 |
| Sep 21 | Student holiday | no class | |
| Sep 25 | The Riemann integral III | 1.24-27 | 2.1-4 |
| Sep 27 | Limits and continuity I | 3.1-6, 7.15 | |
| Sep 28 | Presentations | | |
| Oct 2 | Limits and continuity II | 3.7-11 | |
| Oct 4 | Limits and continuity III | 3.12-16 | Notes G from OCW |
| Oct 5 | Limits and continuity IV | 3.17-20 | Notes H from OCW |
| Oct 9 | Columbus Day | no class | |
| Oct 11 | Differentiation | 4.1-9 | |
| Oct 12 | The fundamental theorem of calculus | 5.1-6 | Notes K from OCW |
| Oct 16 | Exam I | all of the above | |
| Oct 18 | Mean value theorem | 4.10-15 | |
| Oct 19 | Presentations | | |
| Oct 23 | Techniques and applications of differentiation | 4.16-21 | |
| Oct 25 | Logarithm and exponent | Notes M from OCW, 6.1-19 | |
| Oct 26 | Trigonometric functions | Notes L from OCW, Lecture notes | 2.5-8 |
| Oct 30 | Inverse functions & more | 6.20-22... | |
| Nov 1 | Techniques of integration | 6.23-26,5.7-11 | |
| Nov 2 | Presentations | | |
| Nov 6 | Taylor polynomials | 7.1-6 | Notes O from OCW |
| Nov 8 | L'Hôspital's rule | 7.12-17, Notes P from OCW | |
| Nov 9 | Sequences and series I | 10.1-9 | |
| Nov 13 | Sequences and series II | 10.11-16, Notes Q from OCW | |
| Nov 15 | Sequences and series III | 10.17-20,23 | |
| Nov 16 | Presentations | | |
| Nov 20 | Uniform convergence I | 11.1-5 | |
| Nov 22 | Thanksgiving | no class | |
| Nov 23 | Thanksgiving | no class | |
| Nov 27 | Exam II | everything since Exam I | |
| Nov 29 | Uniform convergence II | | |
| Nov 30 | Presentations | | |
| Dec 4 | Power, Taylor series | 11.6-11, Notes R and S from OCW | |
| Dec 6 | Fourier series | Notes T from OCW | |
| Dec 7 | Leftovers and review | | |
| Dec 11 | Review | | |

Рис. Е.1. Екранна копія розкладу занять
(http://math.mit.edu/~vstojanoska/18014/schedule.html)

*Ж Використання інформаційно-комунікаційних технологій у вищій освіті в Україні: поточний стан, проблеми, перспективи розвитку*

У звіті директора Українського інституту інформаційних технологій І. Г. Малюкової [254], виконаного на замовлення Інституту ЮНЕСКО з інформаційних технологій в освіті, зазначено, що значну роль у впровадженні ІКТ в освітню сферу зіграв Закон України «Про Національну програму інформатизації» від 13.09.2001р. № 74/98-ВР [262], в рамках якого було реалізовано кілька проектів інформатизації навчальних закладів. Важливе значення у виборі напрямів і завдань розвитку навчання за допомогою Інтернет і мультимедіа (дистанційного) в Україну мала «Програма розвитку системи дистанційного навчання на 2004-2006 роки», затверджена Постановою Кабінету Міністрів України від 23.09.2003 року № 1494.

На сучасному етапі найбільший вплив на розвиток ІКТ у вищій освіті мають [254; 255; 256; 257; 258; 260; 261; 262; 263]:

– Закон України «Про основні положення розвитку інформаційного суспільства в Україні на 2007-2015 роки» від 09.01.2007 р. № 537-V, Розпорядження Кабінету Міністрів від 15 серпня 2008 року № 653-р «Про затвердження плану заходів з виконання завдань, передбачених Законом України «Про основні положення розвитку інформаційного суспільства в Україні на 2007-2015 роки», які містять положення про ефективне впровадження ІКТ в сфері освіти, в тому числі вищої;

– Державна програма «Інформаційні та комунікаційні технології в освіті і науці» на 2006-2010 роки [256], затверджена Постановою Кабінету Міністрів України від 07.12.2005 року № 1153, безпосередньо визначає план дій з розвитку засобів ІКТ для освітньої сфери, в тому числі, для вищої освіти .

Для реалізації програм, націлених на широкомасштабне та ефективне впровадження ІКТ в систему вищої освіти, було здійснено ряд організаційних заходів – як з боку державних органів влади, так і освітньо-наукового співтовариства:

– при Верховній Раді Україні створено та функціонує Консультативна рада з питань інформатизації;

– при Кабінеті Міністрів України створено Міжгалузеву раду з питань розвитку інформаційного суспільства (Постанова КМУ від 14.03.2009 р. № 4);

– при Міністерстві освіти і науки молоді та спорту Україні створені: Український інститут інформаційних технологій в освіті, УІІТО (на базі Національного технічного університету України «Київський політехнічний інститут», НТУУ «КПІ»); Координаційна рада з питань дистанційного навчання (Науково-технічна рада Державної програми «Інформаційні та комунікаційні технології в освіті і науці» на 2006-2010 роки) [254].

Нормативно-правовим забезпеченням використання ІКТ у вищій освіті регламентуються:

– *Укази Президента України*: «Про заходи щодо розвитку національної складової глобальної інформаційної мережі Інтернет» від 31.07.2007 р. № 928; «Про додаткові заходи щодо забезпечення відкритості у діяльності органів державної влади» від 01.08.2002 р. № 683; «Про невідкладні заходи щодо забезпечення функціонування та розвитку освіти в Україні» від 04.07.2005 р. № 1013; «Про першочергові завдання впровадження новітніх інформаційних технологій» від 20.10.2005 р. № 1497.

– *Закони України*: «Про Концепцію національної програми інформатизації» від 04.02.1998 р. № 75/98-ВР; «Про національну програму інформатизації» від 04.02.1998 р. № 74/98-ВР; «Порядок локалізації програмних продуктів (програмних засобів) для виконання Національної програми інформатизації» від 16.11.1998 р. № 1815; «Про пріоритетні напрями інноваційної діяльності в Україні» від 11.07.2001 р. № 2623-ІІІ; «Про вищу освіту» від 17.01.2002 р. № 2984-ІІІ; «Про пріоритетні напрями розвитку науки і техніки в Україні» від 16.01.2003 р. № 433-IV; «Про електронні документи та електронний документообіг» від 22.05.2003 р. № 851-IV; «Про електронний цифровий підпис» від 22.05.2003 р. № 852-IV; «Про телекомунікації» від 18.11.2003 р. № 1280-IV; «Порядок легалізації комп’ютерних програм в органах виконавчої влади» від 04.03.2004 р. № 253; «Про державні цільові програми» від 18.03.2004 р. № 1621-IV; «Про наукову і науково-технічну діяльність» від 20.11.2003 р. № 1316-IV; «Про захист інформації в інформаційно-телекомунікаційних системах» від 31.05.2005 р. № 2594-IV; «Про основні положення розвитку інформаційного суспільства в Україні на 2007-2015 роки» від 09.01.2007 р. № 537-V.

– *Постанови Верховної Ради України:* «Про затвердження Завдань Національної програми інформатизації на 2006-2008 роки» від 04.11.2005 р. № 3075-IV; «Про Рекомендації парламентських слухань з питань розвитку інформаційного суспільства в Україні» від 01.12.2005 р. № 3175-IV; «Про Рекомендації парламентських слухань на тему: "Створення в Україні сприятливих умов для розвитку індустрії програмного забезпечення"» від 15.03.2012 № 4538-VI;

– *Постанови Кабінету Міністрів України:* «Про Порядок оприлюднення у мережі Інтернет інформації про діяльність органів виконавчої влади» від 04.01.2002 р. № 3; «Про затвердження Порядку підключення до глобальних мереж передачі даних» від 12.04.2002 р. № 522; «Про затвердження Порядку проведення експертизи Національної програми інформатизації та окремих її завдань (проектів)» від 25.07.2002 р. № 1048; «Про заходи щодо подальшого забезпечення діяльності органів виконавчої влади» від 29.08.2002 р. № 1302; «Про затвердження Порядку взаємодії органів виконавчої влади з питань захисту державних інформаційних ресурсів в інформаційних та телекомунікаційних системах» від 16.11.2002 р. № 1772; «Про заходи щодо створення електронної інформаційної системи «Електронний Уряд» від 24.02.2003 р. № 208; «Концепція формування системи національних електронних інформаційних ресурсів» від 05.05.2003 р. № 259-р; «Про затвердження Програми розвитку системи дистанційного навчання на 2004-2006 роки» від 23.09.2003 р. № 1494; «Про затвердження Порядку використання комп'ютерних програм в органах виконавчої влади» від 10.09.2003 р. № 1433; «Про затвердження Державної програми розвитку і функціонування української мови на 2004-2010 роки» від 02.10.2003 р. № 1546; «Про затвердження Положення про Національний реєстр електронних інформаційних ресурсів» від 17.03.2004 р. № 326; «Про затвердження Комплексної програми забезпечення загальноосвітніх, професійно-технічних та вищих навчальних закладів сучасними технічними засобами навчання з природничих, математичних і технологічних дисциплін» від 13.08.2004 р. № 905; «Про затвердження Державної програми інформатизації та комп'ютеризації вищих навчальних закладів І-ІІ рівня акредитації на 2005-2008 роки» від 08.09.2004 р. № 1182; «Про затвердження Положення про Реєстр інформаційних, телекомунікаційних та інформаційно-телекомунікаційних систем органів виконавчої влади, а також підприємств, установ та організацій, що належать до сфери їх управління» від 03.08.2005 р. № 688; «Про затвердження Державної програми« Інформаційні та комунікаційні технології в освіті і науці на 2006-2010 роки» від 07.12.2005 р. № 1153; «Про затвердження Порядку використання у 2006 році коштів, передбачених у Державному бюджеті для інформатизації та комп'ютеризації професійно-технічних та вищих навчальних закладів, забезпечення їх сучасними технічними засобами навчання з природничих, математичних і технологічних дисциплін» від 24.05.2006 р. № 712; «Про схвалення Концепції Державної цільової програми впровадження у навчально-виховний процес загальноосвітніх навчальних закладів ІКТ «Сто відсотків» на період до 2015 року» від 27.08.10 р. № 1722-р;

– *Накази Міністерства освіти і науки України:* «Про створення Українського центру дистанційної освіти» від 07.07.2000 р. № 293; «Про

створення Координаційної ради Міністерства освіти і науки України з питань дистанційного навчання» від 26.02.2001 р. № 91; «Про створення Українського інституту інформаційних технологій в освіті Національного технічного університету України «Київський політехнічний інститут» від 24.11.2004 р. № 880; «Про реалізацію спільного проекту МОН, Представництва ООН в Україні та Всеукраїнської асоціації комп'ютерних клубів. Розвиток доступу до сучасних інформаційно-комунікаційних технологій населення на основі партнерства між школами і комп'ютерними клубами» від 13.12.2004 р. № 935; «Про проведення апробації електронних засобів навчального та загального призначення для загальноосвітніх навчальних закладів» від 20.03.2006 р. № 213; «Про затвердження вимог до специфікації навчальних комп'ютерних комплексів для оснащення кабінетів інформатики та інформаційно-комунікаційних технологій навчальних закладів системи загальної середньої освіти» від 11.05.2006 р. № 363; «Про затвердження тимчасових вимог до педагогічних програмних засобів» від 15.05.2006 р. № 369; «Про затвердження тимчасових рекомендацій визначення трудомісткості створення педагогічних програмних засобів» від 05.06.2006 р. № 432; «Про створення Центру розвитку інформаційного суспільства Національного технічного університету України «Київський політехнічний інститут» від 05.06.2006 р. № 429; «Про забезпечення функціонування інформаційної системи. «Конкурс» від 11.06.2008 р. № 514; «Про Інформаційно-пошукову систему «Конкурс» від 14.01.2009 р. № 16; «Про продовження Всеукраїнського експерименту щодо навчання вчителів ефективному використанню інформаційно-комунікаційних технологій у навчальному процесі та підвищення кваліфікації педагогічних працівників за програмою Intel ® «Навчання для майбутнього» від 24.03.2009 р. № 271; «Про продовження Всеукраїнського експерименту щодо навчання вчителів ефективному використанню інформаційно-комунікаційних технологій у навчальному процесі та підвищення кваліфікації педагогічних працівників за програмою Intel ® «Навчання для майбутнього» від 24.03.2009 р. № 271; «Про впровадження моделі навчання "1 учень 1 комп'ютер"» від 11.03.2010 р. №196; «Про затвердження Положення про дистанційне навчання» від 01.06.2013 р. № 660.

Вимоги до необхідного рівня комп'ютерного і програмного забезпечення, доступу викладачів і студентів до Інтернету (і іншим комунікаційним мережам), наявності електронної бібліотеки у ВНЗ визначаються документом «Порядок ліцензування діяльності вищого навчального закладу з надання освітніх послуг вищої освіти» [254].

Одним з важливих показників рівня впровадження ІКТ у навчальний процес і процес управління ВНЗ є забезпечення доступу викладачам і студентам до телекомунікаційних мереж: локальним (Інтернет), корпоративної в науково-освітній сфері (Українська науково-освітня телекомунікаційна мережа УРАН), глобальної мережі (Інтернет) [254]. Аналіз даних показав, що майже кожен вищий навчальний заклад має локальну мережу. Всі ВНЗ мають підключення до Інтернет. При цьому середня кількість провайдерських каналів для підключення одного ВНЗ до Інтернету складає ~ 1,9. Пропускна здатність каналів, в середньому, становить 150 Мбіт/с. Кількість користувачів електронної пошти в одному ВНЗ, в середньому, становить близько 1000.

Особливий інтерес представляють дані про використання інформаційних засобів мережі УРАН. Українська науково-освітня телекомунікаційна мережа УРАН (Мережа УРАН) створена за рішенням Міністерства освіти України та Національної Академії наук України за підтримки університетів, інститутів Міністерства освіти і Національної Академії наук згідно Спільному Постановою Президії Національної Академії наук України і Колегії Міністерства освіти України від 20 червня 1997 року.

Експлуатація та подальший розвиток Мережі УРАН здійснюється Асоціацією УРАН згідно Концепції Національної програми інформатизації та Державної Програми «Інформаційні та комунікаційні технології в освіті і науці» на 2006-2010 роки.

Діяльність Асоціації є неприбутковою, а розвиток мережної інфраструктури забезпечується, в основному, за рахунок цільового державного фінансування або міжнародних грантів. Асоціація УРАН налічує 80 вищих навчальних закладів та наукових установ.

Головним призначенням Мережі УРАН є забезпечення установ, організацій та фізичних осіб у сфері освіти, науки і культури України інформаційними послугами на основі Інтернет-технологій для реалізації професійних потреб і розвитку зазначених областей. Такі послуги передбачають, зокрема, оперативний доступ до даних, обмін ними, їх розповсюдження, накопичення та обробку для проведення наукових досліджень, навчання за допомогою Інтернет і мультимедіа, електронного тестування, використання методів телематики, функціонування електронних бібліотек, віртуальних лабораторій, проведення телеконференцій, реалізацію дистанційних методів моніторингу тощо [254].

Мережа УРАН будується за ієрархічним принципом: у кожному місті України, який є регіональним центром наукової та освітньої діяльності, створюється регіональний вузол мережі на базі університету або наукової установи.

Базовою організацією Головного центру управління Мережею УРАН є Міністерство освіти і науки в Києві. Головний центр управління Мережі УРАН забезпечує основний інформаційний сервіс мережі та функціонування її бекбону. Крім того, Головний центр управління забезпечує функції регіонального вузла для користувачів Київського регіону.

Розвиток міських волоконно-оптичних сегментів було здійснено протягом 1997-2007 років у рамках інфраструктурних грантів НАТО та державного замовлення з боку Міністерства освіти і науки.

Сьогодні мережа УРАН фізично об’єднує понад 80 науково-дослідницьких та освітніх закладів (180 точок підключення) та експлуатує власні волоконно-оптичні мережі в 15 містах України загальною довжиною близько 230 км і міжнародну волоконно-оптичну лінію зв’язку довжиною 80 км Львів – державний українсько-польський кордон. Топологія Мережі УРАН наведена на рис. Ж.1.

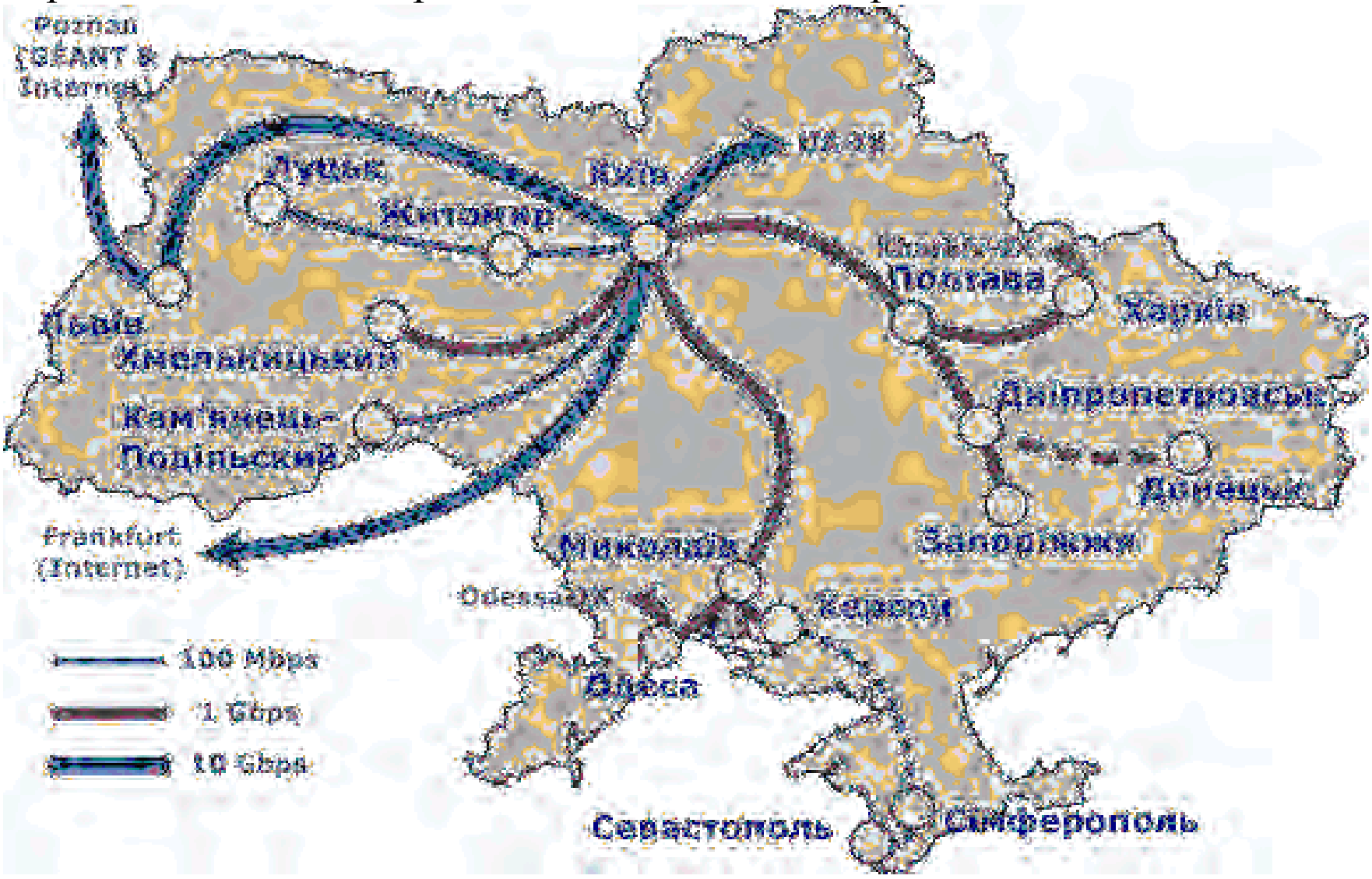

Рис. Ж.1. Топологія Мережі УРАН

У рамках реалізації Державної програми «Інформаційні та комунікаційні технології в освіті і науці» на 2006-2010 рр. в 2007 році був підписаний договір про підключення Мережі УРАН до пан-європейської науково-освітньої мережі GEANT2 і було організовано взаємоз’єднання мереж УРАН і GEANT2 в Польщі через канал 155 Мбіт/с, організований провідним оператором зв’язку для GEANT2 – компанією Memorex Telecommunication (Австрія).

GEANT2 – це високошвидкісна мережа Європи, що об’єднує каналами пропускної здатності 10-40 Гбіт/с національні наукові мережі європейських країн. Крім європейських країн GEANT2 пропонує глобальні зв’язку з повністю інтегрованим сервісом з національними науковими мережами в Північній (Internet2) і Південної (ALICE) Америці, Азії (TEIN2), Середземномор’ї (EUMEDCONNECT), Африці.

Інтеграція до європейських науково-освітніми мережами в рамках GEANT відкриває Україні нові можливості доступу до наукових та освітніх інформаційних ресурсів, зокрема, до віддалених центрів суперкомп’ютерних обчислень і наукових даних, електронних бібліотек, баз даних і знань, інформаційних пошукових систем, ресурсів дистанційного навчання тощо.

Для забезпечення якості навчання, відповідного сучасним вимогам, у вищих навчальних закладах процес напрацювання інформаційних ресурсів відбувається за рахунок власних можливостей, в тому числі фінансових; значна кількість таких ресурсів ініціативно створюється педагогами, науковцями, інженерами та студентами. У більшості ВНЗ акумуляторами напрацьованих інформаційних ресурсів є електронні бібліотеки, де вони накопичуються, в більшості своїй, у вигляді файлів різного формату. На рівні держави також приділяється увага цій проблемі. Так, в 2007-2008 роках з Державного бюджету в рамках Державної програми «Інформаційні та комунікаційні технології в освіті і науці» на 2006-2010 роки профінансовані:

– створення електронних підручників;

– створення пілотного проекту банку атестованих курсів дистанційного навчання для навчальних закладів всіх рівнів освіти;

– створення дистанційних курсів для вищих навчальних закладів;

– створення пілотного проекту типової електронної бібліотеки вищого навчального закладу;

– створення пілотного проекту типового програмно-апаратного комплексу системи архівації та зберігання контенту електронної наукової бібліотеки;

– створення дистанційних курсів для вищих навчальних закладів;

– створення пілотного проекту типової електронної бібліотеки вищого навчального закладу;

– створення пілотного проекту типового програмно-апаратного комплексу системи архівації та зберігання контенту електронної наукової бібліотеки;

– створення і наповнення повнотекстовими документами електронних бібліотек вищих навчальних закладів та порталу Кримської міжвузівської електронної бібліотеки;

– створення типового проекту абонентського бездротового доступу навчального закладу до інформаційних ресурсів;

– побудова освітньої національної GRID-інфраструктури для забезпечення наукових досліджень.

І. Г. Малюкова зазначає, що електронні інформаційні ресурси навчального призначення, створені за рахунок бюджетного фінансування в рамках програми «Інформаційні та комунікаційні технології в освіті і науці», знаходяться у вільному доступі для використання всіма державними вузами а також особами, які бажають їх використовувати для самоосвіти [254].

Крім цього, багато авторів електронних підручників та дистанційних курсів, популярність розробки яких зростає у середовищі викладачів, також надають вільний доступ до своїх інформаційних ресурсів. При цьому особливу перевагу автори віддають розробці дистанційних курсів в модульному об’єктно-орієнтованому динамічному навчальному середовищі (платформі) Moodle.

Український інститут інформаційних технологій в освіті НТУУ «КПІ» активно сприяє подальшій популяризації Moodle і навчання роботі в цій платформі, для чого на своєму сайті (http://www.udec.ntu-kpi.kiev.ua) відкрив вільний доступ викладачам університету, а також усім бажаючим до самої платформі та спеціального дистанційного курсу по роботі в Moodle.

Існує велика диспропорція між кількістю викладачів, які використовують та не використовують у навчальному процесі новітні технології, які пройшли і не пройшли перепідготовку або підвищення кваліфікації за цим напрямом.

Метою Державної програми «Інформаційно-комунікаційні технології в освіті і науці» у 2006-2010 році [256] було створення умов для розвитку освіти і науки, підвищення ефективності державного управління шляхом впровадження ІКТ, забезпечення реалізації прав на вільний пошук, одержання, передачу, виробництво і поширення даних, здійснення підготовки необхідних спеціалістів і кваліфікованих користувачів, сприяння розвитку вітчизняного виробництва високотехнологічної продукції і насамперед – конкурентоспроможних комп’ютерних програм як найважливішої складової інформаційних та комунікаційних технологій сприяння переходу економіки на інноваційний шлях розвитку.

Основними завданнями Програми було виконання комплексу завдань, що сприятимуть забезпеченню:

– підвищення загальної інформаційної грамотності населення;

– оснащення навчальних закладів сучасним комп’ютерним та

телекомунікаційним обладнанням;

– впровадження інформаційних та комунікаційних технологій у навчальний процес і проведення наукових досліджень, забезпечення доступу до національних і світових інформаційних ресурсів;

– розроблення, впровадження та легалізація програмного забезпечення;

– залучення мережевих технічних ресурсів для забезпечення підключення наукових установ та навчальних закладів до Інтернет;

– розвитку технологій дистанційного навчання і використання їх для запровадження в Україні системи навчання протягом усього життя;

– захисту прав інтелектуальної власності (авторів та розробників);

– підвищення кваліфікації та перепідготовка кадрів;

– розбудови інфраструктури науково-освітньої телекомунікаційної мережі (УРАН), підключення до неї наукових установ, наукових бібліотек, центрів науково-технічних даних за допомогою каналів передачі даних, інтеграцію їх з європейською науково-дослідницькою мережею (GEANT);

– розширення мережі електронних бібліотек навчальних закладів та наукових установ;

– розроблення систем забезпечення інформаційної безпеки функціонування мереж та інформаційних ресурсів.

Виконання завдань Програми здійснювалась з урахуванням стратегії соціально-економічного розвитку регіонів, стану та перспектив розвитку інформаційних і комунікаційних технологій, новітніх досягнень в інформаційній сфері.

Серед заходів, спрямованих на розвиток інформаційних та комунікаційних технологій в освіті і науці за 2006-2010 роки було виконано:

– оснащення універсальними навчально-комп’ютерними комплексами: загальноосвітніх та позашкільних навчальних закладів, наукових та науково-методичних установ, інститутів системи післядипломної педагогічної освіти;

– оснащення мобільними навчально-комп’ютерними комплексами: наукових та науково-методичних установ, загальноосвітніх навчальних закладів І-ІІІ ступеня;

– оснащення професійно-технічних навчальних закладів навчально-комп’ютерними комплексами

– створення мережі регіональних, базових та локальних центрів системи дистанційного навчання, оснащення їх програмно-технічними засобами;

– забезпечення вищих навчальних закладів, наукових та науково-

методичних установ технічними засобами і мережевим обладнанням;

– створення і модернізація локальних мереж у навчальних закладах, наукових та науково-методичних установах;

– підключення загальноосвітніх, позашкільних, професійно-технічних та вищих навчальних закладів до Інтернет;

– впровадження бездротових технологій, створення типових проектів абонентського доступу до ресурсів даних;

– розроблення та впровадження мікрохвильових систем широкосмугового радіодоступу до Інтернет у важкодоступних районах;

– розроблення та впровадження технології xDSL для доступу до Інтернет дротовою мережею зв’язку;

– оснащення ліцензійними програмами;

– створення електронних підручників та енциклопедій навчального призначення;

– створення фільмотеки навчальних фільмів на електронних носіях (40 відео-фрагментів);

– створення банку електронних документів нормативно-правового, науково-методичного, психолого-педагогічного, організаційного, програмно-технологічного та інформаційного забезпечення дистанційного навчання;

– створення та впровадження програмних засобів пілотної системи поточного і підсумкового контролю знань студентів у вищих навчальних закладах;

– створення та впровадження програмних засобів для уніфікованої системи дистанційного навчання;

– створення банку атестованих курсів дистанційного навчання для загальноосвітніх, професійно-технічних, вищих навчальних закладів та закладів післядипломної освіти;

– розроблення елементів штучного інтелекту та засобів і технологій для індивідуалізації навчального процесу та їх впровадження в систему дистанційного навчання;

– створення системи дистанційного навчання для перепідготовки та підвищення кваліфікації педагогічних і науково-педагогічних працівників загальноосвітніх, професійно-технічних та вищих навчальних закладів;

– створення центру розроблення та впровадження програмних засобів навчального призначення;

– створення GRID-інфраструктури для забезпечення наукових досліджень;

– розбудова інфраструктури національної науково-освітньої телекомунікаційної мережі (УРАН);

– створення Інтернет-порталу: загальної середньої та професійно-технічної освіти, дистанційного навчання, інформаційних ресурсів освіти і інноваційної діяльності;

– розроблення технології створення віртуальних навчальних інформаційних ресурсів за освітньо-кваліфікаційними рівнями, реалізація проекту підготовки бакалаврів за спеціальністю «Телекомунікації»;

– створення програмного та інформаційного забезпечення для електронних наукових бібліотек і архівів;

– забезпечення функціонування української мови в інформаційному комп'ютерному середовищі;

– створення українського сегмента міжнародної лінгвістичної системи;

– створення електронних бібліотек вищих навчальних закладів;

– створення системи електронних класифікаторів і нормативних документів для забезпечення дистанційної освіти;

– створення та впровадження типових макетів і шаблонів електронних документів для використання у вищих навчальних закладах;

– створення віртуального університету, розроблення та підтримка його інформаційних ресурсів;

– сертифікація та атестація програмних засобів та курсів дистанційного навчання;

– підвищення кваліфікації та перепідготовка кадрів: підготовка науково-педагогічних працівників вищих навчальних закладів та їх сертифікація для роботи з програмними засобами навчального призначення та інформаційними і комунікаційними технологіями;

– створення програмно-методичного комплексу та електронних ресурсів для підвищення кваліфікації працівників загальноосвітніх, професійно-технічних та вищих навчальних закладів, викладачів та наукових працівників у галузі інформаційних та комунікаційних технологій;

– створення державної системи оцінки знань та вмінь в галузі інформаційних та комунікаційних технологій із системою сертифікації, що відповідає міжнародному стандарту ECDL;

– створення комплексу державних стандартів України в галузі інформаційних та комунікаційних технологій;

– удосконалення нормативно-правової бази в сфері інтелектуальної власності щодо захисту розробок у галузі інформаційних та комунікаційних технологій;

– нормативно-правове забезпечення функціонування загальнодержавного реєстру інформаційних ресурсів науково-технічних

та освітніх даних;

– удосконалення нормативно-правової бази системи дистанційного навчання;

– створення системи моніторингу, планування та прогнозування діяльності навчальних закладів;

– розроблення типової автоматизованої системи управління загальноосвітнім навчальним закладом;

– розроблення програмних засобів системи незалежного тестування знань;

– створення та впровадження автоматизованої системи обліку підручників у загальноосвітніх навчальних закладах;

– створення галузевої системи сертифікації програмних засобів наукового і навчального призначення;

– створення інформаційної системи моніторингу результатів наукових досліджень;

– розроблення програмно-технічних систем забезпечення захисту інформаційних ресурсів від несанкціонованого доступу;

– розроблення та впровадження системи управління Програмою.

Удосконалення системи підготовки ІТ-фахівців має супроводжуватися:

– змінами в підходах до розробки освітніх державних стандартів, які повинні враховувати високу швидкість змін на ринку ІКТ;

– створенням (удосконаленням) науково-виробничих комплексів вузів ІТ-компаніями;

– формуванням спільного кадрового складу з викладачів, наукових працівників і фахівців компаній;

– реформуванням системи оплати праці в вузах для ІТ-спеціальностей з урахуванням залучення висококваліфікованих фахівців до навчального процесу.

Підвищення ефективності використання ІКТ у вищій освіті має супроводжуватися моніторингом досягнень у цій сфері, який буде спиратися на міжнародні методики та індикатори, а також на позитивний практичний досвід інших країн, у тому числі США.

*К LMS, що застосовуються у вищій освіті в Україні*

Розглянемо деякі LMS, що використовують для організації та підтримки процесу навчання в технічних ВНЗ України та опишемо їх характеристики.

**Агапа** (http://agapa.com.ua/) – це система, у якій поєднуються навчальний комплекс, комунікаційне середовище та система управління даними. Розробка даного програмного засобу була розпочата у 2003 році компанією «АВ-Консальтінг» і була першим програмним продуктом в Україні, що був зареєстрований як «система дистанційного навчання».

Становлення СПН «Агапа» відбувалося у три етапи [311]:

– на *першому етапі* було розроблено такі функціональні модулі, як «Обліковий запис», «Авторизація», «Реєстрація», «Особисті повідомлення», «Новини», «Зображення», «Особисті сторінки», «Нотатки» та «Друзі». Водночас відбулося і перше експериментальне використання системи. Після успішної апробації системи для побудови сайтів навчального призначення було розпочато роботи над функціональними модулями «Курси», «Глосарій», «Тестування», «Індивідуальні роботи», «Конкурси», «Анкети»;

– на *другому етапі* розвитку система набуває рис, характерних для LCMS;

– на *третьому етапі* розвитку системи було реалізовано функціональні модулі, притаманні LMS: «Журнали», «Планувальник», «Робочий стіл», «Форум», «Чат», що надало можливість ефективно підтримувати процес навчання у групах, які формуються з урахуванням організаційної структури навчального закладу і склад яких може бути як фіксованим, так і змінним;

– на *четвертому етапі* розвитку системи як BLMS реалізуються модулі, що надають можливість керувати значними об'ємами різних за природою та джерелами даних («Репозиторій», «Електронна бібліотека», «Конференції», «Файловий менеджер»), створюючи умови для побудови на базі системи освітні портали вищих навчальних закладів. Відбувається доопрацювання таких модулів як «Планувальник», «Тестування» та «Журнали успішності», а також створюються модулі «Навчальні плани», «Шкали оцінювання», «Імпорт/експорт» курсів та новин.

У СПН «Агапа» можна розміщувати лекції, проводити практичні, лабораторні та семінарські заняття, організовувати лабораторно-обчислювальні практикуми, підтримувати індивідуальну роботу студентів, а також проводити консультування студентів, яке може бути як дистанційне так і мобільне.

Систему «Агапа» у процесі навчання використовують Київський національний економічний університет ім. Вадима Гетьмана,

Криворізький національний університет, Луганський національний університет ім. Тараса Шевченка, Тернопільський національний економічний університет.

**Moodle** (модульне об'єктно-орієнтоване динамічне навчальне середовище) – вільнопоширювальна за ліцензією GNU General Public License система управління навчанням, що поширюється за ліцензією (http://moodle.org).

Система реалізує філософію «соціального конструктивізму» і орієнтована на організацію індивідуальної роботи студентів, що підтримується та курується викладачем.

Система Moodle має зручний інтерфейс провідника; що легко інсталюється майже на всі платформи, які підтримують PHP; потребує лише одну базу даних (і може ділитися нею); повна абстракція баз даних підтримує всі види баз даних (окрім визначення початкової таблиці); у списку курсів висвітлює опис для кожного курсу, включаючи можливість перегляду гостями; існує можливість категоризації та пошуку курсів; добре продумана система безпеки. Форми перевіряються, дані підтверджуються, cookies зашифровуються; більшість текстових областей вводу (ресурси, поштові відправлення) можуть бути відредаговані за допомогою вбудованого WYSIWYG HTML редактора [236].

Особливістю використання у процесі навчання СПН Moodle є те, що нові знання студент може отримати за умови власного індивідуального досвіду і найефективніше ці знання засвоюються лише при умові, що ці знання можуть бути передані іншим учасникам процесу навчання. Таку систему зручно використовувати для підтримки комбінованого навчання вищої математики.

Система Moodle відповідає всім основним критеріям, що висуваються до систем навчання за допомогою Інтернет і мультимедіа, зокрема таким, як [296]:

– *функціональність* – наявність набору функцій різного рівня (форуми, чати, аналіз активності слухачів (студентів), управління курсами та навчальними групами тощо);

– *надійність* – зручність адміністрування та управління навчанням, простота оновлення контенту на базі існуючих шаблонів, захист користувачів від зовнішніх дій тощо;

– *стабільність* – високий рівень стійкості роботи системи стосовно різних режимів роботи та активності користувачів;

– *вартість* – сама система безкоштовна, витрати на її впровадження, розробку курсів і супровід – мінімальні;

– відсутність обмежень за кількістю ліцензій на слухачів (студентів);

– *модульність* – наявність в навчальних курсах набору блоків матеріалу, які можуть бути використані в інших курсах;

– наявність вбудованих засобів розробки та редагування навчального контента, інтеграції різноманітних освітніх матеріалів різного призначення;

– підтримка міжнародного стандарту SCORM (Sharable Content Object Reference Model) – основи обміну електронними курсами, що забезпечує перенесення ресурсів до інших систем;

– наявність системи перевірки та оцінювання знань слухачів у режимі он-лайн (тести, завдання, контроль активності на форумах);

– зручність і простота використання та навігації – інтуїтивно зрозуміла технологія навчання (можливість легко знайти меню допомоги, простота переходу від одного розділу до іншого, спілкування з викладачем-тьютором тощо).

Систему підтримки навчання Moodle в Україні використовують такі навчальні заклади, як Києво-Могилянська академія, Національний технічний університет України «Київський політехнічний університет», Національний педагогічний університет ім. М. П. Драгоманова, Черкаський державний технологічний університет, Національний технічний університет «Харківський політехнічний інститут», ДВНЗ «Криворізький національний університет» та інші ВНЗ.

**TrainingWare Class** (http://www.ksob.ru/) перша російська СПН з відкритими кодами доступу, що розроблена з урахуванням вимог вітчизняної освіти. Розробником TrainingWare Class є ЗАТ «Корпоративні системи навчання».

TrainingWare Class – це технологічна платформа для автоматизації процесів навчання та атестації користувачів. Вона забезпечує взаємодію між викладачем і студентами в процесі навчання, що надає можливість організовувати процес навчання вищої математики за моделлю комбінованого навчання. В цій системі можна розташовувати навчальні відомості курсу, розробляти тести та проводити автоматизовану атестацію користувачів; створювати комплексні системи автоматизації навчальних процесів і системи моніторингу навчання на рівні кафедри, факультету чи університету; формувати єдині бібліотеки навчально-методичних матеріалів, створені учасниками педагогічних спільнот та соціальних мереж.

**Claroline LMS** (http://www.claroline.net/) – це платформа для навчання за допомогою Інтернет і мультимедіа (e-Learning) та електронної діяльності (e-Working), що надає можливість викладачам створювати ефективні онлайн-курси та керувати процесом навчання та спільними діями на основі Web-технологій.

Claroline LMS була створена в інституті педагогіки і мультимедіа католицького університету в Лувені. Продукт є безкоштовним і доступним багатьом навчальним закладам, оскільки одночасно може навчати близько 20000 студентів. Для використання даної СПН необхідно установити PHP / MySQL / Apache.

Використання СПН Claroline надає можливість створювати та адмініструвати курси в режимі онлайн, редагувати їх вміст, управляти ним. СПН включає генератор вікторин, форуми, календар, функцію розмежування доступу до документів, каталог посилань, систему контролю за успіхами студентів, модуль авторизації.

Кожен курс містить ряд інструментів, що надають викладачеві можливість:

– вказати опис курсу;

– опублікувати навчальні матеріали в будь-якому форматі (текст, PDF, HTML, відео);

– адмініструвати публічні та приватні форуми;

– реалізовувати технології індивідуалізованого навчання;

– розробляти для студентів онлайн-курси;

– керувати процесом навчання;

– виконувати розсилки по електронній почті;

– проводити чати;

– переглядати статистику активності користувачів;

– використовувати технологію wiki для спільної роботи над проектами.

Використовується адаптована платформа навчання Claroline на кафедрі медичної інформатики у Львівському національному медичному університеті імені Данила Галицького.

**Прометей** (www.prometeus.ru) – програмна оболонка, що забезпечує процес навчання та тестування студентів, а також надає змогу керувати всією діяльністю навчального закладу. У системі реалізовані такі автоматизовані функції: управління навчальним процесом; розподіл прав доступу до навчальних ресурсів і засобів управління системою; розмежування взаємодії учасників процесу навчання; ведення журналів активності користувачів навчального закладу; навчання та оцінка знань в середовищі Інтернет, в корпоративних та локальних мережах. Дана система має модульну структуру, тому легко розширюється, модернізується та масштабується.

СПН «Прометей» складається з наступних модулів:

– типовий web-вузол – набір HTML-сторінок, що надають відомості про навчальний заклад, список курсів та дисциплін, список фасилітаторів;

– модуль «Адміністратор» забезпечує виконання адміністратором

управління системою, реєструє нових користувачів, формує навчальні групи;

– модуль «Планування навчального процесу» надає можливість створювати план-графіки вивчення курсу, що включають заходи різних типів;

– модуль «Бібліотека» надає можливість зберігати навчальні посібники в довільному файловому форматі, закріплювати їх за певними курсами, виконувати повнотекстовий пошук, передавати по мережі відео- та аудіо файли, зберігати статистику звернень слухачів;

– модуль «Тестування» реалізує перевірку навчальних досягнень студентів в режимах самоперевірки, тренінгу та екзамену;

– модуль «Спілкування» забезпечує різноманітні засоби спілкування між учасниками процесу навчання.

Розробники СПН «Прометей» визначають такі переваги платформи [280]:

– дружній інтерфейс, простота опанування та експлуатації;

– відсутність ліцензій на клієнтські місця;

– можливість використання методики онлайн-навчання, що базується на командній роботі;

– висока продуктивність і масштабованість в міру збільшення числа користувачів і навантаження;

– 10 видів тестів, а також можливість використання графіки та мультимедіа в тестах;

– можливість побудови додаткових звітів;

– можливість об’єднання декількох систем в єдине освітнє середовище;

– помірні вимоги до ресурсів сервера та клієнтських місць.

У системі освіти України платформу «Прометей» використовують Національна академія державного управління при Президентові України, Київський обласний інститут післядипломної освіти педагогічних кадрів.

*Л Анкета з проблем та перспектив вищої математичної освіти*

1. Яких математичних дисциплін Ви навчаєте?

| *Дисципліна* | *Так / Ні* |
|---|---|
| Вища математика | |
| Теорія ймовірностей | |
| Диференціальні рівняння | |
| Математичне програмування | |
| Математична статистика | |
| Аналітична геометрія | |
| Математичний аналіз | |
| Дискретна математика | |
| Лінійна алгебра | |
| Чисельні методи | |
| Математична логіка | |
| Дослідження операцій | |
| Методи оптимізації | |
| Інші дисципліни | |

2. У чому, на Вашу думку, полягає мета вищої математичної освіти? (оцініть рівень вагомості: 1 – дуже низький, 2 – низький, 3 – достатній, 4 – середній, 5 – високий)

| *Мета* | *Рівень вагомості* |
|---|---|
| Формування наукового світогляду | |
| Формування наукового підходу до розв’язування реальних задач | |
| Формування загальнолюдської культури | |
| Формування математичної культури | |
| Формування інформаційної культури | |
| Інтелектуальний розвиток особистості | |
| Підготовка до життя у суспільстві | |
| Формування соціально зрілої, творчої особистості | |
| Забезпечити прикладну спрямованість змісту вищої математичної освіти | |
| Ваш варіант | |

3. Що Ви відповідаєте студентам на питання: «Навіщо мені потрібна математика?» (оцініть рівень вагомості: 1 – дуже низький, 2 – низький, 3 – достатній, 4 – середній, 5 – високий)

| *Відповідь* | *Рівень вагомості* |
|---|---|
| Розвиває логічне мислення | |
| Підвищує рівень вміння розробляти алгоритми розв’язування практичних задач | |
| Підвищує рівень математичної культури | |
| Підвищує рівень культури взагалі | |
| Розвиває вміння навчатися | |
| Допомагає систематизувати знання з теорії і методів розв’язування практичних задач | |
| Формує практичні навички щодо побудови математичних моделей реальних задач | |
| Це необхідний інструмент для розв’язування задач в різних сферах діяльності фахівців з вищою освітою | |
| Одержані знання і вміння з вищої математики сприяють вивченню профільних дисциплін. | |
| Одержати загальну фундаментальну освіту з одночасною умовою формування базового рівня професійної компетентності майбутнього спеціаліста. | |
| Ваш варіант | |

4. Чи влаштовує Вас зміст математичних дисциплін, який визначений стандартами вищої освіти (проектами стандартів)? (підкресліть).

**Влаштовує** **Влаштовує частково** **Не влаштовує**

5. Вкажіть, що Ви пропонуєте для покращення змісту математичних дисциплін?

| *Пропозиція* | *Так / Ні* |
|---|---|
| Дозволити працювати викладачам за авторськими програмами | |
| Збільшити варіативну частину в програмі дисципліни | |
| Зміст потрібно оновити за рахунок останніх наукових досягнень у математиці та її застосуваннях | |
| Більше приділяти уваги практичному застосуванню математичних дисциплін | |
| Ваш варіант | |

6. Які види навчальної діяльності Ви використовуєте при навчанні математичних дисциплін?

| *Вид навчальної діяльності* | *Так / Ні* |
|---|---|
| Лекційні заняття | |
| Практичні заняття | |
| Лабораторні заняття із використанням комп’ютера | |
| Індивідуальні заняття | |
| Організація самостійної роботи студентів | |
| Консультації | |
| Курсові роботи | |
| Наукова робота в проблемних групах | |
| Ваш варіант | |

7. Які форми поточного контролю Ви використовуєте при навчанні математичних дисциплін?

| *Форма контролю* | *Так / Ні* |
|---|---|
| Усне опитування | |
| Математичний диктант | |
| Самостійна робота | |
| Контрольна робота | |
| Модульний контроль | |
| Автоматизований контроль | |
| Колоквіум | |
| Тестування | |
| Комп’ютерне тестування | |
| Взагалі не використовую | |
| Ваш варіант | |

8. Якій системі оцінювання навчальних досягнень студентів Ви надаєте перевагу?

| *Система оцінювання* | *Так / Ні* |
|---|---|
| Чотирибальна | |
| Дванадцятибальна | |
| Багатобальна | |
| Рейтингова | |
| Ваш варіант | |

9. Якій формі підсумкового контролю з математичних дисциплін Ви надаєте перевагу?

| *Форма підсумкового контролю* | *Так / Ні* |
|---|---|
| Усний екзамен | |
| Письмовий екзамен | |
| Залік | |
| Диференційований залік | |
| Екзамен і залік | |
| За результатами поточного контролю (за рейтингом, модульним контролем тощо) | |
| Ваш варіант | |

10. Ваше ставлення до модульно-рейтингової системи організації та оцінювання навчальної діяльності студентів при навчанні математичних дисциплін?

| *Варіант відповіді* | *Так / Ні* |
|---|---|
| Ніколи не чув про таку систему | |
| Знаю, але не застосовую | |
| Вважаю, що дана система оцінювання не потрібна при навчанні математичних дисциплін | |
| Знаю і застосовую частково (для окремих дисциплін) | |
| Знаю і застосовую частково (для окремих студентських груп і спеціальностей) | |
| Працюю за цією системою | |
| Ваш варіант | |

11. Які проблеми при навчанні математичних дисциплін виникають у Вашій професійній діяльності найчастіше?

| *Проблема* | *Так / Ні* |
|---|---|
| Низький рівень підготовки студентів зі шкільної математики | |
| Недостатній рівень практичних умінь та навичок щодо використання теоретичних знань | |
| Низька мотивація студентів при вивченні предметів математичного циклу | |
| Недостатній рівень навчально-пізнавальної активності студентів | |
| Невміння студентів самостійно працювати з навчальним матеріалом | |
| Невміння студентів застосовувати математичні знання для формалізації практичних задач та їх розв'язування | |
| Ваш варіант | |

12. Які шляхи подолання зазначених у п. 10 проблем Ви пропонуєте і впроваджуєте в своїй роботі? (оберіть варіант чи запропонуйте свій)

| *Шлях подолання проблем при навчанні математичних дисциплін* | *Упроваджую* |
|---|---|
| Стимулювання мотивації, підвищення інтересу до навчання | |
| Розвиток мислення, інтелектуальних здібностей студентів | |
| Індивідуалізація та диференціація навчання | |
| Розвиток самостійності | |
| Надання переваги активним методам навчання і діяльнісному підході | |
| Підвищення наочності навчання | |
| Збільшення арсеналу засобів пізнавальної діяльності, опанування сучасними методами наукового пізнання, пов’язаними із застосуванням інформаційних і комунікаційних технологій (ІКТ) | |
| Проведення лабораторних робіт при навчанні математичних дисциплін з використанням ІКТ | |
| Створення методичних і дидактичних матеріалів, зокрема мультимедійних | |
| Розширення доступу до освітніх та наукових інформаційних ресурсів через Internet | |
| Застосування інноваційних педагогічних технологій | |
| Ваш варіант | |

13. Чи потрібно під час навчання вищої математики формувати елементи інформаційної культури студентів? (підкресліть)

**Так** **Ні**

14. Чи використовуєте Ви інформаційно-комунікаційні технології у своїй професійній діяльності? (оберіть варіант відповіді)

| *Так (вкажіть, для чого)* | | *Ні (вкажіть причини)* | |
|---|---|---|---|
| Для створення методичних та дидактичних матеріалів з дисципліни, зокрема мультимедійних | | Не маю доступу до комп’ютера | |
| Як джерело відомостей через Internet | | Не вмію працювати з комп’ютером | |
| Для вимірювання навчальних досягнень студентів (комп’ютерне тестування, автоматизований контроль) | | Не вважаю, що комп’ютер може допомогти при навчанні математичних | |

| *Так (вкажіть, для чого)* | | *Ні (вкажіть причини)* | |
|---|---|---|---|
| На заняттях з математичних дисциплін, як інструмент розв'язування задач | | Потрібне підвищення кваліфікації з використання ІКТ у навчальному процесі | |
| Як засіб дистанційного навчання математичних дисциплін | | У студентів низький рівень інформаційної культури | |
| Для активізації самостійної роботи студентів | | Відсутні умови для використання ІКТ у навчальному процесі | |
| Ваш варіант | | Ваш варіант | |

15. Оцініть рівень Вашого володіння деякими системами комп'ютерної математики та опрацювання даних, а також використання їх у навчальному процесі (0 – невідома, 1 – знаю, але не застосовую в роботі, 2 – знаю і застосовую для своєї роботи, 3 – знаю і застосовую на заняттях з математичних дисциплін).

| *Системи комп'ютерної математики та опрацювання даних* | *Рівень ознайомлення і використання* |
|---|---|
| Mathcad | |
| Matlab | |
| Maple | |
| Mathematica | |
| GAUSS | |
| MuPAD | |
| Derive | |
| Gran1, Gran-2D, Gran-3D | |
| Advanced Grapher | |
| Dynamic Geometry (DG) | |
| Cabri | |
| Excel | |
| Statistica | |
| SPSS | |
| Ваш варіант | |

16. Чи використовуєте Ви програмно-педагогічні засоби у процесі навчання студентів вищої математики? Вкажіть рівень Вашого використання їх у навчальному процесі (0 – невідомо, 1 – знаю, але не застосовую в роботі, 2 – знаю і застосовую для своєї роботи, 3 – знаю і застосовую на заняттях з математичних дисциплін).

| *Програмно-педагогічні засоби* | *Рівень* |
|---|---|
| електронні програмно-методичні комплекси | |
| електронні підручники | |
| електронні довідники | |
| електронні задачники | |
| електронні тренажери | |

17. Чи використовуєте Ви системи підтримки навчання у процесі навчання студентів вищої математики? Вкажіть рівень Вашого використання їх у навчальному процесі (0 – невідомі, 1 – знаю, але не застосовую в роботі, 2 – знаю і застосовую для своєї роботи, 3 – знаю і застосовую при навчанні математичних дисциплін).

| *Системи підтримки навчання* | *Рівень* |
|---|---|
| Moodle | |
| TrainingWare Class | |
| Claroline LMS | |
| Прометей | |
| Агапа | |
| ATutor | |
| Blackboard | |
| CCNet | |
| Chamilo | |
| Claroline | |
| Desire2Learn | |
| eFront | |
| ILIAS | |
| metacoon | |
| OLAT | |
| Sakai Projec | |
| WebCT | |
| SharePointLMS | |
| JoomlaLMS | |
| Pass-port | |
| Yacapaca | |
| CampusCE | |
| Ваш варіант | |

18. Чи використовуєте Ви мобільні математичні середовища у процесі навчання студентів вищої математики? Вкажіть рівень Вашого використання їх у навчальному процесі (0 – невідомі, 1 – знаю, але не застосовую в роботі, 2 – знаю і застосовую для своєї роботи, 3 – знаю і застосовую при навчанні математичних дисциплін).

| *Мобільне математичне середовище* | *Рівень* |
|---|---|
| Лекційні демонстрації | |
| Динамічні моделі | |
| Тренажери | |
| Навчальні експертні системи | |
| Ваш варіант | |

19. Якщо Ви використовуєте інші засоби при викладанні вищої математики, вкажіть які.

*М Створення сайту за допомогою Google*

Перед створенням сайту необхідно створити обліковий запис Google на https://www.google.com/accounts/NewAccount. Щоб створити свій сайт, треба зайти на sites.google.com, вибрати акаунт, під яким створюється сайт, і натиснути кнопку «Створити сайт» (рис. М.1).

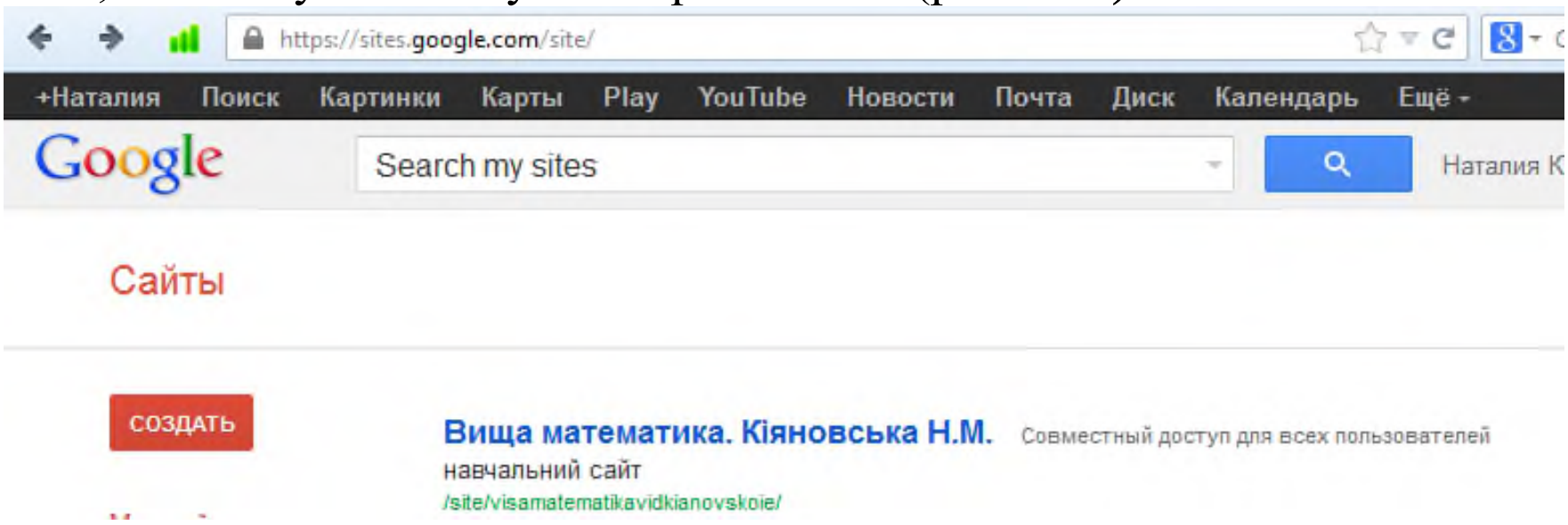

Рис. М.1. Створення сайту на sites.google.com

Після цього необхідно дати своєму сайту ім’я та вибрати рівень доступу до сайту. Вибравши тему сайту, вводиться код в поле для коду та натискається кнопка «Створити» (рис. М.2).

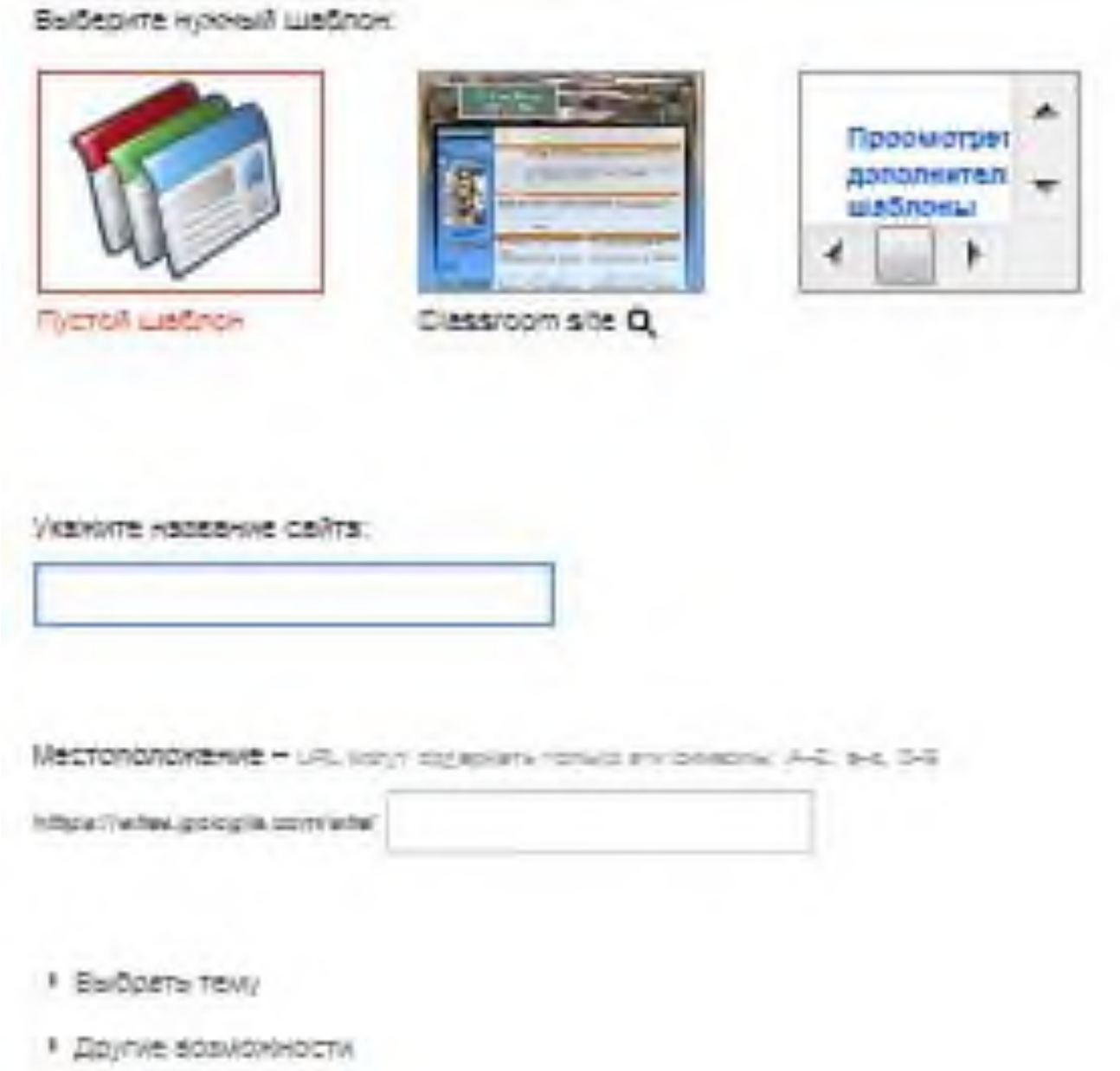

Рис. М.2. Вибір параметрів сайту

Після створення сайту можна додавати текст, документи, відео або фотографії. Для цього треба натиснути «Редагувати сторінку» – «Вставити», в меню вибрати необхідний елемент для вставки: фотографія, календар, документ, презентація, Picasa слайдшоу або відео. Після натискання кнопки «Вибрати», з'являється нове вікно, в якому можна вказати параметри документа, такі як висота, ширина і заголовок. Натиснувши «Зберегти», з'являється порожній шаблон документа. Натиснувши на нього, з'явиться можливість вирівняти документ. Тут же можна включити і обтікання текстом. Після натискання «Зберегти» вгорі сторінки, з'являється вибраний документ на сайті.

Якщо на сторінці сайту виникає необхідність у створенні *посилань*, треба натиснути «Редагувати сторінку» в лівому верхньому кутку, виділити слово, до якого додається посилання, та натисніть кнопку «Посилання» на панелі редагування сторінки. В якості посилання можна вибрати вже існуючу сторінку сайту або будь-яку іншу Web-сторінку.

Для створення *нової сторінки* сайту, треба натиснути кнопку «Створити нову сторінку», вибрати тип сторінки і її рівень.

На сайтах викладачів бажано забезпечити інтерактивне спілкування зі студентами, для цього необхідно вбудувати *форум* на створений сайт. Для того, щоб вбудувати обговорення з груп Google на свій сайт на Сайтах Google за допомогою гаджета форуму необхідно виконати такі дії:

– створити або відкрити вбудований форум груп Google і сайт Google, де буде вбудовано цей форум;

– перейти на сторінку, що буде містити форум;

– натиснути кнопку «Редагувати сторінку»;

– в меню «Вставити» вибрати пункт «Додаткові гаджети»;

– у вікні «Налаштувати гаджет» вибрати пункт «Додати гаджет з URL»;

– у полі «Введіть URL гаджета, який необхідно додати» вказати адресу http://www.gstatic.com/sites-gadgets/forum/forum_content.xml;

– натиснути кнопку «Додати»;

– у вікні «Налаштувати гаджет» виконайте такі дії: 1) у полі «Форум за замовчуванням» ввести назву форуму. (Наприклад, якщо URL форуму http://groups.google.com/group/samplegroup, то назвою групи є samplegroup); 2) не заповнювати поле «Домен, зміст якого необхідно відобразити»; 3) у розділі «Відобразити» задати висоту і ширину гаджета, а також параметри його відображення;

– натиснути кнопку «Переглянути гаджет»;

– якщо влаштовує вигляд гаджета форуму, натиснути кнопку «ОК». Щоб внести зміни, натиснути кнопку «Назад до налаштувань».

Якщо є декілька гаджетів форуму і необхідно, щоб користувачі могли переглядати їх одночасно, можна скористатися гаджетом змісту.

*Н Аналіз, особливості та вимоги до методів навчання: методу проектів, різнорівневого навчання, кейс-методу, методу навчання у групах, навчання у співробітництві та методу портфоліо*

На заняттях із вищої математики неможливо обійтись без репродуктивних методів навчання. Використання цих методів допомагає сформувати у студентів усталених умінь та навичок, необхідних їм у подальшій професійній діяльності. Вміння відтворювати повідомлені викладачем відомості та виконувати операції за зразком сприяють формуванню у студентів базових, фундаментальних знань.

Але завдання вищої технічної освіти полягає у підготовці творчої особистості, здатної працювати у швидкозмінних умовах. Тому використання продуктивних методів навчання: методу проектів, різнорівневого навчання, кейс-методу, методу навчання у групах, навчання у співробітництві, методу портфоліо, методу проблемного навчання, евристичного методу, методу дослідницького навчання тощо сприяє розвитку особистості, здатної до швидкого прийняття правильних рішень.

Розглянемо використання зазначених методів на прикладі розділу вищої математики «Невизначений та визначений інтеграли».

*Метод проектів.* Є. С. Полат [251] в основу методу проектів покладає розвиток пізнавальних навичок студентів, вмінь самостійно будувати свої знання, вміння орієнтуватися в інформаційному просторі, розвиток критичного і творчого мислення. Основою методу проектів є досягнення дидактичної мети через детальну розробку проблеми (технологію), що повинна завершитися цілком реальним, відчутним практичним результатом, оформленим зазначеним способом. В основу методу проектів покладена ідея, що складає суть поняття «проект», його прагматична спрямованість на результат, який можна отримати при вирішенні тієї чи іншої практично або теоретично значущої проблеми. Цей результат можна побачити, осмислити, застосувати в реальній практичній діяльності. Для досягнення такого результату, необхідно навчити студентів самостійно мислити, знаходити і вирішувати проблеми, залучаючи для цієї мети знання з різних областей, уміння прогнозувати результати і можливі наслідки різних варіантів розв'язання, уміння встановлювати причинно-наслідкові зв'язки.

Метод проектів завжди орієнтований на самостійну діяльність студентів – індивідуальну, парну, групову, яку студенти виконують протягом певного відрізка часу. Цей метод органічно поєднується з груповими (collaborative or cooperative learning) методами. Метод проектів завжди припускає розв'язання якоїсь проблеми. Розв'язання проблеми передбачає, з одного боку, використання сукупності, різноманітних методів, засобів навчання, а з іншого, передбачає необхідність інтегрувати знання та уміння з різних галузей науки, техніки, технології, творчих областей з метою отримання нового результату діяльності. Якщо говорити про метод проектів як про педагогічну технологію, то ця технологія передбачає сукупність дослідницьких, пошукових, проблемних методів, творчих за самою своєю суттю.

Основні вимоги до використання методу проектів є [251]:

1) наявність значущої в дослідницькому, творчому плані проблеми, що вимагає інтегрованого знання, дослідницького пошуку для її вирішення;

2) практична, теоретична, пізнавальна значущість передбачуваних результатів;

3) самостійна (індивідуальна, парна, групова) діяльність студентів;

4) структурування змістовної частини проекту (із зазначенням поетапних результатів);

5) використання дослідницьких методів, що передбачають певну послідовність дій: визначення проблеми та похідних від неї завдань дослідження; висунення гіпотез їх вирішення; обговорення методів дослідження; обговорення способів оформлення кінцевих результатів; збір, систематизація та аналіз отриманих даних; підведення підсумків, оформлення результатів, їх презентація; висновки, висування нових проблем дослідження.

На заняттях з вищої математики при вивченні теми «Застосування визначеного інтегралу» групу студентів доцільно розбити на декілька підгруп і кожній підгрупі видати завдання: охарактеризувати задану лінію (або декілька ліній: лемніскату Бернуллі, спіраль Архімеда, логарифмічну спіраль, гіперболічну спіраль, кардіоїду, ін.), знайти довжину, дослідити як буде змінюватися довжина, якщо змінювати деякі параметри (рис. Н.1), знайти об’єм, якщо лінію почати обертати навколо вказаної вісі. Результат роботи групи подати у вигляді презентації.

Під *різнорівневим навчанням* М. Ю. Бухаркіна [153] розуміє таку організацію навчально-виховного процесу, при якій кожен студент має можливість опановувати навчальний матеріал на різному рівні, не нижче базового, залежно від його здібностей та індивідуальних особливостей особистості, при цьому за критерій оцінки діяльності студента приймаються його зусилля з оволодіння цим матеріалом, творче його застосування.

В. Ю. Колісник [217] зазначає, що різнорівневий підхід передбачає роботу в межах однієї групи, але процес навчання здійснюється за допомогою завдань різного рівня складності. Основні принципи цього підходу: 1) попереднє визначення психолого-педагогічної готовності до навчання; 2) залучення всіх студентів окремої групи; 3) можливість працювати невеликими групами; 4) подолання особистісних проблем у процесі навчання.

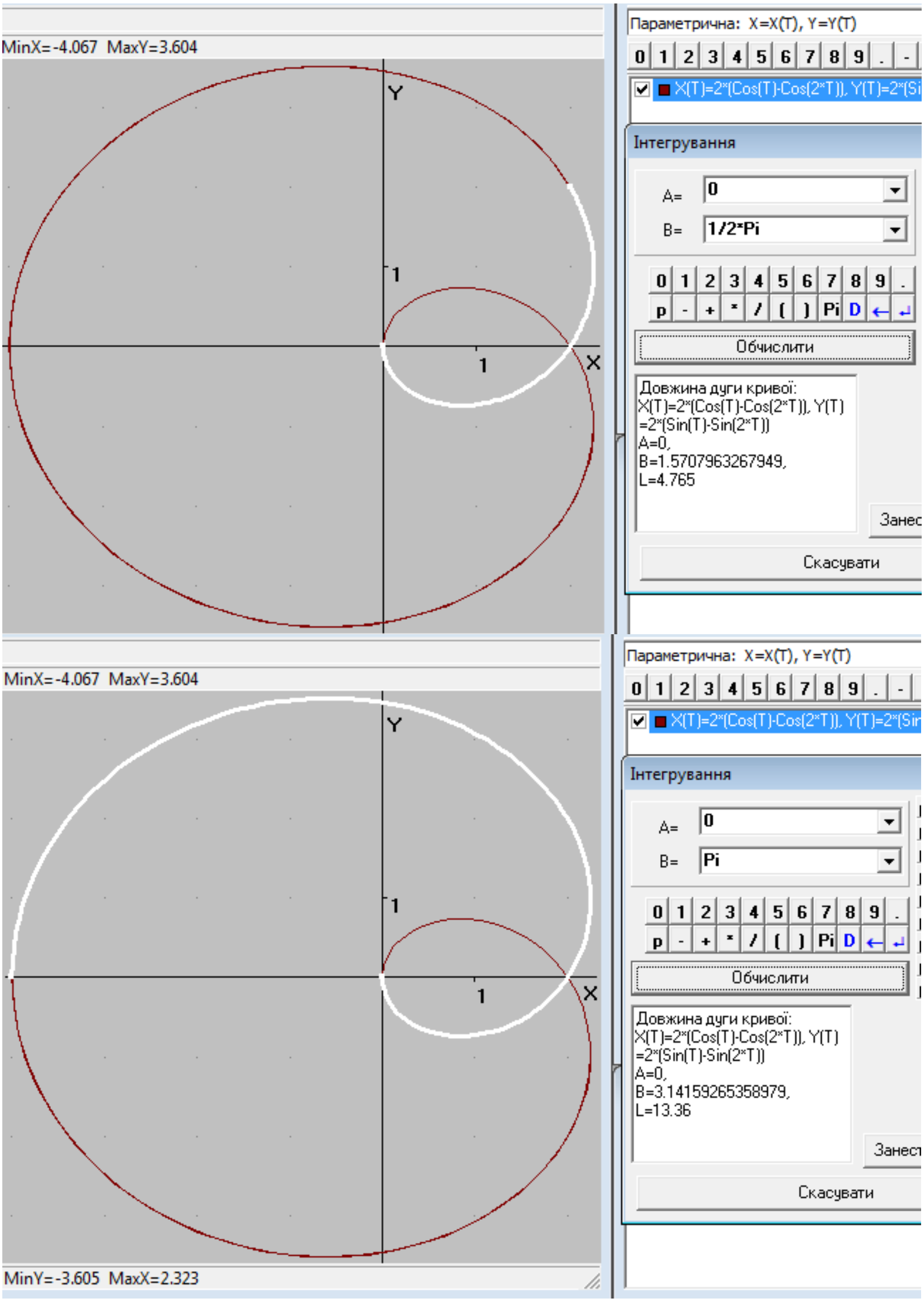

Рис. Н.1. Дослідження зміни довжини кардіоїди із зміною параметру із використанням GRAN 1

При вивченні теми «Застосування визначеного інтегралу» групу доцільно розділити на три частина: сильна підгрупа, середня та проблемна підгрупи. Теоретичний матеріал доцільно представити у вигляді презентації, після чого сильній підгрупі видати завдання на застосування визначеного інтегралу, а з іншою частиною студентів продовжувати працювати далі. Разом із середньою та проблемною підгрупами необхідно розглянути приклади з теми і після цього середній підгрупі видати завдання для самостійного виконання. Із проблемною підгрупою продовжувати розв'язувати завдання із чітким поясненням ходу розв'язання задач.

*Кейс-метод.* Кейс (з англ. «case» – випадок, ситуація) – це розбір ситуації або конкретного випадку, ділова гра. Він може бути названий технологією аналізу конкретних ситуацій, «окремого випадку». Суть технології полягає в тому, що в основі його використовуються описи конкретних ситуацій або випадків. Представлений для аналізу випадок повинен:

– відображати реальну життєву ситуацію;

– в описі має бути присутня проблема або ряд прямих або непрямих утруднень, протиріч, прихованих завдань для вирішення дослідником;

– потрібне оволодіння попереднім комплексом теоретичних знань для залучення їх у вирішенні конкретної проблеми або ряду проблем.

У процесі роботи над кейсом часто потрібні додаткові матеріали. В результаті студенти знаходять власні висновки, розв'язки проблемної ситуації, і часто, у вигляді неоднозначних множинних рішень [142].

О. К. Ільїна [181] вважає, що одним з основних завдань викладача, який використовує кейс-метод, є залучення студентів до аналізу, обговорення та вирішення проблеми. Для цього важливо виконання двох умов: матеріал кейса повинен представляти для студентів професійний інтерес і передбачати можливість особистого внеску студента в свою освіту і в освіту своєї «команди». Цікавий матеріал і можливість застосування професійних знань стимулює участь у дискусії. Бажання вирішити проблему спонукає студентів не просто прочитати кейс, але ретельно його вивчити, оволодіти фактами і деталями.

Процес створення кейсу складається з декількох етапів [278, 46-47]:

– в першу чергу необхідно сформувати цілі кейса, визначити основні проблеми та питання, які будуть перебувати в центрі уваги кейса;

– після визначення загального напряму наступає етап побудови програмної карти кейса, що складається з основних тез, які необхідно втілити;

– збір даних щодо тез програмної карти кейса. Побудова або вибір моделі ситуації;

– далі необхідно визначитися із загальною структурою кейса. Данні в ньому може подаватися в певній послідовності або згідно деякої моделі або схеми;

– написання тексту кейсу;

– діагностика правильності та ефективності кейса; проведення навчального експерименту, побудованого за тією чи іншою схемою, для з'ясування ефективності даного кейса;

– підготовка остаточного варіанту кейса;

– впровадження кейса в практику навчання, його застосування при проведенні навчальних занять, а також його публікацію з метою поширення у викладацькому співтоваристві;

– підготовка методичних рекомендацій з використання кейса: розробка завдання для студентів та можливих питань, описання передбачуваних дій студентів і викладача в момент обговорення кейсу.

На заняттях з вищої математики студентам можна запропонувати згадати залежність роботи сталої сили від величини сили та довжини шляху. Застосовуючи ці поняття до заданої ситуації, отримують: $A = \int_a^b F(x)dx$.

Ці теоретичні відомості використовуються ними надалі під час розв'язування задач.

Наприклад, задача: Обчислити роботу, що витрачається на стискання газу в циліндрі з розмірами $R$ та $H$.

Під час обговорення студенти приходять до висновку: для обчислення роботи, що витрачається на стискання газу в циліндрі радіуса $R$ та висотою $H$, потрібно скористатися рівнянням стану газу $P_1V_1 = P_2V_2$.

На початку процесу тиск в циліндрі $P_0$. Об'єм циліндра $V_0 = \pi R^2H$. Позначаючи $P(x)$ тиск газу в циліндрі при переміщені поршня на відстань $x$ від початкового положення, студенти знаходять об'єм частини циліндру з газом $V(x) = \pi R^2(H - x)$. Тому $P_0V_0 = P(x)V(x)$, $P_0\pi R^2H = P(x)\pi R^2(H - x)$.

З останнього рівняння отримують силу тиску $P(x) = P_0H/(H - x)$.

Отже, студенти приходять до висновку, що, за умови переміщення поршня на відстань $dx$, витрачається елементарна робота $dA = P(x)dx = P_0H/(H - x)dx$. З огляду на це, роботу, що потрібно витратити при стисканні газу, можна обчислити, інтегруючи останню рівність: $A = \int_0^a \frac{P_0H}{H-x}dx = P_0H\int_0^a \frac{dx}{H-x}$, де $a$ – відстань, на яку перемістився поршень.

Наведений приклад охоплює не одну, а кілька фізичних підзадач, що дають можливість студентам засвоїти методи використання визначеного інтеграла в розв'язуванні задач електротехніки, теоретичної механіки та інших спецдисциплін.

Перевірку обчислень доцільно проводити із використанням СКМ. Аналогічну задачу можна запропонувати студентам для самостійної роботи.

*Метод навчання у групах.* Мета цієї технології полягає у формуванні вмінь ефективно працювати спільно в тимчасових командах і групах та отримувати якісні результати. Це така організація занять, в ході яких у студентів формуються інформаційно-комунікативні компетентності, розвиваються розумові здібності в результаті вирішення проблемної ситуації, підготовленої викладачем. Робота студентів будується навколо ключових проблем, виділених викладачем.

Навчаючись у групах, студенти розвивають здібності організовувати спільну діяльність, засновану на принципах співпраці. При цьому у них формуються такі особистісні якості, як толерантність до різних точок зору і поведінки, відповідальність за результати роботи, формується вміння поважати чужу точку зору, слухати партнера, вести ділове обговорення, досягати згоди в конфліктних ситуаціях та спірних питаннях [278, 49].

Виділивши підгрупи із 5-7 студентів, викладач кожній підгрупі видає завдання з теми. Так при вивченні теми «Методи інтегрування невизначених інтегралів», кожній підгрупі видається по 10-14 прикладів на різні методи. Завдання студентів полягає у тому, щоб за певний час розподілити між собою приклади, розв’язати їх та зробити перевірку із використанням СКМ або Web-СКМ (рис. Н.2). Та група студентів, що найшвидше розв’язала всі приклади правильно, отримує додаткові бали на модульній контрольній роботі.

*Навчання у співробітництві.* Суть використання особистісно-орієнтованого підходу і педагогічної технології «навчання у співробітництві» полягає в особистій участі кожного студента у виконанні спільного завдання в залежності від його можливостей і особистих уподобань. Використання педагогічної технології «навчання у співробітництві» у процесі навчання полягає у підвищенні мотивації і стимулюванні інтересу студентів до навчання із використанням різноманітних можливостей проектної технології.

Використання методу навчання у співробітництві сприяє кращій підготовці студентів до модульної роботи. На такому занятті групу необхідно розділити на підгрупи (групу з 20 студентів максимум на 4 підгруп). В кожній групі обрати доповідача, опонента (робить зауваження до доповіді), рецензента (знаходить позитивні моменти доповіді), контролюючого (робить перевірку відповіді за допомогою СКМ (рис. Н.3)). Кожній групі пропонуються приклади з розділу «Визначений інтеграл», в кожному з яких містяться різні методи інтегрування (це підсумкове заняття з цього розділу і завдання повинні бути ускладнені).

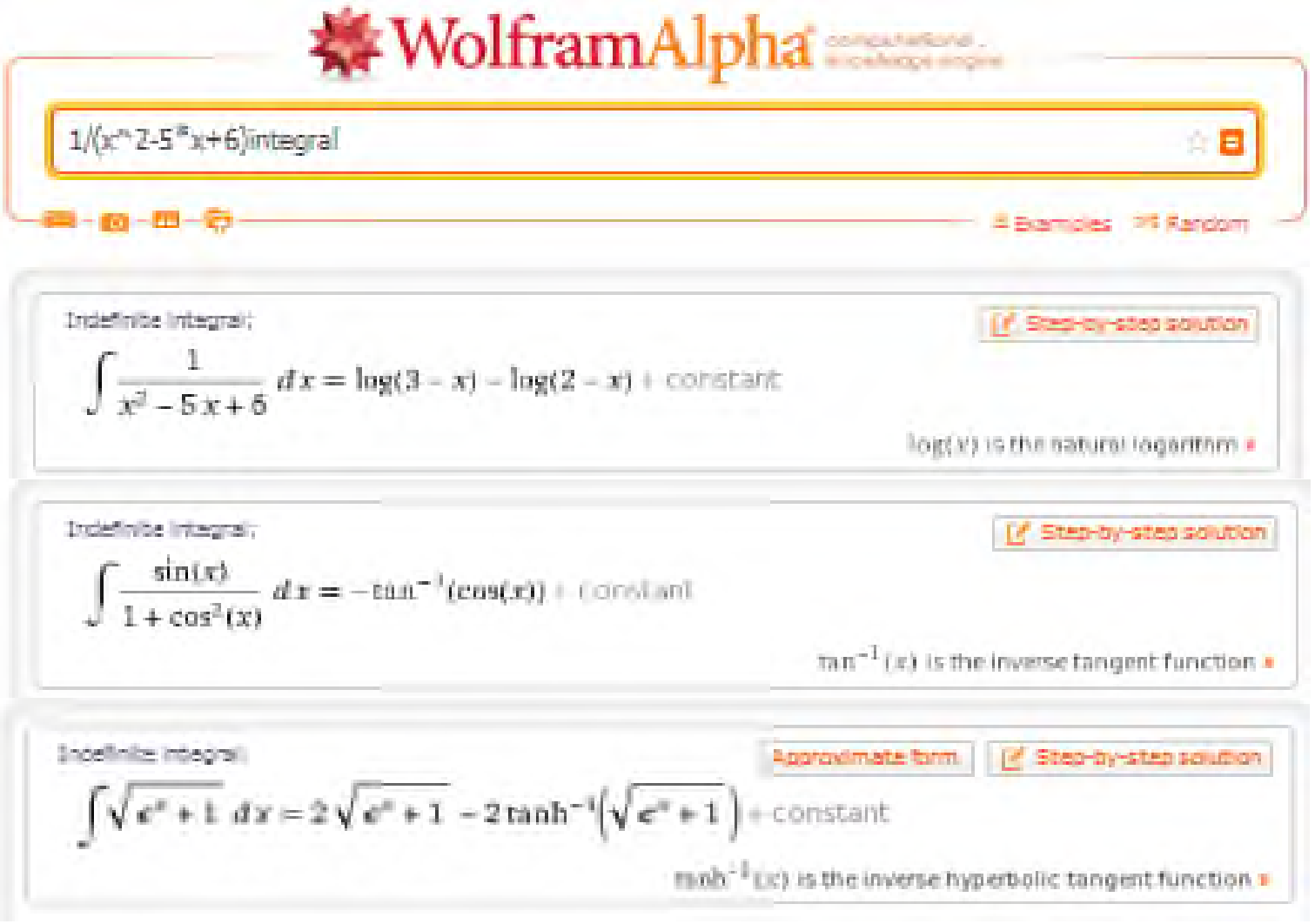

Рис. Н.2. Обчислення невизначених інтегралів із використанням Wolfram|Alpha

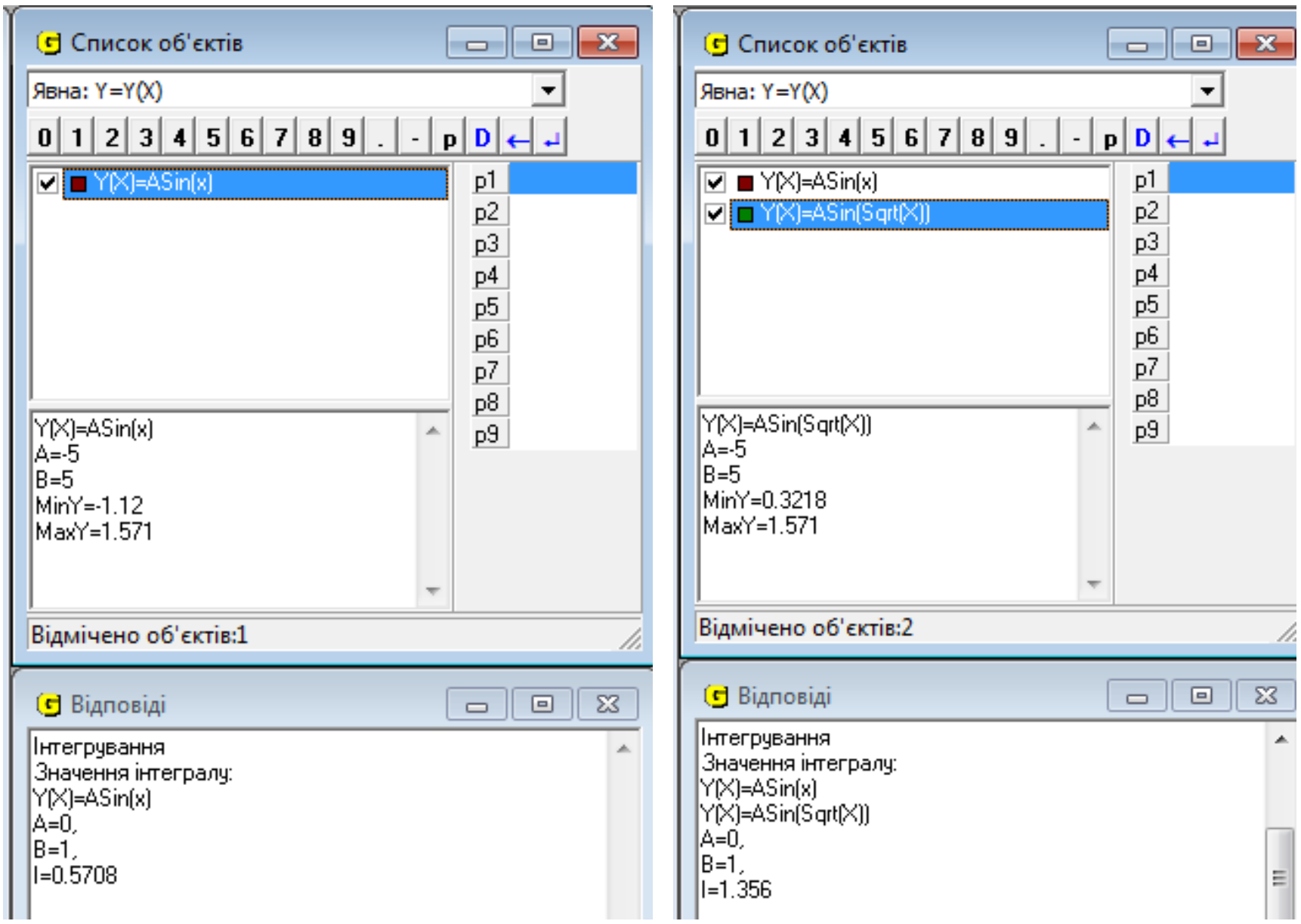

Рис. Н.3. Обчислення визначеного інтегралу за допомогою GRAN 1

Завдання всім підгрупам видаються однакові, дається час подумати, після цього виходить до дошки студент однієї з груп і розв'язує приклад, а його опоненти починають задавати питання.

Подальший розвиток особистісно-орієнтованого підходу до навчання і педагогічної технології «навчання у співробітництві» як його складової призвів до появи так званого портфеля. Головною метою такого навчання став розвиток інтелектуальних і творчих здібностей студентів, готових до самореалізації, самостійного мислення, прийняття важливих для себе рішень [249].

*Метод портфоліо.* Портфоліо з італійської – портфель. О. Г. Смолянінова [301] визначає метод портфоліо як сукупність окремих завдань, виконаних студентами, зібраних для певної мети, наприклад, файл закінчених мультимедіа-проектів. Оцінка за допомогою портфоліо справжня і рефлексивна. Студент збирає матеріали або дані з різних джерел, використовуючи різні методи, в різних часових рамках. Це означає, що зміст портфоліо може включати матеріали у вигляді рисунків, фотографій, відео або звуків, рукописні матеріали або інші зразки робіт, комп'ютерних дисків і копій, стандартизованих або визначених програмою випробувань. Крім того, портфоліо може містити характеристики студента та його самооцінку.

Етапи створення портфоліо [301]:

– *накопичення, збір даних.* Студенти зберігають тільки ті дані (підтвердження компетентності), що представляють їхні успіхи день у день в процесі навчання;

– *відбір даних.* Перегляд і оцінка зібраних даних і визначення тих, що демонструють досягнення згідно вибраним критеріям і зовнішнім вимогам;

– *рефлексія.* Оцінка власного зростання, а також недоліків власного розвитку;

– *перспективна оцінка.* Порівняння власних позицій із зовнішніми вимогами і визначення цілей на майбутнє. Це сприяє професійному зростанню і стимулює потребу в безперервній освіті;

– *презентація.* Демонстрація портфоліо одноліткам, батькам, роботодавцям;

– *рефлексія після презентації портфоліо*, в результаті якої робиться висновок про необхідність його розвитку та покращення. На цьому етапі проводиться самостійна перевірка відповідності отриманих результатів власним очікуванням, формуються навички самооцінки.

Метод портфоліо надає можливість накопичувати матеріал, що свідчить про розвиток інформаційної і комунікативної компетентностей студентів, і моделювати ситуації професійної діяльності. Крім того, метод портфоліо допомагає студенту навчитися адекватно оцінювати власні досягнення і можливості, робити висновки про необхідність виправляти помилки і самовдосконалюватися.

Незважаючи на існування багатьох підходів до формування портфоліо, для будь-яких з них від студентів потрібно [301]:

– накопичувати матеріали для портфоліо;

– відбирати дані (підтвердження компетентності);

– рефлексувати;

– зберігати, представляти та каталогізувати елементи портфоліо.

На заняттях з вищої математики метод портфоліо доцільно використовувати при розв’язуванні задач на застосування визначених інтегралів. В аудиторному занятті доцільно навести ряд задач прикладного характеру з даної теми, а на наступне заняття дати завдання зробити презентацію розв’язання самостійно підібраної задачі.

*Евристичний метод*. Сутність частинно-пошукового (евристичного) методу навчання виражається в таких основних його ознаках:

– навчальні відомості студенти отримують самостійно, спираючись на свій досвід;

– викладач не пояснює новий матеріал, а спонукає студентів до самостійного його виведення;

– студенти самостійно розмірковують, розв’язують завдання, створюють і розв’язують проблемні ситуації, аналізують, порівнюють, роблять висновки, спираючись на чітке та коротке керівництво викладача.

Так при вивченні теми «Невизначений інтеграл» студентам надається таблиця інтегралів, що складається із основних інтегралів, отриманих за означенням первісної та таблиці похідних. Після вивчення методів інтегрування, а саме методу заміни змінної, студентам пропонується отримати самостійно таблицю інтегралів з описанням заміни змінної $t = kx + b$. Студенти також самостійно доповнюють таблицю інтегралів від таких тригонометричних функцій, як $\operatorname{tg} x$ та $\operatorname{ctg} x$. Перевірку отриманих результатів доцільно проводити із використанням СКМ.

*Метод дослідницького навчання* передбачає творче застосування набутих знань, оволодіння методами наукового пізнання, формування досвіду самостійного наукового пошуку [324].

Характерні ознаки цього методу такі:

– викладач разом зі студентами формулює проблему;

– нові знання не повідомляють, а студенти повинні самостійно здобути їх у процесі дослідження проблеми, порівняти різні варіанти відповідей, а також визначити основні засоби досягнення результатів;

– основною метою діяльності викладача є оперативне управління

процесом розв’язання проблемних завдань;

– навчання характеризується високою активністю, підвищеним інтересом студентів, а набуті знання є більш глибокими.

Оволодіння навчальним матеріалом може здійснюватись у процесі спостереження, пошуку висновків, під час роботи з книгою, письмової вправи з доведенням закономірності, практичних і лабораторних робіт.

Виконання дослідницького завдання передбачає такі етапи:

1) спостереження і вивчення фактів, виявлення суперечностей у предметі дослідження (постановка проблеми);

2) формулювання гіпотези щодо розв’язання проблеми;

3) побудова плану дослідження та його реалізація;

4) аналіз і систематизація одержаних результатів, формулювання висновків.

На заняттях з вищої математики метод дослідницького навчання доцільно використовувати при розв’язуванні задач прикладного характеру. Так при вивченні теми «Застосування визначених інтегралів» перед студентами можна поставити проблему обчислити об’єм кар’єру, що має форму еліптичного параболоїду. Студентам надається рівняння параболи і за допомогою визначеного інтегралу обчислити об’єм тіла обертання і результати порівняти із результатами, одержаними із використанням СКМ (рис. Н.4).

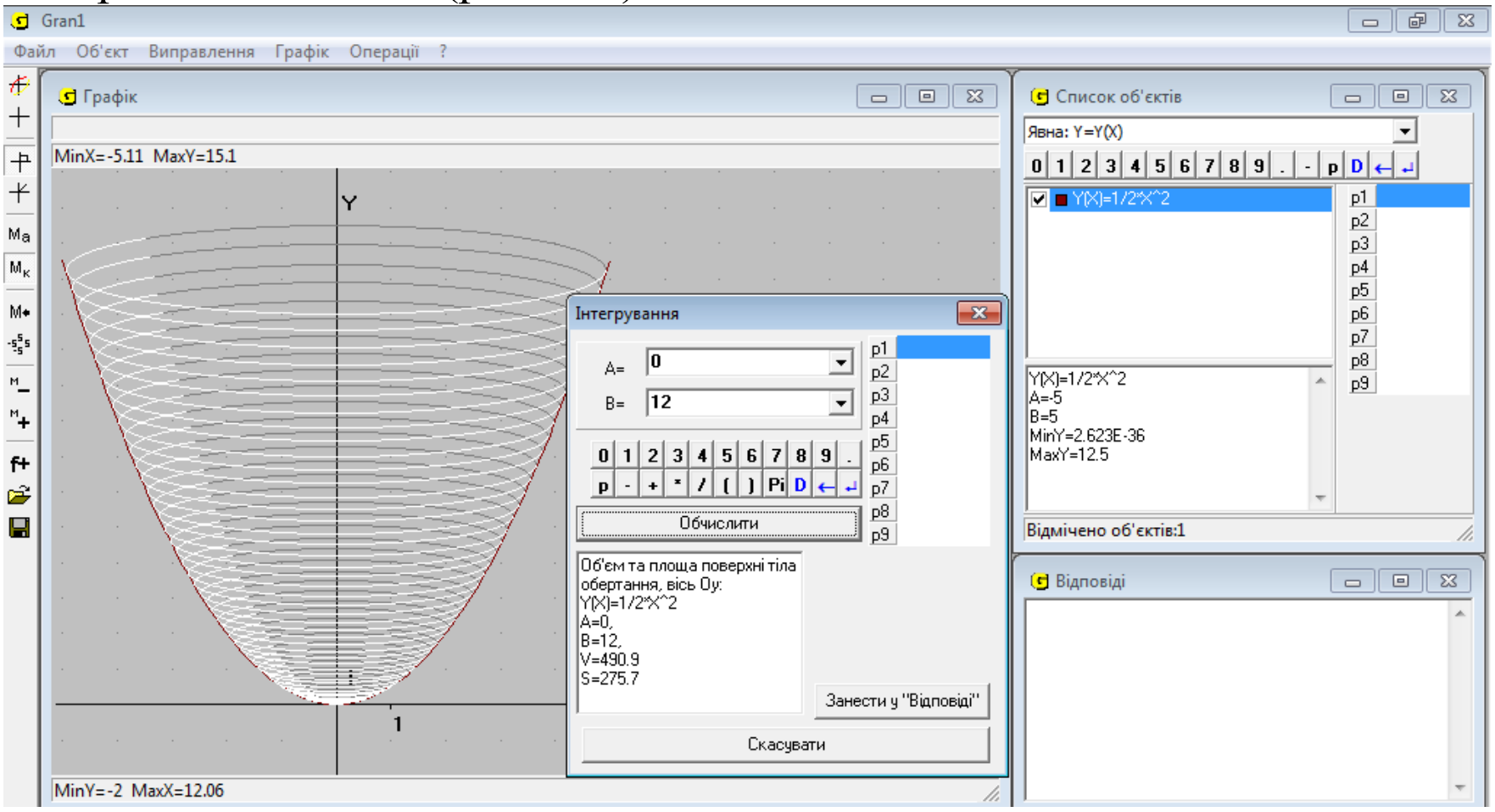

Рис. Н.4. Обчислення об’єму тіла обертання із використанням GRAN 1

*Метод проблемного навчання*. Під час вивчення розділу «Кратні інтеграли» поставити проблему: за допомогою відомого матеріалу з розділу «Визначений інтеграл» самостійно вдома отримати означення подвійного інтегралу, сформулювати основні властивості подвійного інтегралу, з'ясувати геометричний та фізичний змісти, визначити сферу його застосування, а отримані результати оформити у вигляді презентацій та надіслати на пошту викладача або винести на обговорення в Piazza. Після вивчення теми «Подвійний інтеграл» доцільно план лекції з теми «Потрійний інтеграл» надати на сайті, а лекцію провести в аудиторії за допомогою презентацій студентів, що були обрані як найкращі або просто порівнюючи їх в аудиторії.

Зазначені методи можуть бути використані при організації процесу навчання за моделлю змішаного навчання як окремо, так і в поєднанні один з одним. В реальних умовах одні й ті ж самі методи викладач може використовувати по-різному, спрямовуючи діяльність студентів або на відтворення набутих раніше знань (репродуктивна діяльність), або на самостійне розв'язання нових навчальних завдань (творча діяльність) [303].

**МІНІСТЕРСТВО ОСВІТИ І НАУКИ УКРАЇНИ**
**ДВНЗ «КРИВОРІЗЬКИЙ НАЦІОНАЛЬНИЙ УНІВЕРСИТЕТ»**
**КАФЕДРА ІНЖЕНЕРНОЇ МАТЕМАТИКИ**
**КАФЕДРА ФУНДАМЕНТАЛЬНИХ**
**І СОЦІАЛЬНО-ГУМАНІТАРНИХ ДИСЦИПЛІН**

"ЗАТВЕРДЖУЮ"
Проректор
________/В. І. Вербицький/
________20__ року

# НАВЧАЛЬНА ПРОГРАМА З ДИСЦИПЛІНИ

«ІНФОРМАЦІЙНО-КОМУНІКАЦІЙНІ ТЕХНОЛОГІЇ
НАВЧАННЯ ВИЩОЇ МАТЕМАТИКИ
СТУДЕНТІВ ІНЖЕНЕРНИХ СПЕЦІАЛЬНОСТЕЙ»
для підвищення кваліфікацій викладачів математичних дисциплін технічних ВНЗ та підготовки магістрів математики

| Кількість навчальних тижнів | Кількість годин на тиждень | Разом | Кількість кредитів за ECTS | Аудиторні | | | | | | | Самостійна робота | Підсумковий захід |
|---|---|---|---|---|---|---|---|---|---|---|---|---|
| | | | | Лекції | Практичні | Семінари | Лабораторні | Модуль | Консультації | Разом | | |
| 18 | 3 | 54 | 1,5 | 2 | 14 | | | | | 16 | 36 | 2 |
| Разом | | 54 | 1,5 | 2 | 14 | | | | | 16 | 36 | 2 |

Кривий Ріг – 2013

**Розробник: Кіяновська Н. М.,** асистент кафедри інженерної математики

Робоча програма затверджена на засіданні кафедри фундаментальних і соціально-гуманітарних дисциплін

Протокол від «22» лютого 2013 року № 6

Завідувач кафедри фундаментальних і соціально-гуманітарних дисциплін

«22» лютого 2013 року ____________________ (Семеріков С. О.)

(підпис) (прізвище та ініціали)

## Пояснювальна записка

Сучасний етап інформатизації освіти пов'язаний із широким впровадженням у систему освіти методів і засобів ІКТ, створенням на цій основі комп'ютерно орієнтованого інформаційно-комунікаційного середовища, з наповненням цього середовища електронними науковими, освітніми та управлінськими інформаційними ресурсами, з наданням можливостей суб'єктам освітнього процесу здійснювати доступ до ресурсів середовища, використовувати його засоби і сервіси при розв'язуванні різних завдань.

На сучасному етапі розвитку суспільства і освіти головною метою інформатизації освіти є підготовка тих, хто навчаються, до активної і плідної життєдіяльності в інформаційному суспільстві, забезпечення підвищення якості, доступності та ефективності освіти, створення освітніх умов для широких верств населення щодо здійснення ними навчання протягом усього життя за рахунок широкого впровадження в освітню практику методів і засобів ІКТ та комп'ютерно орієнтованих технологій підтримки діяльності людей.

Створення за рахунок і на основі впровадження ІКТ нових і додаткових умов підвищення якості освіти досягається шляхом:

– розробки і широкого впровадження в практику освіти нових особистісно орієнтованих технологій навчання і учіння;

– диференціації та демократизації навчально-виховного процесу для найбільш повного розвитку схильностей і здібностей людини, задоволення її запитів і потреб, розкриття її творчого потенціалу;

– організації ефективної колективної навчальної діяльності, в тому числі екстериторіальної та спільної міжнародної (освіта без кордонів);

– розширення простору і підвищення ефективності вільного доступу до інформаційних (в тому числі міжнародних) освітніх ресурсів, баз даних і знань, розвиток засобів формування, зберігання, пошуку і представлення інформаційних освітніх матеріалів, створення автоматизованих бібліотечних систем;

– створення нового покоління комп'ютерно орієнтованих засобів навчання, в тому числі, комп'ютерних програмних засобів навчального призначення;

– розвитку засобів оцінювання результатів навчальних досягнень учнів, впливу педагогічних інновацій на результати навчальної діяльності, засобів управління навчанням.

Ефективність процесу інформатизації освіти та його результативність залежить від багатьох чинників, але передусім, від людини, від тих, хто створює комп'ютерно орієнтовані системи навчання й освіти, забезпечують їх впровадження і розвиток в освітній практиці і, безумовно, від якості управління і обсягів ресурсного забезпечення цього процесу. Для забезпечення успіху повинна здійснюватись відповідна високоякісна підготовка і перепідготовка викладацького складу, кадрів управління навчальними закладами і освітою – основної рушійної сили модернізації національної системи освіти на сучасному етапі її розвитку, зокрема її інформатизацією.

Педагоги у всіх країнах світу дуже добре усвідомлюють переваги, що надає методично обґрунтоване використання сучасних інформаційних і комунікаційних технологій у сфері освіти. Використання ІКТ допомагає вирішувати проблеми в тих галузях, де суттєве значення мають знання і комунікація. Сюди входять: вдосконалення процесів учіння / навчання, підвищення рівня навчальних досягнень студентів та їх навчальної мотивації, покращення взаємодії викладачів зі студентами, спілкування в мережі і виконання спільних проектів, вдосконалення організації та управління освітою та навчанням. У зв'язку з цим перед викладачами постає завдання бути обізнаними в останніх досягненнях комп'ютерно орієнтованих технологій, розвивати свою ІКТ-компетентність.

Базуючись на дослідженнях О. М. Спіріна та О. В. Овчарук, під *ІКТ-компетентністю* будемо розуміти підтверджені знання, вміння, ставлення та здатність особистості автономно і відповідально використовувати на практиці ІКТ для задоволення власних індивідуальних потреб і розв'язування суспільно значущих, зокрема професійних, задач у певній предметній галузі або виді діяльності.

Інформаційно-комунікаційна компетентність включає свідоме та критичне застосування технологій інформаційного суспільства для роботи, навчання, відпочинку та спілкування. Вона побудована на застосуванні базових інформаційно-комунікаційних навичок: використання засобів ІКТ для доступу, накопичення, вироблення, подання та обміну даними і відомостями та для спілкування, участі в спільнотах через мережу Інтернет.

Основні знання, вміння та ставлення, що відносяться до цієї компетентності:

– інформаційно-комунікаційна компетентність вимагає свідомого *розуміння* та *знання* природи, ролі та можливостей технологій інформаційного суспільства в особистісному та соціальному житті, навчанні та роботі. Це включає використання комп'ютерних технологій (як, наприклад, текстових редакторів, електронних таблиць, баз даних, масивів даних локального та хмарного зберігання), розуміння можливостей та потенціальних ризиків Інтернету та спілкування через електронні медіа для роботи, навчання, відпочинку, обміну даними і відомостями та колаборативного мережного спілкування, навчання та досліджень;

– особистості повинні також *усвідомлювати*, як технології інформаційного суспільства можуть підтримувати творчість та інноваційність, бути обізнаними про валідність та відповідність даних і відомостей, що на етичних та правових принципах є доступними та залучають до їх використання;

– *уміння* передбачають здатність знаходити, збирати та опрацьовувати дані, відомості і повідомлення та використовувати її систематично та критично, відповідно до реального та віртуального середовища. Особистості повинні володіти вмінням використовувати засоби для розробки, подання та усвідомлення комплексу даних та здатністю до доступу, пошуку та використання сервісів мережі Інтернет;

– також особистості повинні бути *здатними* використовувати засоби ІКТ для підтримки критичного мислення та відповідного ставлення до доступних даних і відомостей та відповідально використовувати медіа. Ця компетентність передбачає здатність входження до соціальних, культурних, професійних спільнот та мереж.

Особистості також повинні бути здатними використовувати ІКТ для підтримки не лише критичного мислення, а й творчості та інновацій:

1) ІКТ-бачення: розуміння та усвідомлення ролі та значення ІКТ для роботи та навчання впродовж життя;

2) ІКТ-культура: спосіб розуміння, конструювання, світоглядного бачення цифрових технологій для життя та діяльності в інформаційному суспільстві;

3) ІКТ-знання: фактичні та теоретичні знання, що відображають галузь ІКТ для навчання та практичної діяльності;

4) ІКТ-практика: практика застосування знань, умінь, навичок у галузі ІКТ для особистих та суспільних професійних та навчальних цілей;

5) ІКТ-удосконалення: здатність удосконалювати, розвивати, генерувати нове у сфері ІКТ та засобами ІКТ для навчання, професійної діяльності, особистого розвитку;

6) ІКТ-громадянськість: підтверджена якість особистості демонструвати свідоме ставлення через дію, пов’язану із застосуванням ІКТ для відповідальної соціальної взаємодії та поведінки.

Доцільно, щоб ці характеристики були максимально можливо відображені наскрізно і на всіх рівнях ІКТ-компетентностей у процесі їх набуття.

Існує безліч причин, які заважають студентам і викладачам в повній мірі використовувати можливості, що з’являються із використанням ІКТ. Це і брак коштів на закупівлю обладнання, і обмежений доступ в Інтернет, і відсутність цифрових освітніх ресурсів на рідній мові. Але головна причина в тому, що викладачі не завжди знають, як ефективно використовувати ІКТ.

Під *ІКТ-компетентністю викладача* будемо розуміти професійно значуще особистісне утворення – здатність педагогічно виважено та методично обґрунтовано виокремлювати, добирати, досліджувати, проектувати та використовувати ІКТ з метою усебічного забезпечення процесу навчання.

У рекомендаціях ЮНЕСКО щодо інформаційно-комунікаційної компетентності викладачів виділено такі напрями:

– *розуміння ролі ІКТ в освіті:* педагоги повинні бути знайомі з освітньою політикою і вміти пояснити на професійній мові, чому їх педагогічні практики відповідають освітній політиці і як її реалізують;

– *навчальна програма і оцінювання:* педагоги повинні відмінно знати освітні стандарти та вимоги щодо оцінювання навчальних досягнень зі свого навчального предмету. Крім того, педагоги повинні бути здатні включити використання засобів ІКТ у свою навчальну програму;

– *педагогічні практики*: педагоги повинні знати передовий досвід теорії та методики використання ІКТ у навчальній роботі та для подання навчального матеріалу;

– *технічні і програмні засоби ІКТ:* педагоги повинні знати базові прийоми роботи з технічними і програмними засобами; програмні засоби, що підвищують продуктивність роботи; Web-браузери; комунікаційні програмні засоби; засоби презентаційної графіки; програми для вирішення завдань управління;

– *організація і управління освітнім процесом:* педагоги повинні вміти використовувати засоби ІКТ для роботи з усім класом, у малих групах, а також для індивідуальної роботи. Вони повинні надавати всім студентам рівний доступ до використання ІКТ;

– *професійний розвиток:* педагоги повинні володіти навичками роботи з ІКТ і знати Web-ресурси, щоб отримувати додаткові навчально-методичні матеріали, необхідні для їх професійного розвитку.

Рівень розвитку інформаційно-комунікаційної складової педагогічної компетентності викладача визначається:

– мотивацією: прагненням до якісного виконання інформаційної діяльності, ефективного використання ІКТ, самовдосконалення тощо;

– інформаційно-науковими знаннями з оперуванням поняттями, актами, властивостями, закономірностями, методами, алгоритмами та ін.;

– уміннями і навичками досвіду роботи з інформаційними джерелами;

– інформаційним світоглядом, інноваційним мисленням та інтуїцією, мобільністю, ціннісними орієнтаціями щодо доцільності використання

інформаційно-комунікаційних технологій.

Спираючись на запропоноване О. М. Спіріним описання рівнів ІКТ-компетентностей, пропонуємо три *рівні ІКТ-компетентності викладача*:

*І рівень, базовий.* Систематично використовувати стандартні засоби ІКТ для підтримки навчання. Самостійно добирати засоби ІКТ для реалізації цілей навчання. Правильно добирати і використовувати ІКТ для розв'язування основних навчальних задач.

*ІІ рівень, поглиблений.* Проводити проектування процесу навчання із використанням ІКТ. Створювати предметно орієнтоване навчальне середовище, сприяти розвитку персональних навчальних середовищ. Застосувати ІКТ для комбінування форм організації, методів та засобів навчання. Уміти розв'язувати професійні задачі підвищеної складності з використанням ІКТ, адаптувати засоби ІКТ для розв'язування основних професійних задач, зокрема бути здатним проектувати, конструювати і вносити інновації до елементів наявних ІКТ навчання.

*ІІІ рівень, дослідницький.* Засвоїти та демонструвати повне володіння методикою використання ІКТ у предметній галузі. Досліджувати, добирати та проектувати засоби ІКТ організації навчального процесу. Зробити оригінальний вклад у розвиток теорії та методики використання ІКТ у процесі навчання, розробляти інноваційні ІКТ навчання.

Запропоновані рівні було визначено, виходячи з вимог державної цільової програми впровадження у навчально-виховний процес загальноосвітніх навчальних закладів ІКТ на період до 2015 року «Сто відсотків».

У визначенні компетентностей викладача вищої математики доцільно скористатися результатами дослідження С. А. Ракова, де вказується на необхідність формування:

– процедурної компетентності як умінь розв'язувати типові математичні та інформатичні задачі;

– логічної компетентності як володіння дедуктивним методом доведення та спростування тверджень;

– технологічної компетентності як умінь застосовувати у професійній діяльності засобів ІКТ;

– дослідницької компетентності як володіння методами дослідження соціально та індивідуально значущих задач математичними методами;

– методологічної компетентності як умінь оцінювати доцільність використання математичних методів для розв'язування індивідуально і суспільно значущих задач.

Викладач з вищої математики технічного ВНЗ має володіти такими **компетентностями**:

– *інформаційна* – здатність викладача до проведення критичного

аналізу джерел даних, пошуку необхідних ресурсів, синтезу, узагальненню та структуруванню опрацьованих відомостей;

– *технічна* – здатність та готовність викладача до ефективного використання та опанування апаратних та програмних засобів ІКТ;

– *технологічна* – здатність та готовність викладача до інформаційно-технологічної діяльності, а саме: постановка цілей створення електронних освітніх ресурсів, використанню існуючої або розробки нової технології для створення електронного освітнього ресурсу, тестуванню створеного продукту на відповідність до певних вимог тощо;

– *педагогічна* – здатність та готовність викладача до педагогічного проектування, змістового наповнення та використання електронних освітніх ресурсів у власній професійній діяльності;

– *мережна та телекомунікаційна* – здатність до опанування основними принципами побудови і використання локальних мереж та глобальної мережі Інтернет;

– *дослідницька* – здатність проводити дослідження доступними засобами ІКТ;

– *в питаннях інформаційної безпеки* – здатність запобігти можливим інформаційним атакам у комп'ютерних системах, володіти знаннями з принципів захисту даних, вміти проводити апаратні та програмні методи захисту даних.

Відповідно до означених компетентностей викладачі з вищої математики технічного ВНЗ повинні знати про: засоби ІКТ для навчання вищої математики у інженерній освіті; можливості та призначення засобів ІКТ; техніку безпечного користування засобами ІКТ; структуру мережі Інтернет та її значення для освіти.

У викладачів вищої математики технічного ВНЗ повинні бути сформовані такі основні **вміння** використання ІКТ:

– працювати з електронною поштою;

– працювати з Web-браузерами;

– використовувати математичні пакети для проведення обчислень та моделювання;

– застосовувати мережні засоби для підтримки спілкування;

– працювати з науковими текстовими процесорами (зокрема LaTeX, Kile);

– використовувати системи відображення документів;

– користуватися програмами автоматизації роботи з даними;

– працювати з периферійним комп'ютерним обладнанням (принтер, сканер, модем, Web-камера тощо);

– проектувати та створювати нові засоби навчання математики.

Ураховуючи, що ефективною та дієвою формою розвитку ІКТ-

компетентності викладачів вищої математики технічних ВНЗ є підвищення кваліфікації.

Розглянемо програму навчального спецкурсу «Інформаційно-комунікаційні технології навчання вищої математики студентів інженерних спеціальностей».

### 1. Мета та завдання навчальної дисципліни

**Метою спецкурсу** є розвиток інформаційно-комунікаційної компетентності викладачів вищої математики технічних ВНЗ та магістрів математики: навчити створювати навчальні текстові документи, таблиці, малюнки, діаграми, презентації, використовуючи Інтернет-технології та бази даних; здійснювати анкетування, діагностування, тестування, пошук даних; розробляти власні електронні продукти та комбінувати готові.

**Завдання спецкурсу:**

– сприяти формуванню стійкої педагогічної позиції щодо необхідності використання засобів ІКТ у власній професійній діяльності;

– систематизувати технічні, педагогічні, інформаційні, технологічні, телекомунікаційні, дослідницькі знання та вміння викладачів, знання викладачів з питань інформаційної безпеки для створення та використання електронних освітніх ресурсів в професійній діяльності;

– удосконалити інформаційно-технологічні навички роботи з апаратними та програмними засобами ІКТ;

– розвинути та удосконалити інформаційно-технологічні вміння щодо створення та використання освітніх інформатичних продуктів в професійній діяльності.

Особливості реалізації мети та завдань спецкурсу:

– навчання у співпраці та через практичну діяльність;

– спрямованість кожного на досягнення особистісно значущого результату;

– використання ІКТ для реалізації педагогічних ідей;

– використання ІКТ для підтримки комунікації зі студентами;

– постійна зміна видів діяльності;

– спрямованість на успіх у всіх видах діяльності;

– постійне обговорення власних думок з колегами;

– проектування всіх видів діяльності з використанням ІКТ;

– можливість самостійного виконання завдання за комп’ютером;

– використання методу демонстраційних прикладів.

У контексті даного спецкурсу ***слухачі курсу повинні знати*** про можливості використання створених або готових ІКТ в організаційно-методичній діяльності або навчальному процесі.

**Очікуваний результат спецкурсу, *слухачі курсу повинні вміти*:**

ефективно використовувати комп'ютерно орієнтовані методи навчання на заняттях з вищої математики.

На заняттях передбачено такі **форми і методи діяльності слухачів курсу**: тематичні лекції, самостійна робота за комп'ютером; виконання практичних завдань; колективне обговорення сучасних проблем стосовно впровадження ІКТ у навчальний процес; використання Інтернет-ресурсів для розроблення дидактичних матеріалів; виконання навчальних завдань: створення комп'ютерної презентації, публікації, Web-сайту, форуму; оцінювання навчальних досягнень за допомогою тестових систем тощо.

Курс призначений для підвищення кваліфікації викладачів вищої математики та магістрів математики, які опановують нові технології та засоби ІКТ з метою забезпечення зростання ефективності освітнього процесу та якості навчання вищої математики у вищих технічних закладах. Зміст курсу може бути адаптований для будь-якої категорії слухачів очної або очно-дистанційної форми організації навчання. На вивчення спецкурсу відводиться 54 годин / 1,5 кредити ECTS.

## 2. Програма навчальної дисципліни

**Змістовий модуль 1**. Інформаційно-комунікаційні технології навчання вищої математики студентів інженерних спеціальностей.

**Тема 1**. Інформаційно-комунікаційні технології навчання вищої математики студентів інженерних спеціальностей – класифікація та дидактичні можливості.

*Інформаційні технології. Інформаційно-комунікаційні технології. Використання термінів «інформаційні технології» та «інформаційні та комунікаційні технології» у наукових джерелах. Етапи розвитку теорії та методики використання ІКТ в інженерній освіті. ІКТ, що можуть бути використані у процесі навчання математики: електронно-навчальні ресурси, системи підтримки навчання, мобільне математичне середовище, системи комп'ютерної математики, системи тестування.*

**Тема 2.** Наочні засоби навчання: засоби мультимедіа, мультимедійна дошка, Web-конференції.

*Методика застосування сучасних інформаційно-комунікаційних засобів навчання. Мультимедіа. Вебінари у Dropbox. Мультимедійні дошки. Онлайн навчання у WiZiQ.*

**Тема 3.** Використання систем комп'ютерної математики, бази знань Wolfram|Alpha, математичних онлайн ресурсів на заняттях вищої математики. Розробка та використання тренажерів.

*Місце систем комп'ютерної математики на заняттях з вищої математики. Розробка та використання тренажерів. Математичні*

*онлайн ресурси на заняттях вищої математики. Використання бази знань Wolfram|Alpha.*

**Тема 4.** Використання систем управління навчанням (LMS) для розробки навчальних електронних курсів.

*Робота з Moodle. Про користувачів. Навчальні матеріали курсу. Робота з глосаріями. Робота з мультимедійними матеріалами. Ресурс Книга. Типи тестів у Moodle. Створення та використання тестів. Завдання на курсі. Модуль Урок.*

**Тема 5.** Розробка та проведення онлайн тестування з використанням тестових систем.

*Сервіс Майстер-тест. Сервіс Банк Тестов.RU. Сайт Твой тест. Система електронного тестування Tests Online.*

**Тема 6.** Використання Інтернет-ресурсів для інтерактивної позааудиторної взаємодії викладачів та студентів: платформа Piazza, Skype, Google+.

*Створення навчального класу в Piazza, використання наукових текстових редакторів LaTeX (середовище для створення математичних текстів). Реєстрація та проведення конференцій в Skype. Створення форуму в Google+.*

**Тема 7.** Використання сервісів Google у методичній роботі викладача: таблиці, презентації, документи, Writely та Geogebra форми для проведення анкетування. Сервіси Microsoft: Office Web Apps.

*Використання таблиць Google для ведення обліку успішності студентів. Проведення презентацій на заняттях та онлайн. Використання документів та рисунків Google. Особливості анкетування в Google. Динамічні документи Geogebra. Сервіси Microsoft: Office Web Apps.*

**Тема 8.** Проектування та розробка персональних педагогічних Web-ресурсів – авторські сайти викладачів.

*Служба Сайти Google. Створення посилань на сторінці сайту. Створення нової сторінки сайту. Вбудовування форуму у сайт. Служби блогів.*

### 3. Структура навчальної дисципліни

| Назви змістових модулів і тем | Кількість годин | | | | |
|---|---|---|---|---|---|
| | усього | у тому числі | | | |
| | | лек. | практ. | інд. | с. р. |
| **Модуль 1** | | | | | |
| **Тема 1**. Інформаційно-комунікаційні технології навчання вищої математики студентів інженерних спеціальностей – | | 2 | | | 4 |

| Назви змістових модулів і тем | Кількість годин | | | | |
|---|---|---|---|---|---|
| | усього | у тому числі | | | |
| | | лек. | практ. | інд. | с. р. |
| класифікація та дидактичні можливості. | | | | | |
| **Тема 2.** Наочні засоби навчання: засоби мультимедіа, мультимедійна дошка, Web-конференції. | | | 2 | | 4 |
| **Тема 3.** Використання систем комп'ютерної математики, бази знань Wolfram\|Alpha, математичних онлайн ресурсів на заняттях вищої математики. Розробка та використання тренажерів. | | | 2 | | 4 |
| **Тема 4.** Використання систем управління навчанням (LMS) для розробки навчальних електронних курсів. | | | 2 | | 4 |
| **Тема 5.** Розробка та проведення онлайн тестування з використанням тестових систем. | | | 2 | | 6 |
| **Тема 6.** Використання Інтернет-ресурсів для інтерактивної позааудиторної взаємодії викладачів та студентів: платформа Piazza, Skype, Google+. | | | 2 | | 4 |
| **Тема 7.** Використання сервісів Google у методичній роботі викладача: таблиці, презентації, документи, Writely та Geogebra, форми для проведення анкетування. Сервіси Microsoft: Office Web Apps. | | | 2 | | 4 |
| **Тема 8.** Проектування та розробка персональних педагогічних Web-ресурсів – авторські сайти викладачів. | | | 2 | | 6 |
| Підсумковий тестовий контроль. | | | | | 2 |
| Разом за змістовим модулем 1 | 54 | 2 | 14 | | 38 |
| **Усього годин** | 54 | 2 | 14 | | 38 |

| | Аудит. | Самост. | Загальна |
|---|---|---|---|
| Всього годин | 16 | 38 | 54 |
| Всього кредитів за ECTS | 1,5 | | |

## 4. Теми та завдання практичних занять

| № з/п | Назва теми | Кількість годин |
|---|---|---|
| 1 | Наочні засоби навчання: засоби мультимедіа, мультимедійна дошка, Web-конференції.<br>**Практичне завдання № 1.**<br>Організувати проведення Інтернет-конференції.<br>**Практичне завдання № 2.**<br>Прийняти участь в Інтернет-конференції. | 2 |
| 2 | Використання систем комп'ютерної математики, бази знань Wolfram\|Alpha, математичних онлайн ресурсів на заняттях вищої математики. Розробка та використання тренажерів.<br>**Практичне завдання № 1.**<br>Виконати запропоновані індивідуальні завдання з використанням СКМ, бази знань Wolfram\|Alpha, математичних онлайн ресурсів.<br>**Практичне завдання № 2.**<br>Протестувати роботу розроблених тренажерів. | 2 |
| 3 | Використання систем управління навчанням (LMS) для розробки навчальних електронних курсів.<br>**Практичне завдання № 1.**<br>Ознайомитись з можливостями роботи в Moodle.<br>**Практичне завдання № 2.**<br>Розробити навчальний електронний курс в Moodle. | 2 |
| 4 | Розробка та проведення онлайн тестування з використанням тестових систем.<br>**Практичне завдання № 1.**<br>Розробити в тестових онлайн системах тест.<br>**Практичне завдання № 2.**<br>Розробити тест за допомогою Moodle. | 2 |
| 5 | Використання Інтернет-ресурсів для інтерактивної позааудиторної взаємодії викладачів та студентів: платформа Piazza, Skype, Google+.<br>**Практичне завдання № 1.**<br>Створити власний навчальний клас в Piazza.<br>**Практичне завдання № 2.**<br>Зареєструватися та провести конференцію в Skype.<br>**Практичне завдання № 3.**<br>Створити власний форум в Google. | 2 |
| 6 | Використання сервісів Google у методичній роботі | 2 |

| №<br>з/п | Назва теми | Кількість<br>годин |
|---|---|---|
| | викладача: таблиці, презентації, документи, Writely та Geogebra, форми для проведення анкетування. Сервіси Microsoft: Office Web Apps.<br>**Практичне завдання № 1.**<br>Створити електронний журнал викладача.<br>**Практичне завдання № 2.**<br>Створити презентацію вибраної теми.<br>**Практичне завдання № 3.**<br>Розробити онлайн анкетування.<br>**Практичне завдання № 4.**<br>Використовуючи сервіси Google, створити документ. | |
| 7 | Проектування та розробка персональних педагогічних Web-ресурсів – авторські сайти викладачів.<br>**Практичне завдання № 1.**<br>Розробити власний Web-сайт.<br>**Практичне завдання № 2.**<br>Доповнити розроблений сайт підтримкою онлайн спілкування. | 2 |
| | Разом | 14 |

## 5. Самостійна робота

| №<br>з/п | Назва теми | Кількість<br>годин |
|---|---|---|
| 1 | Інформаційно-комунікаційні технології навчання вищої математики студентів інженерних спеціальностей – класифікація та дидактичні можливості.<br>**Самостійна робота.**<br>Опрацювати конспект лекції. | 4 |
| 2 | Наочні засоби навчання: засоби мультимедіа, мультимедійна дошка, Web-конференції.<br>**Самостійна робота.**<br>Підготувати тези доповіді на конференції. | 4 |
| 3 | Використання систем комп’ютерної математики, бази знань Wolfram\|Alpha, математичних онлайн ресурсів на заняттях вищої математики. Розробка та використання тренажерів.<br>**Самостійна робота.**<br>Підібрати приклади для розв’язання в СКМ. | 4 |
| 4 | Використання систем управління навчанням (LMS) для | 4 |

| № з/п | Назва теми | Кількість годин |
|---|---|---|
| | розробки навчальних електронних курсів. **Самостійна робота.** Підібрати методичні матеріали з вибраного розділу з вищої математики. | |
| 5 | Розробка та проведення онлайн тестування з використанням тестових систем. **Самостійна робота.** Розробити текст тесту з розділу з вищої математики. | 4 |
| 6 | Використання Інтернет-ресурсів для інтерактивної позааудиторної взаємодії викладачів та студентів: платформа Piazza, Skype, Google+. **Самостійна робота.** Конспект лекції та практичного з обраної теми з вищої математики. | 4 |
| 7 | Використання сервісів Google у методичній роботі викладача: таблиці, презентації, документи, Writely та Geogebra, форми для проведення анкетування. Сервіси Microsoft: Office Web Apps. **Самостійна робота.** Розробити текст анкети. Підібрати матеріали для презентації. | 6 |
| 8 | Проектування та розробка персональних педагогічних Web-ресурсів – авторські сайти викладачів. **Самостійна робота.** Підібрати методичні матеріали з вибраного розділу з вищої математики. | 6 |
| 9 | Підсумковий контроль. | 2 |
| Разом | | 38 |

**Оформлення звіту про самостійну роботу**

Підсумком самостійної роботи над вивченням дисципліни «Інформаційно-комунікаційні технології навчання вищої математики студентів інженерних спеціальностей» є система електронних ресурсів, що була створена для певного розділу вищої математики із використанням вивчених засобів.

## 6. Індивідуальні завдання

| № з/п | Назва теми | Кількість годин |
|---|---|---|
| 1 | Інформаційно-комунікаційні технології навчання вищої математики студентів інженерних спеціальностей – класифікація та дидактичні можливості. | 4 |
| 2 | Наочні засоби навчання: засоби мультимедіа, мультимедійна дошка, Web-конференції. | 4 |
| 3 | Використання систем комп’ютерної математики, бази знань Wolfram\|Alpha, математичних онлайн ресурсів на заняттях вищої математики. Розробка та використання тренажерів. | 4 |
| 4 | Використання систем управління навчанням (LMS)для розробки навчальних електронних курсів. | 4 |
| 5 | Розробка та проведення онлайн тестування з використанням тестових систем. | 6 |
| 6 | Використання Інтернет-ресурсів для інтерактивної позааудиторної взаємодії викладачів та студентів: платформа Piazza, Skype, Google+. | 6 |
| 7 | Використання сервісів Google у методичній роботі викладача: таблиці, презентації, документи, Writely та Geogebra, форми для проведення анкетування. Сервіси Microsoft: Office Web Apps. | 4 |
| 8 | Проектування та розробка персональних педагогічних Web-ресурсів – авторські сайти викладачів. | 4 |
| 9 | Підсумковий тестовий контроль. | 2 |
| | Разом | 38 |

Звіт про виконану роботу – це є портфоліо по курсу та виконується кожним слухачем курсу згідно індивідуально отриманих завдань. У звіті мають бути показані вміння застосовувати отримані теоретичні знання та навички у впроваджені ІКТ у процес навчання вищої математики.

## 7. Методи контролю
## Критерії оцінювання знань слухачів курсу

| Рівні навчальних досягнень | Критерії оцінювання рівня навчальних досягнень |
|---|---|
| І. Базовий | Слухач курсу описує поняття ІКТ, знає про навчальне призначення стандартних засобів ІКТ. Слухач курсу самостійно добирає засоби ІКТ для реалізації цілей навчання. Правильно добирає і |

| Рівні навчальних досягнень | Критерії оцінювання рівня навчальних досягнень |
|---|---|
| | використовує ІКТ для розв’язування основних навчальних задач. |
| ІІ. Поглиблений | Слухач курсу вміє проектувати процес навчання на основі ІКТ, створювати предметно орієнтоване навчальне середовище, сприяти розвитку персональних навчальних середовищ. Вміє застосувати ІКТ для комбінування форм організації, методів та засобів навчання. Уміє розв’язувати професійні задачі підвищеної складності з використанням ІКТ, адаптувати засоби ІКТ для розв’язування основних професійних задач, зокрема бути здатним проектувати, конструювати і вносити інновації до елементів наявних ІКТ навчання. |
| ІІІ. Дослідницький | Слухач курсу засвоїв та демонструє повне володіння методикою використання інформаційно-комунікаційних технологій у предметній галузі. Досліджує, добирає та проектує засоби ІКТ організації навчального процесу. Може зробити оригінальний вклад у розвиток теорії та методики використання ІКТ у процесі навчання, вміє розробляти інноваційні ІКТ навчання. |

## 8. Розподіл балів, які отримують студенти
## Шкала оцінювання: національна та ECTS

| Сума балів | Оцінка ECTS | Оцінка за національною шкалою | |
|---|---|---|---|
| | | **для екзамену** | **для заліку** |
| 90 – 100 | **A** | відмінно | зараховано |
| 82 – 89 | **B** | добре | |
| 74 – 81 | **C** | | |
| 64 – 73 | **D** | задовільно | |
| 60 – 63 | **E** | | |
| 35 – 59 | **FX** | незадовільно з можливістю повторного складання | не зараховано з можливістю повторного складання |
| 0 – 34 | **F** | незадовільно з обов’язковим повторним вивченням дисципліни | не зараховано з обов’язковим повторним вивченням дисципліни |

## 9. Методичне забезпечення

**Орієнтований перелік програмного забезпечення**

| Тип програмного забезпечення | Приклади програм |
|---|---|
| Операційна система з графічним інтерфейсом | Windows, Linux |
| Web-браузер | Chrome, Internet Explorer, Opera, Firefox |
| Програми для створення електронних презентацій | MS PowerPoint, OO Impress |
| Текстовий редактор | MS Word, OO Writer, LaTeX |
| Програвачі, що надають можливість відтворювати різноманітні мультимедійні файли | Winamp, Mplayer |

**Для навчально-методичного забезпечення спецкурсу необхідні такі технічні засоби:**

– сучасна комп'ютерна техніка, що об'єднано у локальну мережу;

– підключення до глобальної мережі Internet.

**Під час вивчення матеріалу розглядаються такі мережні сервіси:**

1. Платформа для проведення онлайн семінарів (www.wiziq.com);
2. Платформа Piazza (piazza.com);
3. Сервіси Google (drive.google.com);
4. Сервіси Microsoft (Office Web Apps и т.п.);
5. Сайти Google (sites.google.com);
6. Тестові системи Tests Online (tests-online.ru), Твой тест (www.make-test.ru), Банк Тестов.RU (www.banktestov.ru), Майстер-тест (www.master-test.net/uk).

## 10. Рекомендована література

**Базова**

**Допоміжна**

*Р Рекомендації щодо вибору методів, форм та засобів навчання розділів вищої математики студентів інженерних спеціальностей у технічних ВНЗ*

*Таблиця* Р. 1

**Методи, форми та засоби навчання відповідно до розділів вищої математики, що вивчаються студентами інженерних спеціальностей у технічних ВНЗ**

| № з/п | Розділи | Форми | Методи | Засоби | |
|---|---|---|---|---|---|
| | | | | спільні | специфічні |
| 1. | Лінійна алгебра | лекція; комп'ютерно-орієнтовані практичні заняття; виконання індивідуальних завдань; онлайн-консультації | проблемний | Сайти Google, GoogleDocs, DropBox, Piazza, Moodle, Wolfram\|Alpha, WiZiQ, YouTube, Google+ | Sage, MathCad |
| 2. | Векторна алгебра | лекція; практичне заняття; виконання індивідуальних завдань; контрольні заходи; онлайн-консультації | інформаційно-рецептивний; графічні роботи | | GeoGebra, Wolfram Demonstrations Project |
| 3. | Аналітична геометрія | лекція; комп'ютерно-орієнтовані практичні заняття; самостійна робота студентів; комп'ютерний контроль якості знань; онлайн-консультації | інформаційно-рецептивний; демонстрація (у динаміці) | | Smath Studio, Maple |
| 4. | Границя | лекція; практичне заняття; виконання індивідуальних завдань; онлайн-консультації | частково-пошуковий (евристичний) | | Maple, GRAN, MathCad, онлайн сервіс FooPlot (fooplot.com) |
| 5. | Диференціальне числення функції однієї змінної | лекція; комп'ютерно-орієнтовані практичні заняття; самостійна робота студентів; онлайн-консультації | інформаційно-рецептивний; демонстрація (у динаміці) | | Maple, GRAN, Wolfram Demonstrations Project, FooPlot (fooplot.com) |

| № з/п | Розділи | Форми | Методи | Засоби | |
|---|---|---|---|---|---|
| | | | | спільні | специфічні |
| 6. | Невизначений та визначений інтеграли | лекція; практичне заняття; виконання індивідуальних завдань; комп'ютерний контроль якості знань; онлайн-консультації | пошуковий (дослідний) | | Sage, GRAN, MathCad |
| 7. | Функції багатьох змінних | лекція; комп'ютерно-орієнтовані практичні заняття; контрольні заходи; самостійна робота студентів; онлайн-консультації | інформаційно-рецептивний; демонстрація (у динаміці) | | Maple, GeoGebra |
| 8. | Диференціальні рівняння | лекція; комп'ютерно-орієнтовані практичні заняття; виконання індивідуальних завдань; контрольні заходи; онлайн-консультації | пошуковий (дослідний) | | Smath Studio, MathCad |
| 9. | Кратні інтеграли | лекція; комп'ютерно-орієнтовані практичні заняття; контрольні заходи; самостійна робота студентів; онлайн-консультації | проблемний | | MathCad |
| 10. | Ряди | лекція; практичне заняття; виконання індивідуальних завдань; комп'ютерний контроль якості знань; онлайн-консультації | інформаційно-рецептивний | | MathCad |
| 11. | Теорія ймовірностей | лекція; практичне заняття; контрольні заходи; самостійна робота студентів; онлайн-консультації. | проблемний | | Voki, GRAN, MathCad |

| № з/п | Розділи | Форми | Методи | Засоби | |
|---|---|---|---|---|---|
| | | | | спільні | специфічні |
| 12. | Основи математичної статистики | лекція; комп'ютерно-орієнтовані практичні заняття; практична підготовка; контрольні заходи; онлайн-консультації | частково-пошуковий (евристичний) | | GRAN, Statistica |

Розглянемо детальніше розділ «Функції багатьох змінних».

Для надання основних відомостей з теми (текст лекцій у друкованому виді, приклади завдань до практичних занять, завдання для індивідуальних робіт) доцільно використовувати власний сайт викладача, де в зручному та зрозумілому для студентів порядку можна викласти весь навчальний матеріал. Основною перевагою власного сайту є те, що можна надати відкритий доступ до навчальних матеріалів без необхідності реєстрації (рис. Р.1).

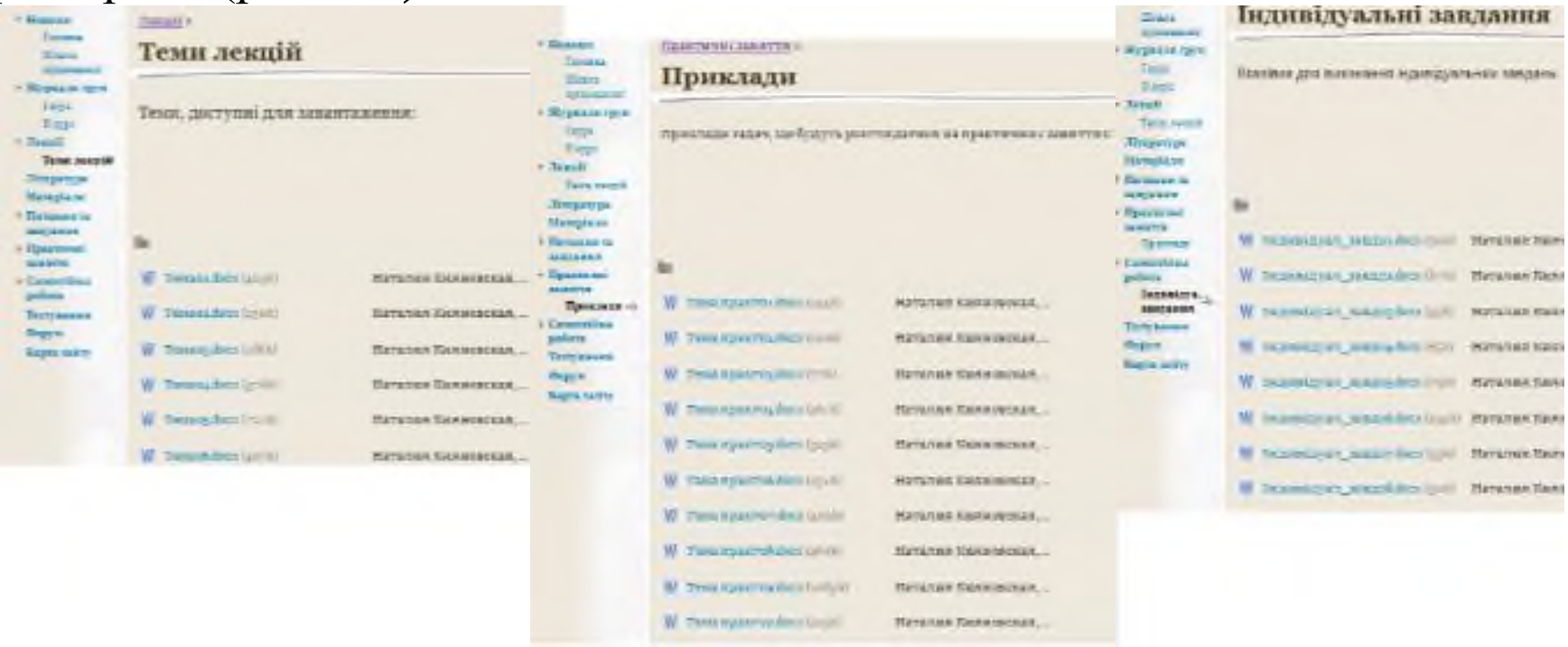

Рис. Р.1. Матеріали, доступні для завантаження

Проведення аудиторних занять та консультацій необхідно доповнити онлайн консультаціями в платформі Piazza. Використання Piazza надає можливість провести додатково роз'яснення навчального матеріалу та з'ясувати рівень його засвоєння. Крім того, в Piazza студенти можуть оперативно вирішити проблемні питання (рис. Р.2).

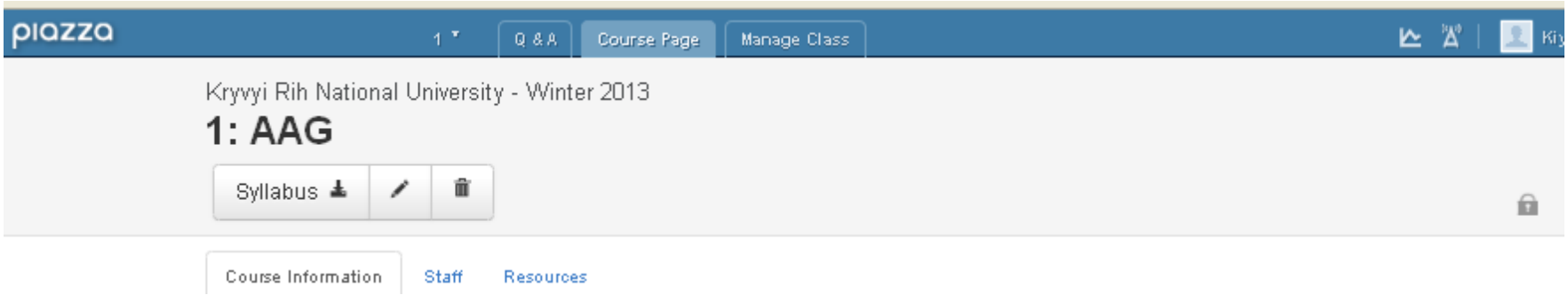

Рис. Р.2. Платформа Piazza для позааудиторних консультацій

Інтегрування відеоматеріалу з теми, розміщеного на YouTube, на авторський сайт та в систему дистанційного навчання, зокрема Moodle, допоможе урізноманітнити подання навчального матеріалу та надасть можливість студентам обрати кращій спосіб для засвоєння нового матеріалу (рис. Р.3).

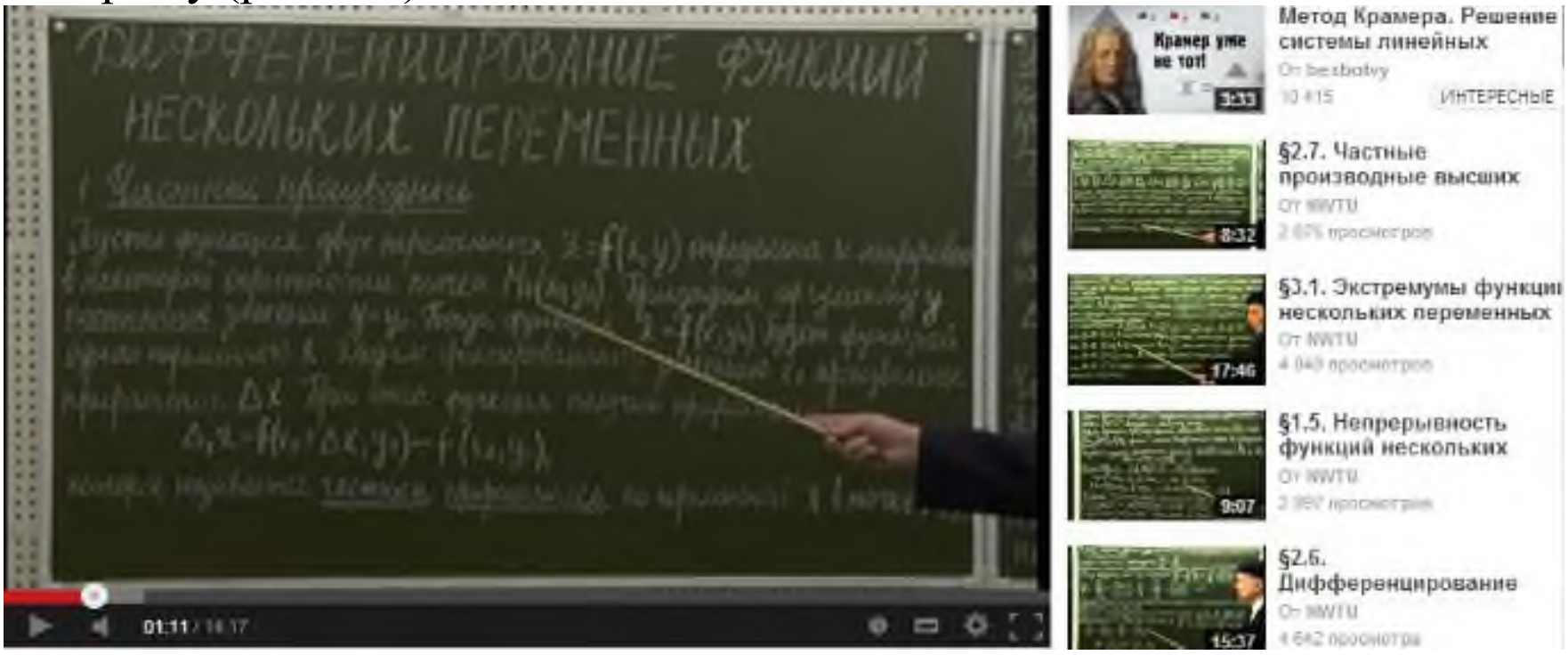

Рис. Р.3. Матеріал з теми, розміщений на YouTube

Проведення лекційних демонстрацій та візуалізація навчального матеріалу на практичних заняттях можлива із використанням СДГ GeoGebra. Так, при необхідності накреслити лінії рівня заданої функції, наприклад, $z = (x - y)^2$, необхідно провести побудову функцій наступного виду: $z = C$, тобто $C = x^2 - 2xy + y^2$, де $C \in R$. Побудову ліній рівня заданої функції було проведено у вікні GeoGebra, надаючи $C$ довільні значення з множини дійсних чисел (рис. Р.4).

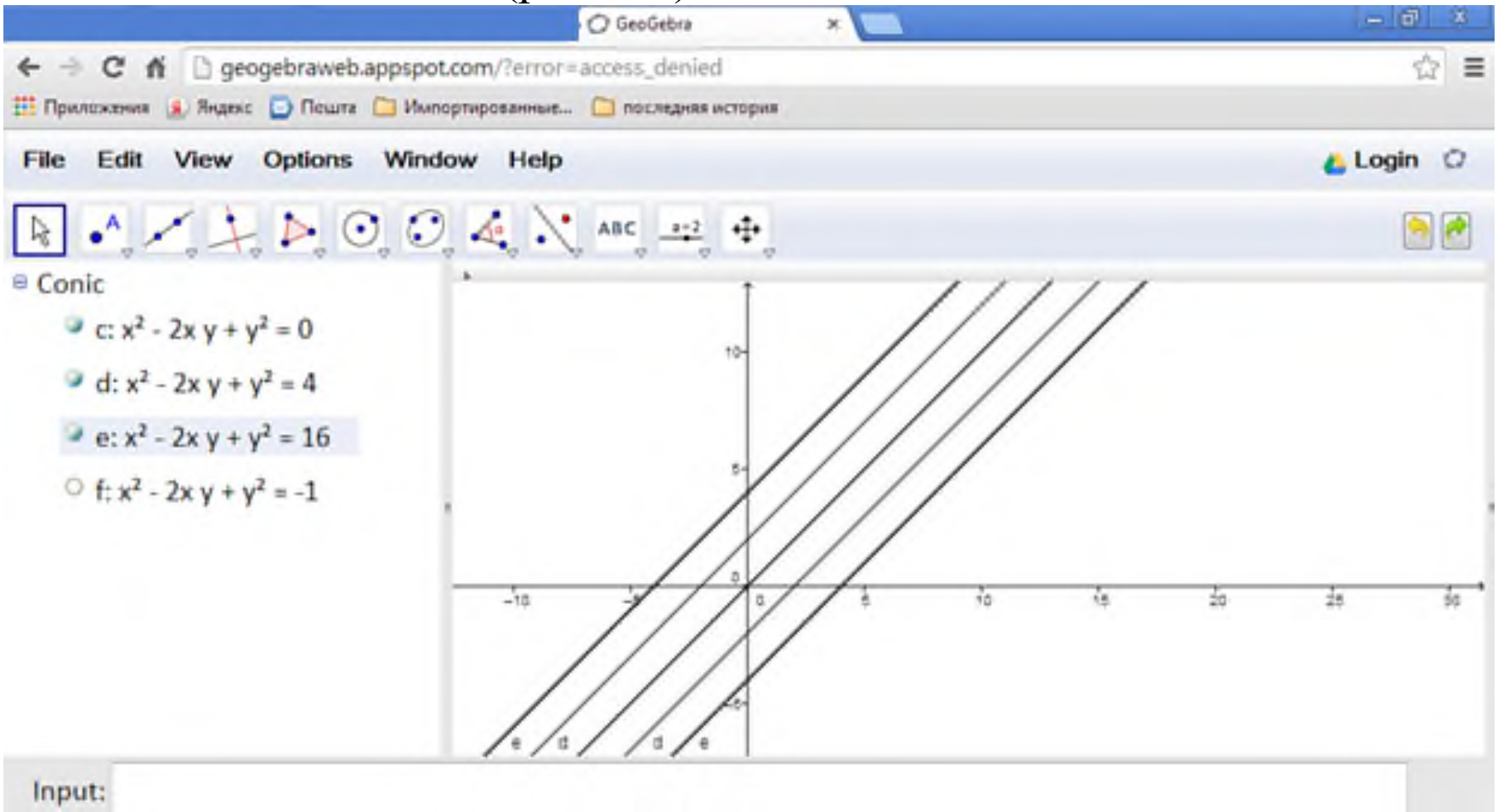

Рис. Р.4. Побудова ліній рівня заданої функції у вікні GeoGebra

При розв’язуванні завдань на знаходження екстремумів функції перевірку відповідей зручно проводити у Wolfram|Alpha, де крім побудованого графіка функції, наведено значення частинних похідних та значення локального екстремуму (рис. Р.5).

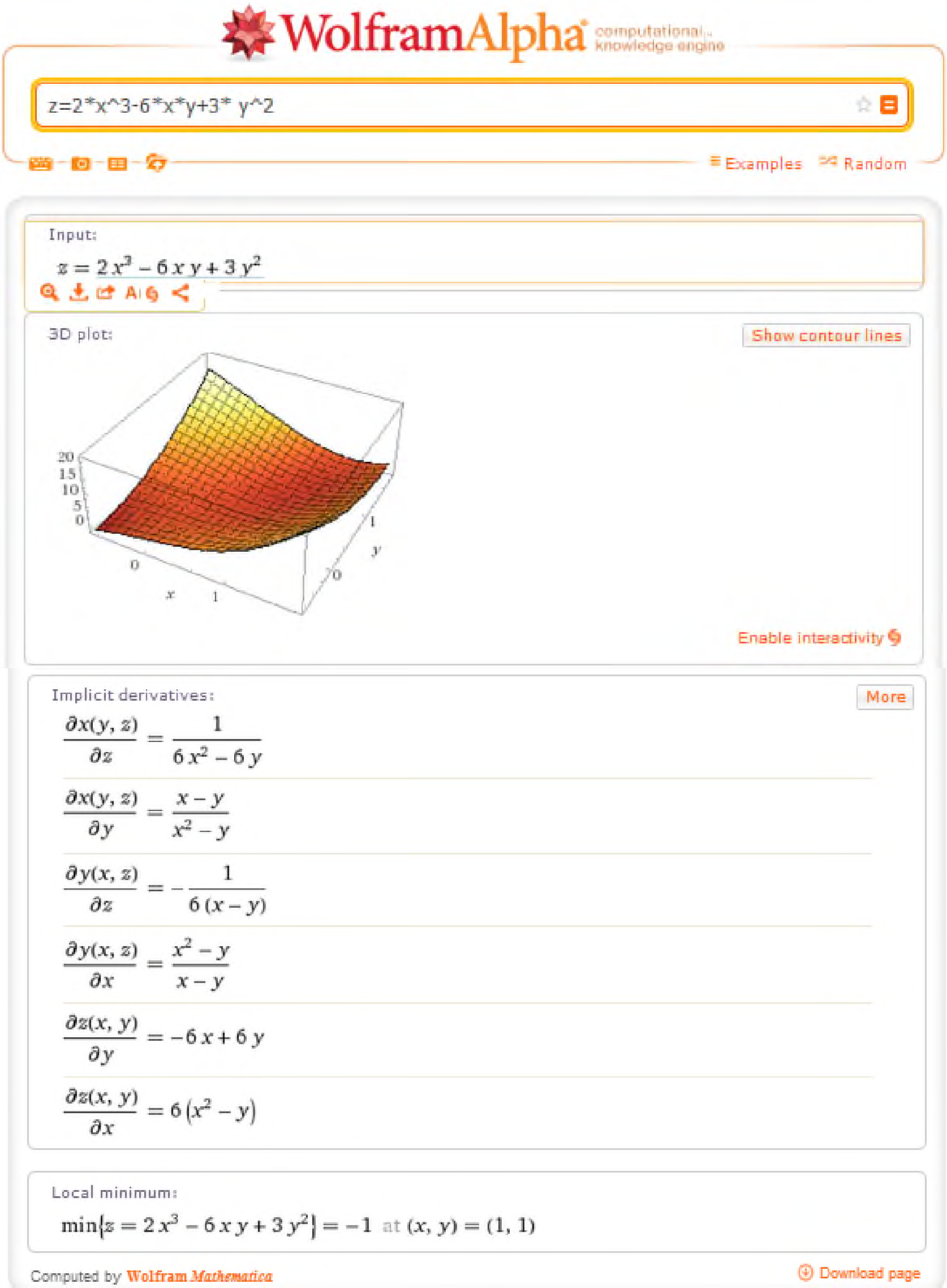

Рис. Р.5. Використання Wolfram|Alpha при вивченні теми «Функції багатьох змінних»

Проводячи побудови поверхонь у СКМ, студенти отримують зображення поверхонь у тривимірному просторі з можливістю зміни кута, під яким поверхня розташована у просторі, що сприяє кращому розумінню навчального матеріалу та розвитку просторової уяви (із використанням СКМ Maple на рис. Р.6).

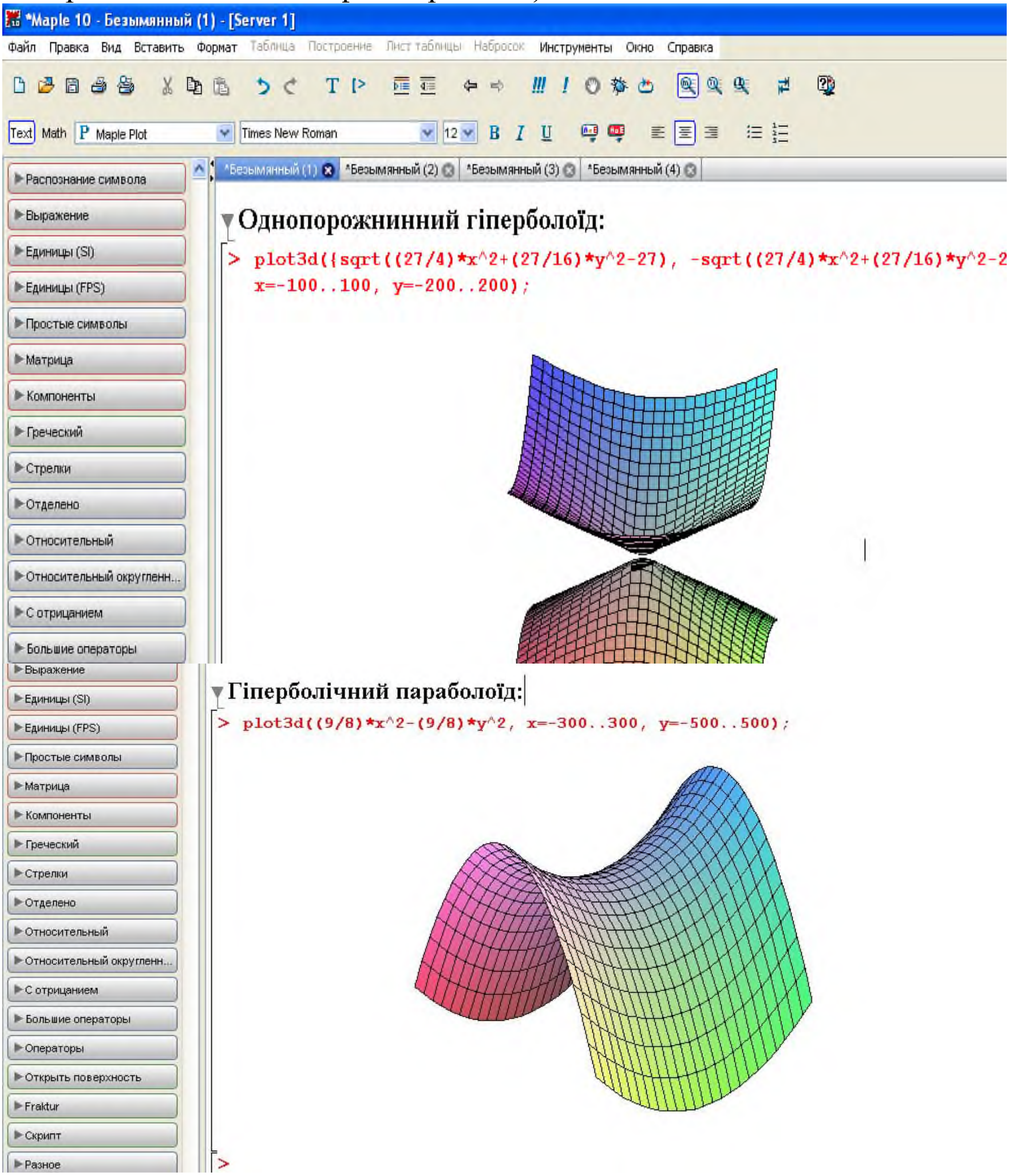

Рис. Р.6. Використання Maple для побудови поверхонь

| | |
|---|---|
| engineering education | інженерна освіта |
| Accreditation Board for Engineering and Technology (ABET) | Рада з акредитації інженерних розробок і технологій |
| Engineering Criteria 2000 – EC2000 | Інженерні критерії 2000 року |
| STEM – science, technology, engineering, mathematics | Наука, технологія, інженерія та математика |
| National Academy of Engineering | Національна академія інженерії |
| Massachusetts Institute of Technology | Массачусетський технологічний інститут |
| Fall 2012 | Осінній семестр 2012 |
| Course | напрям підготовки |
| Differential Equations | курс диференціальних рівнянь |
| MIT School of Engineering | Школа інженерії МТІ |
| Civil and Environmental Engineering | цивільна інженерія та охорона навколишнього середовища |
| Mechanical Engineering | машинобудування |
| Materials Science and Engineering | матеріалознавство |
| Electrical Engineering and Computer Science | електротехніка та комп'ютерні науки |
| Chemical Engineering | хімічна інженерія |
| Aeronautics and Astronautics | аеронавтика та астронавтика |
| Biological Engineering | біотехнології |
| Nuclear Science and Engineering | ядерна інженерія |
| The MIT Flexible Engineering Degree Program | Гнучка програма підготовки інженерів МТІ |
| Freshman | першокурсники |
| core curriculum | набір обов'язкових навчальних дисциплін |
| General Institute Requirements – GIRs | загальноінститутські вимоги |
| Calculus I | аналіз функції однієї змінної, Числення однієї змінної |
| Calculus II | аналіз функції багатьох змінних, Числення багатьох змінних |
| Linear Algebra | лінійна алгебра |
| Computer Engineering | комп'ютерна інженерія |
| Information and communications | Інформаційно-комунікаційні технології |
| technology (ICT) | (ІКТ) |
| National Education Technology Plan, Transforming American Education: Learning Powered by Technology | Трансформація американської освіти: навчання за допомогою технологій |
| Virtual learning environments (VLE) | віртуальні навчальні середовища |
| Advanced Distributed Learning (ADL) | Розвиток розподіленого навчання |
| SCORM (Sharable Content Object Reference Model) | Зразкова модель об'єкта вмісту для спільного використання |
| LMS (Learning Management System) | Система управління навчанням |
| mobile learning, M-Learning | мобільне навчання |
| MIT Computer Science and Artificial Intelligence Laboratory | Лабораторія комп'ютерних наук та штучного інтелекту MTI |
| Math Learning Center | студентський Центр навчання математики |
| Math Academic Services – MAS | «Математичні академічні послуги» |
| The National Assessment of Educational Progress – NAEP | «Національна оцінка освітніх досягнень» |
| Massive open online course – MOOC | Масові відкриті дистанційні курси МВДК |